\documentclass{acta_deli}

\usepackage{xcolor}

\usepackage{amsfonts}
\usepackage{amsmath}
\usepackage{amssymb}
\usepackage{bm}
\usepackage{graphicx}
\usepackage{harvard_acta}
\usepackage{tikz}

\usetikzlibrary{arrows,calc}
\numberwithin{figure}{section}

\def\abso#1{{|#1|}}
\def\alpb{{\boldsymbol\alpha}}
\def\alpbd{{\alpb_\del}}
\def\alpbs{{\alpb_s}}
\def\bal{\begin{aligned}}
\def\ballx{{B_\del(\xb)}}

\def\ibt{\begin{itemize}}

\def\ccc{{c}}
\def\CCC{{C}}
\def\CCd{{C_\del}}
\def\cds{{C_{\ddd,s}}}
\def\cdotsp{{\,\cdot\,}}

\def\cds{{C_{\ddd,s}}}
\def\Dbf{{\textbf{\em D}}}
\def\ddd{{d}}

\def\del{{\delta}}
\def\dxb{{\D \xb}} %GS%
\def\dyb{{\D \yb}} %GS%
\def\dydx{{\D \yb \D \xb}} %GS%
\def\eal{\end{aligned}}
\def\eit{\end{itemize}}
\def\eps{{\epsilon}}
\def\fb{{\bm f}}
\def\fv{{\lrang{f}{v}}}
\def\fvh{{\lrang{f}{\vh}}}
\def\gam{{\gamma}}
\def\Gam{{\Gamma}}
\def\Gamp{{\Gamma_\phi}}
\def\Gamt{{\Gamma_\theta}}
\def\gamd{{\gamma_\del}}

\def\gams{{\gam_{s}}}
\def\GGG{{G}}

\def\hb{{\bm h}}
\def\HHH{{H}}
\def\HHHh{{\widehat\HHH}}

\def\Hsoi{{H^s(\omgomgi)}}
\def\ib{{\bm i}}
\def\Ibf{{\textbf{\em I}}}
\def\intbx{{\int_{{\ballx}}}}

\def\intoi{{\int_\omgomgi}}

\def\jb{{\bm j}}

\def\laa{{\lambda}}
\def\lrang#1#2{\langle{#1},{#2}\rangle}

\def\Mbzd{{{\mathbb Z}^\ddd}}
\def\mcA{{\mathcal A}}
\def\mcB{{\mathcal B}}
\def\mcC{{\mathcal C}}

\def\mcD{{\mathcal D}}
\def\mcDd{{{\mathcal D}_\del}}
\def\mcDf{{{\mathcal D}_s}}
\def\mcDs{{{\mathcal D}^\ast}}
\def\mcDsd{{{\mathcal D}_\del^\ast}}
\def\mcDsf{{{\mathcal D}^\ast_s}}
\def\mcDdsw{{{\mcD^\ast_{\del,\omega}}}}
\def\mcE{{\mathcal E}}
\def\mcF{{\mathcal F}}
\def\mcG{{\mathcal G}}

\def\mcI{{\mathcal I}}
\def\mcJ{{\mathcal J}}

\def\mcL{{\mathcal L}}
\def\mcLd{{\mcL_\del}}
\def\mcLdh{{\mcL_{\del,h}}}

\def\mcLz{{\mcL_0}}
\def\mcM{{\mathcal M}}
\def\mcN{{\mathcal N}}
\def\mcNd{{\mcN_\del}}
\def\mcQ{{\mathcal Q}}
\def\mcR{{\mathcal R}}
\def\mcT{{\mathcal T}}
\def\mcTh{{\mcT_h}}
\def\mcV{{\mathcal V}}
\def\mcX{{\mathcal X}}
\def\mcXd{{\mcX_{\ballx}(\yb)}}
\def\midtilde{{\raisebox{-0.25\baselineskip}{\textasciitilde}}}
\def\mmm{{m}}
\def\mumu{{\mu}}
\def\nb{{\bm n}}
\def\NNh{{N_h}}
\def\norm#1{{\|#1\|}}
\def\nub{{\boldsymbol\nu}}
\def\nubd{{\nub_\del}}
\def\omg{{\Omega}}
\def\omgi{{\omg_{\mcI}}}
\def\omgid{{\omg_{{\mcI}_\del}}}
\def\omgii{{\omg_{{\mcI}_\infty}}}
\def\omgomgi{{\omg\cup\omg_{{\mcI}_\del}}}
\def\omgomgip{{(\omgomgi)\times(\omgomgi)}}
\def\pb{{{\bm p}}}
\def\pbp{{{\bm p}_\phi}}
\def\pbt{{{\bm p}_\theta}}
\def\phid{{\phi_\del}}

\def\phiv{{\varphi}}
\def\Rd{{{\mathbb R}^\ddd}}
\def\Rdc{{\Rd\setminus\omg}}

\def\rhob{{\boldsymbol\rho}}
\def\RNp{{{\mathbb R}^{N_{\phi}}}}
\def\RNt{{{\mathbb R}^{N_{\theta}}}}
\def\Ro{{\mathbb R}}
\def\rrr{{\rho}}

\def\Theb{{\boldsymbol\Theta}}
\def\Thebd{{{\boldsymbol\Theta}_\del}}
\def\Thebs{{{\boldsymbol\Theta}_s}}
\def\thet{{\vartheta}}
\def\thed{{\theta_\del}}
\def\thes{{\theta_s}}
\def\ub{{\bm u}}

\def\udhp{{u_{\del,h,\pb}}}
\def\uh{{u_h}}

\def\UUU{{U(\omg)}}
\def\vb{{\bm v}}

\def\vh{{v_h}}
\def\VVc{{V_c(\omgomgi)}}
\def\VVcd{{V_{c,\del}(\omgomgi)}}
\def\VVcdn{{V^n_{c,\del}(\omgomgi)}}
\def\VVch{{V_c^h(\omgomgi)}}
\def\VVd{{V_d(\omg)}}

\def\VVt{{V_t(\omgid)}}
\def\VVV{{V(\omgomgi)}}
\def\VVVh{{V^h(\omgomgi)}}
\def\wb{{\bm w}}
\def\www{{w}}
\def\WWW{{W}}
\def\xb{{\bm x}}

\def\xib{{\boldsymbol\xi}}
\def\xy{{\xb,\yb}}
\def\xyp{{(\xb,\yb)}}
\def\yb{{\bm y}}
\def\ymx{{\yb-\xb}}

\def\zb{{\bm z}}

\def\zzb{{\bm 0}}

\newtheorem{rem}{Remark}[section]
\def\gh{{g_h}}

\def\mub{{\boldsymbol\mu}}

\def\mcK{{\mathcal K}}
\def\mcP{{\mathcal P}}
\def\mcPd{{\mcP_\del}}

\def\VVpd{{V_{p,\del}}}
\def\mcAd{{\mcA_\del}}

\def\abs#1{|#1|} %%\def\abs#1{\left|#1\right|}
\def\mat#1{\textbf{\em #1}}
\def\pair#1#2{\langle #1,#2\rangle} %%\def\pair#1#2{\left\langle #1,#2\right\rangle}

\def\edoc{\end{document}}

\pagerange{\pageref{firstpage}--\pageref{lastpage}}
\allowdisplaybreaks
\doi{XXXXXXXX}  %GES% 8 figures, starting with 1
\usepackage[UKenglish]{babel} %GS%
	
\newcommand{\ignore}[1]{}

\newcommand{\citeb}[1]{\citename*{#1}\ \citeyear*{#1}}

\usepackage[twoside,
	paperheight=247mm,
	paperwidth=174mm,
	marginratio={3:5,2:3},
	left=18mm,
	top=20mm]{geometry}
\usepackage{subfigure}
\newcommand{\ro}[1]{{#1}} %%{\textcolor{blue}{#1}}
\newcommand{\roo}[1]{{#1}} %%{\textcolor{red}{#1}}
\newcommand{\fo}[1]{{#1}} %%{\textcolor{blue}{#1}}
\newcommand{\bo}[1]{{#1}} %%{\textcolor{blue}{#1}}
\newcommand{\boo}[1]{{#1}} %%{\textcolor{blue}{#1}}
\newcommand*{\ie}{{\it i.e.}}
\newcommand*{\tie}{{that~is}}

\newcommand*{\eg}{{\it e.g.}}
\newcommand*{\feeg}{{for example}}
\newcommand*{\fieg}{{for instance}}

\newcommand*{\scf}{{see}} 
\numberwithin{figure}{section}
\numberwithin{table}{section}
\newcommand{\actaqed}{\hfill $\actabox$}
 
\newcommand{\noi}{\noindent}

\newcommand{\gfrac}[2]{{{#1}/{#2}}}

\newcommand{\ffrac}{\displaystyle \frac}
\newcommand{\ds}{\displaystyle}
\newcommand{\fns}{\footnotesize}

\newcommand{\rme}{\mathrm{e}}
\newcommand{\rmi}{{\mathrm{i}}}
\newcommand{\dd}{\mathrm{d}}

\newcommand{\Bigoh}{{O}}
\newcommand{\qfa}{\quad\text{{for all} }}

\newcommand{\hspz}{\hspace{-6pt}}
\newcommand{\diam}{\operatorname{diam}}
\newcommand{\trip}[1]{\vert\!\vert\!\vert #1 \vert\!\vert\!\vert}
\newcommand{\zz}{\hspace{-1.4pt}}
\newcommand{\zzz}{\hspace{2.5pt}}
\newcommand{\vst}{\vspace{-3pt}}

\usepackage{draftwatermark}
\SetWatermarkText{\textsf{preprint}}
\SetWatermarkScale{6}

\begin{document}

\title[Numerical methods for nonlocal and fractional models]{Numerical methods for nonlocal\\ and fractional models} %%\thanks{{{Colour online for monochrome figures available at \textsf{journals.cambridge.org/anu}.}}}}

\author[M.\zzz D'Elia, Q.\zzz Du, C.\zzz Glusa, M.\zzz Gunzburger, X.\zzz Tian and Z.\zzz Zhou] 
       {Marta D'Elia{$^1$}, Qiang Du{$^2$}, Christian Glusa{$^1$}, \cr
         Max Gunzburger{$^3$}, Xiaochuan Tian{$^4$} and Zhi Zhou{$^5$}\\
{$^1$} Center for Computing Research,\rule{0pt}{24pt}\\
Sandia National Laboratories, \\
Albuquerque, New Mexico 87185, USA\\
E-mail: {\sf mdelia@sandia.gov}, {\sf caglusa@sandia.gov}\\
{\sf https://sites.google.com/site/martadeliawebsite}
\\
{$^2$} Department of Applied Physics and Applied Mathematics\rule{0pt}{24pt}\\
and Data Science Institute,\\ 
Columbia University, \\
New York, NY 10027, USA\\
E-mail: {\sf qd2125@columbia.edu}\\
{\sf http://www.columbia.edu/\midtilde{qd2125}}
\\
{$^3$} Department of Scientific Computing,\rule{0pt}{24pt}\\ 
Florida State University, \\
Tallahassee, Florida 32306, USA\\
E-mail: {\sf mgunzburger@fsu.edu}\\
{\sf https://www.sc.fsu.edu/gunzburg}
\\
{$^4$} Department of Mathematics,\rule{0pt}{24pt}\\
University of Texas, \\
Austin, Texas 78751, USA\\
E-mail: {\sf xtian@math.utexas.edu}\\
{\sf https://web.ma.utexas.edu/users/xtian}
\\
{$^5$} Department of Applied Mathematics,\rule{0pt}{24pt}\\
The Hong Kong Polytechnic University, \\
Kowloon, Hong Kong, China\\
E-mail: {\sf zhizhou@polyu.edu.hk}\\
{\sf https://sites.google.com/site/zhizhou0125}
}

\maketitle

\label{firstpage}

\newpage

\begin{abstract}
  Partial differential equations (PDEs) are {used with huge success to model phenomena across} all scientific and engineering disciplines. However, across an equally wide swath, there exist situations in which {PDEs} fail to adequately model observed phenomena, or are not the best available model for that purpose. On the other hand, in many situations, {\em nonlocal models} that account for interaction occurring at a distance have been shown to more faithfully and effectively model observed phenomena that involve possible singularities and other anomalies. In this article we consider a generic nonlocal model, beginning with a short review of its definition, the properties of its solution, its mathematical analysis and of specific concrete examples. We then provide extensive discussions about numerical methods, including finite element, finite difference and spectral methods, for determining approximate solutions of the nonlocal models considered. In that discussion, we pay particular attention to a special class of nonlocal models that are the most widely studied in the literature, namely those involving fractional derivatives. The article ends with brief considerations of several modelling and algorithmic extensions, which serve to show the wide applicability of nonlocal modelling.
\end{abstract}

%%\tableofcontents
\begin{center}\parbox{0.78\textwidth}{
  \section*{CONTENTS}
\contentsline {part}{PART 1: Nonlocal diffusion models, including\\ fractional models}{3}
\contentsline {section}{\numberline {1}General models for nonlocal diffusion}{4}
\contentsline {subsection}{\numberline {1.1}A brief review of a nonlocal vector calculus}{7}
\contentsline {section}{\numberline {2}Weak formulations of nonlocal models}{11}
\contentsline {subsection}{\numberline {2.1}Function spaces, weak formulations and \leavevmode {\color {red}well-posedness}}{11}
\contentsline {subsubsection}{\numberline {2.1.1}Function spaces and norms}{12}
\contentsline {subsubsection}{\numberline {2.1.2}Weak formulations and \leavevmode {\color {red}well-posedness}}{14}
\contentsline {subsubsection}{\numberline {2.1.3}Relation to an energy minimization principle}{15}
\contentsline {subsection}{\numberline {2.2}Kernel choices and the corresponding energy spaces}{16}
\contentsline {subsubsection}{\numberline {2.2.1}A list of kernel functions in common use}{17}
\contentsline {paragraph}{General integrable kernel functions.}{17}
\contentsline {paragraph}{`Critical' kernel functions.}{17}
\contentsline {paragraph}{`Peridynamic' kernel functions.}{18}
\contentsline {paragraph}{Fractional kernel functions.}{18}
\contentsline {section}{\numberline {3}Fractional diffusion models}{19}
\contentsline {subsection}{\numberline {3.1}Integral fractional Laplacian models for diffusion on bounded domains}{20}
\contentsline {subsubsection}{\numberline {3.1.1}Truncated interaction domains}{23}
\contentsline {subsection}{\numberline {3.2}Spectral fractional Laplacian models for diffusion on bounded domains}{24}
\contentsline {subsection}{\numberline {3.3}Additional considerations about fractional Laplacian models}{25}
\contentsline {subsubsection}{\numberline {3.3.1}Other integral representations of the fractional Laplacian}{25}
\contentsline {subsubsection}{\numberline {3.3.2}Extension representations}{25}
\contentsline {subsubsection}{\numberline {3.3.3}Regularity of solutions}{26}
\contentsline {subsubsection}{\numberline {3.3.4}Analytic solutions in one and two dimensions}{27}
\contentsline {subsection}{\numberline {3.4}Inhomogeneous constitutive functions in fractional models}{28}
\contentsline {part}{PART 2: Numerical methods for nonlocal and\\ fractional models}{29}
\contentsline {section}{\numberline {4}Introductory remarks}{29}
\contentsline {subsection}{\numberline {4.1}Asymptotically compatible schemes}{30}
\contentsline {section}{\numberline {5}Finite element methods for nonlocal models}{31}
\contentsline {subsection}{\numberline {5.1}{Asymptotically compatible conforming} finite element methods for {nonlocal diffusion}}{33}
\contentsline {subsection}{\numberline {5.2}Nonconforming and DG FEMs for nonlocal models with sufficiently singular kernels}{35}
\contentsline {paragraph}{Compactness results.}{36}
\contentsline {subsection}{\numberline {5.3}Adaptive mesh refinement for nonlocal models}{37}
\contentsline {subsection}{\numberline {5.4}Approximations of fractional models as limit of nonlocal models with a finite range of interactions}{39}
\contentsline {section}{\numberline {6}Finite element methods for the integral fractional\\ Laplacian}{41}
\contentsline {subsection}{\numberline {6.1}Quadrature rules}{42}
\contentsline {section}{\numberline {7}{Finite element methods for the spectral fractional\\ Laplacian}}{43}
\contentsline {subsection}{\numberline {7.1}Truncation in the extended direction}{44}
\contentsline {subsection}{\numberline {7.2}Spectral method in the extended direction}{45}
\contentsline {section}{\numberline {8}Spectral-Galerkin methods for nonlocal diffusion}{46}
\contentsline {section}{\numberline {9}Spectral-Galerkin methods for fractional diffusion}{51}
\contentsline {subsection}{\numberline {9.1}Spectral-Galerkin methods in unbounded domains}{51}
\contentsline {subsubsection}{\numberline {9.1.1}Approximation by Hermite polynomials}{52}
\contentsline {subsubsection}{\numberline {9.1.2}Approximation by modified Gegenbauer polynomials}{53}
\contentsline {subsubsection}{\numberline {9.1.3}Approximation by modified Chebyshev polynomials}{55}
\contentsline {subsection}{\numberline {9.2}Spectral-Galerkin methods in bounded domains}{57}
\contentsline {section}{\numberline {10}Finite difference methods for the strong form of\\ nonlocal diffusion}{60}
\contentsline {paragraph}{Quadratic exactness.}{61}
\contentsline {paragraph}{Uniform truncation error.}{61}
\contentsline {section}{\numberline {11}Numerical methods for the strong form of\\ fractional diffusion}{64}
\contentsline {subsection}{\numberline {11.1}Quadrature-rule based finite difference methods}{64}
\contentsline {subsection}{\numberline {11.2}Monte Carlo method by Feynman--Kac formula}{67}
\contentsline {subsection}{\numberline {11.3}Radial basis function methods}{69}
\contentsline {section}{\numberline {12}Conditioning and fast solvers}{71}
\contentsline {subsection}{\numberline {12.1}Fast algorithm for kernels with non-smooth truncation}{71}
\contentsline {subsection}{\numberline {12.2}Conditioning and solvers for finite element discretizations of the integral fractional Laplacian model}{73}

}\end{center}

\begin{center}\parbox{0.78\textwidth}{
\contentsline {part}{PART 3: Selected extensions}{74}
\contentsline {section}{\numberline {13}Weakly coercive, indefinite and non-self-adjoint\\ problems}{74}
\contentsline {subsection}{\numberline {13.1}Indefinite and non-self-adjoint problems}{74}
\contentsline {subsection}{\numberline {13.2}Asymptotically compatible schemes}{76}
\contentsline {subsection}{\numberline {13.3}Mixed finite element methods}{76}
\contentsline {section}{\numberline {14}Nonlocal convection--diffusion problems}{77}
\contentsline {subsection}{\numberline {14.1}Non-symmetric kernels and nonlocal convection--diffusion operators}{78}
\contentsline {subsubsection}{\numberline {14.1.1}Steady-state nonlocal convection--diffusion problems}{80}
\contentsline {section}{\numberline {15}Time-dependent nonlocal problems}{82}
\contentsline {subsection}{\numberline {15.1}Time-dependent nonlocal diffusion}{96}
\contentsline {subsection}{\numberline {15.2}Time-dependent convection--diffusion problems}{97}
\contentsline {subsection}{\numberline {15.3}Nonlocal wave equations}{97}
\contentsline {section}{\numberline {16}Inverse problems}{85}
\contentsline {subsection}{\numberline {16.1}Inverse problems for nonlocal diffusion}{83}
\contentsline {subsubsection}{\numberline {16.1.1}Distributed optimal control for a matching functional}{83}
\contentsline {subsubsection}{\numberline {16.1.2}Coefficient identification}{84}
\contentsline {subsection}{\numberline {16.2}Inverse problems for fractional operators}{86}
\contentsline {subsubsection}{\numberline {16.2.1}Distributed optimal control for a matching functional}{87}
\contentsline {section}{\numberline {17}Variational inequalities and obstacle problems}{93}
\contentsline {section}{\numberline {18}Reduced-order modelling}{96}
\contentsline {section}{\numberline {19}A turbulent flow application}{99}
\contentsline {paragraph}{Stage 1.}{102}
\contentsline {paragraph}{Stage 2.}{102}
\contentsline {section}{\numberline {20}Peridynamics models for solid mechanics}{101}
\contentsline {section}{\numberline {21}Image denoising}{105}
\contentsline {subsection}{\numberline {21.1}Image deblurring via minimization of a nonlocal functional}{108}
\contentsline {subsubsection}{\numberline {21.1.1}Optimization of the denoising parameters}{110}
\contentsline {refsection}{References}{110}

}\end{center}

\vspace{20pt} %\vspace{5pt}
\part[Nonlocal diffusion models, including\\ fractional models]{Nonlocal diffusion models, including fractional models}\label{part1}

{In order to make the article as self-contained as possible,} in this part we introduce the model equations {and discuss solutions of those equations used in the development and analysis of the numerical methods described in Part~\ref{part2}}. We focus on a general class of {\em nonlocal models} that are characterized by interactions at a distance {via} integral equation formulations as opposed to partial differential equations. For the same {applications,} the nonlocal models considered {provide different representations of physics}\footnote{Throughout, for the sake of economy of exposition, we refer to physics as the setting for the generic models we consider. Of course, those models also arise in all other physical science settings, as well as in engineering and the biological and social sciences.} {compared to partial} differential equation models. {\em Fractional derivative models} are an example of the general {class we} consider. Because there exists a {large mathematical} and computational literature devoted to fractional derivative models, that {class is} given special attention throughout the article. However, also highlighted are the opportunities afforded by the more general nonlocal models we consider that are not available through the use of fractional models.
 
{We use the following notational conventions throughout. Exceptions should} not cause confusion. Please note that the adherence to these conventions {becomes} less strict as one moves from Part~1 to Part~2 to Part~3.

\begin{center}
  \renewcommand{\tabcolsep}{8pt} 
  \begin{tabular*}{\textwidth}{ll}
    \hline\hline
    Roman letters & functions depending on a single point,
\\&\eg\ $u(\xb)$, $\vb(\xb)$, $\Dbf(\xb)$ 
\\[3pt]
Greek letters & {functions} depending on two points,
\\
&\eg\ $\eta\xyp$, $\nub\xyp$, $\Theb\xyp$
\\[3pt]
plain font & scalar-valued functions, \eg\ $u(\xb)$, $\eta\xyp$
\\[3pt]
boldface font & vector-valued functions, \eg\ $\ub(\xb)$, $\nub\xyp$
\\[3pt]
upper-case boldface font & tensor-valued functions, \eg\ $\Dbf(\xb)$, $\Theb\xyp$
\\[3pt]
calligraphy font & operators, functionals, bilinear forms, 
\\
&\eg\ $\mcA$, $\mcD$, $\mcL$
\\
\hline 
\end{tabular*}
\end{center}
In choosing these notational conventions, our goal is to be as consistent {as possible} or practical throughout the article. Unfortunately, because {of the} different notations {adopted in the literature}, the notations used in the article may differ from those used in some of the cited books and papers.

\vspace{5pt}
\section{General models for nonlocal diffusion}\label{sec:mod-strong}

We consider an integral equation model that is a nonlocal analogue of the classical Poisson problem 
\begin{equation}\label{str-pde}
\begin{cases}
   -\mcLz u :=
    -\nabla\cdotsp (\Dbf\nabla u) = f(\xb) & \text{for all }  \xb\in\omg,
\\
        \mcB u = g(\xb) & \text{for all } \xb\in\partial\omg,
\end{cases}
\end{equation}
where $\omg\in\Rd$ denotes a bounded, open domain having boundary $\partial\omg$, $f(\xb)$ and $g(\xb)$ denote given functions defined on $\omg$ and $\partial\omg$, respectively, and $\Dbf(\xb)$ denotes a given symmetric, positive definite $\ddd\times\ddd$ matrix. For the boundary conditions operator $\mcB$, we have the choices
\begin{equation}\label{str-bc}
\mcB u = \begin{cases}
u  & \Leftarrow\mbox{Dirichlet boundary condition},
\\
(\Dbf\nabla u)\cdot\nb & \Leftarrow\mbox{Neumann boundary condition},
\\
(\Dbf\nabla u)\cdot\nb + r(\xb) u& \Leftarrow\mbox{Robin boundary condition},
\end{cases}
\end{equation}
\bo{where $r(\xb)$ denotes a given function}, as well as mixed boundary conditions, \feeg\ Dirichlet and Neumann boundary {conditions} applied on \bo{disjoint, covering parts} of the boundary. {Problem} \eqref{str-pde} is a model {for steady-state diffusion,} \feeg.

As will be immediately clear, the nonlocal model we consider is not an integral or boundary integral reformulation {of problem} \eqref{str-pde} based on, \feeg, the use of Green's functions, but rather {it models different physics from}~\eqref{str-pde}.

\begin{figure} %%fig1.1
\centering
\subfigure[]{\includegraphics[width=150pt,viewport=0 0 322 230,clip]{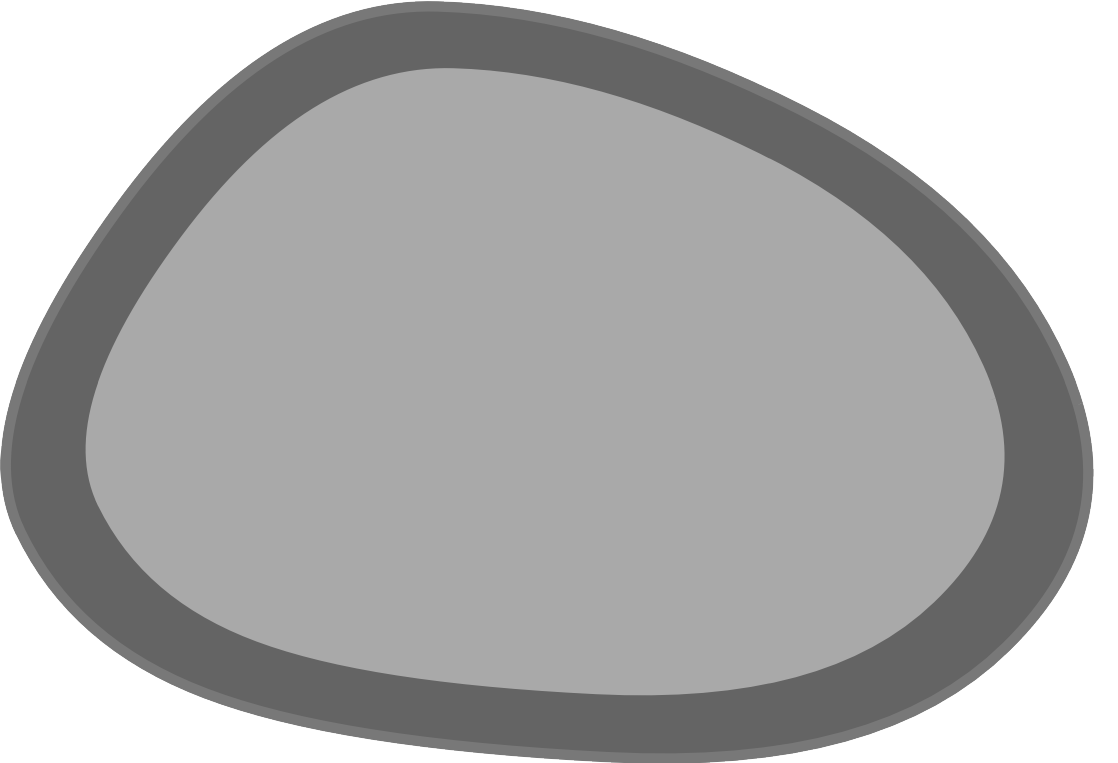}}\hspace{18pt} %%interaction
\subfigure[]{\includegraphics[width=150pt,viewport=342 0 664 230,clip]{Figures/fig1_1.pdf}}\\ %%interaction
\caption{$\omg$ {is light} grey whereas $\omgid$ {is dark} grey and has thickness $\del$. {(a)}~$\del$ smaller than the diameter of $\omg$, {(b)}~$\del$ larger than the diameter of $\omg$.}
\label{fig:interaction}
\end{figure}

\pagebreak %%20200331
Given the bounded, open domain $\omg\in\Rd$ and given a constant $\del>0$, we define the {\em interaction domain} corresponding to $\omg$ as
\begin{equation}\label{str-inter}
\omgid := \{
\text{$\yb\in\Rdc $ such that $\yb\in\ballx $ for some $ \xb\in\omg$} \},
\end{equation}
where $\ballx$ denotes the ball of radius $\del$ centred at $\xb$. We refer to $\del$ as the {\em interaction radius} or {\em horizon}. The nomenclature `interaction' used here is appropriate because $\omgid$ contains all points in the complement domain $\Rdc$ that are within a distance $\del$ of some point in $\omg$. Note that the four cases 
\begin{equation}\label{str-dcases}
0<\del\ll \text{diam\,}\omg, \quad \text{diam\,}\omg \ll \del < \infty,\quad
\del \approx \text{diam\,}\omg,\quad \del = \infty
\end{equation}
are all of interest. An illustration of the first two cases is given in \fo{Figure~\ref{fig:interaction}}. For the fourth case we have that the interaction domain $\omgii=\Rdc$. Also, the case of $\omg=\Rd$ is of considerable interest.

For $\del>0$, we consider the nonlocal problem for a scalar-valued function $u(\xb)$ defined on $\omgomgi$, given by 
\begin{equation}\label{str-vcp}
\begin{cases}
  - \mcLd u =  f(\xb) & \text{for all } \xb\in\omg,
\\
  \mcV u = g(\xb) &\text{for all } \xb\in\omgid.
\end{cases}
\end{equation}
Note that {\em the constraint $\mcV u = g$ is applied on the domain $\omgid$ having nonzero volume in $\Rd$}, in contrast to the constraint $\mcB u = g$ in \bo{\eqref{str-pde}}, {which} is applied on the boundary surface $\partial\omg$. For this reason we refer to $\mcV u = g$ as being a {\em volume constraint} and to \eqref{str-vcp} as a being a {\em volume-constrained problem} {\cite{du12sirev}.}

In \eqref{str-vcp}, we have that $f(\xb)$ and $g(\xb)$ denote given scalar-valued functions defined on $\omg$ and $\omgid$, respectively, and
\begin{equation}\label{str-lap}
\mcLd u {(\xb)} := 2 \intoi (u(\yb)-u(\xb)) \gamd\xyp \dyb
\quad\text{for all } \xb\in\omg,
\end{equation}
where $\gamd\xyp$ is a symmetric function, \tie,
\begin{equation}\label{str-gamsym}
\gamd\xyp = \gamd(\yb,\xb),
\end{equation}
and, for any $\xb$, 
\begin{equation}\label{str-gamz}
\mbox{supp}(\gamd\xyp) = \ballx,
\end{equation}
\tie, $\gamd\xyp=0$ whenever $|\ymx|>\del$. We note that in some of the specific cases we consider below, the integral in \eqref{str-lap} has to be viewed in the {principal value} sense. 
{For {notational} simplicity, we omit `p.v.' in front of the  integral sign even in such cases.}
The operator $\mcLd$ is a nonlocal analogue of the {partial differential equation (PDE)} operator $\nabla\cdotsp (\Dbf\nabla u)$. That connection is made explicit in Section~\ref{sec:mod-nlvc}. The nonlocal model \eqref{str-vcp} is a nonlocal analogue of the local PDE diffusion model \eqref{str-pde} model.

Analogous to \eqref{str-bc}, we have the choice of volume constraints 
\begin{equation}\label{str-vc}
\mcV u = \begin{cases}
u  & \Leftarrow\mbox{Dirichlet volume constraint},
\\
\mcNd u & \Leftarrow\mbox{Neumann volume constraint},
\\
\mcNd u + r(\xb) u& \Leftarrow\mbox{Robin volume constraint},
\end{cases}
\end{equation}
as well as mixed types of volume constraints, where \bo{the linear operator $\mcNd$ can be defined in several ways, for example,} 
\begin{equation}\label{str-flux}
\mcNd u := -2 \intoi (u(\yb)-u(\xb)) \gamd\xyp \dyb
\quad\text{for all } \xb\in\omgid.
\end{equation}
 It is obvious that the definitions of the operators $\mcLd$ and $\mcNd$ involve exactly the same integrand but have different domains of definition, \ie\ $\omg$ for the former and $\omgid$ for the latter. So, clearly, they can be combined into a single operator defined over $\omgomgi$. However,  it is convenient to keep using the two operators.

Because points $\xb$ only interact with points $\yb\in\ballx$, \eqref{str-lap} and \eqref{str-flux} can be equivalently written as
\begin{equation}\label{str-lap1}
   \mcLd u := 2 \intbx (u(\yb)-u(\xb)) \gamd\xyp \dyb  \quad\text{for all } \xb\in\omg
\end{equation}
and
\begin{equation}\label{str-flux1}
   \mcNd u := -2 \int_{(\omgomgi)\cap\ballx} (u(\yb)-u(\xb)) \gamd\xyp \dyb \quad\text{for all } \xb\in\omgid,
\end{equation}
respectively, where we have used the facts that due to the definition of the interaction domain $\omgid$, 
 $(\omgomgi)\cap\ballx=\ballx$ for $\xb\in\omg$ but $(\omgomgi)\cap\ballx\subset\ballx$ for $\xb\in\omgid$.

\subsection{A brief review of a nonlocal vector calculus} \label{sec:mod-nlvc}

To make some sense of how the model \eqref{str-vcp} arises and why we refer to it as a nonlocal analogue of the model \eqref{str-pde}, we need to first introduce some elements of a nonlocal vector calculus. Note that in this article we also use elements of the {\em fractional calculus.} However, because there are many excellent references about the fractional calculus and its applications (see \eg\ \citeb{baleanu}, \citeb{mainardi}, \citeb{Meerschaert2012book}), we do not discuss them here. 

The classical vector calculus provides a set of tools that are in ubiquitous use for the modelling, analysis, discretization and numerical analysis of {PDE} models, an obvious example being the use of Green's first identity to transform the strong formulation of a PDE into a weak formulation. The {foundations} of that calculus are the familiar divergence, gradient and curl differential operators, upon which an edifice is built that includes vector identities (\eg\ div curl $\vb$ = 0), integral theorems (\eg\ {Gauss's} theorem) and much, much more.

A nonlocal vector calculus has been developed \cite{dt-GuLe10,du12sirev,du13m3as,alalinlvc,MeDu16,du19cbms,deliahelmholtz} to deal with nonlocal models such as \eqref{str-vcp} in much the same way as the classical vector calculus is used to deal with PDE models such as \eqref{str-pde}. Here, mostly following \citeasnoun{du13m3as}, we provide a brief introduction {to} the nonlocal vector {calculus}, including notions that are used in the rest of the article. We remark that certain elements of the nonlocal vector calculus have previously appeared in
\citeasnoun{cagr:98}, \citeasnoun{zhsc:05}, \citeasnoun{Gilboa2008}, \citeasnoun{lete:10}, \citeasnoun{dt-LZOB10} {and} \citeasnoun{jlyy:11}, \feeg. However, the discussions in those papers, compared to that in \citeasnoun{du13m3as}, are limited in scope and in application and provide only a partial development of a nonlocal vector calculus that mimics the classical vector calculus.

The foundation of the nonlocal vector calculus {is} {\em integral} operators that mimic the three differential operators upon which the classical vector calculus is built. 

Given the vector-valued functions $\nubd\xyp \colon \omgomgip\to\Rd$ and $\alpbd\xyp\colon\omgomgip\to \Rd$, the action of the {\em nonlocal divergence operator} $\mcDd\colon \omgomgi \to \Ro$ on $\nub\xyp$ is defined as
\begin{equation}\label{str-div}
(\mcDd\nub)(\xb) 
:= \int_{\omgomgi} (\nub(\yb,\xb)+\nub\xyp)\cdot\alpbd\xyp \mcXd\dyb 
\end{equation}
{for all $\xb\in\omgomgi$,}
where $\alpbd\xyp$ denotes an antisymmetric function, \tie, $\alpbd(\yb,\xb)=-\alpbd\xyp$ for all $\xb,\yb\in\omgomgi$, and $\mcX_{(\cdotsp)}(\cdotsp)$ denotes the indicator function. Note that 
unlike its differential counterpart $\nabla\cdot$, 
the nonlocal operator $\mcDd$ is not uniquely defined, \tie, not only do we have an unspecified parameter {$\del>0$ but also we} have said nothing about $\alpbd\xyp$ other than it is an antisymmetric function. The choices one makes for $\del$ and especially $\alpbd\xyp$ can result in operators having very different properties. Thus $\del$ and $\alpbd\xyp$ are modelling choices dictated by the specific application one considers. We have already mentioned, in \eqref{str-dcases}, the wide choices of $\del$ that are of interest; several choices for $\alpbd\xyp$ are considered in Section~\ref{sec:klist}.

Simple manipulations show that, {under suitable regularity assumptions that we do not dwell on here,}
\begin{equation}\label{str-adjdef}
      \int_{\omgomgi} v \mcDd \nub  \D \xb = \int_{\omgomgi}\int_{\omgomgi} \nub \cdotsp \mcDsd v \dydx,
\end{equation}
where
\begin{equation}\label{str-divs}
(\mcDsd u)\xyp := - (u(\yb) -u(\xb))  \alpbd\xyp \mcXd  
\end{equation}
{for all $\xy\in\omgomgi$.}
From \eqref{str-adjdef}, it is natural to refer to operator $\mcDsd$ as being  the {\em adjoint operator} corresponding to $\mcDd$. With $\mcDsd$ being the adjoint of $\mcDd$, {one may formally} refer to $-\mcDsd$ as being a {\em nonlocal gradient operator.}

Let $\Thebd\xyp\colon (\omgomgi)\times$ $(\omgomgi) \to \Ro^{\ddd\times\ddd}$ denote a tensor-valued function that is symmetric in the function sense, \ie\ \bo{$\Thebd\xyp=\Thebd (\yb, \xb)$}, {and symmetric} and \bo{non-negative} definite in the matrix sense. Then, from the definitions of the nonlocal divergence and gradient operators, we have that  
\begin{align}\label{str-com}
&-\mcDd (\Thebd\mcDsd u) 
\\&\quad = 2\int_{\omgomgi} (u(\yb)-u(\xb)) \alpbd\xyp\cdotsp ( \Thebd\xyp\alpbd\xyp )\mcXd \dyb. \notag
\end{align}

For all $\xy\in\omgomgi$, we define the {\em kernel} $\gamd\xyp$ as
\begin{equation}\label{str-gamd}
\gamd\xyp := \alpbd\xyp\cdotsp (\Thebd\xyp\alpbd\xyp)\mcXd. 
\end{equation}
Note that $\gamd\xyp$ defined in this way \roo{satisfies the} symmetry condition \eqref{str-gamsym} and the support condition \eqref{str-gamz}, where the former follows from the fact that $\mcX_{\ballx}(\yb)$ is itself a symmetric function, \tie, if $\yb\in\ballx$ then necessarily $\xb\in B_\del(\yb)$.

We then have that \eqref{str-com} can be expressed as
\begin{align}\label{str-com1}
\mcDd (\Thebd\mcDsd u)
&=
-2\int_{\omgomgi} (u(\yb)-u(\xb)) \gamd\xyp  \dyb \notag \\*
& = -2\int_{\ballx} (u(\yb)-u(\xb)) \gamd\xyp \dyb,
\end{align}
where $\gamd\xyp$ is given by \eqref{str-gamd} and $\mcDd$ and $\mcDsd$ are given by \eqref{str-div} and \eqref{str-divs}, respectively. 

Of course, we recognize that 
\begin{equation}\label{str-com2}
\mcDd (\Thebd\mcDsd u)=-\mcLd u,
\end{equation}
where $\mcLd u$ is defined in \eqref{str-lap} or \eqref{str-lap1}. The fact that the operator $\mcLd$ can be written as a composition of nonlocal divergence and gradient operators justifies referring to \eqref{str-vcp} as a \bo{(variable coefficient)} nonlocal Poisson problem, \tie, $\mcLd u$ 
{may be viewed as}
a nonlocal analogue of $\nabla \cdotsp ( \Dbf \nabla u)$ {for suitable choices of $\alpbd\xyp$}, and, if $\Thebd$ is the identity tensor, $\mcLd u$ is indeed a nonlocal analogue of the~$\Delta u$.

\begin{rem}[variable coefficients in nonlocal model]\label{rem:var-coef}\hspace{-8pt}% 
An advantage accruing from the nonlocal calculus is that, through compositions such as \eqref{str-com} or \eqref{str-com1}, it provides a natural means for generalizing nonlocal constant coefficient operators such as the fractional Laplace operator {to variable coefficient settings}
{in a manner completely analogous} to how one defines variable coefficient operators in the PDE case. As an example, {in Section~\ref{sec:flvarcoef}, we describe} one way {in which} this can be done for fractional diffusion models.
\hfill\actaqed
\end{rem}
 
The final bit of the nonlocal vector calculus that we need  is a (generalized) nonlocal Green's first identity. In the context of this section, that identity is given by \cite{du13m3as}
\begin{equation}\label{str-gfi}
\int_{\omgomgi} v\mcDd (\Thebd\mcDs u) \dyb =
\int_{\omgomgi}\int_{\omgomgi} \mcDsd v \cdotsp (\Thebd\mcDsd)u \dydx.
\end{equation}
We {have}
\begin{equation}\label{str-poi2}
 \mcDsd v \cdotsp (\Thebd\mcDsd)u = (v(\yb)-v(\xb))(u(\yb)-u(\xb)) \gamd\xyp,
\end{equation}
{so  substituting} \eqref{str-com1}, \eqref{str-flux} and \eqref{str-poi2} into \eqref{str-gfi} results in the equivalent form of the nonlocal Green's first identity given by
\begin{align}\label{str-gfi1}
& \int_{\omgomgi}\int_{\omgomgi} (v(\yb)-v(\xb))(u(\yb)-u(\xb)) \gamd\xyp
\dydx \notag  \\*
&\quad  =-\int_{\omg} v(\xb) \mcLd u (\xb)\dxb + \int_{\omgid} v(\xb) \mcNd u (\xb)\dxb.
\end{align}
 \bo{We} recognize \eqref{str-gfi1} as a nonlocal analogue of the (generalized) local Green's first identity 
\[
\int_\omg \nabla v\cdotsp(\Dbf\nabla u) \dxb = -\int_\omg   v\nabla\cdotsp(\Dbf\nabla u)  \dxb + \int_{\partial\omg} v \nb \cdotsp(\Dbf\nabla u)  \dxb.
\]
Derivations in a more rigorous functional analytic setting of nonlocal integral identities  \bo{similar to} \eqref{str-gfi1} and the connection to their local analogues can be found in \citeasnoun{MeDu16} {and} \citeasnoun{du19cbms}.

\setcounter{rem}{2}
\medskip\noi{\bf Remark 1.2 ({roles} played in the kernel by operator definitions and constitutive functions).}
It is important to differentiate between the roles {of the functions $\alpbd\xyp$ and $\Thebd\xyp$ in} the nonlocal models we are considering. It is clear that $\alpbd\xyp$ {\em serves to define operators} (\ie\ $\mcDd$, $\mcDsd$ and $\mcNd$ given by \eqref{str-div}, \eqref{str-divs} and \eqref{str-flux}, respectively), {\em irrespective of how those operators are used.} On the other hand, $\Thebd\xyp$ {\em serves as a constitutive function.} Thus both $\alpbd\xyp$ and $\Thebd\xyp$ are modelling choices. Both influence the properties of solutions of nonlocal models such as their regularity.
Of course, the situation is much the same in the local PDE case, \tie, the operators $\nabla\cdot$ and $\nabla$ are defined irrespective of how they are used and $\Dbf(\xb)$ denotes a constitutive tensor. \hfill \actaqed
\medskip

\setcounter{rem}{3}
\noi{\bf Remark 1.3 (the choices that define  a nonlocal diffusion model).}
%%\begin{rem}[\bo{the choices that define} a nonlocal diffusion model]\label{rem:modchoise} 
  Recapitulating, 
  to define a specific nonlocal diffusion model, one must make three modelling choices:
  \renewcommand{\leftmargini}{20pt}
  \begin{itemize}\setlength\itemsep{4pt}
\item the horizon $\del$ that defines the extent of nonlocal interactions,

\item the antisymmetric function $\alpbd\xyp$ that defines the nonlocal divergence operator $\mcDd$ and nonlocal gradient operator $-\mcDsd$,

\item the constitutive tensor $\Thebd\xyp$ that defines the properties of the media.
  \end{itemize}
  These three choices are all that enter into the definition \eqref{str-gamd} of the kernel $\gamd\xyp$, {so specifying} them uniquely {defines} the operator $\mcLd$. \hfill\actaqed 
\medskip %%\end{rem}

\setcounter{rem}{4}
\noi{\bf Remark 1.4 (additional operators of the nonlocal vector calculus).}
The composition of the nonlocal divergence operator $\mcDd$ and its adjoint operator $\mcDsd$ are the only nonlocal operators needed to define the nonlocal operator $\mcLd$ that operates on scalar-valued functions. However, 
other aspects of the nonlocal vector calculus, such as nonlocal vector identities and nonlocal operators acting on vector-valued functions,  make use of additional nonlocal operators. 

Thus, in addition to $\mcDd$ and $\mcDsd$, we have the nonlocal gradient operator $\mcG_\del$ defined by, for a scalar-valued function $\nu\xyp$,
\begin{equation}\label{str-grad}
 (\mcG_\del\nu)(\xb) := \int_{\omgomgi} (\nu(\yb,\xb)+\nu\xyp)\alpbd\xyp \mcXd \dyb
\end{equation}
for all $\xb\in\omgomgi$ and its adjoint operator $\mcG^\ast_\del$ given by, for a vector-valued function $\ub(\xb)$,
\begin{equation}\label{str-grads}
(\mcG^\ast_\del \ub)\xyp = - (\ub(\yb)-\ub(\xb)) \cdotsp \alpbd\xyp\mcXd  
\end{equation}
for all $\xy\in\omgomgi$. Similarly, we have the nonlocal curl operator $\mcC_\del$ defined by, for a vector-valued function $\nub\xyp$,
\begin{equation}\label{str-curl}
 (\mcC_\del\nub)(\xb) := \int_{\Rd} (\nub(\yb,\xb)+\nub\xyp)\times\alpbd\xyp\mcXd \dyb 
\end{equation}
for all $\xb\in\omgomgi$ and its adjoint operator $\mcC^\ast_\del$ defined by, for a vector-valued function $\ub(\xb)$,
\begin{equation}\label{str-curls}
\bo{(\mcC^\ast_\del\ub)\xyp = - \alpbd\xyp\times({\ub}(\yb)-{\ub}(\xb))\mcXd  }
\end{equation}
for all $\xy\in\omgomgi$.
Thus, in the nonlocal vector calculus we have pairs of divergence operators $\mcDd$ and $-\mcG^\ast_\del$, gradient operators $\mcG_\del$ and $-\mcD^\ast_\del$, and curl operators $\mcC_\del$ and $\mcC^\ast_\del$, with $\mcDd$, $\mcG_\del$, $\mcC_\del$ acting on functions of two points, \ie\ $\nub\xyp$ and $\nu\xyp$, and $\mcDsd$, $\mcG^\ast_\del$, $\mcC^\ast_\del$ acting on functions of one point, \ie\ $\ub(\xb)$ and $u(\xb)$.

With these operators in hand, we have the nonlocal vector identities
\begin{equation}\label{str-vlap}
\begin{cases}
  \mcDd(\mcC^\ast_\del \ub) = 0,
  \qquad
  \mcC_\del(\mcG^\ast_\del u) = 0,
\\
\mcDd(\mcDsd \ub)  =  \mcG_\del(\mcG^\ast_\del \ub) + \mcC_\del(\mcC^\ast_\del \ub),
\end{cases}
\end{equation}
{which} are nonlocal analogues of $\nabla\cdotsp(\nabla\times\ub)=0$, $\nabla\times(\nabla u)=0$ and $\Delta \ub = \nabla \cdotsp (\nabla\ub) = \nabla( \nabla\cdotsp \ub) - \nabla\times (\nabla\times\ub)$, respectively. In particular, we can view $-\mcDd(\mcDsd \ub)  =  -\mcG_\del(\mcG^\ast_\del \ub) - \mcC_\del (\mcC^\ast_\del \ub)$ as a nonlocal vector Laplacian.

In addition, the operators introduced in this remark are used in \citeasnoun{deliahelmholtz} to define the nonlocal Hodge--Helmholtz decomposition of vector-valued {functions} that depend on two points $\xb$ and $\yb$.
Note that the operators $\mcDd$, $\mcDsd$, $\mcG_\del$ and  $\mcG^\ast_\del$, as operators between functions of two points $\xb$ and $\yb$ and functions of one point $\xb$ or $\yb$, do not recover their classical counterparts in the local limit $\del\to0$. 
Thus, in {parallel,}  additional examples of nonlocal operators acting on vector-valued functions of {one point only} are introduced in \citeasnoun{du13m3as} and further analysed in \citeasnoun{MeDu16}, \citeasnoun{du19cbms}, \citeasnoun{dtty19ima} {and} \citeasnoun{ld19ng}. 
In fact, a nonlocal Helmholtz decomposition of a vector-valued function  $\ub(\xb)$ can be established which recovers the classical one in the local limit; {see} \citeasnoun{ld19ng}. {These operators can also be used to construct a nonlocal analogue of the Stokes equation, which was used to provide a foundation for the analysis of smoothed particle hydrodynamics} \cite{DuTi19fcm,ld19ng}. Some of these developments can be found in Section~\ref{sec:mod-peri}.

\vspace{5pt}
\section{Weak formulations of nonlocal models}\label{sec:mod-weak}

The model \eqref{str-vcp} is a strong formulation of a nonlocal problem that is required to hold pointwise on $\omgomgi$. 
\bo{For nonlocal problems subject
to various inhomogeneous data, different weak formulations can be considered.}
In this section we consider 
\bo{one of the possible weak formulations} 
corresponding to \eqref{str-vcp} and then relate that formulation to a minimization principle. We conclude the section by providing specific examples of the kernels $\gamd\xyp$ and their associated energy spaces.

\subsection{Function spaces, weak formulations and {well-posedness}}\label{weak-wposed}

\subsubsection{Function spaces and norms}\label{sec:spanorm}

We define  \bo{an} {\em `energy' space} 
\begin{equation}\label{weak-espace}
   \VVV =  \{  v(\xb)\in L^2(\omgomgi) \colon \| v \|_{\VVV} <\infty\},
\end{equation}
where $\| \cdot \|_{\VVV} $ denotes the {\em `energy' norm} defined as
\begin{align}\label{weak-enorm}
  & \| v \|_{\VVV} \\*
  &\quad  = \biggl( 
\int_{\omgomgi}\int_{\omgomgi}  (v(\yb) -v(\xb)) ^2  \gamd\xyp  \dydx 
 + \|v\|^2_{L^2(\omgomgi) }\biggr)^{1/2}. \notag
\end{align}
The space $\VVV$ is by definition a subspace of $L^2(\omgomgi)$
and 
can be shown, for suitably chosen kernels, to be a Hilbert space that comes equipped with the inner product
\begin{align}\label{weak-ip}
  &(u,v)_{\VVV} \notag \\*
  &\quad =  
  \int_{\omgomgi}\int_{\omgomgi}  (v(\yb)-v(\xb))  (u(\yb)-u(\xb))   \gamd\xyp  \dydx \notag \\*
  &\quad\quad\,  + (u,v)_{L^2(\omgomgi) }\quad\text{for all } u,v\in\VVV.
\end{align}

We also have the {\em constrained energy space}
\begin{equation}\label{weak-espacec}
   \VVc =  \{  v(\xb)\in \VVV \colon  v(\xb) = 0 \ \mbox{for $\xb\in\omgid$} \}.
\end{equation}
Under suitable conditions on the kernels, we can prove a nonlocal Poincar\'e inequality \bo{\cite{du12sirev,du13dcdsb,du19cbms}}. Thus, we have that the seminorm
\begin{equation}\label{weak-enormc}
| v |_{\VVc} = \biggl( 
\int_{\omgomgi}\int_{\omgomgi}  (v(\yb)-v(\xb)) ^2  \gamd\xyp
 \dydx \biggr)^{1/2}
\end{equation}
is a norm on $\VVc$ equivalent \roo{to the} norm \eqref{weak-enorm} with equivalence constants that are independent of $\del$ \bo{as $\delta\to 0$}. Correspondingly, $\VVc$ is a Hilbert space equipped with the inner product
\begin{align}\label{weak-ip1}
  (u,v)_{\VVc} 
  =  
\int_{\omgomgi}\int_{\omgomgi}  (v(\yb)-v(\xb))  (u(\yb)-u(\xb))   \gamd\xyp
 \dydx
\end{align}
for all $u,v\in\VVc$.
For more detailed discussions of the conditions on the kernels and rigorous proofs of the properties on these spaces, we refer to
  {\citename{MeDu16} \citeyear{md15non,MeDu16} and \citeasnoun{du19cbms}}.

\bo{To allow inhomogeneous data, we define the}
space $\VVt$,
\begin{equation}\label{weak-trace}
\bo{\VVt =  \{ v = w|_{\omgid} \; \text{for some}\; w\in\VVV \},}
\end{equation} 
\bo{which involves restrictions to the domain $\omgid$ having finite volume in $\Rd$.} 
\bo{A norm on $\VVt$ can be  defined by}
\begin{align}\label{weak-normt}
\roo{\| v \|_{\VVt} 
  = \inf\{ \|w\|_{\VVV}\mid w\in\VVV,  w|_{\omgid}=v\}.}
\end{align}
\bo{One can view $\VVt$, induced by $\VVV$, as a nonlocal analogue of the trace space $H^{1/2}(\partial\omg)$ induced by $H^1(\omg)$.}

As mentioned in Remark~1.2,
specific function spaces are defined by making specific choices for $\alpbd\xyp$. For some choices, the function spaces so defined can be shown to be equivalent to standard function spaces such as Sobolev spaces. See Section~\ref{sec:klist} for some examples.

\subsubsection{Weak formulations and {well-posedness}}

A weak formulation {of problem} \eqref{str-vcp} with $\mcV u=u$ is defined as follows{:} 
\begin{equation}\label{weak-weakf}
  \begin{minipage}[b]{23pc}
{Given $\gamd\xyp$ defined in \eqref{str-gamd} and given $g(\xb)\in\VVt$
and $f(\xb)\in\VVd$,
find $u(\xb)\in\VVV$ such that}
\[
\begin{cases}
{\mcAd}(u,v) = \fv  &\text{for all } v\in\VVc ,
\\
u(\xb) = g(\xb) &\text{for all } \xb\in\omgid,
\end{cases}
\]
\end{minipage}
\end{equation}
where, for all $u(\xb),v(\xb)\in\VVV$, we have the bilinear form
\begin{align}\label{weak-bform}
  &{\mcAd}(u,v) 
  =\int_{\omgomgi}\int_{\omgomgi}  (v(\yb)-v(\xb))  (u(\yb)-u(\xb))  \gamd\xyp
 \dydx
\end{align}
and, for all $f(\xb)\in\VVd$, 
 \bo{the {\em dual space} of $\VVc$ defined by extending the $L^2$ inner product as the duality pairing and with the induced norm $\|\cdot\|_{\VVd}$,}
we have the linear functional
\begin{equation}\label{weak-lfunc}
\fv=\int_{\omg} f(\xb)v(\xb) \dxb.   
\end{equation}
That \eqref{weak-weakf}--\eqref{weak-lfunc} are indeed a weak formulation corresponding to \eqref{str-vcp} is an immediate consequence of setting $v(\xb)=0$ on $\omgid$ in the nonlocal Green's first identity \eqref{str-gfi1} so that the term involving $\mcNd u$ vanishes. Equation \eqref{weak-weakf} is a nonlocal analogue of the weak formulation {of problem} \eqref{str-pde} with $\mcB u=u$ 
given by 
$$\bo{\int_\omg \nabla v\cdotsp(\Dbf \nabla u)\dxb = \fv,}$$
for all  $v\in H^1_0(\omg)$ along with $u=g$ on~$\partial\omg$.

\bo{Under suitable conditions on the kernel $\gamd\xyp$, it} can be shown
that the bilinear form ${\mcAd}(u,v)$ is symmetric and is continuous and coercive with respect to $\VVc$ so that, by the Lax--Milgram theorem, {problem} \eqref{weak-weakf} is {well-posed}.
In particular, we have the {\it a~priori} estimate with respect to the energy norm for the solution $u(\xb)$  of \eqref{weak-weakf} given by
\begin{equation}\label{weak-est}
 \| u \|_\VVV \le \CCd  ( \|f\|_\VVd + \| g\|_\VVt) ,
\end{equation}
where the constant $\CCd$ \bo{does not} depend on $\del$.
\bo{Moreover, for suitably chosen kernels, one can consider its local ($\delta\to 0$) and fractional ($\delta\to \infty$) limits, as discussed in Part \ref{part2}.}
\medskip

\setcounter{rem}{1}
\noi{\bf Remark 2.1 (a nonlocal Neumann model).}
%%\begin{rem}[\bo{a nonlocal Neumann model}]\label{rem:neumann}\hspz%  
We {can of course also} define a weak formulation of \eqref{str-vcp} for $\mcV u = \mcN_{{\del}} u$, \tie, for a nonlocal Neumann problem. In this case a weak formulation is given by: \bo{find $u\in \VVV\setminus \Ro $ such that }
\begin{align}\label{weak-weakfn}
\begin{cases}
\ds {\mcAd}(u,v) = 
\fv + \int_{\omgid} g(\xb)v(\xb)\dxb \qfa  v\in\VVV\setminus\Ro,
\\
\ds 
\int_{\omg} f(\xb) \dxb + \int_{\omgid} g(\xb)\dxb =0, 
\end{cases}
\end{align}
\bo{where $f\in V_d(\Omega)$ and for the Neumann problem $g(\xb)$ belongs to the dual space of $\VVV|_{\Omega_{{\mathcal I}_\del}}$. In \eqref{weak-weakfn}, the quotient space $\VVV\setminus\Ro$ is used to ensure uniqueness of solutions.  The second equation in \eqref{weak-weakfn} is a compatibility condition on the data $f$ and $g$ that is necessary to ensure existence.}

Equation \eqref{weak-weakfn} is a nonlocal analogue of the weak formulation corresponding to 
\eqref{str-pde} with $\mcB u=(\Dbf \nabla u)\cdot\nb$
given by
\[
\int_\omg \nabla v\cdotsp(\Dbf \nabla u)\dxb = \int_\omg fv\dxb + \int_{\partial\omg}gv\dxb
\]
for all $v\in H^1(\omg)\setminus\Ro$ \bo{and a pair of
prescribed data, also denoted by $f$ and $g$, that satisfy}
 the compatibility condition 
\[
{\int_\omg f\dxb + \int_{\partial\omg}g\dxb=0.}
\] 
We note that other ways of defining a nonlocal Neumann problem 
\bo{can also be used.}
For additional
discussions about nonlocal Neumann problems, see \eg\  \citeasnoun{cortazar}, \citeasnoun{du19cbms},
{\citename{MeDu16} \citeyear{md15non,MeDu16}},
\citeasnoun{ttd17amc} {and} \citeasnoun{delianeumann}. 
\bo{One  can also find a fractional version of \eqref{str-vcp} with $\mcV u = \mcN_{{\del}} u$ in \roo{\citeasnoun{Dengw:2018}}
for a nonlocal Neumann problem of the fractional
Laplacian in the bounded domain. In addition, for specialized power-like kernels, one can also relate this to problems associated with problems related to the regional Laplacian defined by \eqref{intflr}.}
\hfill\actaqed
%%\end{rem}
\medskip

For the sake of economy of exposition,  {we will} focus on \eqref{weak-weakf} in which $\mcV u = u$ and, at times, we {will} even focus on the homogeneous version of \eqref{weak-weakf} for which $g(\xb)=0$.

\subsubsection{Relation to an energy minimization principle}\label{sec:weak-energy}

The weak formulation \eqref{weak-weakf} can also be derived as the Euler--Lagrange equation for an energy minimization principle. Consider the energy functional
\begin{align}\label{weak-energy}
\mcE_\del(u;f,\gamd)& =
\ffrac12\int_{\omgomgi}\int_{\omgomgi}  (u(\yb)-u(\xb)) ^2 \gamd\xyp 
\dydx \notag \\*
&\quad\, -
\int_{\omg} f(\xb)u(\xb) \dxb 
\end{align}
{and the following minimization principle:} 
\begin{equation}\label{weak-princ}
  \begin{minipage}[b]{23pc}
{Given $\gamd\xyp$ defined in \eqref{str-gamd} and given $g(\xb)\in\VVt$
and $f(\xb)\in\VVd$,
find $u(\xb)\in\VVV$ such that}
\[
\mcE_\del(u;f,\gamd)\ \text{is minimized}
\]
subject to $u(\xb) = g(\xb)$ for $\xb\in\omgid$.
\end{minipage}
\end{equation}
Using standard techniques from the calculus of variations, one easily sees that \eqref{weak-weakf} is  indeed the Euler--Lagrange equation corresponding to the energy minimization principle \eqref{weak-princ} \bo{under suitable conditions on the kernel $\gamd\xyp$ and data $f$ and $g$}.

\subsection{Kernel choices and the corresponding energy spaces}\label{sec:mod-kernels} 

We assume that $\gamd\xyp$ given by \eqref{str-gamd} can be written in the form
\begin{equation}\label{gamd1}
\gamd\xyp= \phid\xyp \thed\xyp \mcX_{\ballx}(\yb) ,
\end{equation}
where $\thed\xyp$ and $\phid\xyp$ denote \bo{non-negative,} symmetric, scalar-valued functions. This form is very general with respect to both operators and constitutive functions. In fact, by setting
\begin{equation}\label{ealpbd}
\alpbd\xyp = \ffrac{\ymx}{|\ymx|}  \sqrt{\phid\xyp}
\end{equation}
and
\begin{equation}\label{ethed}
\thed\xyp :=\ffrac{\ymx}{|\ymx|} \cdotsp  \biggl(
     \Thebd\xyp    \ffrac{\ymx}{|\ymx|}\biggr)
\end{equation}
in \eqref{str-gamd}, we arrive at \eqref{gamd1}. Thus the only assumption made is that $\alpbd\xyp$ is directed along the vector $\ymx$. Note that $\alpbd\xyp$ given by \eqref{ealpbd} is indeed an antisymmetric function, whereas $\thed\xyp$ given by \eqref{ethed} is indeed a 
\bo{non-negative,} symmetric function if $\Thebd\xyp$ is a symmetric in the function sense and symmetric and  \bo{non-negative} definite in the matrix sense.

We do not combine $\thed\xyp$ and $\phid\xyp$ into a single function because we want to keep separate the operator definition and constitutive roles played by $\gamd\xyp$; note that $\phid$ is only related to $\alpbd$ whereas $\thed$ is only related to $\Thebd$. {Thus $\phid\xyp$ now} takes over the role of $\alpbd\xyp$ in this regard, and $\thed\xyp$ can be viewed as a constitutive function. Note also that $\alpbd\xyp$ as given by \eqref{ealpbd} can be used to define the `first-order' nonlocal operators $\mcDd$ and $\mcDsd$, whereas $\gamd\xyp$ is used to define the `second-order' nonlocal operator $\mcLd$ and the `Neumann' nonlocal operator $\mcNd$. {We will} refer to $\gamd\xyp$ as the {\em kernel}, $\phid\xyp$ as the {\em kernel function} and $\thed\xyp$ as a {\em constitutive function}.

\begin{rem}[scalar constitutive tensors]\hspz%  
The case of $\Thebd\xyp$ being a scalar tensor (\ie\ if $\Thebd\xyp = \thed\xyp \Ibf$, where $\Ibf$ denotes the identity tensor and $\thed\xyp$ denotes a scalar-valued function) is by far the most common case considered, just as it is for a scalar diffusivity tensor $\Dbf = a(\xb)\Ibf$ in the local PDE case. In this case, \eqref{ethed} becomes a tautology. Also, in this case, \eqref{str-gamd} simplifies to 
\begin{align*}
\gamd\xyp&= \thed\xyp \alpbd\xyp\cdot\alpbd\xyp\mcX_{\ballx}(\yb) \\*
&= \thed\xyp |\alpbd\xyp|^2\mcX_{\ballx}(\yb). 
\end{align*}
\end{rem}

\setcounter{rem}{3}
\noi{\bf Remark 2.3 (homogeneous and inhomogeneous nonlocal operators).}
If $\Thebd\xyp=\Thebd(|\ymx|)$ and $\phid\xyp=\phid(|\ymx|)$, \tie, $\Thebd$ and $\phid$ are {\em radial functions}, then \eqref{str-vcp} is a model for a homogeneous medium. This is not surprising because $\Thebd$ and $\phid$ being  radial functions means that the nature of the interaction between a point $\xb$ and another point $\yb$ is independent of the location of the point $\xb$; that interaction only depends on the distance $|\ymx|$ between $\xb$ and $\yb$. In this case, \eqref{str-vcp} is a nonlocal analogue of the PDE $-\nabla\cdot(\Dbf\nabla u)=f$ in which $\Dbf$ is a constant tensor. If in addition $\Thebd$ is a scalar tensor, \eqref{str-vcp} is a nonlocal analogue of the PDE $-\kappa\Delta u=f$ in which $\kappa$ is a constant. Thus, to obtain an inhomogeneous
model, both $\Thebd\xyp$ and $\phid\xyp$ cannot be radial functions. Note that $\phid\xyp$ is indeed often chosen to be a radial function, in which case $\Thebd\xyp$ has to be non-radial for $\mcLd$ to be {an} inhomogeneous operator. {It is important to note that a constant coefficient nonlocal problem does not mean that the kernel and constitutive functions are themselves constant functions; they merely need to be radial functions.}

\subsubsection{A list of kernel functions in common use}\label{sec:klist}

Because $\thed\xyp$ is a constitutive function which is not specific even within a single application, we focus on choices for the function $\phid\xyp$ that determines generic properties of a model. However, we again note that because we have assumed that $\Thebd\xyp$ is a symmetric,  \bo{non-negative} definite tensor, it follows that $\thed\xyp$ is a  \bo{non-negative} function for all $\xb,\yb\in\omgomgi$.

We also assume, as is usually the case, that for all $\xb,\yb\in\omgomgi$, the constitutive function $\thed\xyp$ is bounded from above and below by \bo{positive} constants whose values do not depend on $\del$. In fact, often $\thed\xyp$ does not depend on $\del$. Thus we focus on a list of kernel functions $\phid\xyp$. In \fo{Table~\ref{ker-prop}} we provide the energy spaces corresponding to the kernel functions we list.

\begin{table}
\caption{Kernel functions $\phid\xyp$ and the corresponding energy\newline spaces $\VVV$.}\label{ker-prop}
\begin{tabular}{lcl}
\hline\hline
Type & Definition & Energy space {$\VVV$}
\\
\hline
\bo{Translation-invariant,} integrable & \eqref{intker} & $L^2(\omgomgi)$
\\
Critical & \eqref{dddker} & $\subsetneq L^2(\omgomgi)$
\\
Peridynamic & \eqref{perker} & $L^2(\omgomgi)$ for $\ddd=2,3$
\\
Fractional & \eqref{fraker} & $\Hsoi$
\\
\hline
\end{tabular}
\end{table}

\paragraph{\bo{Translation-invariant,} integrable kernel functions.}
The kernel function \bo{is translation-invariant, that is, $\phid\xyp=\phid(\ymx)$, and} satisfies   \bo{for some positive constant $C>0$,}
\begin{equation}\label{intker}
    \bo{C\leq \int_{\omgomgi} \phid (\ymx)  \dyb < \infty \quad\text{for all } \xb\in\omgomgi.}
\end{equation}
In the design of numerical methods, one may need to differentiate between cases for which
  \renewcommand{\leftmargini}{20pt}
\begin{itemize}\setlength\itemsep{4pt}
\item \bo{$\phid {(\zb)}$} is a smooth, bounded function of  \bo{${\zb}$},

\item \bo{$\phid {(\zb)}$} is a bounded function of  \bo{${\zb}$} but is not smooth, \eg\ a piecewise constant function,

\item \bo{$\phid {(\zb)}$} is a singular but integrable function, \eg\  \bo{$\phid(\zb)\propto1/|\zb|^r$ with $r<\ddd$}.
\end{itemize}
\bo{A special case of translation-invariant kernel functions \roo{is that of} {\em radial kernel functions} for which $\phid(\zb)=\phid(|\zb|)$.}

\paragraph{`Critical' kernel functions.}
The kernel function satisfies 
\begin{equation}\label{dddker}
    \phid\xyp \propto \ffrac1{|\ymx|^\ddd}.
\end{equation}
We refer to this kernel as {\em critical} because it `just misses being integrable', \tie, $1/{|\ymx|^{\ddd-\eps}}$ is integrable for any $\eps>0$ and for any $\ddd$, but is not integrable for any $\ddd$ if $\eps=0$. For the other kernels in this list, the corresponding energy space is precisely known; see \fo{Table~\ref{ker-prop}}. For critical kernels, that is not the case; what is known is that energy space is a Hilbert space and is a strict subspace of $L^2(\omgomgi)$.

\paragraph{`Peridynamic' kernel functions.}
The kernel function satisfies 
\begin{equation}\label{perker}
    \phid\xyp \propto \ffrac1{|\ymx|}.
\end{equation}
We refer to this kernel as {\em peridynamic} because it has the same singular behaviour as {one of} the  commonly used kernels in peridynamic models \bo{for solid mechanics}. Note that for $\ddd=1$ this kernel function is a member of the critical class, whereas for $\ddd=2$ it is integrable and for $\ddd=3$ it is {square-integrable}.

\paragraph{Fractional kernel functions.}
The kernel function satisfies 
\begin{equation}\label{fraker}
    \phid\xyp \propto \ffrac{1}{|\ymx|^{\ddd+2s}}.
\end{equation}
We refer to this kernel as {\em fractional} because the singularity at $\yb=\xb$ is that of the kernel for the fractional Laplace operator. This kernel plays a central role in Section~\ref{sec:mod-frac}.

%%Tab 2.1 original position

\begin{rem} \bo{We note that many {\em square-integrable} kernel functions of the more general form $\phid(\xb,\yb)$ also have $\VVV=L^2(\omgomgi)$ as the energy space. \roo{However, we do not} consider this class of kernels in the later sections.} \hfill\actaqed
\end{rem}

\begin{rem}[solution smoothness in nonlocal diffusion]\hspz% 
  For \bo{trans\-lation-invariant} integrable kernel functions such as \bo{
  \eqref{intker},}  the {\em solution \bo{may not be}  smoother than the data.} For example, if $g(\xb)=0$ and $f(\xb)\in L^2(\omg)$, then \bo{generically,} the solution $u(\xb)=({-}\mcLd)^{-1}f\in L^2(\omgomgi)$. Because $L^2(\omg)$ includes, \feeg, functions with jump discontinuities \bo{possibly over some \roo{co-dimension~$1$} interfaces}, this means that if the data $f(\xb)$ has a jump discontinuity, so \bo{may}   the solution $u(\xb)$. For the fractional kernel function, we could have $2s$ derivatives of smoothing, \tie, if the data $f(\xb)\in H^{-s}(\omg)$, then $u(\xb)\in H^{s}(\omgomgi)$, {where $H^{s}(\omgomgi)$} denotes a fractional Sobolev space and $H^{-s}(\omg)$ its dual space. Note that for $s<\gfrac12$, $H^{s}(\omgomgi)$ also contains functions having jump discontinuities. \hfill\actaqed
\end{rem}

\begin{rem}[scaling constants and limit behaviour] 
  Often kernel\linebreak
  functions take the form
\[
  \phid\xyp = C \phiv\xyp,
\]
{where $C$} denotes a scaling constant that, depending on the specific instance, may depend on $\del$ and also on other parameters appearing in the definition of $\phid\xyp$. One approach towards defining the scaling constant $C$ is to choose it so that in the limit $\del\to0$, a nonlocal model and its solution {reduce} to a local PDE model and its solution, respectively. Another approach is to choose $C$ so that in the limit $\del\to\infty$, the nonlocal model and its solution reduce to a  model that is posed on $\Rd$. {Below, we} will have occasion to follow {both approaches}. In any case, the scaling constant is a modelling choice so that, for example, if neither type of limiting behaviour is of interest, then the scaling constant should be chosen based on the physics being modelled.~\hfill\actaqed
\end{rem}

\vspace{5pt}
\section{Fractional diffusion models}\label{sec:mod-frac}

In this section we consider a widely used nonlocal diffusion model that involves the fractional Laplace operator {$-(-\Delta)^s$} with $s\in(0,1)$. In a real sense, fractional diffusion models can be thought of as being a special case of the general nonlocal diffusion models discussed in Sections~\ref{sec:mod-strong} and~\ref{sec:mod-weak}. After all, fractional diffusion models are defined in terms of a specific choice for the kernel function; see \eqref{fraker}. On the other hand, the scope of fractional diffusion models and the literature devoted to their analysis, approximation and application hugely dwarf what is available for other nonlocal diffusion models. {Thus a more detailed consideration of fractional models is warranted here.} 

The {\em integral definition of the fractional Laplacian}, {for  $u(\xb) \colon  \Rd\to \Rd$ with $\ddd\in\mathbb{N}_{+}$, is given by}
\begin{align}\label{intfl}
 (-\Delta)^{s} u ({\xb})  & =  \int_{\Rd} ( u( {\xb})-u( {\yb}) ) \gams{\xyp} \dyb
  \quad\text{for all } \xb\in\Rd
  \end{align}
{with}
\begin{equation}\label{gams}
 \gams{\xyp}=\ffrac{\cds}{|\ymx|^{\ddd+2s}}.
\end{equation}
Note that this is a case in which the integral in \eqref{intfl} should be interpreted in the {principal value} sense. In \eqref{intfl}, \bo{the constant $\cds$ is given by}
\begin{equation}\label{cds}
  \cds = \ffrac{2^{2s}s\Gamma (s+\gfrac{\ddd}{2})}{\pi^{\ddd/2}\Gamma (1-s)},
\end{equation}
where $\Gamma$ denotes the {gamma} function. 

The integral representation of the fractional Laplacian given in \eqref{intfl} 
can be viewed as a special case of the nonlocal  Laplacian given by \eqref{str-lap}. It is also equivalent to and can be directly derived from the Fourier representation \cite{Valdinoci:2009}
\begin{align}\label{E:definition_spectral}
(-\Delta)^{s} u(\xb) &= \ffrac{1}{(2\pi)^\ddd} \int_{\Rd} |\xib|^{2s}  ( u , \rme^{-\rmi \xib \cdot \xb} ) \,\rme^{\rmi \xib \cdot \xb} \D \xib \notag   \\*
  &= \mathcal{F}^{-1} ( |\xib|^{2s} \mathcal{F}\{u\}(\xib) )(\xb),
\end{align}
where $\mathcal{F}$ denotes the Fourier transform.
{If $s=1$, the spectral operator \eqref{E:definition_spectral} 
coincides with the usual PDE Laplacian $-\Delta$, whereas it reproduces the identity operator when $s=0$.  In fact, it is {well known} that  $(-\Delta)^s u(x) \rightarrow u(x)$ as $s\to0^+$ and $(-\Delta)^s u(x) \rightarrow -\Delta u(x)$ as $s\to 1^-$ for a regular function $u$; see \eg\ \citeasnoun[Theorems~3 and~4]{Stinga:2019} and \citeasnoun[Proposition~4.4]{FracGuide:2012}}. 

The fractional Laplacian \eqref{intfl} is often used to model superdiffusion for which the mean-squared particle displacement grows faster than that for PDE models of diffusion. At microscopic scales, in contrast to standard diffusion that is described by Brownian motion, superdiffusion can be described by L\'evy flights in which the length of particle jumps follows a {heavy-tailed} distribution, reflecting the long-range interactions between particles.
See \eg\ \citeasnoun{MetzlerKlafter:2000} {and} \citeasnoun{Sokolov:2022} for discussions of the physical background and practical applications of anomalous diffusion.
 
\subsection{Integral fractional Laplacian models for diffusion on bounded domains}

Letting $\omg\subset\Rd$ denote a bounded Lipschitz domain, we define the \emph{integral fractional Laplacian} $-(-\Delta)^{s}$ to be the restriction of the full-space operator to functions satisfying a volume constraint on ${\omgii}:=\Rd\setminus\omg$, \ie\ on the interaction domain corresponding to $\omg$. The \emph{fractional Poisson problem} is then given by
\begin{equation}\label{eq:fracPoisson}
\begin{cases}
 (-\Delta)^{s}u = f &\text{for all }\xb\in\omg,
\\
  u=g&\text{for all }\xb\in{\omgii} ,
\end{cases}
\end{equation}
where we have a given source term $f(\xb)$ and Dirichlet volume constraint data $g(\xb)$ defined for $\xb\in\omg$ and $\xb\in\omgii$, respectively.
{Problem} \eqref{eq:fracPoisson} is a nonlocal analogue of the local Poisson problem for the PDE Laplacian $\Delta$. {It is known (see \eg\ \citeb{Biccari:2018}) that as $s\rightarrow 1^-$, the solution of integral fractional diffusion model \eqref{eq:fracPoisson} strongly converges to the solution of the local diffusion problem in $H^{1-\epsilon}(\Omega)$.} 

A related operator on $\omg$ is the {\em regional fractional Laplacian} \cite{BogdanBurdzyChen:2003,ChenKim:2002}
\begin{equation}\label{intflr}
 -  (-\Delta) _{\text{regional}}^{s} u ({\xb})  = \cds~  \int_{\omg} \ffrac{u( {\yb})-u( {\xb})}{|\ymx|^{\ddd+2s}}\dyb,
\end{equation}
{which} is used in one of several approaches for generalizing the PDE Poisson problem with a homogeneous Neumann boundary condition {to the fractional Laplacian case \cite{Dipierro:2017}}. Note that in \eqref{eq:fracPoisson} the operator has not changed, \tie, the operator defined in \eqref{intfl} is used; what has changed is that domain of the operator is changed from $\Rd$ to the bounded domain $\omg$. On the other hand, in using the operator in \eqref{intflr}, not only is the domain changed in the same manner, but the operator itself has changed, \tie, \eqref{intfl} involves an integral over $\Rd$ whereas \eqref{intflr} involves an integral over $\omg$. 
 
To discuss the variational form of  the fractional Poisson problem \eqref{eq:fracPoisson}, we use the standard fractional Sobolev space $H^{s} (\Rd) $ defined via the Fourier transform as
\[  
H^s(\Rd)  = \biggl\{ u\in L^2(\Rd)  \colon   \int_{\Rd} (1+ |\xib|^{2s}) |\mcF u(\xib)|^2 \D \xib <\infty \biggr\}. 
\]
If $\omg\subset \Rd$ is a bounded domain, we define the Sobolev space $H^{s}(\omg)$ as \cite{Mclean:2000}
\begin{align*}
   H^{s} (\omg)&:=\{u\in L^{2} (\omg)   \colon  \norm{u}_{H^{s} (\omg)} < \infty\}
\end{align*}
equipped with the norm
\[
  \norm{u}_{H^{s} (\omg)}^{2}= \abso{u}^{2}_{H^{s} (\omg)} + \norm{u}_{L^{2} (\omg)}^{2},
\]
where we have the seminorm
\[
  \abso{u}^{2}_{H^{s} (\omg)}=\int_{\omg}  \int_{\omg}  \ffrac{ (u( {\yb})-u( {\xb})) ^{2}}{\abso{\ymx}^{\ddd+2s}} \dydx.
\]
Moreover, when imposing a homogeneous Dirichlet volume constraints, \eg\ $g=0$ in \eqref{eq:fracPoisson}, we use the space
\[
  {H}^{s}_c(\omg):=\{u\in H^{s}(\Rd)  \colon  u=0 \ \text{for all } \xb\in {\omgii}\}
\]
that is equipped with the norm
\begin{align*}
  \norm{u}_{{H}^{s}_c (\omg)}^{2} = \norm{u}_{H^{s} (\Rd) }^{2}=\norm{u}_{L^{2} (\omg)}^{2} + \abso{u}^{2}_{H^{s} (\Rd) }.
\end{align*}
For $s>1/2$, ${H}^{s}_c(\omg)$ coincides with the space $H_{0}^{s}(\omg)$ that is the closure of $C_{0}^{\infty}(\omg)$ with respect to the $H^{s}(\omg)$-norm, whereas for $s<1/2$, ${H}^{s}_c(\omg)$ is identical to $H^{s}(\omg)$.
In the critical case $s=1/2$, ${H}^{s}_c(\omg)\subsetneq H^{s}_{0}(\omg)$. See \eg\ \citeasnoun[Chapter~3]{Mclean:2000} for a detailed discussion.

A variational form of \eqref{eq:fracPoisson} is derived starting from the {integration-by-parts} formula (\bo{\eg\ the nonlocal Green's first identity \eqref{str-gfi1}}) 
\begin{align*}  
&\ffrac12\int_{\mathbb{R}^d}\int_{\mathbb{R}^d}
(u({\yb})-u({\xb}))  (v({\yb})-v({\xb}))  \gams\xyp  \dydx \\*
&\quad = \int_{\omg} ((-\Delta)^{s}u({\xb})) v({\xb}) \dxb  + \int_{\omgii} (\mcN_{s}u({\xb}))  v({\xb}) \dxb ,
\end{align*}
where, as in \eqref{str-flux1}, we have the nonlocal  operator given by
  \[
  \mathcal{N}_{s}u({\xb})= \int_{\mathbb{R}^d} (u({\xb})-u({\yb}))  \gams\xyp  {\dyb}  \quad\text{for all } \xb\in{\omgii}.
  \]
\bo{An alternative definition of
$\mathcal{N}_s$ has  been presented in  \citeasnoun{Dipierro:2017}, corresponding to a nonlocal operator $\mathcal{N}_\delta$ 
that is defined on $\Omega_\delta$ as an integral over the domain $\Omega$ only, instead of $\omgomgi$ as in the case of
\eqref{str-flux}.}
Substituting the volume constraint in \eqref{eq:fracPoisson} and $v(\xb)=0$ on $\omgii$, we then have the {following weak formulation:} 
\begin{equation}\label{eq:fracPoissonVariational}
  \begin{minipage}[b]{25pc}
  {Find} $ u(\xb)\in {H}^{s}(\omg)$  such that $ u(\xb)=g(\xb) $ for all $ \xb\in\omgii $  and
  \[
\bo{  \mcA_s(u,v)=\fv + (g,v)_{{\omg,\omgii}}\qfa v\in {H}^{s}_c(\omg),}
  \]
  \end{minipage}
\end{equation}
where
\vst %%20200406
\begin{equation}\label{aaaaaa}
\bal
 \mcA_s(u,v)
  &=\underbrace{\ffrac{1}{2} \int_{\omg}  \int_{\omg}\;   (u({\yb})-u({\xb}))  (v({\yb}) -v({\xb}) )  \gams\xyp \dydx}_{\mcA_{\omg,\omg}(u,v)} \\
    &\quad\, + \underbrace{\int_{\omg}u({\xb})v({\xb})\int_{{\omgii}} \gams\xyp \dydx}_{\mcA_{\omg,\omgii}(u,v)}
\eal
\end{equation}
and
\vst %%20200406
\begin{equation}\label{rhsfl}
(g,v)_{\omg,\omgii} = \int_\omg v(\xb) \int_\omgii g(\yb)\gams\xyp \dydx.
\end{equation}
The bilinear form $\mcA_s(\cdotsp,\cdotsp)$ is ${H}^{s}_c(\omg)$-coercive and continuous so that as long as the right-hand side of \eqref{eq:fracPoissonVariational} is bounded, the {well-posedness} of that problem follows from the Lax--Milgram theorem. \bo{\roo{If} $g=0$, $\mcA_s(\cdotsp,\cdotsp)$ is identical to the bilinear form $\mcA_\delta(\cdotsp,\cdotsp)$ defined in \eqref{weak-bform} with $\delta=\infty$.}
\boo{On the other hand, for $g\neq 0$, $\mcA_s(\cdotsp,\cdotsp)$ can be seen as the $\delta=\infty$ limit of the bilinear form corresponding to some inhomogeneously constrained nonlocal diffusion problems discussed in \citeasnoun{du19cbms}, instead of the bilinear form $\mcA_\delta(\cdotsp,\cdotsp)$ given by \eqref{weak-bform}.}
Note that the bilinear form \bo{$\mcA_{{\omg,\omg}}(\cdotsp,\cdotsp)$} can be seen to correspond to the regional fractional Laplacian.

Note that
\vst %%20200406
\[
  \mcA_{\omg,\omgii}(u,v) = \cds \int_{\omg}  u({\xb}) v({\xb})\bigg( \int_{\omgii}   \ffrac{1}{\abso{\ymx}^{\ddd+2s}}\dyb\bigg) \dxb.
\]
The identity
\vst %%20200406
\[
\bo{  \ffrac{1}{\abso{\ymx}^{\ddd+2s}} =  - \ffrac{1}{2s} \nabla_{{\yb}}\cdot \ffrac{\ymx}{\abso{\ymx}^{\ddd+2s}}}
\]
and {Gauss's} theorem results in
\[
\mcA_{\omg,\omgii}(u,v) = \ffrac{\cds}{2s} \int_{\omg}  u({\xb}) v({\xb})\bigg( \int_{\partial\omg}   \ffrac{(\ymx)\cdot{\bm n}_\yb}{\abso{\ymx}^{\ddd+2s}}\dyb\bigg) \dxb,
\]
where ${\bm n}_\yb$ denotes the unit outer normal to $\partial\omg$ at ${\yb}$. We then have that 
\begin{equation}\label{aass}
\mcA_s(u,v) = \mcA_{{\omg,\omg}}(u,v) + 
\ffrac{\cds}{2s} \int_{\omg}  u({\xb}) v({\xb})\bigg( \int_{\partial\omg}   \ffrac{(\ymx)\cdot{\bm n}_\yb}{\abso{\ymx}^{\ddd+2s}}\dyb\bigg) \dxb
\end{equation}
so that an integral over the unbounded domain $\omgii$ can be avoided in any computation involving $\mcA_s(u,v)$. Of course, for a homogeneous volume constraint $g(\xb)=0$ on $\omgii$, the right-hand side in \eqref{eq:fracPoissonVariational} also only involves integrals over $\omg$. For $g(\xb)\ne 0$, one does have to evaluate, for $\xb\in\omg$, the {\em data integral} $\int_\omgii g(\yb)\gams\xyp  \D \yb$ that, in principle, can be approximated using a quadrature rule. 

Above we have incorporated the inhomogeneous volume constraint directly into the weak form. An alternative is to enforce that condition via Lagrange multipliers \cite{Acosta:2019}.

\subsubsection{Truncated interaction domains}\label{truncted-flap}

An alternative to the approach based on the bilinear form \eqref{aass} is to   {\em truncate the integration domain} $\omgii$ as is considered in \citeasnoun{MG:2013}. For $\del> {\text{diam}} (\omg)$, we pose the truncated problem
\begin{equation}\label{eq:fracPoissont}
\begin{cases}
 (-\Delta)^{s}_\del u = f &\text{for all }\xb\in\omg,
\\
  u=g&\text{for all }\xb\in{\omgid} ,
\end{cases} 
\end{equation}
where
\[
 (-\Delta) ^{s}_\del u ({\xb})  := \int_{\omgomgi}  ( u( {\xb})-u( {\yb}))  \gams{\xyp} \dyb.
\]
{\citeasnoun{MG:2013} and \citeasnoun{burko1} have shown} that if $u_\infty$ and $u_\del$ denote the solutions of \eqref{eq:fracPoisson} and \eqref{eq:fracPoissont}, respectively, then
\begin{equation}\label{inflimit}
   \| u_\del - u_\infty\|_{H^s(\omg)} \le \ffrac{C}{\del^{2s}}
   \quad\mbox{as } \del\to\infty,
\end{equation}
where $C$ depends on norms of $f$ and $g$.

\begin{rem}\label{truncfl}\hspace{-7.7pt}% 
Clearly, \eqref{inflimit} implies that the solution of \eqref{eq:fracPoissont} can be viewed as an approximation of the solution of \eqref{eq:fracPoisson}. However, \eqref{eq:fracPoissont} and its solution are useful in their own right. In \eqref{eq:fracPoisson}, the horizon $\del$ is infinite, but, in practice, although $\del$ may be large, it is not likely to be infinite. In \eqref{eq:fracPoissont}, $\del$ is assumed large but finite so that it may provide a better model for practical applications that feature large horizons. In this case one can reverse roles and interpret the solution of
\eqref{eq:fracPoisson} to be an approximation of the solution of \eqref{eq:fracPoissont}, with {\eqref{inflimit} now} telling us something about the error incurred by replacing the perhaps more useful model \eqref{eq:fracPoissont} {with} the much more studied model \eqref{eq:fracPoisson}. \hfill\actaqed
\end{rem}

\subsection{Spectral fractional Laplacian models for diffusion on bounded domains}

An {alternative} definition of a fractional-order Laplacian on a bounded domain $\omg$ makes use of spectral information about the PDE Laplacian. Here we focus on the case of homogeneous constraints.

Let $0<\lambda_{0}\leq\lambda_{1}\leq\cdots$ and $\phiv_{0},\phiv_{1},\ldots$ denote the eigenvalues and the corresponding eigenfunctions of the PDE Laplacian $-\Delta$ with a homogeneous Dirichlet {\em boundary} condition, \tie, we have
\begin{equation}\label{eq:Eig}
\begin{cases}
-\Delta\phiv_{m}({\xb}) =\lambda_{m}\phiv_{m}({\xb})& \text{for all }{\xb}\in\omg ,\\
    \phiv_{m} ({\xb})  =0 & \text{for all } {\xb}\in\partial\omg,
\end{cases}
\end{equation}
where the orthogonal eigenfunctions are normalized so that $\norm{\phiv_{m}}_{L^{2}(\omg)}=1$. Then the eigenfunctions $\{\phiv_{m}\}_{m=0}^{\infty}$ form a complete orthonormal basis for $L^{2}(\omg)$. As a result, any function $u\in L^{2}(\omg)$ can be expanded as
\[
  u=\sum_{m=0}^{\infty}u_{m}\phiv_{m} \quad\text{with } u_{m}= ({u},{\phiv_{m}}), 
\]
where $(\cdot,\cdot)$ denotes the inner product of $L^2(\omg)$. We then have 
\[
 (-\Delta) u({\xb})= \sum_{m=0}^{\infty}u_{m} \lambda_m \phiv_{m} ({\xb}) ,
\]
and the {\em spectral fractional Laplacian} of order $s\in(0,1)$ with homogeneous boundary condition is given by
\[
   ({-\Delta}|_{\omg,0}) ^{s}u ({\xb}) 
  = \sum_{m=0}^{\infty} u_{m}\lambda_{m}^{s} \phiv_{m}({\xb}).
\]
As $s\rightarrow 0$, the identity operator is recovered, whereas the integer order Laplacian is recovered as $s\rightarrow1$. 

For a given source term $f$, the spectral fractional Laplacian Poisson problem for $u(\xb)$, $\xb\in\omg$, is then given by
\begin{align}\label{eqn:poisson-spectral}
   ({-\Delta}|_{\omg,0}) ^{s}u ({\xb})  = f(\xb) \quad\text{for all }\xb\in\omg.
\end{align}

Using the heat kernel $p_{\omg,0}({\xb},{\yb},t)$, the spectral fractional Laplacian can be rewritten (see \eg\ \citeb{Abatangelo:2017}) in the form of the nonlocal operator
\[
   ({-\Delta}|_{\omg,0}) ^{s}u({\xb})
  = \int_{\omg} \gamma_{\omg,0;s}({\xb},{\yb}) (u({\yb})-u({\xb}))  \dyb + \kappa_{\omg,0;s}({\xb})u({\xb}),
\]
\roo{where}
\begin{equation*} \begin{aligned} 
  \gamma_{\omg,0;s}\xyp &= \ffrac{s}{\Gamma(1-s)} \int_{0}^{\infty}  p_{\omg,0}(t,{\xb},{\yb}) \ffrac{1}{t^{1+s}}{\D t},
  \\*
  \kappa_{\omg,0;s}({\xb}) &= \ffrac{s}{\Gamma(1-s)} \int_{0}^{\infty} \biggl(1-\int_{\omg}p_{\omg,0}(t,{\xb},{\yb}) \D \yb\biggr) \ffrac{1}{t^{1+s}} \D t.
  \end{aligned}
  \pagebreak %%20200331
\end{equation*} 
 If $\omg=\Rd$, the integral and the spectral definitions of the fractional Laplacian coincide, but they are different for bounded domains \cite{Servadei:2014}. For example, in the particular case of $\omg=\Rd_{+}$, \ie\ the half-plane, the kernel of the spectral fractional Laplacian is given by
\[
  \gamma_{\Rd_{+},0;s}\xyp = \gamma_{s}\xyp-\gamma_{s}({\xb},-{\xb}) \neq \gamma_{s}\xyp \quad \text{for all } {\xb},{\yb}\in\Rd_{+}.
\]

Note that the spectral definition can allow for the treatment of non-homogeneous boundary conditions of various types; see \eg\ \citeasnoun{Antil:2018} {and} \citeasnoun{Cusimano:2018}.

\subsection{Additional considerations about fractional Laplacian models}

\subsubsection{Other integral representations of the fractional Laplacian}\label{dunftay}

The {\em Dunford--Taylor integral} is a powerful tool in the numerical analysis of the fractional diffusion problem  \cite{Bonito:2015,Bonito:2017,Bonito:2019}. In particular, the solution $u$ of the Poisson problem \eqref{eqn:poisson-spectral} involving the spectral fractional Laplacian can be given \bo{by} the Dunford--Taylor representation 
\[
  u= ({-\Delta}|_{\omg,0})^{{-s}} f = \ffrac{\sin (s\pi) }{\pi} \int_{0}^{\infty} {\mu^{-s} (\mu-\Delta) ^{-1}f  \D \mu}.
\]
A combination of using a sinc quadrature rule and a {finite element} discretization of the reaction--diffusion-type term in the integrand was explored {by \citename{Bonito:2015} \citeyear{Bonito:2015,Bonito:2017}}.

Using similar techniques, the bilinear form $\mcA_s(\cdot,\cdot)$ associated with the integral fractional Laplacian can be rewritten as \cite{Bonito:2019} 
\begin{equation}\label{dtbf}
 \mcA_s(u,v) = \frac{2 \sin (\pi s)}{\pi} \int_0^\infty  {\mu^{1-2s}} \int_{\Rd}
 ((-\Delta)(I-\mu^2 \Delta)^{-1} u )(\xb)\,  v(\xb) \dxb \D\mu ,
\end{equation}
\bo{for any $u,v \in H_c^s(\Omega)$.}

\subsubsection{Extension representations}

An {\em extension representation} of the fractional Laplacian was introduced in \citeasnoun{Molchanov:1969} and popularized by \citeasnoun{CaffarelliSilvestre:2007}, {who showed} that the fractional Poisson problem $\eqref{eq:fracPoisson}$ posed on $\Rd$ can be recast as a Neumann-to-Dirichlet mapping over the extended domain $\Rd\times[0,\infty)$, \tie, we have
\begin{equation}\label{eq:extensionProblem}
\begin{cases}
    -\nabla\cdotsp z^{\beta} \nabla U ({\xb},z ) = 0 &  
    \text{for all }({\xb},z)\in\Rd\times[0,\infty) , \\[7pt]
    \ffrac{\partial U}{\partial \nu^{\beta}} ({\xb}) = 2^{1-2s}\ffrac{\Gamma(1-s)}{\Gamma(s)} f({\xb})& \text{for all }{\xb}\in \Rd,
\end{cases}
\end{equation}
where $\beta = 1-2s$ and
\[
  \ffrac{\partial U}{\partial \nu^{\beta}}({\xb})= -\lim_{{z}\rightarrow 0^{+}}z^{\beta}\ffrac{\partial U}{\partial z} ({\xb},{z}) ,
\]
with the solution to \eqref{eq:fracPoisson} recovered by taking the trace of $U$ on $\Rd$.

The spectral fractional Laplacian can be recovered, by restricting the
extension domain from $\Rd\times[0,\infty)$ to $\omg\times[0,\infty)$ and imposing a
homogeneous Dirichlet boundary condition on the lateral surface \cite{Stinga:2010}.

The apparent advantage of {the} extension problem is that the nonlocal problem is replaced {with} a classical, integer-order local problem.
This comes at the price of having to deal with a singular weight function and an additional spatial dimension.

\subsubsection{Regularity of solutions}\label{fracreg}

In the classical integer-order case, smoother domains and smoother right-hand sides result in smoother solutions of the Poisson problem; see \eg\  \citeasnoun{taylor}. A lifting property of this type does not hold for the fractional-order Poisson problems \eqref{eq:fracPoisson} and \eqref{eqn:poisson-spectral}.
The typical solution behaviour of \eqref{eq:fracPoisson} close to the boundary can be characterized (see \citeb{ros2014dirichlet}) as
\begin{equation}\label{bdybeha}
  u({\xb})\approx \operatorname{dist}(\xb,\partial\omg)^{s}.
\end{equation}

The Sobolev regularity of the solution of \eqref{eq:fracPoisson} was studied in \citeasnoun{Acosta:shortFEM} as a special case of a result in \citeasnoun{Grubb:2015}.
Suppose that $\partial\omg\in C^{\infty}$ and $f\in H^{r} (\omg)$ for $r\geq-s$ and let $u\in{H}^{s}_c (\omg)$ denote the solution of the fractional Poisson problem \eqref{eq:fracPoissonVariational}. Then {we obtain} the regularity estimate
\[
    u\in
    \begin{cases}
      H^{2s+r} (\omg) & \text{if } 0<s+r<1/2 , \\
      H^{s+1/2-\varepsilon} (\omg)\ \text{for all}\ \varepsilon>0 & \text{if } 1/2\leq s+r .
    \end{cases}
\]

This result is extended to the case of non-homogeneous volume condition in \citeasnoun{Acosta:2019}. Additional regularity results with respect to $s$ and, for the truncated fractional Laplacian, with respect to $\del$, are derived in \citeasnoun{burko1}.

The solution of the spectral fractional Poisson problem \eqref{eqn:poisson-spectral} displays different properties.
Its behaviour close to the boundary is given by (see \citeb{CaffarelliStinga2016_FractionalEllipticEquationsCaccioppoliEstimatesRegularity})
\begin{align*}
  u(\xb)\approx \operatorname{dist}(\xb,\partial\omg)^{\min\{2s,1\}}, \quad s\neq 1/2.
\end{align*}
{\citeasnoun{Grubb2015_RegularitySpectralFractionalDirichletNeumannProblems} showed} that if $f\in H_c^{r} (\omg)$, $r\geq-s$, then the weak solution of the spectral fractional Poisson problem \eqref{eqn:poisson-spectral} satisfies $u\in H_c^{r+2s} (\omg)$.

More detailed regularity results for the spectral fractional Laplacian can be found in \citeasnoun{Grubb2015_RegularitySpectralFractionalDirichletNeumannProblems}.

\subsubsection{Analytic solutions in one and two dimensions}

Closed-form solutions of the fractional Poisson problem \eqref{eq:fracPoisson} posed on the unit ball are available; see \citeasnoun{DydaKuznetsovEtAl2016_FractionalLaplaceOperatorMeijerGFunction} for detailed derivations. These solutions are useful to have in hand, \eg\ for verifying numerical results.

In $\ddd=1$ dimension, for the source term
\[
  f_{k,0}^{1D}= 2^{2s}\Gamma (1+s) ^{2}\binom{s+k-1/2}{s}\binom{s+k}{s} P_{k}^{(s,-1/2)} (2x^{2}-1) 
\]
with $k\geq0$, the solution is given by
\[
  u_{k,0}^{1D}= P_{k}^{(s,-1/2)} (2x^{2}-1)   (1-x^{2}) ^{s}_{+},
\]
where
$$\binom{a}{b}=\ffrac{\Gamma (a+1) }{\Gamma (b+1) \Gamma (a-b+1) }$$
denotes a generalized binomial coefficient, $P_{k}^{(a,b)}$ denote the Jacobi polynomials, and $a_{+}=\max\{0,a\}$. Moreover, for
\[
  f_{k,1}^{1D} = 2^{2s}\Gamma (1+s) ^{2}\binom{s+k+1/2}{s}\binom{s+k}{s} x P_{k}^{(s,1/2)} (2x^{2}-1) 
\]
with $k\geq0$, the solution is given by
\[
  u_{k,1}^{1D}= x P_{k}^{(s,1/2)} (2x^{2}-1)   (1-x^{2}) ^{s}_{+}.
\]

Turning to $\ddd=2$ dimensions, for
\[
  f_{k,\ell}^{2D}
  = 2^{2s}\Gamma (1+s) ^{2}\binom{s+k+\ell}{s}\binom{s+k}{s} r^{\ell}\cos (\ell\phiv)  P_{k}^{(s,\ell)} (2r^{2}-1) 
\]
with $\ell,k\geq0$ and $(r,\phiv)$ denoting polar coordinates, the solution is given by
\[
  u_{k,\ell}^{2D} = r^{\ell}\cos (\ell\phiv)  P_{k}^{(s,\ell)} (2r^{2}-1)  (1-r^{2}) ^{s}_{+}.
\]

\begin{figure} %%fig3.1
\centering
\includegraphics[width=160pt,viewport=0 0 200 200,clip]{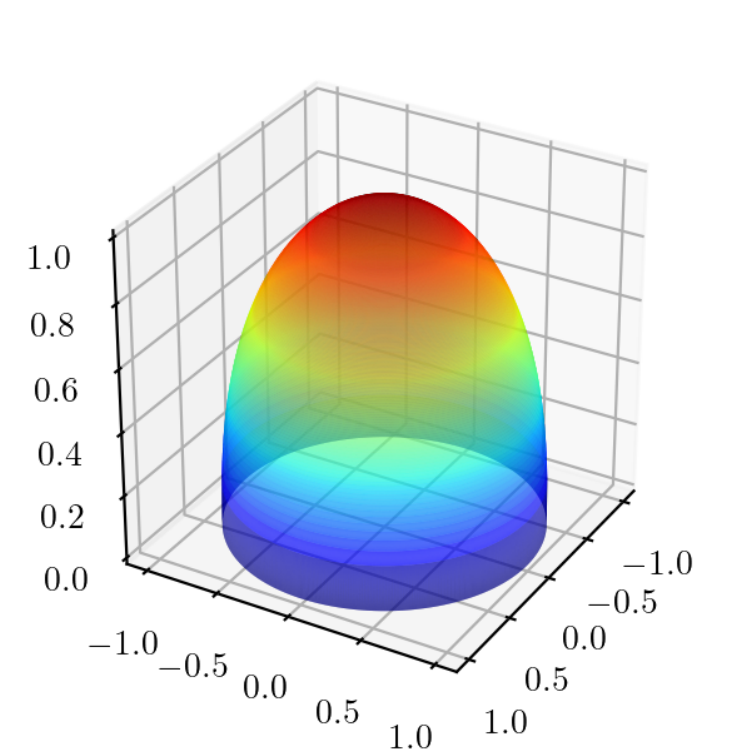}\hspace{24pt} %%glusa1
\includegraphics[width=160pt,viewport=0 0 200 200,clip]{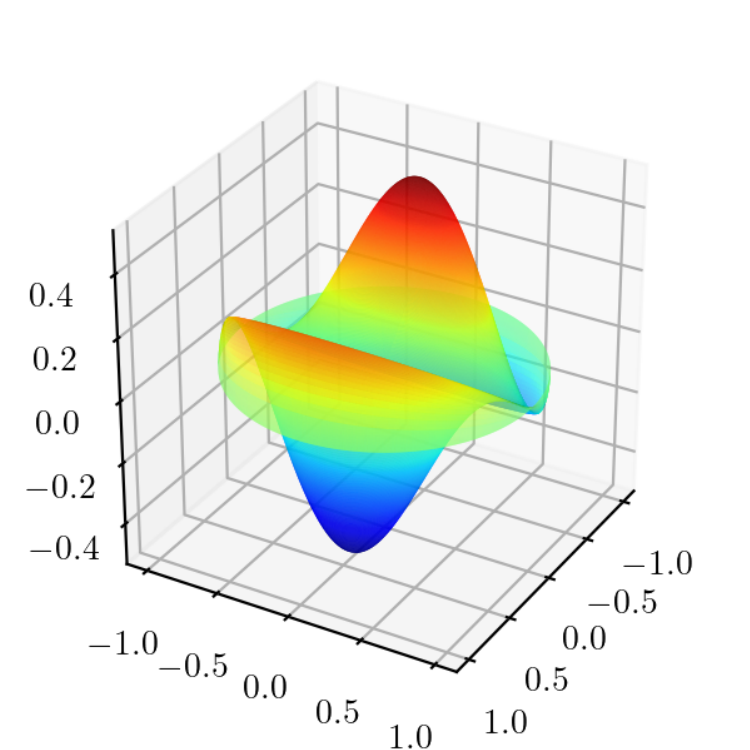}\\ %%glusa2
\caption{Analytic solutions $u_{0,0}^{2D}$ for $s=0.4$ and $u_{1,1}^{2D}$ for $s=0.6$. The behaviour \eqref{bdybeha} close to the boundary is apparent. 
}\label{fig:integralFractionalSolutions}
\end{figure}

Two example solutions are shown in \fo{Figure~\ref{fig:integralFractionalSolutions}}. 
Observe that the solutions {contain a term that} behaves like $\delta (\xb) ^{s}$, where $\delta (\xb) $ is the distance from $\xb\in\omg$ to the boundary.

\subsection{Inhomogeneous constitutive functions in fractional models}\label{sec:flvarcoef}

So {far in this section we have only considered the fractional Laplacian and its variants,} which is analogous to considering local diffusion problems having a diffusivity that is not only {constant but equal} to one. Of course, in the local case, it is an easy matter to treat inhomogeneous diffusivities, \tie, {we merely replace $\Delta u$ with} $\nabla\cdot(\Dbf \nabla u)$, where {$\Dbf(\xb)$ generally denotes an} inhomogeneous constitutive tensor. Using the nonlocal vector calculus, it \roo{is} almost as easy to define integral fractional models in which the constitutive properties of the media considered are {inhomogeneous even} if the nonlocal constitutive function is a tensor. Here we only consider generalizations of the integral fractional Laplacian \eqref{intfl}; generalizations of the other integral operators discussed in this section follow immediately from that for \eqref{intfl}.

The inhomogeneous fractional Laplacian is simply defined by choosing the kernel in \eqref{intfl} to now be
\begin{equation}\label{gamd1f}
\gams\xyp= \cds\ffrac{\thes\xyp}{|\ymx|^{\ddd+2s}}  \quad\text{for all } \xb,\yb\in\Rd,
\end{equation}
where
\begin{equation}\label{ethedf}
\thes\xyp :=\ffrac{\ymx}{|\ymx|} \cdotsp \biggl(
     \Thebs\xyp    \ffrac{\ymx}{|\ymx|}\biggr)
\end{equation}
with $\Thebs\xyp$ denoting a constitutive tensor. Thus we have defined the {\em inhomogeneous integral fractional Poisson operator}
\begin{equation}\label{intfi}
(-\Delta)^{s}_{\Thebs} u ({\xb}) :=  \int_{\Rd}  \cds\thes\xyp \ffrac{u( {\xb})-u( {\yb})}{|\ymx|^{\ddd+2s}} \dyb
 \quad\text{for all } \xb\in\Rd.
\end{equation}
Note that we could have arrived {at} \eqref{intfi} by defining 
\[
(-\Delta)^{s}_{\Thebs} u ({\xb}) := \frac12 \mcDf (\Thebs \mcDsf u),
\]
where $\mcDf$ and $\mcDsf$ are defined by \eqref{str-div} and \eqref{str-divs}, respectively, with
\begin{equation}\label{ealpbdf}
\alpbs\xyp = \ffrac{\ymx}{|\ymx|}  \sqrt{\ffrac{\cds}{|\ymx|^{\ddd+2s}}}.
\end{equation}

If $\Thebs\xyp=\Thebs(|\ymx|)$, \ie\ $\Thebs$ is a {\em radial function}, 
then \eqref{intfl} is a model for a homogeneous medium, \tie, in this case, \eqref{intfl} is a fractional analogue of the PDE $-\nabla\cdot(\Dbf\nabla u)=f$ in which $\Dbf$ is a constant tensor. If in addition $\Thebs$ is a scalar tensor,
{then} \eqref{intfl} is a fractional analogue of the PDE $-\kappa\Delta u=f$ in which $\kappa$ is a constant. Thus, to obtain an inhomogeneous integral fractional model, $\Thebs\xyp$ cannot be a radial function. See Remark~2.3.

\vspace{5pt}
\part[Numerical methods for nonlocal and\\ fractional models]{Numerical methods for nonlocal and fractional models}\label{part2}

In this part we consider the approximation, via finite element, finite difference and spectral methods, of solutions of both weak and strong formulations of both nonlocal and fractional diffusion models. In {so doing}, we encounter, as we did in Part~1, problems posed on $\Rd$ or on bounded domains $\omg\subset\Rd$, and also encounter problems for which the horizon $\del=\infty$ or $\diam (\omg) < \del < \infty$ or $\del < \diam (\omg)$.

\vspace{5pt}
\section{Introductory remarks}

In most cases we consider, two parameters will appear in the design of discretization algorithms, namely the horizon $\del$ and a discretization parameter such as a grid size parameter $h$ for finite element and finite difference methods or the dimension $N$ of a basis for spectral methods. The limiting behaviours of continuous and discretized models are thus of interest. Four types of limits can arise: for continuous, \ie\ un-discretized, models, 
  \renewcommand{\leftmargini}{20pt}
\begin{itemize}
\item the limits $\del\to0$ and $\del\to\infty$,
\end{itemize}
 and for discretized problems,
\begin{itemize}\setlength\itemsep{4pt}
\item for fixed $\del$, the limit $h\to0$ or $N\to\infty$,

\item for fixed $h$ or $N$, the limits $\del\to0$ and $\del\to\infty$,

\item simultaneous limits such as both $\del\to0$ and $h\to0$.
\end{itemize}
\bo{For the first of these, one can find discussions \roo{of} the local and global (fractional) limits of nonlocal continuum  models in, for example, \roo{\citeasnoun{md15non}, \citeasnoun{tdg16acm} and \citeasnoun{du19cbms}}.}
 For the second of these, one must be cognizant of how the constants appearing in error estimates depend on $\del$, whereas for the third, the same can be said for the dependence of constants on $h$ or $N$.

It turns out that for some algorithms, the order in which limits are taken matters, \tie, the limits obtained are different if the limits are taken in a different order, \eg\ $\del\to0$ and then $h\to0$, or $h\to0$ and then $\del\to0$, or if $\del$ and $h$ are related in some way so that they simultaneously tend to zero. Presumably, \bo{it is possible that for some algorithms} at least one of the limits is `wrong' in some sense. Ideally, unless one is totally uninterested in limiting behaviours, one would {prefer the} way in which limits are taken {not to} affect the limit obtained. In the remainder of our introductory remarks, we expand on this concept.

\subsection{Asymptotically compatible schemes}\label{ACschemes}

It is known in practice that to obtain consistency between nonlocal models and corresponding local PDE models, the mesh or quadrature point spacing may have to be reduced at a faster pace than the reduction of the horizon parameter \cite{BYASSX09,BoDu10,dt-ChGu11}. Otherwise there could potentially be complications, most notably inconsistent limiting solutions when the horizon parameter is coupled proportionally to the discretization parameter \cite{tian17thesis,TiDu13,TiDu14,TiDu19}. {\em Asymptotically compatible} (AC) schemes, motivated by the findings in \citeasnoun{TiDu13} and formally introduced in \citeasnoun{TiDu14}, are numerical discretizations of nonlocal models that converge to nonlocal continuum models for a fixed horizon parameter and to the local discrete schemes as the horizon vanishes for both discrete schemes with a fixed numerical resolution and for continuum models with increasing numerical resolution.

\begin{figure} %%fig4.1
  \vspace{4pt}
  \centering
       \begin{tikzpicture}[scale=1.3]
     \tikzset{to/.style={->,>=stealth',line width=.8pt}}   
 \node(v1) at (0,3.5) {\textcolor{purple}{$u_{\delta,h}$}};
  \node (v2) at (5.5,3.5) {\textcolor{purple}{$u_{0,h}$}};
   \node (v3) at (0,0) {\textcolor{purple}{$u_\delta$}};
    \node (v4) at (5.5,0) {\textcolor{purple}{$u_0$}};
     \draw[to] (v1.east) --  node[midway,below] {\footnotesize{\textcolor{red}{$\delta\to 0$}}}      
     (v2.west);
     \draw[to] (v1.south) -- node[midway,right] {\footnotesize{\textcolor{blue}{$h\to 0$}}} (v3.north);
       \draw[to] (v3.east) -- node[midway,above] {\footnotesize{\textcolor{red}{$\delta\to 0$}}} (v4.west);
       \draw[to] (v2.south) -- node[midway,left] {\footnotesize{\textcolor{blue}{$h\to0$}}}(v4.north);
     \draw[to] (v1.south east) to[out = 2, in = 180, looseness = 1.1] node[midway] {\footnotesize{\textcolor{red}{ $\delta \to 0$}\ \ \ \ \textcolor{blue}{$h\to0\;$}}} (v4.north west);
      \end{tikzpicture}\\
 \parbox{225pt}{\caption{An illustrative diagram for AC schemes.}
      \label{dt-fig-diagram}}
\end{figure}
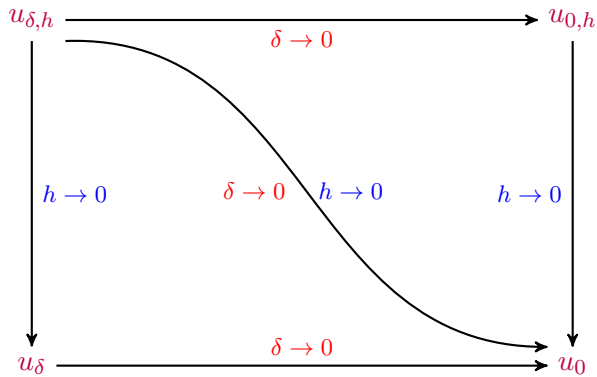

Let $h>0$ denote a {grid size} parameter (or particle spacing) and {let} $\del$ denote  the horizon parameter or even a more generic model parameter. Instead of $h$, we could use the dimension $N$ of a spectral basis. The implications of the AC property are illustrated in \fo{Figure~\ref{dt-fig-diagram}}. There, $u_{\del}$ denotes the solution of the continuous nonlocal problem with $\del>0$, $u_0$ the solution of the corresponding continuous local problem, {$u_{\del,h}$} the solution of the discretized nonlocal problem, and {$u_{0,h}$} the solution of the discretized local problem. \fo{Figure~\ref{dt-fig-diagram}} is meant to illustrate the basic property of AC schemes, namely that for such schemes it does not matter in which order the limits are taken, or even if $\del$ and $h$ are related in some way so that they simultaneously tend to zero.
Note that figures similar to \fo{Figure~\ref{dt-fig-diagram}} can be drawn for $N\to\infty$ and/or $\del\to\infty$.

As seen from detailed studies given in \citeasnoun{TiDu13}, some popular discretization schemes for nonlocal peridynamics fail to be AC. In particular,  if $\del$ is taken to be proportional to $h$, then as $h\to 0$, piecewise constant conforming finite element  solutions actually converge to the incorrect limit, similarly to those based on simple Riemann sum quadrature approximations to nonlocal operators. Similar discussions can be found in \citeasnoun{BYASSX09} {and} \citeasnoun{dt-ChGu11} for diffusion models and \citeasnoun{du16cmame} for systems of peridynamic equations. 

Asymptotically compatible (AC) schemes, {such as} conforming Galerkin-type approximations of weak forms \cite{dt-ChGu11,TiDu14,du16mms}, nonconforming discontinuous Galerkin approximations \cite{djlt19camc,djl19mc,dy19jsc} or {collocation- or quadrature-based} approximations of strong forms \cite{du14asme,du16cmame,dtty19ima,collo16,collo2},
offer the potential to solve for approximations of a model of interest with different choices of parameters to gain efficiency and to avoid the pitfall of reaching inconsistent limits.

\vspace{5pt}
\section{Finite element methods for nonlocal models}\label{femnl}

Given weak formulations of nonlocal models, it is natural to consider finite element approximations \cite{dt-ChGu11,TiDu13,TiDu14,xu2016multiscale,du16mms,TiDu19,du19cbms,jha2019finite}.
To derive a finite element discretization, one possible way is to follow the same recipe as that used for the local PDE setting. We assume that $\omg\subset\Rd$ is {a polytope, so the} first step is to construct a regular subdivision of $\omgomgi$ into finite elements, {\eg\ triangles or quadrilaterals for $\ddd=2$}. Based on the grid, {we then define} a finite element space $\VVVh$, {usually consisting of piecewise polynomial functions with respect to the grid, and then {choose} a basis for that space containing functions whose support extends} over a few contiguous elements. We do not dwell on the construction of finite element spaces because it is the same as that for the local PDE setting \cite{brenner,ciarlet,ErnGuermond:book}, 
  except perhaps that it is prudent to have the grid contain a regular subdivision of $\omg$ itself so that the boundary of $\omg$ is subdivided into $(\ddd-1)$-dimensional faces. 
{For} {\em conforming finite element methods}, \tie, the finite element space $\VVVh$ is a subset of the energy space $\VVV$ for the continuous problem,
we define the constrained finite element space 
\[
\VVch= \{ v(\xb)\in\VVVh \colon  v(\xb) = 0 \ \text{for all } \xb\in\omgid\}.
\]
Note that $\VVch\subset\VVc$. We {let $\gh(\xb)$ denote} an approximation to $g(\xb)$, \feeg, $\gh(\xb)$ could be the interpolant of $g(\xb)$ in $\VVch$ restricted to $\omgid$, or, if $g(\xb)$ is not smooth enough to have pointwise values, a {least-squares} approximation could be used instead.
A conforming finite element approximation $\uh(\xb)\in\VVVh$ of the solution of $u(\xb)\in\VVV$ of \eqref{weak-weakf} is then determined by solving the following problem{:} 
\begin{equation}\label{weak-weakffe}
  \begin{minipage}[b]{22pc}
{Given $\gamd\xyp$ defined in \eqref{str-gamd} and given $f(\xb)\in\VVd$,
find $\uh(\xb)\in\VVVh$ such that}
\[
\begin{cases}
{\mcAd}(\uh,\vh) = \fvh  &\text{for all } v\in\VVch ,
\\
\uh(\xb) = \gh(\xb) &\text{for all } \xb\in\omgid,  
\end{cases}
\]
\end{minipage}
\end{equation}
where $\gh(\xb)$ is defined as discussed above. Because the bilinear form ${\mcAd}(\cdot,\cdot)$ is continuous and coercive with respect to $\VVc$, it is likewise continuous and coercive with respect to $\VVch$, {so the well-posedness of problem} \eqref{weak-weakffe} follows immediately from the Lax--Milgram theorem. {For problems involving singular nonlocal interaction kernels such as the fractional type given in \fo{Table~\ref{ker-prop}}, special numerical quadrature is needed for the evaluation of stiffness matrices. We refer to Section~\ref{sec:quadrature} for related discussions.}

\setcounter{rem}{1}
\medskip\noi{\bf Remark 5.1 ({DG} methods in the nonlocal setting).}
One of the major differences between finite element methods for local PDEs and nonlocal models occurs in the use of discontinuous Galerkin (DG) methods. For elliptic PDEs, DG methods are not conforming, \tie, a finite element space containing functions with jump discontinuities cannot be a subspace of the space $H^1(\omg)$ in which the PDE is {well-posed}. As a result, the use of DG methods requires an accounting for fluxes across element faces; otherwise, elements would be uncoupled. On the other hand, for {some}  kernel choices, finite element functions with jump {discontinuities belong} to the energy space $\VVV$ in which the nonlocal problem is {well-posed}. This is the case for \bo{kernels that are both radial and integrable and for} {the fractional kernel with $s<\gfrac12$, as}
listed in \fo{Table~\ref{ker-prop}}.
For such {kernels 
DG methods are conforming,} so that, at least for $\del> h$, no explicit accounting for fluxes across element faces is needed. The coupling between elements is taken care of by nonlocality, \tie, even though the basis functions live only on single elements, any specific element is coupled to all other elements that \bo{contain points that are within a distance $\del$ to some point in the specific element.}  \bo{For more singular kernels where DG methods are not conforming, \eg\ the fractional kernels with $s>1/2$, we refer to the discussions in Section~\ref{nonconformingDG}.}

\subsection{{Asymptotically compatible conforming} finite element methods for {nonlocal diffusion}}\label{hskernel}

To delve deeper into finite element methods for nonlocal diffusion {with a finite range of interactions}, we {first} specialize to a more specific setting \bo{for which the local  diffusion limit has been rigorously established}. Consider the model given by \eqref{weak-weakf} with $g=0$. We also choose the kernel $\gamd\xyp$ as a rescaled translation-invariant kernel, \tie, 
\begin{equation}\label{eq:rescale}
 \gamd\xyp  =\ffrac{1}{\del^{d+2}}\gamma\biggl(\ffrac{|\ymx|}{\del}\biggr)
\end{equation}
for $\del\in(0,1]$ and some kernel $\gamma( {|\cdot |})$ that is {assumed to be radial}, non\-negative, compactly supported in $B_1({\bf 0})$ (the unit ball centred at the origin),  and has a bounded second-order moment, \tie,
\begin{equation}\label{eq:kernel}
\widehat \gamma (|\bm\xi |)=|\bm\xi|^2\gamma(|\bm\xi|)\in L^1_{\text{loc}}(\Rd)\quad \text{and}\quad  \int_{{B_1({\bf 0})} } \widehat\gamma(|\bm\xi|) \D \bm\xi=d.
\end{equation}
For kernels satisfying these conditions, we refer to the discussions in Section~\ref{sec:mod-kernels}. For such kernels,
\bo{the local $\delta \to  0$ limit becomes the homogeneous Dirichlet boundary value  problem of the Poisson equation. Let us}  denote the energy space $\VVcd=\VVc$ as defined in \eqref{weak-espacec} {to highlight the dependence on $\del$}.

Now, for any fixed $\del\in(0,1]$, we introduce {conforming}  finite element spaces \bo{$\{{V_{\del, h}} \} \subset  \VVcd$}
associated with the
triangulation $\tau_h=\{ K \}$ of the domain {$\omgomgi$} (or {$\omg\cup\omg_{\mcI_1}$ {that contains the domain $\omgomgi$ for any $\del<1$).} 
We set
\[
V_{\del, h} := \{ v\in  \VVcd \colon v|_K \in P_p(K) \ \text{for all}\ K \in \tau_h \},
\]
where $P_p(K)$ denotes the space of  polynomials on $K\in\tau_h$ of degree less than or equal to $p$. Again, for different $\del$, in order to have the finite element functions defined
on a common spatial domain,  we also assume, {as for} $ \VVcd$, that any element in $V_{\del, h}$
vanishes outside $\omg$.

 As $h\to0$,  {we assume that} \bo{$\{{V_{\del, h}} \}$} is dense in $ \VVcd$, \tie, for any $v\in  \VVcd$, there exists a sequence
\bo{$\{v_h\in {V_{\del, h}}\}$} such that, for a given $\del>0$,
\begin{equation}\label{FEsq}
\| v_h-v\|_{\VVcd}\to0\quad\text{as } h\to 0.
\end{equation}
These properties are easily satisfied by standard finite element spaces.

The Galerkin approximation is defined by replacing $\VVcd$ {with} $V_{\del, h}$ in \eqref{weak-weakf}: 
\begin{equation}\label{eq:wfscalarapprox}
\text{{Find} $u_{h,\del}\in {V_{\del, h}}$   such that 
$\mcA_\del(u_{h,\del}, v_h) = \fvh$  for all $ v_h\in  {V_{\del, h}}$.}
\end{equation}

{The analysis of Galerkin approximations of parametrized nonlocal variational problems can be formulated within the general framework of AC schemes \cite{TiDu14}. Other than the necessary properties of functions and operators that guarantee $u_\del\to u_0$ as $\del\to 0$, where $u_0$ denotes the solution of the corresponding PDE, and standard density properties such as \eqref{FEsq}, one of the key ingredients for the asymptotic compatibility of the scheme is the {\em asymptotic density} property of the finite element space. Rather than a general definition, we provide a specific instance, adapted to the models studied here}, {as follows.}
A family of finite-dimensional spaces
$\{ V_{\del,h}\subset  \VVcd$, $\del\in (0,1)$, $h\in(0, h_0]\}$
is {\em asymptotically dense} in $H_0^1(\Omega)$, {if}, {for all $ v\in H_0^1(\Omega)$},
there exists a sequence $\{v_n\in V_{\del_n, h_n}\}_{h_n\to0, \del_n\to 0}$ as $n\to\infty$  such that
\begin{equation}\label{eq:asdense}
\| v-v_n\|_{H^1(\omg)}\to 0 \quad \text{as }  n\to\infty.
\end{equation}

{\citeasnoun{TiDu14} have shown rigorously} that, for scalar nonlocal diffusion equations such as \eqref{weak-weakf}, all conforming Galerkin approximations of the nonlocal models containing continuous piecewise linear functions are automatically AC {in any space dimension}. This means that they can recover the correct local limit as long as both $\del$ and $h$ are {decreasing, even if the nonlocal parameter} $\del$ is reduced at a much faster pace than the mesh spacing $h$. {Even though} the analysis of the above AC property is highly technical, an intuitive explanation is that even with $h$ larger than $\del$, {the nonlocal features that ensure the correct local limit are still encoded in the stiffness matrices} thanks to higher-order (than constant) basis functions. In fact, even though, for the class of kernels we are considering, discontinuous piecewise constant finite element spaces are conforming, they do not generally result in AC approximations. This was first noticed numerically {by} \citeasnoun{dt-ChGu11}.

For general nonlocal systems in any dimension,
{\citename{TiDu14}  \citeyear{TiDu14,TiDu19} have shown}  that as long as the  condition $h=o(\del)$ is met as $\del\to0$, then {we obtain} the correct local limit even for discontinuous piecewise constant finite element approximations when they are of the conforming type. Practically speaking, this implies that a mild growth of the bandwidth in the finite element stiffness matrix is needed as the mesh is refined in order to recover the correct local limit for piecewise constant finite element
schemes. In fact,
{\citename{TiDu13} \citeyear{TiDu13,du15lec} have shown}
  that if a constant bandwidth is kept as the mesh is refined, the approximations may converge to an incorrect local limit.
 
Naturally, schemes using higher-order basis functions tend to provide higher-order accuracy in the nonlocal setting as well,
{should the solutions enjoy sufficient regularity}. At the moment, the theory on AC schemes does not offer any estimate of the order of convergence with respect to different couplings of $h$ and $\del$. Preliminary numerical experiments in \citeasnoun{TiDu14} offer some insight about the balance of modelling and discretization errors, but additional theoretical analyses need to be carried out, except for the case of Fourier spectral approximations of nonlocal models with periodic boundary conditions, for which precise error estimates can be found in {\citename{dy16sinum} \citeyear{dy16sinum,DY:2017} and} \citeasnoun{smd18jcp}.

\subsection{Nonconforming and DG FEMs for nonlocal models with sufficiently singular kernels} \label{nonconformingDG}

The framework of AC schemes is very general. For example, the theory also guided the development of nonconforming discontinuous {Galerkin approximations} \cite{td15sinum} for nonlocal diffusion and nonlocal 
peridynamic models with sufficiently singular interaction kernels.
Consider the scalar model given by \eqref{weak-weakf} with $g= 0$. If the nonlocal interaction kernel satisfies  
\begin{equation}\label{eq:kernel2}
\int_{|\xb|<\epsilon} |\xb|\gamma_\delta(|\xb|)  \D \xb=\infty\quad\text{for all } \epsilon\in (0,\epsilon_0],
\end{equation}
then it may be the case that discontinuous finite element {solutions do not belong to the associated energy space and alternative formulations have to be developed.}
For technical reasons,  {\citeasnoun{td15sinum} also} assumed  that
\begin{equation}\label{eq:kernel21}
\lim_{\epsilon\to 0} \epsilon^2 \biggl( \int_{|\xb|<\epsilon} |\xb|^2\gamma_\delta(|\xb|)  \D \xb \biggr)^{-1} = 0.
\end{equation}
Note that for kernels $\gamma_\delta(r)$ that behave like $1/r^{d+2s}$ as $r\to 0$ with $s\in (1/2,1)$, conditions \eqref{eq:kernel2} and \eqref{eq:kernel21} are satisfied.

{\citeasnoun{td15sinum} introduced} a nonconforming DG scheme based on the removal of the singularity in the nonlocal interaction kernel in a sufficiently small neighbourhood of the origin parametrized by a cut-off level $n$. In other words, $\gamma_\delta(r)$ in \eqref{weak-weakf} is replaced {with}
\begin{equation}\label{modifiedkernel}
\gamma_\delta^n (r) :=
\begin{cases}
\gamma_\delta(r)  & \text{if } \gamma_\delta(r)\leq n ,\\ 
n  & \text{if } \gamma_\delta(r)> n.
 \end{cases}
\end{equation}
 For nonlocal problems with the regularized kernel $\gamma_\delta^n (r)$, discontinuous element spaces, such as conventional discontinuous finite element spaces, can be used as conforming Galerkin approximations, leading to a discrete solution $u_{h,\delta,n}$. Naturally, there are other ways to define the cut-off. The essential requirement is that the resulting modified kernel is \bo{both radial and} integrable (for a given $n$) and it converges pointwise to the original kernel.

The goal is then to demonstrate that $u_{h,\delta,n}$
converges to the solution $u_\delta$ of the original continuous nonlocal model with a singular kernel as $n\to \infty$ and $h\to 0$. The convergence theory can be established by generalizing the relevant
compactness results given in \citeasnoun{BBM01} as $n \to \infty$ for a given $\delta>0$. Indeed, the following generalization is made in \citeasnoun{td15sinum}.

\paragraph{Compactness results.}
Given the kernels $\gamma_\delta^n$ and $\gamma_\delta$ and the corresponding energy spaces $\VVcdn$ and $\VVcd$, assume 
that the energy norms of $\{v_n\in\VVcdn\}$ have
the uniform bound
\[
{\text{supp}_n}  \int_\omgomgi\int_\omgomgi \gamma_\delta^n(|\yb-\xb|) (v_n(\yb)-v_n(\xb))^2  \D \yb  \D \xb \leq C_0.
\]
Then $\{v_n\}$ is relatively compact in $L^2(\omgomgi)$ and any limit function $v$ belongs to $\VVcd$ with
\begin{equation}\label{limitbound}
\int_\omgomgi \int_\omgomgi \gamma_\delta(|\yb-\xb|) (v(\yb)-v(\xb))^2 \D \xb  \D \yb \leq C_0 .
\end{equation}

The classical Bourgain--Brezis--Mironescu compactness result established  in \citeasnoun{BBM01} can be seen as the local limit of the new compactness result for nonlocal spaces in \citeasnoun{td15sinum} with a  finite $\delta$. 
{By applying} the framework of asymptotically compatible discretization \cite{TiDu14,TiDu19}, reviewed earlier
 with respect to the
mesh parameter $h$ and the cut-off level $n$, we can obtain the convergence of $u_{h,\delta,n}$ to $u_\delta$ unconditionally
as $n\to \infty$ and $h\to 0$ if the underlying finite element space contains the continuous piecewise linear finite element space.
Moreover, {in this case, if $n\to \infty$, we expect that for any given $\delta$ and $h$,  $u_{h,\delta,n} \to u_{h,\del}$ as
$n\to \infty$, where $u_{h,\del}$ denotes the conforming finite element approximation of $u_\del$ in the space ${V_{\delta,h}}\cap \VVcd$.}

If piecewise constant finite elements are used, a conditional convergence theorem has been established in \citeasnoun{td15sinum}, {provided the definition of the cut-off is suitably modified.}  For instance, consider a kernel of the type
\begin{equation}\label{frackernel}
\ffrac{\gamma_\star}{|\yb-\xb|^{d+2s}}\leq \gamma_\delta(|\xb -\yb|)\leq\ffrac{\gamma^\star}{|\yb-\xb|^{d+2s}} \quad \text{ for } \xb,\, \yb\in\omgomgi ,
\end{equation}
for some $s\in (1/2,1)$ and positive constants $\gamma_\star$ and $\gamma^\star$.
Then one may let
\begin{equation}\label{modifiedkernel2}
\gamma_\delta^n(r) :=
\begin{cases}
\gamma_\delta(r) & \text{for } r\geq1/n ,\\ 
\gamma_\delta(1/n)  & \text{for } 0<r<1/n.
\end{cases}
\end{equation}
For the approximate solution $u_{h,\delta,n}$ defined as the Galerkin approximation to the nonlocal problem with kernel $\gamma_\delta^n$ using piecewise constants, we have that $\| u_{h,\del,n}-u_\del \|_{L^2}\to0$ if $h=o(1/n)$ as $n\to\infty$.

We note that the nonconforming approximations discussed in \citeasnoun{td15sinum} in the local limit do not yield a standard nonconforming finite element approximation {or} a DG approximation to the local problem. One may construct other alternative formulations that
can give rise to the conventional nonconforming and DG discretization of local PDEs{;} see \eg\ \citeasnoun{djlt19camc} for a study based on the DG with penalty formulation.

\subsection{Adaptive mesh refinement for nonlocal models}

Due {to reduced sparsity, \fieg,} nonlocal models generally incur greater computational costs {than their} local PDE-based counterparts. Thus, designing effective adaptive methodologies is important, and {it} is an area worthy {of attention. So far}  the work
has been limited to nonlocal models with a finite range of interactions.

{\citeasnoun{du13mc} provided} an {\it a~posteriori} error analysis of conforming finite element methods for solving linear nonlocal diffusion and peridynamic models. The approach adopted is a {residual-type} error estimator in the {$L^2$-norm}, {\eg\ of the form $\|{-}\mathcal{L}_\delta(u_{\del,h})-f\|_{L^2}$, which} remains well-defined and
 can be easily computed for \bo{kernel functions that are both radial and integrable} from  element-wise contributions without worrying about flux jumps across element boundaries. This is in sharp contrast to the case of second-order elliptic PDEs. The theory of {\it a~posteriori} error analysis has been rigorously derived for nonlocal {volume-constrained} problems associated with scalar equations. The reliability and  efficiency of the estimators are {proved}, and relationships between nonlocal and classical local {\it a~posteriori} error estimates are also studied. 

{\citeasnoun{du13sinumb} have also developed} a convergent adaptive finite element algorithm for the numerical solution of scalar nonlocal models. For problems involving certain \bo{radial but} non-integrable kernel functions, the convergence of the adaptive algorithm is rigorously derived with the help of several basic ingredients, such as an upper bound on the estimator, the estimator reduction, and the orthogonality property. How these estimators and methods work in the local limit and for general time-dependent and nonlinear peridynamic models remains to be investigated.

{For} nonlocal problems having solutions with jump discontinuities, 
an adaptive finite element method  is given in \citeasnoun{xu2016multiscale}. There, an algorithm is developed that first detects the location of the discontinuity and then refines the grid near the discontinuity. To preserve the $h^2$ accuracy possible with the use of {piecewise linear} elements even when the exact solution contains jump discontinuities at unknown locations, the elements surrounding the discontinuity should have thickness of $\Bigoh(h^4)$ across the discontinuity. This was already observed for the one-dimensional case in \citeasnoun{dt-ChGu11} {and} \citeasnoun{du16mms}. In higher dimensions, a naive refinement strategy that results in small, well-shaped elements in the vicinity of a $(\ddd-1)$-dimensional surface can thus result in an excessive number of degrees of freedom. The adaptive refinement strategy of \cite{xu2016multiscale} instead results in elongated elements having thickness $\Bigoh (h^4)$ across the {discontinuity but} length $\Bigoh (h)$ along the discontinuity, as illustrated in \fo{Figure~\ref{fig:stretch}}. The presence of elongated elements is not harmful to the error because of the anisotropic behaviour of the solution, \tie, it is smooth along the discontinuity but discontinuous across the discontinuity. {Whereas robust meshing algorithms for this {type} of anisotropic refinement and accurate predictions of solution jump discontinuities remain computationally challenging in \bo{higher space}  dimensions,} the numerical examples in \citeasnoun{xu2016multiscale} illustrate that the adaptive strategy developed there does indeed result in $h^2$ convergence, as is also the case for the one-dimensional numerical results in \citeasnoun{dt-ChGu11} {and} \citeasnoun{du16mms}.

\begin{figure} %%fig5.1
\centering
\includegraphics[width=120pt,viewport=0 0 362 382,clip]{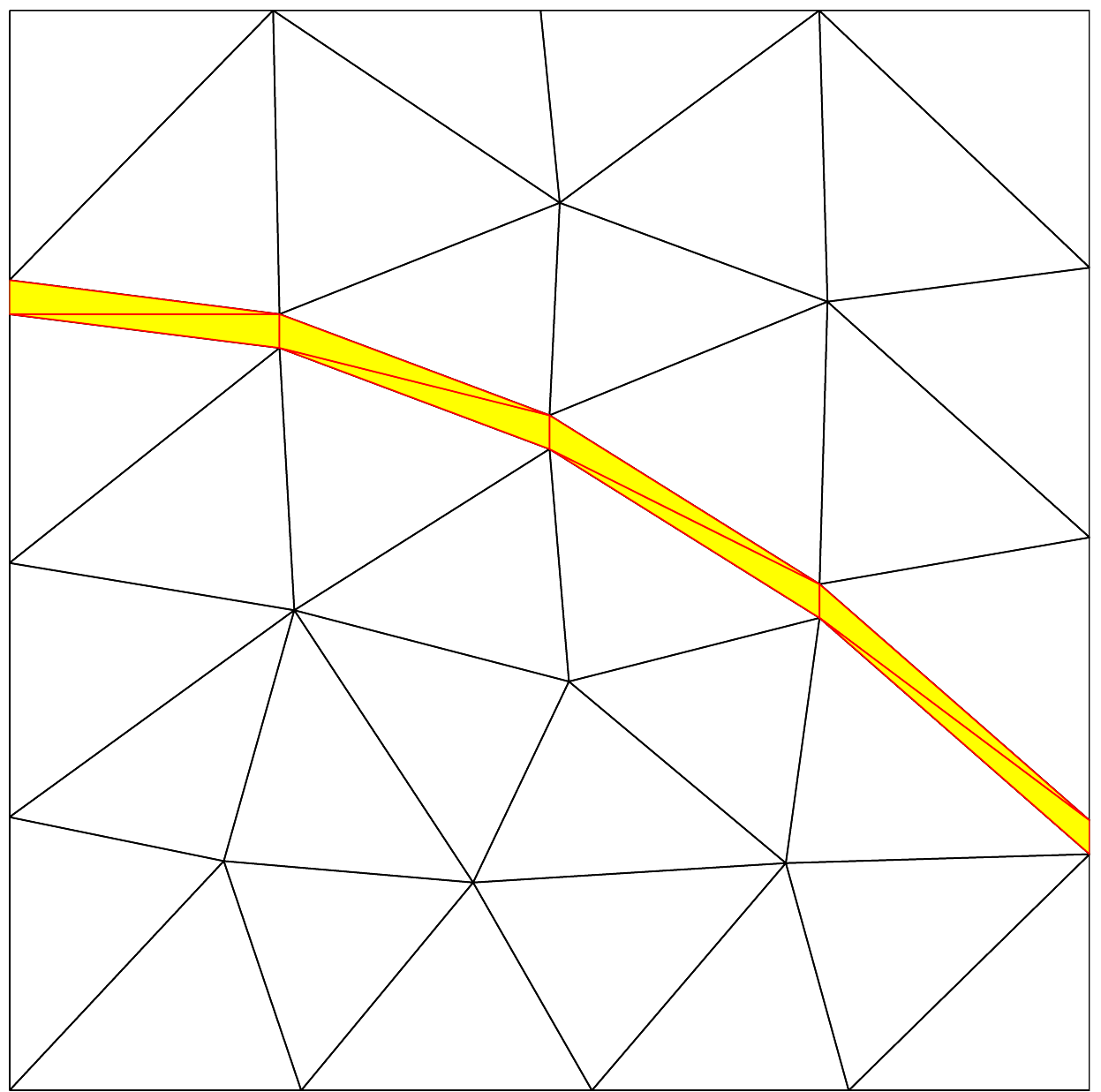}\\ %%cvt10
\caption{{A} coarse-grid illustration of a mesh resulting from the adaptive strategy discussed in {this} section. The elongated yellow elements surround the discontinuity in the solution.}\label{fig:stretch}
\end{figure}
 
Another useful feature of the adaptive strategy developed in \citeasnoun{xu2016multiscale} is that the transition between the elongated elements along discontinuities and well-shaped elements with sides of length $\Bigoh (h)$ away from the discontinuity can be abrupt, \tie, there is no need to have a transition zone in which elements gradually grow in size. This is also illustrated in \fo{Figure~\ref{fig:stretch}}. This feature is already present in the methods developed in \citeasnoun{dt-ChGu11} {and} \citeasnoun{du16mms} for the one-dimensional case. The importance of this feature is that it also helps to keep down the number of degrees of freedom.

The next means of savings in the number of degrees of freedom is to use discontinuous basis functions only in the elongated elements, and continuous basis functions in elements in which the exact solution is smooth. The final means is to switch from the nonlocal model to the corresponding PDE model in all elements that do not interact with the elongated elements that surround the discontinuity, \ie\ all elements whose vertices are at distance greater than $\del$ from the vertices of the elongated elements.

The four degrees of freedom-saving features (\ie\ elongated elements, abrupt grid transitions, the use of discontinuous basis function only in the elongated elements, and the switch to a local PDE model away from the discontinuity), result in tremendous savings in the number of degrees of freedom. In fact, the resulting number of degrees of freedom is comparable to that for a local PDE model using a regular grid of elements having sides of $\Bigoh (h)$. Moreover, the relative savings become greater as the grid size is reduced and as the dimension $\ddd$ is increased. See \citeasnoun{xu2016multiscale} for several numerical illustrations supporting these conclusions.

{With regard to more general coupling strategies of local and nonlocal models for computational effectiveness,  relevant works are discussed in a recent review on the subject \cite{du19cbms}.
AC schemes have also been studied in the context of multiscale problems modelled by local--nonlocal coupling formulations with a spatially heterogeneous horizon \cite{dt18iccm,ttd19mms}, based on a new trace theorem for nonlocal function spaces that provides a stronger version of the classical local counterparts \cite{td17sima}.}
}

\subsection{Approximations of fractional models as limit of nonlocal models with a finite range of interactions}

One of the key messages we want to get across is to show that many popularly studied fractional PDEs are either specialized nonlocal models or can be treated as limiting cases.

Naturally, for fractional differential operators defined in integral form, by truncating the fractional kernel in the fractional operators to a finite range measured by the horizon parameter $\delta$, we end up with a nonlocal model parametrized by $\delta$. Hence we see from Section~\ref{truncted-flap} that on the continuum level, the fractional models may then be viewed as the infinite horizon limit of  nonlocal models with a finite $\delta$, after a suitable scaling of the kernel.

\begin{figure} %%fig5.2
  \vspace{5pt}
\centering
\begin{tikzpicture}[thick,scale=1, every node/.style={scale=1.}]
       \node[draw, circle, minimum width=22mm] at ({0}:0cm) {};
         \node(v3) at (0,-0.6)  {\fns{$\Omega$}}; 
       \node(v1) at (0.,-1.4) {\textcolor{red}{\fns{$\partial \Omega$}}};
    \node(v2) at (0, 0) {\textcolor{red}{\fns{$-\mathcal{L}_{0,h} u_{0,h}$}}\fns{=$\,{b}_h$}};
             \node(v4) at (0,0.6)  {\textcolor{blue}{{Local}}}; 
      \node(v5) at (0,1.35) {\textcolor{red}{\fns{$ {{u}_{0,h}}$}}\fns{=$\,{0}$}};
       \node(v1) at (1.8,0.4) {\textcolor{red}{\fns{$0 \leftarrow \delta$}}}; 
  \draw[red, thick, <-]  (1.2,0.) -- (2.4,0); 
      \node[draw, circle, minimum width=34mm] at ({0}:43mm) {};
 \node[draw, circle, minimum width=23mm] at ({0}:43mm) {};
    \node(v3) at (4.3,-0.6) {\fns{$\Omega$}}; 
  \node(v1) at (4.3,-1.4) {\textcolor{red}{\fns{${\Omega_{\mathcal{I}_\delta}}$}}};
    \node(v2) at (4.3,0.) {\textcolor{red}{\fns{$- \mathcal{L}_{\delta,h} {u}_{\delta,h}$}}\fns{=$\,{b}_h$}};
                 \node(v4) at (4.3,0.6)  {\textcolor{blue}{Nonlocal}}; 
        \node(v5) at (4.3,1.4) {\textcolor{red}{\fns{$ {u}_{\delta,h}$}}\fns{=$\,{0}$}};
\end{tikzpicture}\\[2mm]
\begin{tikzpicture}[thick,scale=1, every node/.style={scale=1.}]
      \node[draw, circle, minimum width=39mm] at ({0}:0cm) {};
 \node[draw, circle, minimum width=28mm] at ({0}:0cm) {};
    \node(v3) at (0.0,-0.7) {\fns{$\Omega$}}; 
  \node(v1) at (0.0,-1.7) {\textcolor{red}{\fns{${\Omega_{\mathcal{I}_\delta}}$}}};
    \node(v2) at (0.0,0.) {\fns{$-$}\textcolor{red}{\fns{$\mathcal{L}_{\delta,h} {u}_{\delta,h}$}}\fns{=$\,{b}_h$}};
                 \node(v4) at (-.0,0.6)  {\textcolor{blue}{Nonlocal}}; 
        \node(v5) at (-.,1.65) {\textcolor{red}{\fns{$ {u}_{\delta,h}$}}\fns{=$\,{0}$}};
 \node(v1) at (2.7,0.4) {\textcolor{red}{\fns{$\delta \rightarrow \infty$}}}; 
  \draw[red, thick, ->]  (2.2,0.0) -- (3.2,0); 
    \node[draw, circle, minimum width=30mm] at ({0}:49mm) {};
         \node(v3) at (4.9,-0.7)  {\fns{$\Omega$}}; 
                \node(v1) at (5.0,-1.8) {\textcolor{red}{\fns{$\omgii$}}};
        \node(v2) at (4.9,0) {\fns{$-$}\textcolor{red}{\fns{$\mathcal{L}_{\infty,h} {u}_{\infty,h}$}}\fns{={$\,{b}_h$}
        }};
                \node(v4) at (4.9,0.6)  {\textcolor{blue}{{Fractional}}}; 
      \node(v5) at (4.9,1.75) {\textcolor{red}{\fns{$ {{u}_{\infty,h}}$}}\fns{=$\,{0}$} };
\end{tikzpicture}\\
\caption{\bo{Limits of numerical solutions of nonlocal models with finite interaction radius $\delta$ as $\delta\to 0$ (local limit, top) and as $\delta\to\infty$ (fractional limit, bottom).}}
\label{fig:diagram3}
\end{figure}
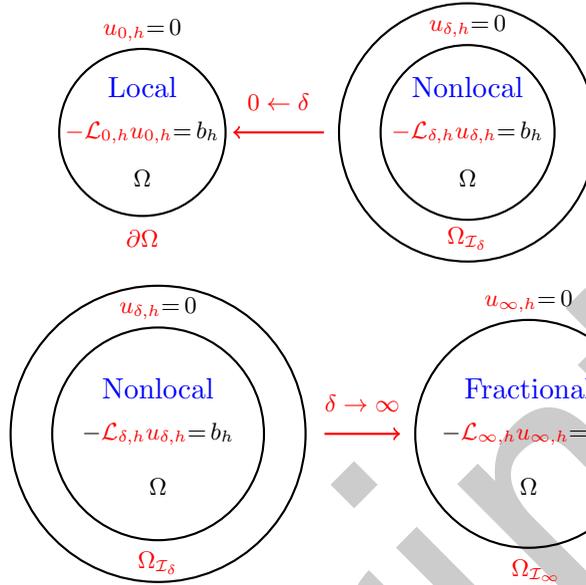

{\citeasnoun{MG:2013} and \citeasnoun{burko1} have shown} that as $\delta\to \infty$, the Galerkin approximations of the nonlocal model with the parameter $\delta$ can converge to the solution of the fractional equations, provided $h$ is taken to be sufficiently small. 
One can actually apply 
the framework of the AC schemes 
to show that all conforming Galerkin approximations are AC in the limit $\delta\to 0$ \cite{tdg16acm}. That is, the discretization of fractional equations associated with the fractional Laplacian  can be viewed as the global (with an infinite nonlocal horizon parameter) limit of nonlocal models with properly scaled {fractional-type} nonlocal interaction kernels  \cite{tdg16acm}. Moreover, the convergence does not require the dependence of $h$ on $\delta$.

Let us give an illustration of the work presented in \citeasnoun{tdg16acm}. Consider the {volume-constrained} problems defined for the fractional Laplacian in $\Omega$ with homogeneous Dirichlet condition in $\omgii=\Rd\setminus\omg$.
By truncations of both the spatial domain and the range of nonlocal interactions, we may
end up with a class of parametrized problems
\begin{equation}\label{chap4:fracdiff_ap}
\begin{cases}
-\mcL_\del u(\xb)=-\int_{\Rd}(u(\yb)-u(\xb))\gamma_\del(|\yb-\xb|)  \D \yb=f &\text{for all } \xb\in \Omega ,\\
u=0 &\text{for all } \xb\in{\omgid}
\end{cases}
\end{equation}
with
\begin{equation}\label{chap4:fractionalkernel}
  \gamma_\del(|\yb-\xb|)=
  \begin{cases}
\ffrac{{C_{d,s,\delta}}}{
|\yb-\xb|^{d+2s} } & \yb\in B_\del(\xb) , \\[9pt]
0 & \yb\in \Rd \backslash B_\del(\xb)  .
\end{cases}
\end{equation}
Here,  \bo{ $C_{d,s,\delta}=\cds$ is given as in  \eqref{cds} for the limiting case of $\delta\to \infty$ so that $\mcL_\infty = (-\Delta)^{s}$  as in \eqref{intfl}.  Note that in the limit of $\delta\to 0$, we take a $C_{d,s,\delta}$ to be a different scaling factor so that its local limit is  $\mcL_0 =-\Delta$.}

Let $h$ denote the discretization parameter associated with
Galerkin approximations of \eqref{chap4:fracdiff_ap}.
\pagebreak %%20200331
For example, $h$ could be the mesh parameter for the finite element discretization or the reciprocal of the number of 
basis functions in spectral approximations.
We can then apply the AC framework   to obtain the convergence of the Galerkin approximations of \eqref{chap4:fracdiff_ap} 
to the solution of the fractional equation as 
$h\to 0$ and $\del\to \infty$ as long as the finite-dimensional approximations spaces are subspaces of the underlying energy space, \tie, as long as we adopt conforming discretizations.  

As pointed out in \citeasnoun{tdg16acm} {and} \citeasnoun{TiDu19}, analyses of the fractional and local limits ($\del\to \infty$ and $\del\to 0$ respectively) of nonlocal models parametrized by a finite $\delta$ and with fractional-type kernels demonstrate that nonlocal models are more general than their fractional and local counterparts and they also serve as a bridge between fractional and local models{;} see \fo{Figure~\ref{fig:diagram3}} and, again, \fo{Figure~\ref{dt-fig-diagram}} and similar diagrams given in \citeasnoun{du19cbms} {and} \citeasnoun{TiDu19}.
For AC schemes,  the bridging roles of nonlocal diffusion with a finite range of nonlocal interactions, {presented above} at the continuum level,
remain valid on the discrete level. 

%% fig5.2 original position

\vspace{5pt}
\section{Finite element methods for the integral fractional Laplacian}\label{sec:femintlap}

 \bo{Let} $\Omega$ denote a polygon \bo{and} let $\mcTh$ denote a family of shape-regular and locally quasi-uniform triangulations of $\Omega$ \cite{ciarlet,brenner,ErnGuermond:book}.
Let $\mathcal{Z}_{h}$ denote the set of vertices of $\mcTh$ \bo{and} {let} $h_{K}$ denote the diameter of the element $K\in\mcTh$.
Moreover, let
\begin{align*}
  h:=\max_{K\in\mcTh}h_{K},\quad 
  h_{\min}:=\min_{K\in\mcTh}h_{K}.
\end{align*}
Let $\varphi_{i}$ denote the usual piecewise linear Lagrange basis function associated with a node ${\zb}_{i}\in\mathcal{Z}_{h}$, satisfying $\varphi_{i} ({\zb}_{j}) =\delta_{ij}$ for ${\zb}_{j}\in\mathcal{Z}_{h}$, and let $X_{h}:=\operatorname{span}\{\varphi_{i} \colon  {\zb}_{i}\in\mathcal{Z}_{h}\}$.
The finite element subspace $V_{h}\subset H_c^{s} (\Omega) $ is given by $V_{h}=X_{h}$ when $s<1/2$ and by
\begin{align*}
  V_{h} = \{v_{h}\in X_{h}  \colon  v_{h}=0 \text{ on }\partial\Omega\} = \operatorname{span}\{\varphi_{i} \colon  {\zb}_{i}\not\in \partial\Omega\}
\end{align*}
when $s\geq 1/2$.
The corresponding cardinality ${N}_h$ of $V_h$ is equal to the number of nodes in $\mathcal{Z}_{h}$ when $s<1/2$ and {is otherwise} the number of interior nodes in $\omg$.

Now, let $u_{h}$ denote the solution of the finite element discrete problem {given} by 
\[
  \text{{Find} $ u_{h}\in V_{h}$  such that \bo{$ {\mcA_s}(u,v) =\langle f,v\rangle $} for all $ v\in V_{h}$,}
\]
where the bilinear form \bo{$\mcA_s(\cdot,\cdot)$ is defined in \eqref{aaaaaa} as a special case of the nonlocal bilinear form \eqref{weak-bform} with $\del=\infty$ and kernel function of the type \eqref{fraker} and \roo{the} right-hand side is defined in the usual manner.} 
The following error estimates are derived in \citeasnoun{AcostaBorthagaray2017_FractionalLaplaceEquation}.

\bo{If the family of triangulations $\mcTh$ is globally quasi-uniform, and if $u\in H^{t} (\Omega) $ for $t\in (1/2,2]$ and if $0<s\leq t$, then
  \begin{align}
    \norm{u-u_{h}}_{H_c^{s} (\Omega) }
    &\leq C h^{t-s}\abs{u}_{H^{t} (\Omega)}.\label{eq:uniformRate}
  \end{align}
The estimate \eqref{eq:uniformRate} implies that for $u\in H^{2} (\Omega) $, the expected rate of convergence on a globally quasi-uniform mesh is $h^{2-s}=\Bigoh(N_h^{(s-2)/d})$.
Because the solutions of \eqref{eq:fracPoisson} generally have limited regularity up to the boundary (see Section~\ref{fracreg}), {it is advisable to use} a mesh that is more refined close to the boundary.
While using such a mesh can restore the optimal rate of convergence $\Bigoh(N_h^{(s-2)/d})$ in $\ddd=1$ dimensions,} this does not hold for higher-dimensional problems, where {shape-regularity} of the elements becomes the limiting factor.
Therefore we can expect no more than $\Bigoh(N_h^{-1/2+\varepsilon})$ 
 rate of convergence in $\ddd=2$ dimensions.

Whereas it may be possible to construct an appropriately graded mesh for simple geometries, the general case will require adaptively refined meshes to resolve the boundary singularity.
{\it A~posteriori} error estimators for the integral fractional Laplacian have been developed in \citeasnoun{NochettoPetersdorffEtAl2010_PosterioriErrorAnalysisClass}, \citeasnoun{FaustmannMelenkEtAl2019_OptimalAdaptivityPreconditioningFractionalLaplacian} {and} \citeasnoun{AG:2017}.

\subsection{Quadrature rules}
\label{sec:quadrature}

The fractional Poisson equation \eqref{eq:fracPoisson} leads to a dense linear algebraic system
\begin{align}
  \mat{A}_{s}\vec{u}=\vec{f}, \label{eq:systemPoisson}
\end{align}
in which the entries in the matrix $\mat{A}_{s}=\{\mathcal{A}_{s} (\varphi_{i},\varphi_{j}) \}_{ij}$ involve singular integrals.
In order to compute these entries, we decompose the integrals into contributions between pairs of elements $K,\widehat{K}\in\mcTh$ and between pairs consisting of elements $K\in\mcTh$ and external edges $e\in\partial\mcTh$ as follows:
\begin{align*}
  \mathcal{A}_{s}(\varphi_{i},\varphi_{i}) &= \sum_{K}\sum_{\widehat{K}} \mathcal{A}_{s}^{K\times\widehat{K}}(\varphi_{i},\varphi_{j}) + \sum_{K}\sum_{e} \mathcal{A}_{s}^{K\times e}(\varphi_{i},\varphi_{j}).
\end{align*}
The individual contributions $\mathcal{A}_{s}^{K\times\widehat{K}}$ and $\mathcal{A}_{s}^{K\times e}$ are given {by}
\begin{align*}
  \mathcal{A}_{s}^{K\times\widehat{K}}(\varphi_{i},\varphi_{j})
  & = \ffrac{C(d,s)}{2} \int_{K}  \D \xb \int_{\widehat{K}}  \D \xb \ffrac{ (\varphi_{i}(\xb)-\varphi_{i}(\yb))  (\varphi_{j}(\xb)-\varphi_{j}(\yb)) }{\abs{\xb-\yb}^{\ddd+2s}} ,
    \\
    \mathcal{A}_{s}^{K\times e}(\varphi_{i},\varphi_{j})
  &= \ffrac{C(d,s)}{2s} \int_{K}  \D \xb \int_{e}  \D \xb \ffrac{\varphi_{i} (\xb)  \varphi_{j} (\xb)  ~ \vec{n}_{e}\cdot (\xb-\yb) }{\abs{\xb-\yb}^{\ddd+2s}}.
  \end{align*}

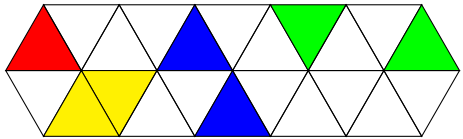
\begin{figure} %%fig6.1
  \vspace{5pt}
  \centering
  \begin{tikzpicture}
    \coordinate (node0) at (0,0);
    \coordinate (node1) at (1,0);
    \coordinate (node2) at (0.5,0.866);
    \coordinate (node3) at (1.5,0.866);
    \coordinate (node4) at (2,0);
    \coordinate (node5) at (0.5,-0.866);
    \coordinate (node6) at (1.5,-0.866);
    \coordinate (node7) at (2.5,0.866);
    \coordinate (node8) at (2.5,-0.866);
    \coordinate (node9) at (3,0);
    \coordinate (node10) at (3.5,0.866);
    \coordinate (node11) at (3.5,-0.866);
    \coordinate (node12) at (4,0);
    \coordinate (node13) at (4.5,0.866);
    \coordinate (node14) at (4.5,-0.866);
    \coordinate (node15) at (5,0);
    \coordinate (node16) at (5.5,0.866);
    \coordinate (node17) at (5.5,-0.866);
    \coordinate (node18) at (6,0);

    \draw[fill=red] (node0) -- (node1) -- (node2) -- cycle;
    \draw (node1) -- (node3) -- (node2) -- cycle;
    \draw (node1) -- (node4) -- (node3) -- cycle;
    \draw (node0) -- (node1) -- (node5) -- cycle;
    \draw[fill=yellow] (node1) -- (node4) -- (node6) -- cycle;
    \draw[fill=yellow] (node1) -- (node5) -- (node6) -- cycle;

    \draw (node4) -- (node3) -- (node7) -- cycle;
    \draw (node4) -- (node6) -- (node8) -- cycle;
    \draw[fill=blue] (node4) -- (node9) -- (node7) -- cycle;
    \draw (node4) -- (node9) -- (node8) -- cycle;

    \draw (node9) -- (node7) -- (node10) -- cycle;
    \draw (node9)[fill=blue] -- (node8) -- (node11) -- cycle;
    \draw (node9) -- (node12) -- (node10) -- cycle;
    \draw (node9) -- (node12) -- (node11) -- cycle;

    \draw[fill=green] (node12) -- (node10) -- (node13) -- cycle;
    \draw (node12) -- (node11) -- (node14) -- cycle;
    \draw (node12) -- (node15) -- (node13) -- cycle;
    \draw (node12) -- (node15) -- (node14) -- cycle;

    \draw (node15) -- (node13) -- (node16) -- cycle;
    \draw (node15) -- (node14) -- (node17) -- cycle;
    \draw[fill=green] (node15) -- (node18) -- (node16) -- cycle;
    \draw (node15) -- (node18) -- (node17) -- cycle;
  \end{tikzpicture}\\
  \caption{Classes of element pairs {of} configurations that need to be handled:
    pairs of identical elements ({red}), element pairs with common edge ({yellow}), with common vertex ({blue}) and separated element pairs ({green}).
  }\label{fig:elementConfigs}
\end{figure}

Contributions from non-disjoint pairs of elements (see \fo{Figure~\ref{fig:elementConfigs}}) are not directly amenable to numerical quadrature, due to their singular nature.
Fortunately, these can be treated by adapting techniques used in the boundary element method (BEM) literature to address similar issues arising from singular kernels \cite{SauterSchwab2010_BoundaryElementMethods}.

However, the fractional Laplacian does pose new difficulties beyond those addressed by the BEM literature, but which can be treated as described {by} \citeasnoun{AinsworthGlusa2018_TowardsEfficientFiniteElement} {and} \citeasnoun{Acosta:shortFEM}.
In particular, {\citename{AinsworthGlusa2018_TowardsEfficientFiniteElement} \citeyear{AG:2017,AinsworthGlusa2018_TowardsEfficientFiniteElement} have developed} non-uniform order Gauss-type quadrature rules.
Alternatively, {one could take the approach proposed by \citeasnoun{ChernovPetersdorffEtAl2011_ExponentialConvergenceHpQuadrature}.} It
allows for transforming quadrature rules given on the unit hypercube $[0,1]^{2d}$ to any pair of elements $K\times\widehat{K}$.
Here, the singularity is taken into account through the choice of the weight in the quadrature rules. 

{\citeasnoun{AG:2017} gave the following result.} Denote the quadrature approximation to the bilinear form $\mathcal{A}_{s} (\cdot,\cdot) $ by $\mathcal{A}_{s}^{Q} (\cdot,\cdot) $.
  \bo{Then, by using \roo{$\Bigoh (\log N_h^{2d})$} quadrature points per element pair $K\times\tilde{K}$, the consistency error due to quadrature $\abs{\mathcal{A}_{s}(u,v)-\mathcal{A}_{s}^{Q}(u,v)}$ is dominated by the discretization error.}

A brief discussion about solvers and condition numbers for finite element discretizations of the integral fractional Laplacian is given in {Section}~\ref{sec:femslcon}.

\vspace{5pt}
\section{{Finite element methods for the spectral fractional Laplacian}}\label{sec:femsl}

The extension problem~\eqref{eq:extensionProblem} has been extensively used to solve equations involving the spectral fractional Laplacian.
For a bounded domain, the extension problem is posed on ${\omg_z}=\Omega\times[0,\infty)$ and takes the form
\begin{align}\label{eq:extensionSpectral}
  \begin{cases}
    -{\nabla\cdot} z^{\beta} \nabla U (\xb,z)  = 0 & \text{for all } (\xb,z) \in{\omg_z} , \\[7pt]
    U(\xb,z) = 0 & \text{for all }  (\xb,z)  \in \partial_{L}{\omg_z} , \\[7pt]
    \ffrac{\partial U}{\partial n^{\beta}} (\xb)  = d_{s}f (\xb)  & \text{for all } \xb\in \Omega,
  \end{cases} 
\end{align}
where $d_s$ is the constant given in \eqref{eq:extensionProblem} and where $\partial_{L}{\omg_z}:=\partial\Omega\times[0,\infty)$ denotes the lateral surface of the semi-infinite cylinder.
The solution $u$ to the spectral fractional Poisson problem \eqref{eqn:poisson-spectral} is then recovered by taking the trace on $\Omega$, \ie\ $u=\operatorname{trace}_{\Omega}U$.

We define the solution space $\mathcal{H}^{1}_{\beta} ({\omg_z}) $ on the semi-infinite cylinder ${\omg_z}$ as
\begin{align*}
  \mathcal{H}_{\beta}^{1} ({\omg_z}) &=\{V\in H^{1}_{z^{\beta}} ({\omg_z})   \colon  V=0 \text{ on } \partial_{L}{\omg_z}\},
\end{align*}
with norm $\norm{V}_{\mathcal{H}^{1}_{\beta}}= \abs{V}_{H^{1}_{z^{\beta}}}$.
The weak form of the extension problem is then {as follows:} 
\begin{align}\label{eq:varExtensionSpectral}
  \text{{Find} $ U\in\mathcal{H}_{\beta}^{1} ({\omg_z})  $  such that}
    \int_{{{\omg_z}}}z^{\beta}\nabla U\cdot \nabla V = d_{s}\pair{f}{\operatorname{tr}_{\Omega}V}   
\end{align}
{for all $ V \in\mathcal{H}_{\beta}^{1} ({\omg_z}) $.}

In the literature, the fact that the domain is unbounded in the $z$-direction has been handled in different ways.
One class of methods exploits the fact that truncation of the domain to ${\omg}_{z_{\text{trun}}}:=\Omega\times[0,z_{\text{trun}}]$ can be shown to be exponentially converging in $z_{\text{trun}}$, and then discretizes the domain ${\omg}_{z_{\text{trun}}}$ using finite elements.
A second class uses a hybrid approach and discretizes using finite elements in the $\xb$-direction and a suitable spectral method in the $z$-direction.
We will describe both these methods below.

\subsection{Truncation in the extended direction}

The truncation approach was pioneered {by} \citeasnoun{NochettoOtarolaEtAl2015_PdeApproachToFractional}.
A first result shows that truncation of the semi-infinite cylinder ${\omg_z}$ only leads to an exponentially small error.

  Let $U\in\mathcal{H}^{1}_{\beta}$ denote the solution of \eqref{eq:varExtensionSpectral} and {let $U_{z_{\text{trun}}}$ be} the extension by zero to ${\omg_z}$ of the solution to \eqref{eq:varExtensionSpectral} posed on the truncated domain ${\omg}_{z_{\text{trun}}}$, $z_{\text{trun}}\geq 1$, with a homogeneous Dirichlet condition enforced at $z=z_{\text{trun}}$.
  Then
  \begin{align*}
    \norm{U-U_{z_{\text{trun}}}}_{\mathcal{H}^{1}_{\beta}({\omg_z})} & \leq C \rme^{-\sqrt{\lambda_{0}}z_{\text{trun}}/4} \norm{f}_{H^{-s}(\Omega)},
  \end{align*}
  where $\lambda_{0}$ is the smallest eigenvalue of the integer-order Laplacian  as given by~\eqref{eq:Eig}.

Now, let $\mcTh$ denote a quasi-uniform mesh of $\Omega$ with mesh size $h$ and define a graded mesh $\mathcal{T}_{z_{\text{trun}}}$ of the interval $[0,z_{\text{trun}}]$ by
\begin{align*}
  z_{k}= (k/M) ^{\gamma}z_{\text{trun}}, \quad k=0,\ldots,M,
\end{align*}
with $M\sim h^{-1}$ and $\gamma>3/(1-\beta)$.
Then define a mesh $\mathcal{T}$ on the extended domain ${\omg}_{z_{\text{trun}}}$ as the tensor product of $\mcTh$ and $\mathcal{T}_{z_{\text{trun}}}$, and let ${V}_{h}$ denote the space of piecewise \bo{linear}  finite element functions on $\mathcal{T}$.
Let $N=\operatorname{dim}{V}_{h}$ be the overall number of degrees of freedom and choose the truncation to be $z_{\text{trun}}\sim \abs{\log N}$.
Then the following convergence result holds.

Let $f\in H_c^{1-s}(\Omega)$, $U\in\mathcal{H}^{1}_{\beta}({\omg_z})$ {be} the solution of \eqref{eq:varExtensionSpectral}, and {let $U_{h}\in{V}_{h}$ be} the solution of the restriction of \eqref{eq:varExtensionSpectral} to ${V}_{h}$.
\pagebreak %%20200331
  Moreover, let $u=\operatorname{trace}_{\Omega}U$ denote the solution to the spectral fractional-order Poisson problem \eqref{eqn:poisson-spectral}, and {let $u_{h}=\operatorname{trace}_{\Omega}U_{h}$ be} its discrete approximation.
  Then
  \begin{align*}
    \norm{u-u_{h}}_{H_c^{s}(\Omega)}
    &\leq C\norm{U-U_{h}}_{\mathcal{H}^{1}_{\beta}({\omg_z})}\\[2pt]
    &\leq C\abs{\log h}^{s} h\norm{f}_{H_c^{1-s}(\Omega)} \\[2pt]
    &\leq C\abs{\log N}^{s} N^{-1/(\ddd+1)}\norm{f}_{H_c^{1-s}(\Omega)}.
  \end{align*}

As observed {by} \citeasnoun{NochettoOtarolaEtAl2015_PdeApproachToFractional}, this shows that the result is optimal in terms of {regularity but} sub-optimal in terms of complexity.
This shortcoming can be addressed using sparse grids \cite{BanjaiMelenkEtAl2019_TensorFemSpectralFractionalDiffusion} or $hp$-finite elements \cite{Meidner:2018,BanjaiMelenkEtAl2019_TensorFemSpectralFractionalDiffusion}. \bo{Using the latter, exponential convergence can be obtained, including for cases with low regularity of the domain and incompatibility of the forcing term \cite{BanjaiMelenkEtAl2019_TensorFemSpectralFractionalDiffusion}.}

\subsection{Spectral method in the extended direction}

A different approach was taken {by} \citeasnoun{AinsworthGlusa2018_HybridFiniteElement}.
Instead of truncating the domain, a spectral method is used for the extended direction, based on the observation that the eigenfunctions of the extension problem~\eqref{eq:extensionSpectral} are given by
\begin{align*}
  \varphi_{m}(\xb)\psi_{m} ({z}) , \quad m\in\mathbb{N},
\end{align*}
where
\begin{align}
  \psi_{m} ({z}) 
  &= c_{s} (\lambda_{m}^{1/2}{z}) ^{s} K_{s} (\lambda_{m}^{1/2}{z}) ,   \label{eq:spectralBasis}
\end{align}
{and where} $c_{s}= 2^{1-s} / \Gamma (s) $,  $ (\varphi_{m},~\lambda_{m}) $ are the eigenpairs of~\eqref{eq:Eig} \bo{and $K_{\nu}$ are the modified Bessel functions of the second kind}.

This motivates the choices
\begin{align*}
  V_{h}&=\{v_{h}\in H^{1}_{0} (\Omega)   \colon  {v_{h}}|_{K}\in \mathbb{P}_{k} (K)  \ \text{for all}\ K\in\mathcal{T}_{h}\}, \\[2pt]
  {V}_{h,M}&=\biggl\{V_{h,M}=\sum_{m=0}^{M-1}v_{h,m} (\xb) \widehat{\psi}_{m} ({z})   \colon  v_{h,m}\in V_{h} \biggr\}\subset \mathcal{H}_{\beta}^{1} ({\omg_z}) ,
\end{align*}
where
\begin{align}
  \widehat{\psi}_{m} ({z}) 
  &:= c_{s} (\widehat{\lambda}_{m}^{1/2}{z}) ^{s} K_{s} (\widehat{\lambda}_{m}^{1/2}{z}) 
\end{align}
are defined using suitable approximations $\widehat{\lambda}_{m}$ instead of the true eigenvalues~$\lambda_{m}$.

It turns out that using the asymptotic law for the eigenvalues of the Laplacian and using a coarse finite element
\pagebreak %%20200331
approximation gives sufficiently good approximations $\widehat{\lambda}_{m}$.
Moreover, by decimation, only $\abs{\log h}^{p}$ for some $p>0$ eigenvalue approximations {is} required.
Overall, this results in the following error estimate.

  Let $f\in H_c^{r}(\Omega)$, $r\geq -s$, {let $U\in\mathcal{H}^{1}_{\beta}({\omg_z})$ be} the solution of \eqref{eq:varExtensionSpectral} and {let $U_{h}\in{V}_{h}$ be} the solution of the restriction of \eqref{eq:varExtensionSpectral} to ${V}_{h}$.
  Moreover, let $u=\operatorname{trace}_{\Omega}U$ {be} the solution to the spectral fractional-order Poisson problem \eqref{eqn:poisson-spectral} and {let $u_{h}=\operatorname{trace}_{\Omega}U_{h}$ be} its discrete approximation.
  Then
  \begin{align*}
    \norm{u-u_{h}}_{H_c^{s}(\Omega)}
    &\leq C \norm{U-U_{h,M}}_{\mathcal{H}^{1}_{\beta}({\omg_z})}\\
    & \leq C \abs{f}_{H_c^{r}(\Omega)}h^{\min\{k,r+s\}}\sqrt{\abs{\log h}} \\
    &= C \abs{f}_{H_c^{r}(\Omega)}N^{-\min\{k,r+s\}/\ddd}\abs{\log N}^{p}.
  \end{align*}

\begin{rem}[solution of the linear system]\label{femslls} 
Provided that the discretization of the extension problem \eqref{eq:extensionSpectral} displays a tensor structure, the solution of the arising linear system reduces to a sequence of discretized integer-order reaction--diffusion-type problems.
These problems are readily solved using classical iterative linear solvers such as conjugate gradient and multigrid; see \eg\ \citeasnoun{ChenNochettoEtAl2016_MultilevelMethodsNonuniformlyElliptic} {and} \citeasnoun{AinsworthGlusa2018_HybridFiniteElement}. \hfill\actaqed
\end{rem}

\vspace{5pt}
\section{Spectral-Galerkin methods for nonlocal diffusion}\label{sec:smn}

Compared with its fractional counterparts discussed {in} Section~\ref{sec:smf}, algorithms and numerical analyses of spectral method for nonlocal diffusion have been much less studied.  Here, we simply consider a Fourier spectral-Galerkin method developed
{by \citename{dy16sinum} \citeyear{dy16sinum,DY:2017}}. 

For simplicity, we focus on the one-dimensional periodic nonlocal diffusion problem given by
\begin{equation}\label{nonloc}
-(\mcLd u)(x)=f(x) \quad\text{for all } x\in(-\pi,\pi)
\end{equation}
with periodic boundary condition, where we express the nonlocal operator $\mcLd$ in the equivalent form 
\begin{equation}\label{nondeffff}
(\mcLd u)(x)=  2 \int_{|z|\leq \del}\gamd(z) (u(x+z)-u(x))   \D z,
\end{equation}
where $f(x)\in L^2_{\text{per}}[-\pi,\pi]$, where $L^2_{\text{per}}[-\pi,\pi]$ denotes the space of periodic functions in $L^2[-\pi,\pi]$ having  zero mean. As  \bo{in the earlier  sections}, $\del$ denotes the horizon. Here, the kernel function \bo{$\gamd(z) =\gamd(|z|)$ is defined by \eqref{eq:rescale} and \eqref{eq:kernel} with $d=1$, which
ensures} that as the horizon $\del$ tends to zero, the nonlocal operator $\mcLd$ reduces to the PDE Laplacian, \ie\ $\dd^2/\dd x^2$ in one dimension.

\pagebreak %%20200331
For any positive integer $n$, \bo{$\rme^{\pm \rmi nx}$ are eigenfunctions of $-\mcLd$} under
periodic boundary conditions with the corresponding eigenvalues\footnote{Note that $-\dd^2/\dd x^2$ has the same eigenfunctions $\rme^{{\pm}\rmi nx}$, as does $-\mcLd$, and has the eigenvalues $n^2$. Note that, in the limit $\del\to0$, $\lambda_{n,\del}\to \lambda_{n,0}:=n^2$.}
\begin{equation}\label{eigendel}
\bo{\lambda_{n,\del}=4\int_{0}^{\del}\gamd(z)(1-\cos(nz)) \D z.}
\end{equation}
Let $X_N^0:=\text{span}\{\rme^{{\pm}\rmi nx}\}_{n=1,2,\ldots, N}$ and let
\[
\bo{v_N(x) =\sum_{n=1}^N\widehat{v}_n \rme^{\rmi nx} + \widehat{v}_{-n} \, \rme^{-\rmi nx} \in X_N^0}
\]
for any set of constants $\{\widehat{v}_n\}$. Note that because, for any \bo{nonzero integer $n$}, $\rme^{\rmi nx}$ has zero mean with respect to any interval of length $2\pi$, so does $v_N(x)$.
We have that
\[
\bo{-\mcLd v_N=\sum_{n=1}^N \lambda_{n,\del} \bigl(\widehat{v}_n \rme^{\rmi nx}+ \widehat{v}_{-n} \,\rme^{-\rmi nx}\bigr) \quad \text{for } \del \geq 0.}
\]

The Fourier spectral-Galerkin scheme of {\citename{dy16sinum} \citeyear{dy16sinum,DY:2017}} 
{is as} follows: 
\begin{equation}\label{eqn:sol-uN}
  \begin{minipage}[b]{22pc}
{Given} $f(x)\in L^2_{\text{per}}[-\pi,\pi]$, seek $u_{N,\del}(x)\in X_N^0$ such that
\[
(-\mcLd u_{N,\del}, v_N) = (f,v_N) \quad\text{for all } v_N \in X_N^0
\]
{or equivalently}
\[
u_{N,\del} = -\mcLd^{-1} P_N f = \sum_{n=1}^N \lambda_{n,\del}^{-1}  \bigl(\widehat f_n \,\rme^{\rmi nx}+ {\widehat f_{-n}\,\rme^{-\rmi nx}}\bigr),
\]
    \end{minipage}
  \end{equation}
{where $(\cdot,\cdot)$} denotes the standard $L^2(-\pi,\pi)$ inner product and $P_N$ denotes the standard spectral Fourier projection onto $X_N:=\text{span}\{\rme^{{\pm}\rmi nx}\}_{n=0,1,\ldots, N}$ of an element in $ L^2[-\pi,\pi]$. 

{\citeasnoun[Lemma~1]{dy16sinum} proved that} 
the numerical scheme \eqref{eqn:sol-uN} is asymptotically compatible. In particular, {they proved} that
\begin{equation}\label{eqn:AC1}
\|  u_{N,\del} - u_{N,0}  \|_{L^2(-\pi,\pi)} \le \CCC \del^2 \| P_N f \|_{L^2(-\pi,\pi)},
\end{equation}
where the constant $\CCC$ is independent of $N$ and $\del$. Here, $u_{N,0}$ denotes the Fourier spectral solution in $X_N^0$ of the local diffusion problem
\begin{equation}\label{eqn:AC2}
  -\mathcal{L}_0 u(x) := -u_{xx}(x) = f(x) \quad\text{for all } x\in (-\pi,\pi)
  \pagebreak %%20200331
\end{equation} 
with the periodic boundary condition. The proof of the estimate \eqref{eqn:AC1} follows from the error representation
\[
\|u_{N,\del} - u_{N,0}\|_{L^2(-\pi,\pi)} ^2 =\sum_{n=1}^N\biggl|\ffrac{1}{\lambda_{n,\del}}-\ffrac{1}{\lambda_{n,0}}\biggr|^2\bigl(|\widehat{f}_n|^2 + {|\widehat{f}_{-n}|^2}\bigr)
\]
and a careful analysis {of} the asymptotic behaviour of $\lambda_{n,\del}$.

Let $u_0(\xb)$ denote the exact solution of the local problem \eqref{eqn:AC2}. Then, using standard error estimation theories for local models, it is easy to show that $u_{N,0}$ converges in {the $L^2$ sense} to $u_0$ at least quadratically with respect to $1/N$,
provided that $f\in L^2(-\pi,\pi)$. As a result, we arrive at
\begin{equation}\label{eqn:AC-spectral}
\|{u_{N,\del}-u_0}\|_{L^2(-\pi,\pi)}\le \CCC \,(\del^2+ N^{-2}) \|  f \|_{L^2(-\pi,\pi)},
\end{equation}
with a generic constant $\CCC$ independent of $\del$, $N$, $f$, $u_0$ and $u_\del$, where $u_\del$ denotes the exact solution of \eqref{nonloc}. This estimate indicates that the Fourier spectral method of
{\citename{dy16sinum} \citeyear{dy16sinum,DY:2017}} 
is asymptotically compatible. Moreover, we obtain a uniform error estimate in the sense that the estimates hold for any  sufficiently small $\del$ and any sufficiently large $N$ without any restriction on the relative sizes of $\del$ and $N$.

These one-dimensional results can be easily extended to higher-dimen\-sional nonlocal diffusion problems in rectangular and rectangular {parallelepiped} domains with periodic boundary conditions. See \citeasnoun{DY:2017} for discussions about two and three-dimensional problem. Furthermore, the analysis was extended to the nonlocal Allen--Cahn equations \cite[Section~4]{dy16sinum}. 

The eigenvalues $\lambda_{n,\del}$ defined by \eqref{eigendel} cannot, in general, be determined analytically. Thus, when implementing the scheme given in $\eqref{eqn:sol-uN}$, those eigenvalues $\lambda_{n,\del}$ have to {be} estimated. The numerical evaluation of the integral in \eqref{eigendel} is challenging because $\cos(nz)$ is highly oscillatory for large $n$. As an example, we focus on a kernel that is singular at the origin, namely,
\begin{equation}\label{kernels-singular}
	\gamd(z)=\ffrac{\CCC_\beta}{\del^{3-\beta}|z|^\beta} \quad\text{for all } z\in[-\del,0)\cup(0,\del]\ \mbox{and for $\beta\in(0,3)$.}
\end{equation}
Then, by Taylor expansion, we have that
\begin{align}\label{eqn:sig-eig}
\lambda_{n,\del} &= - \ffrac{{4C_\beta}}{\del^2}\int_0^1 r^{-\beta}  (\cos(n\del r) -1)  \D r \notag \\*
&= - \ffrac{{4C_\beta}}{\del^2} \sum_{k=1}^\infty \ffrac{(-1)^k|n\del|^{2k}}{(2k)!(2k+1- {\beta})} \notag \\*
& := - \ffrac{{4C_\beta}}{\del^2} K(n\del).
\end{align}
Therefore it suffices to compute $K(n\del)$ accurately. {\citeasnoun{DY:2017} have given a hybrid algorithm} for computing $K(n\del)$. If $n\del$ {is `small' then} the series in \eqref{eqn:sig-eig} converges {fast, so one} may approximate $K(n\del)$ by using a truncation of that series. On the other hand, for $n\del$ `large', {\citeasnoun{DY:2017} observed} that $K(n\del)$ is a solution of the ordinary differential equation
\begin{equation}\label{eqn:ini-K}
 K'(z) = \ffrac{z-1}{z} K(z) + \ffrac{\cos(z)-1}{z} \quad\mbox{on $(0,n\del]$}
\end{equation}
evaluated at $z=n\del$. For a starting condition, $K(1)$ may be used; it can be computed using a truncation of the series in \eqref{eqn:sig-eig}. As a result, $K(n\del)$ can be accurately evaluated for large $n\del$ by using, \feeg, a high-order Runge-Kutta method.

{\citeasnoun{smd18jcp} proposed a Fourier spectral method} to solve nonlocal diffusion models on the unit sphere $\mathcal{S}^2\subset \mathbb{R}^3$,
where the nonlocal Laplace--Beltrami operator is defined by 
\begin{equation}\label{nondefggg}
\mathcal{L}_\delta u(x)= 2\int_{\mathcal{S}^2} (u(x+z)-u(x)) 
\gamma_\delta(z) \D \nu(z),
\end{equation}
where $\nu$ denotes the standard measure on $\mathcal{S}^2$. The basic idea is to apply spherical harmonics \cite{Atkinson:book}:
\[Y_l^m(x) = Y_l^m (\theta,\varphi) = \ffrac{\rme^{\mathrm{i}m\varphi}}{\sqrt{2\pi}}\mathrm{i}^{m+|m|} \sqrt{\biggl(l+\ffrac12\biggr)\ffrac{(l-m)!}{(l+m)!}} P_l^m(\cos\theta)\]
with $l\ge0$ and $-l\le m\le l$, for which it is shown that they are
the eigenfunctions of the nonlocal Laplace--Beltrami operator.  Here $(\theta,\varphi)$ {are the spherical coordinates} of $x\in\mathbb{S}^2$ and $P_l$ 
{denote} the associated Legendre polynomials. Specifically, there {exists} 
$ \mathcal{L}_\delta Y_l^m(x) = \lambda_\delta(l) Y_l^m(x)$ with
\[
\lambda_\delta(l) = 4\pi \int_{-1}^1 ( P_l(y)-1) \gamma_\delta(\sqrt{2(1-y)})\D y.
\]

To compute the eigenvalue $\lambda_\delta(l)$, one  evaluates an integral which is highly oscillatory due to the Legendre polynomials.
 \citeasnoun{smd18jcp} proposed a numerical scheme by using a modified Clenshaw--Curtis quadrature rule
with the computation complexity $O(l^2)$ per eigenvalue $\lambda_\delta(l)$. For sufficiently large $l$, they applied 
{Szeg\H{o}'s} asymptotic formula for Legendre polynomials, which reduces the complexity to $O(l\log l)$ per eigenvalue.
Then a fast spherical harmonic transform 
\cite{Slevinsky:fast} was applied to accelerate synthesis and analysis of the series expansion in spherical harmonics.

\begin{rem}[extension to {phase field} methods]\hspz% 
{The study of spectral approximation to nonlocal problems has been extended to nonlinear models such as phase field equations \cite{dy16sinum,djlq18jcp}.  Conventional phase field (diffuse interface) models take on a free energy of the type \cite{df19phase}
\begin{equation}\label{lcEn}
E^0(u)=\int_\Omega \biggl(\ffrac{\epsilon}{2}|\nabla u(x)|^2 +\ffrac{1}{4\epsilon}(u(x)-1)^2\biggr) \D x.
\end{equation}
The nonlocal version of the free energy can be written as \cite{dy16sinum}
\begin{equation}\label{nonEn}
E^\delta(u)=\int_\Omega \int_\Omega \biggl(\ffrac{\epsilon}{2}\gamma_\delta(|s|)\ffrac{(u(x+s)-u(x))^2}{2} \D s+\ffrac{1}{4\epsilon}(u(x)-1)^2\biggr) \D x,
\end{equation}
which recovers $E^0$ as its {gamma} limit for $\delta\to 0$.
}

A nonlocal Allen--Cahn equation can be defined as the $L^2$ gradient flow of the  energy in
\eqref{nonEn}. The convergence of AC spectral approximations for problems defined on periodic domains has been established in \citeasnoun{dy16sinum}.  With the assumption on smooth minimizers, spectral accuracy can be assured.
{When} nonlocal interaction kernels do not have
sufficiently strong singularities, the minimizers of the nonlocal free energy (subject to a total mass condition) may develop a discontinuous profile, in contrast to the conventional diffuse interface models with smooth phase field functions. {However}, the convergence of spectral approximations can be assured even for such cases.  We note that, historically, \eqref{lcEn} 
is often derived from a nonlocal form \eqref{nonEn} via the {\em Landau expansion}; see the additional discussions on this perspective given in \citeasnoun{du19cbms}. Naturally, one may also consider the $H^{-1}$ gradient flow of the nonlocal energy or gradient flows with respect to a nonlocal space
embedded between $H^{-1}$ and $L^2$; some related numerical analysis can be found
in \citeasnoun{li2019convergence} {and} \citeasnoun{ainsworth2017analysis}. \hfill\actaqed
\end{rem}
 
\begin{rem}\hspz% 
Nonlocal diffusion allows for singular solutions, which pose new challenges in the design of numerical discretizations. Whereas spectral methods are shown to possess many good properties such as asymptotic compatibility and are able to capture the discontinuities in solutions, spurious Gibbs phenomena do appear. Given the development of spectral methods in the context of PDEs that are effective in singularity detection and in high-order recovery of solution information near and away from singularities (see \eg\ \citeb{Tadmor:2007}), 
it is natural to pursue further studies {of} how to adopt and extend those techniques to nonlocal models.

Furthermore, the algorithms introduced in this section focus on periodic boundary conditions and allows for the derivation of uniform convergence rates with respect to both $\del$ and $N$; see \eg\ \eqref{eqn:AC-spectral}. The extension to other types of boundary conditions, \eg\ Dirichlet or Neumann volume constraints, 
remains an open problem. The main difficulty stems from the lack of a closed expression for the eigenvalues analogous to the expression given in \eqref{eigendel}. Hence the difference $\lambda_{n,\del}^{-1} - \lambda_{0,\del}^{-1}$ is difficult to estimate. An investigation of the AC property of spectral method {in such settings is also}  warranted. \hfill\actaqed
\end{rem}

\vspace{5pt}
\section{Spectral-Galerkin methods for fractional diffusion}\label{sec:smf}

In this section we consider spectral methods for determining approximate solutions of fractional diffusion problems. We first address the fractional Poisson problem posed on $\omg=\Rd$. In this case, the solution decays slowly with algebraic rates as $|\xb|\to\infty$, even if the source term is compactly supported. As a result, some traditional discretization strategies used for the local PDE diffusion problem are not applicable here. Subsequently, we consider fractional diffusion in a bounded domain $\omg$ for which the solution possesses a weakly singular layer near the boundary $\partial\omg$, even if the source term is smooth. This property engenders many challenges in the development and analysis of numerical schemes. Compared to, say, finite difference and finite element methods, spectral methods with specially constructed basis functions may approximate the solution of fractional diffusion problems on bounded domains with higher accuracy.

\subsection{Spectral-Galerkin methods in unbounded domains}\label{sec:smfu}

For $s\in(0,1)$, we consider the fractional diffusion problem in $\omg=\Rd$ given~by
\begin{equation}\label{eqn:fde0}
\begin{cases}
 (-\Delta)^s u (\xb) + \ccc u(\xb) = f(\xb) &\text{for all } \xb\in\Rd ,\\
 u(\xb)\to 0  & \text{as } |\xb|\rightarrow\infty,
\end{cases}
\end{equation}
where $\ccc$ denotes a given positive constant and $f$ denotes a given source term. A weak formulation {of problem} \eqref{eqn:fde0} {is as} follows:  
\begin{equation}\label{eqn:var}
  \begin{minipage}[b]{18pc}
{Given} a positive constant $\ccc$ and $f\in H^{-s}(\Rd)$,
seek $u\in H^s(\Rd)$ such that
\[
{\mcA_s}(u,v) =   (f,v ) \quad\text{for all } v\in H^s(\Rd),
\]
  \end{minipage}
\end{equation}
where, for all $u,v\in H^s(\Rd)$,
\[
\mcA_s(u,v) = \ffrac{\cds}{2}\int_{\Rd}\int_{\Rd} \ffrac{ (u(\yb) - u(\xb))   (v(\yb)-v(\xb)) }{|\ymx|^{\ddd+2s}}\dydx + \ccc(u,v) 
\]
with $\cds$ is defined in \eqref{cds}. The {well-posedness of problem} \eqref{eqn:var} follows directly from the Lax--Milgram theorem.

Several \roo{successful} strategies have been developed for determining approximations of solutions of local PDE diffusion problems posed on $\Rd$; see \eg\ the survey paper {by} \citeasnoun{ShenWang:2009} and the references cited therein. One way is to truncate the domain $\Rd$, provided that the solution decays rapidly. However, it is known that solution of the fractional diffusion equation decays slowly with a power law at infinity so that the naive truncation approach results in a relatively large error due to domain truncation.\footnote{However, see Remark~\ref{truncfl}, where it is pointed out that truncated domains may be of interest in their own right.} Another approach is to design an artificial boundary condition imposed on the boundary of a truncated domain that results in the same solution within a bounded region; the use of this approach remains largely an open question in the fractional case. A third approach, which is the approach used here, does not involve domain truncation, but instead approximates the solution in terms of orthogonal functions that potentially can well approximate the solution in the unbounded domain. 

For simplicity, we restrict our discussion to the one-dimensional case, \ie\ $\ddd=1$; higher-dimensional cases can be treated in a similar manner by using tensor products of one-dimensional orthogonal polynomials.

\subsubsection{Approximation by Hermite polynomials}\label{sec:hermite}

The orthonormal Hermite polynomials $\{H_n(x)\}$ are defined by the three-term recurrence relation \cite{aands,MaoShen:2017}
\begin{align*}
H_0 &= \pi^{-\gfrac14},\qquad H_1(x) = \sqrt{2}\pi^{-\gfrac14}x,\\*
H_{n+1} &= x\sqrt{\ffrac{2}{n+1}}H_n(x)- \sqrt{\ffrac{n}{n+1}}H_{n-1}(x)\quad\text{for }  n=1,2,\ldots.
\end{align*}
It is {well known} that those polynomials form an orthonormal basis in the weighted space $L_{\www}^2(\Ro)$, where the weight function is given by $\www(x)=\rme^{-x^2}$.
Let $V_N^H$ denote the space spanned by the Hermite polynomials of degree less than or equal to $N$, and let
\[
X_N^H = \{ v \colon \sqrt{\www(x)}u(x) \ \text{for all } u \in V_N^H \}.
\] 
We define the $L^2_\www(\Ro)$-projection operator $Q^H_N\colon L^2_\www(\Ro)\to X_N^H$ by
\[
\int_{\Ro} u(x) v_N(x)\www(x)\D x= \int_{\Ro} (Q^H_N u)(x) v_N(x)\www(x)\D x
\quad\text{for all } v_N\in X_N^H.
\]
We also define the modified operator $\widehat Q^H_N\colon L^2(\Ro) \rightarrow X_N^H$ by
\[
\widehat Q^H_N u =  \www^{\gfrac12}Q^H_N (u\,\www^{-\gfrac12})  \quad\text{for all } u\in L^2(\Ro).
\]
Then we have the approximation property 
\begin{equation}\label{eqn:approx-Herm}
   \| \widehat Q^H_N u - u    \|_{H^s(\Ro)} \le \CCC N^{(\ell-m)/2}  \| u    \|_{\HHHh^m(\Ro)}
\quad\text{for all } s\in[0,m],
\end{equation}
where $\CCC>0$ does not depend on $u$ or $N$ and where the modified Sobolev space $\HHHh^m(\Ro)$ comes equipped with the norm and seminorm
\begin{align*}
  \|  u  \|_{\HHHh^m(\Ro)}^2 = \sum_{0\le k\le m}  \biggl\| \biggl( {\ffrac{\dd}{\dd x}} + x\biggr)^k u \biggr\|_{L^2(\Ro)}^2 , \\*
|  u  |_{\HHHh^m(\Ro)} =    \biggl\|  \biggl( {\ffrac{\dd}{\dd x}} + x\biggr)^m u \biggr\|_{L^2(\Ro)}^2,  
\end{align*}
respectively. 

The Hermite--Galerkin discretization of the fractional diffusion problem \eqref{eqn:fde0} is {then as} follows: 
\begin{equation}\label{eqn:spectral-Galerkin-Herm}
  \begin{minipage}[b]{18pc}
{Given} a positive constant $\ccc$ and $f\in H^{-s}(\Rd)$,
seek $u_N(x) \in X_N^H$ such that
\[
\mcA_s(u_N, v_N)  =   (I_N^H f, v_N )  \quad\text{for all } v_N\in X_N^H,
\]
    \end{minipage}
\end{equation}
where $I_N^H \colon C(\Ro)\rightarrow  {X_N^H}$ denotes the Gauss--Hermite interpolation operator
with respect to the Gauss--Hermite points \cite{ShenTang:Book}. Then the standard energy estimate {and} the approximation property \eqref{eqn:approx-Herm} result in the error estimate
\begin{equation}\label{eqn:err-Herm}
\|  u - u_N  \|_{H^s(\Ro)} \le \CCC N^{s-m} |  u |_{\HHHh^m(\Ro)}
+  \CCC N^{1/6-\ell/2} | f |_{\HHHh^\ell(\Ro)} ,
\end{equation}
provided that $u\in\HHHh^m(\Ro)$ and $f\in\HHHh^\ell(\Ro)$, where again $\CCC>0$ does not depend on $u$, $f$ or $N$. 
{This estimate indicates that the numerical solution fails to converge exponentially, because either $u$ and or $f$ will decay algebraically as $x\to \infty$; \scf\ \citeasnoun[Section~5]{MaoShen:2017}.}

The Hermite--Galerkin discretization \eqref{eqn:spectral-Galerkin-Herm} is simple to implement because the mass matrix is diagonal, \bo{whereas} the stiffness matrix can be computed efficiently using the Fourier transform.
An implementation of the Hermite--Galerkin method \eqref{eqn:spectral-Galerkin-Herm} using Fourier transforms is provided in \citeasnoun{MaoShen:2017}. See also \citeasnoun{MaoShen:2017} {and} \citeasnoun{TangYuanZhou:2018} for related Hermite collocation methods, whose error estimate also depend on the Fourier transform of the solution $u$ and the source term $f$.

\subsubsection{Approximation by modified Gegenbauer polynomials}\label{sec:gegen}

In \citeasnoun{tang:2019rational}, the approach taken in Section~\ref{sec:hermite} is extended to approximations in which the basis is constructed by using Gegenbauer polynomials that are modified using a nonlinear singular mapping. The resulting basis has {proved} to be better suited {to} approximating functions with algebraic decay rates (see \eg\ \citeb{Boyd:JCP1987}, \citeb{Boyd:2001}) when compared to classical {bases of} orthogonal polynomials such as Hermite or Laguerre polynomials. 

The Gegenbauer polynomials, denoted by $\GGG_{n}^{\laa}(t)$ with $ t\in (-1,1)$ and scaling parameter $\laa>-1/2$, are defined by the three-term recurrence relation
\begin{align*}
\GGG_{0}^{\laa}(t)&=1,\quad \GGG_{n}^{\laa}(t)=2\laa t,\\*
n\GGG_{n}^{\laa}(t)&=2t(n+\laa-1)\GGG_{n-1}^{\laa}(t)-(n+2\laa-2)\GGG_{n-2}^{\laa}(t), \;\; n=2,3,\ldots.
\end{align*}
The Gegenbauer polynomials are orthogonal with respect to the weight function $\www_\laa(t)=(1-t^2)^{\laa-1/2}$, and they are closely related to the hypergeometric functions \cite[{page}~1000]{Gradshteyn:7ed} that are widely used for the analysis of fractional differential equations; see \eg\ \citeasnoun{ervin2018regularity}.

We seek approximate solutions in the space
\[
X_{N}^{G}  =\text{span}\{R_{n}^{\laa}(x) \colon n=0,1,\ldots,N\},
\]
where the mapped Gegenbauer functions $R_{n}^{\laa}(x)$ are defined by
\[
R_{n}^{\laa}(x):= (1+x^2) ^{-\tfrac{\laa+1}{2}}\GGG_{n}^{\laa}\biggl(\ffrac{x}{\sqrt{1+x^2}}\biggr)\quad\text{for all } x\in{\Ro}.
\]
The $L^2(\Ro)$-orthogonal projection operator $Q_{N}^G: L^2(\Ro)\to X_N^G$ is defined by
\[
\int_{\Ro} u(x) v_N(x) \D x = \int_{\Ro} Q_{N}^G u(x) v_N(x) \D x \quad\text{for all } v_N\in X_N^G .
\]
Next, let $B_m^G(\Ro)$ denote the weighted Sobolev space defined by the norm and seminorm
\begin{align*}
\|  u  \|_{B_m^G(\Ro)}^2 & = \sum_{0\le k\le m}
\biggl\| (1+x^2)^{\tfrac 1 4-\tfrac{\laa+m} 2}\biggl( (1+x^2)^{\tfrac32} \ffrac{\dd}{\dd x}\biggr)^k  \Bigl((1+x^2)^{\tfrac{\lambda+1}2}u\Bigr) \biggr\|_{L^2(\Ro)}^2 , \\*
|  u  |_{B_m^G(\Ro)} & =    \biggl\| (1+x^2)^{\tfrac 1 4-\tfrac{\laa+m} 2} \biggl( (1+x^2 )^{\tfrac32} \ffrac{\dd}{\dd x}\biggr)^m  \Bigl((1+x^2)^{\tfrac{\lambda+1}2}u\Bigr)   \biggr\|_{L^2(\Ro)}^2,
\end{align*}
respectively. Then, for any $u\in H^s(\Ro)\cap B_m^G(\Ro)$ with integer $ 1\le m\le N+1$, $ s\in (0,1)$ and $\laa>-1/2$, we have the approximation property
\[
 \| Q_{N}^G u-u \|_{H^s(\Ro)}\le \CCC N^{s-m} |u|_{B_m^G(\Ro)},
\]
where $\CCC>0$ does not depend on $N$ and $u$. 

The mapped Gegenbauer--Galerkin discretization of the fractional diffusion problem \eqref{eqn:fde0} {is as} follows: 
\begin{equation*}
  \begin{minipage}[b]{18pc}
{Given} a positive constant $\ccc$ and $f\in H^{-s}(\Rd)$,
seek $u_N(x) \in X_N^G$ such that
\[
\mcA_s(u_N, v_N)  =   (I_N^G f, v_N )  \quad\text{for all } v_N\in X_N^G,
\]
  \end{minipage}
  \pagebreak %%20200331
\end{equation*}
where $I_N^G\colon C(\Ro)\rightarrow X_N^G$ denotes the mapped Gegenbauer--Gauss interpolation operator based on the Gegenbauer--Gauss quadrature points \cite[Chapter~3]{ShenTang:Book}. {Then, for any $u\in H^s(\mathbb R)\cap {B}_{m}^G(\mathbb R)$ and $f\in  {B}_{\ell}^G(\mathbb R) $ with integers $ 1\le m,\ell\le N+1$ and $\laa>-1/2$, we have} the error estimate
\begin{equation}\label{IntestaS}
 \| u-u_N\|_{H^s(\Ro)}\le \CCC N^{s-m} |u|_{{ B}_{m}^G(\Ro)}+\CCC N^{-\ell} |f|_{{B}_{\ell}^G(\Ro)},
\end{equation}
where $\CCC$ denotes a positive constant having value independent of $N$, $u$ and $f$. The proof relies on the standard energy estimate as well as the approximation properties of $Q_N^G$ and $I_N^G$. Compared to the Hermite--Galerkin spectral method of Section~\ref{sec:hermite}, we have that $(-\Delta)^s R_n^\laa(x)$ is not as compact as for the Hermite polynomials, but can be easily approximated by a finite series of hypergeometric functions.

\subsubsection{Approximation by modified Chebyshev polynomials}\label{sec:rational}

It is noteworthy that $|\xib|^{2s}$, the symbol of {the}
  Fourier transform of the fractional Laplacian, is non-separable so that direct extensions to higher dimensions of methods developed for the one-dimensional setting can be very complicated. To overcome this computational difficulty, {\citeasnoun{sheng:2019fast} proposed a reformulation of the weak formulation \eqref{eqn:var} by the Dunford--Taylor formula (see Theorem~4.1 of \citeb{Bonito:2019},  and also Section~\ref{dunftay}), and constructed a spectral-Galerkin method having} Fourier-like bi-orthogonal mapped Chebyshev polynomials as basis functions. This method leads to a diagonalized system that can be efficiently solved via the fast Fourier transform (FFT), and the approach can be easily extended to higher-dimensional problems via tensor products of the basis functions. 

For any $u,v\in H^{s}(\Ro)$, by the Dunford--Taylor formula, we have (see \eqref{dtbf})
\begin{align*}
&\ffrac{\cds}{2}\int_{\Rd}\int_{\Rd} \ffrac{ (u(\yb) - u(\xb))   (v(\yb)-v(\xb)) }{|\ymx|^{\ddd+2s}} \dydx \\*
& \quad =\ffrac{2 \sin (\pi s)}{\pi} \int_0^\infty t^{1-2s} \int_{\Rd}
  ((-\Delta)( I-t^2 \Delta)^{-1} u ) (\xb)\,  v(\xb) \dxb \D t,
\end{align*}
{where  $I$} denotes the identity operator. Letting 
$w(\xb):= ( I-t^2 \Delta)^{-1} u(\xb)$, 
we can rewrite the weak formulation \eqref{eqn:var} as: seek $u\in H^s(\Rd)$ such that for all $v\in  H^s(\Rd)$
\begin{equation}\label{auveqn2}
\mcA_s(u,v)=\ffrac{2 \sin (\pi s)}{\pi}  \int_0^\infty t^{-1-2s} (u-w, v) \D t+ \ccc\,(u,v) =(f,v),
\end{equation}
where $w=w(\xb;u,t)\in H^{1}(\Rd)$ is, for any $t>0$, the solution of
\begin{equation}\label{coneta}
t^{2} (\nabla w,\nabla \psi) +(w, \psi) =(u, \psi) \quad\text{for all } \psi \in H^{1}(\Rd).
\end{equation}

\pagebreak %%20200331
The spectral-Galerkin approximation to \eqref{auveqn2}--\eqref{coneta} is then given as follows: find  $u_N\in   {X_N^C}$ such that
\begin{align}\label{auveqn23}
 \mcA_s(u_N,v_N)&=\ffrac{2 \sin (\pi s)}{\pi} \int_0^\infty t^{-1-2s} (u_N-w_N, v_N) \,  {\rm d}t+\ccc(u_N,v_N) \notag  \\*
 & =(I_N f,v_N) \quad\text{for all } v_N\in X_{N}^C,
\end{align}
where $w_N(\xb;u_N,t)\in X_{N}^C$ is, for any $t>0$, the solution of 
\begin{equation}\label{bSeqn}
t^{2} (\nabla w_N,\nabla \psi_{{N}}) +(w_N, \psi_{{N}}) =(u_N, \psi_{{N}}) \quad\text{for all } \psi_{{N}} \in X_{N}^C.
\end{equation}

To solve the numerical scheme \eqref{auveqn23}--\eqref{bSeqn}, {\citeasnoun{sheng:2019fast} used} a basis of Fourier-like mapped Chebyshev polynomials. Recall that the classical Chebyshev polynomials are defined by, for $n=0,1,\ldots,$
\[
T_n(y)=\cos(n\, {\rm arc\,cos}(y))\quad\text{with } y\in (-1,1).
\]
Then we define the mapped Chebyshev polynomials by
\begin{equation}\label{modifiedrational}
R_{n}  (x) := (1+x^2) ^{-\gfrac{1}{2}}T_n\biggl(\ffrac{x}{\sqrt{1+x^2}}\biggr),\quad x\in{\mathbb R}.
\end{equation}
It is well known that the mapped Chebyshev functions are orthonormal in $L^2(\mathbb R)$ but are not orthogonal in $H^1(\Rd)$.
Following the spirit of \citeasnoun{ShenWang:2009},  in \citeasnoun{sheng:2019fast} the approximation space is defined by
$
X_N^C    =\text{span}\{R_{n}\colon n=0,1,\ldots,N\}
$
for which a Fourier-like basis that is bi-orthogonal can be constructed, \tie, the basis functions are orthogonal with respect to both the {$L^2$ and $H^1$ inner} products. As a result,  \eqref{bSeqn} is diagonalizable, as is \eqref{auveqn23}, and hence they can be efficiently solved by an FFT related to Chebyshev polynomials.
 
The integration in terms of $t$ in \eqref{auveqn23} can be evaluated exactly by  using  a sinc quadrature scheme \cite{Bonito:2019} whose computational cost is negligible compared with that of an FFT. Furthermore, this fast algorithm can be easily applied in higher-dimensional settings by using tensor products of the one-dimensional basis functions. The resulting complexity of solving \eqref{auveqn23} is then of \bo{$\Bigoh (N(\log N)^\ddd)$, where $N$ is the number of degrees of freedom}. The relevant error estimate is similar to that in \eqref{IntestaS} with $\laa=0$.

\begin{rem}\hspz% 
Although the solution to the fractional diffusion problem \eqref{eqn:fde0} has a simple and closed form in the Fourier domain, all existing spectral methods fail to achieve at an exponential convergence rate. This is because the solution generally decays at an algebraic rate at infinity, even if the source term is smooth and compactly supported. The proposed Hermite and the mapped rational/Chebyshev polynomials can only approximate exponentially decaying functions with spectral accuracy. How one designs a spectral method that accurately captures the algebraic decaying behaviour of the solution to \eqref{eqn:fde0} remains a largely open problem.

All the estimates stated in this subsection are derived with respect to the energy norm, \eg\ $H^s(\Rd)$. There is, of course, interest in obtaining higher approximation rates with respect to the {$L^2(\Rd)$-norm}. However, the approximation properties of orthogonal polynomials are derived in weighted Sobolev spaces, {so the} usual duality arguments are not directly applicable here. Thus optimal $L^2(\Rd)$-norm error estimates are still not well understood and warrant further investigation. \hfill\actaqed
\end{rem}

\subsection{Spectral-Galerkin methods in bounded domains}\label{ssec:spec-frac-bound}

Now we turn to fractional diffusion in a bounded domain. Owing to the low regularity of the solution of the underlying problems involving the fractional Laplacian (see \eg\ \citeb{ervin2018regularity}, \citeb{grubb2015spectral}, \citeb{ros2014dirichlet}), classical spectral methods cannot approximate the solution very well. In this section we focus on the particular fractional Poisson problem \eqref{eq:fracPoisson} with $g(\xb)=0$, \tie, on the problem
\begin{equation}\label{eqn:fde1}
\begin{cases}
 (-\Delta)^{s}u = f &\text{for all }\xb\in\omg ,
\\
  u=0&\text{for all }\xb\in{\omgii}=\Rd\setminus\omg ,
\end{cases}
\end{equation}
where $f(\xb)$ is a smooth function and $\omg$ is a bounded domain with a smooth boundary.

{\citeasnoun{mao2016efficient} considered the one-dimensional setting with} $\omg=(-1,1)$ and $s\in(1/2,1)$. In this case, the fractional Laplacian
is equivalent to the {\em Riesz-fractional differential derivative} that is defined by, for $u\in C_c^\infty(-1,1)$ \bo{and $x\in(0,1)$}
\[
 -{D^{2s}} u(x) := - \ffrac{1}{2\cos((1-s)\pi)}\ffrac{\dd^2}{\dd x^2}\int_{-1}^1 |x-t|^{1-2s} u(t) \D t =(-\Delta)^s u (x) .
\]
{Then problem} \eqref{eqn:fde1} is equivalent to the two-point boundary value problem: find the function $u$ satisfying
\begin{equation}\label{eqn:fde2}
\begin{cases}
-{D^{2s}} u(x) =f(x)   \quad\text{for all } x\in(-1,1),\\
 u(0)=u(1)= 0.
\end{cases}
\end{equation}

The {discussion} in \citeasnoun{mao2016efficient} is motivated by the observation that
\begin{equation}\label{eqn:frac-J}
(-\Delta)^{s}\{(1-x^2)^s_+ P_n^{s,s}(x)  \}
= \ffrac{\Gam(2s+n+1)}{\Gam(n+1)}  P_n^{s,s}(x),
\end{equation}
where $P_n^{a,b}$ with $a, b>-1$ denotes the classical Jacobi polynomials that are orthogonal with respect to the weight function $\www^{a,b}(x) = (1-x)^a(1+x)^b$. This implies that if the source term $f(x)$ can be expanded in terms of the Jacobi polynomials $P_n^{s,s}(x)$, then the solution of the fractional diffusion problem \eqref{eqn:fde1} is a series in $(1-x^2)^s_+ P_n^{s,s}(x)$.

By choosing the test space $V_N = \text{span}\{ P_n^{s,s}(x) \colon  n=0,\ldots,N \}$ and the trial space $X_N = \text{span}\{ (1-x^2)_+^s P_n^{s,s}(x) \colon  n=0,\ldots,N \}$, a Petrov--Galerkin spectral method is given as follows: find $u_N \in X_N$ such that, for all $v_N\in V_N$,
\[
\int_{-1}^1 (-\Delta)^{s} u_N(x) \, v_N(x) (1-x^2)^s \D x = \int_{-1}^1 f(x)  \, v_N(x) (1-x^2)^s \D x.
\]
Using the relation \eqref{eqn:frac-J} and the approximation property of Jacobi polynomials \cite[Theorem~3.35]{ShenTang:Book}, the following error estimate in the weighted $L^2$-norm is derived in \citeasnoun[Theorem~5]{mao2016efficient}:
\[
\|u_N - u \|_{B^0_{-s}} \le \CCC N^{-2s-m} \|  f \|_{B^m_s}  \quad \text{for any integer} ~~m \ge 0,
\]
where the norm $\| \cdot\|_{B^m_\lambda}$ is defined by
\begin{equation*}
 \| v \|_{B^m_\lambda}^2 = \sum_{j=0}^m \int_{-1}^1 | v^{(j)} |^2 (1-x^2)^{j+\lambda} \D x.
\end{equation*}
{If} $m$ is not an integer, the space $B^m_\lambda$ is defined by interpolation \cite{GuoWang:2004}. The argument also applies to the case {when} $s\in(0,1/2]$ and $\omg$ is multi-interval \cite{Acosta:MC2018}.

The above estimate indicates the exponential convergence under the assumption that the source term $f$ is smooth. The results in \citeasnoun{mao2016efficient} extend the studies {of} spectral methods for solving one-sided fractional equations found in
\citeasnoun{ChenShenWang:2015}, \citeasnoun{Zayernouri:2014} {and} \citeasnoun{Zayernouri:2015}; see also
\citeasnoun{Mao:2018two} {and \citename{Saminee:2019a} \citeyear{Saminee:2019a,Saminee:2019b}} for
general two-sided fractional equations, \citeasnoun{ZayernouriAK:2015} for tempered fractional diffusion equations, and \citeasnoun{DengZhao:2019} for the study of superconvergence points.

A similar analysis may be extended to higher-dimensional cases. {\citeasnoun{XuDarve:2018} studied a spectral method} for solving the fractional Poisson problem \eqref{eqn:fde2} with $\omg=B_1({\bf 0})\subset \mathbb{R}^\ddd$, $\ddd=2,3$, \tie, the unit ball in two or three {dimensions}. The basic idea is to construct a spectral basis by using Jacobi polynomials and spherical harmonic polynomials \cite{Dyda:2017}. In particular, spherical  harmonic polynomials of degree $\ell\ge0$ form a finite-dimensional space, having dimension
\[
M_{\ddd,\ell}=\ffrac{\ddd+2\ell-2}{\ddd+\ell-2}
\binom{\ddd +\ell-2}{\ell}.
\]
For each $\ell$ we fix a linear basis for this space, denoted by $V_{\ell,m}$ with $m = 1,\ldots, M_{\ddd,\ell}$, that is orthonormal with respect to the surface measure on the unit sphere. Then the functions
\[  
\phi_{\ell,m,n}(\xb) =  V_{\ell,m}(\xb)P_n^{s,\gfrac{\ddd}2+\ell-1}(2|\xb|^2-1)\quad  \text{with }  \ell\ge0, \ n\ge 0,\ 0\le m\le M_{\ddd,\ell}
\]
form an orthogonal basis in $L^2_\www(\omg)$ with respect to the weight function $\www({\xb})=(1-|\xb|^2)^s$ \cite{Dyda:2017}.
Analogous to \eqref{eqn:frac-J}, we {have}
\[
(-\Delta)^{s}   ((1-|\xb|^2)_+^s \phi_{\ell,m,n}  ) 
= \ffrac{2^{2s}\Gamma(s+n+1)\Gamma(\gfrac{\ddd}{2}+\ell+s+n)}{n!\Gamma(\gfrac{\ddd}{2}+n)}  \phi_{\ell,m,n} (\xb).
\]
{\citeasnoun{XuDarve:2018} have developed and analysed fast and accurate spectral methods based} on this property and the approximation property of the Jacobi polynomials.

Unfortunately, the above argument heavily relies on the relation \eqref{eqn:frac-J} so that, due to the singular behaviour of the solution near the boundary, it cannot be directly applied to the fractional diffusion with lower-order terms. {\citeasnoun{zhang2019error} studied the regularity of the one-dimensional fractional diffusion--reaction model
\begin{equation}\label{eqn:fde3}
\begin{cases} 
(-\Delta)^{s} u(x) + \ccc u(x)=f(x)  & \text{for all } x\in (-1,1),\\
 u(x)= 0  & \text{for all } x\in \Ro\setminus(-1,1)
\end{cases} 
\end{equation}
with the constant $s\in(0,1)$ and $\ccc >0$.} In particular, {\citename{zhang2019error} proved that}
$(1-x^2)^{-s}u$ belongs to $B_{s}^{(2s\wedge m) +2s }$
if $f\in B_{s}^{m}$. Moreover, the regularity index could be slightly improved, \tie,
\[
 (1-x^2)^{-s}u \in B_{s-1}^{(6s-\epsilon)\wedge (1+4s)}\quad \text{provided that $ f \in B_{s-1}^{m}$  for large $m$.}
\] 
The convergence rate of the spectral method certainly deteriorates accordingly, due to a lack of sufficient solution regularity. See also  \citeasnoun{mao2018spectral} for  a spectral element method based on geometric meshes.

\begin{rem}\hspz% 
In one dimension, the fractional Poisson problem without lower-order terms has been comprehensively studied. The analysis of the regularity and numerical approximations \bo{relies} crucially on the simple relation \eqref{eqn:frac-J}, which indicates that the Jacobi polynomials are suitable to approximate the solution. This is largely the main reason why problems without lower-order terms dominate the study of spectral methods for fractional diffusion. Thus it is of significant interest to develop proper techniques for rigorously handling convection and reaction terms.

The numerical analysis of higher-dimensional fractional problems is a nontrivial endeavour. One reason is the difficulty of designing spectral methods for irregular domains. Another is due to the unclear regularity theory of \eqref{eqn:fde1} in higher-dimensional domains, even for convex polygonal domains. Therefore it is of substantial interest to analyse spectral methods for solving fractional diffusion problems in bounded domains in higher dimensions. Many theoretical questions, \eg\ regularity in suitable spaces, appropriate selection of basis functions and optimal convergence rates, remain largely open problems. \hfill\actaqed
\end{rem}

%%\vspace{5pt}
\section[Finite difference methods for the strong form of\\ nonlocal diffusion]{Finite difference methods for the strong form of nonlocal diffusion}\label{sec:FD-nonlocal}

The AC property was first demonstrated for quadrature-based finite difference discretizations of nonlocal models in \citeasnoun{TiDu13}. Additional observations can be found in \citeasnoun{du15lec}. {\citeasnoun{du16grad} studied an AC quadrature difference  discretization together} with superconvergent nonlocal gradient recovery. {\citeasnoun{dtty19ima} considered} quadrature rules for scalar models in multidimensional spaces, resulting in AC schemes. The key to {obtaining} AC non-variational methods in these {works} is to guarantee the uniform truncation errors of the numerical schemes independent of the parameter $\del$. Numerical methods for the strong form of nonlocal diffusion are analysed under the standard framework of truncation error analysis and numerical stability. Quadrature-based finite difference schemes for nonlocal diffusion equation are introduced in \citeasnoun{dtty19ima} {and} \citeasnoun{TiDu13}, and a reproducing kernel (RK) collocation method is studied in \citeasnoun{LTTF19}. Note that although {\citeasnoun{tyyp19cmame} refer to an AC mesh-free scheme,
their numerical scheme converges only to the corresponding local solution but} not to the nonlocal solution with a fixed $\del$.   

We denote a uniform Cartesian grid on $\Rd$ as
\[
\mcTh := \{ \xb_\jb = h\jb  \colon  \jb\in \Mbzd \}. 
\]
We can then rewrite the nonlocal diffusion operator at any node $\xb_\ib \in \mcTh\cap \omg$ as 
\[
\mcLd u(\xb_\ib)
= 2\int_{\mathcal{B}_\del({\bm 0})}\dfrac{u(\zb+\xb_\ib)-u(\xb_\ib)}{W(\zb)}W(\zb)\gamma_\del(|\zb | ) \D \zb,
\]
where $W(\zb)$ is a weight function that is crucial {to guaranteeing} the AC property. 
One-dimensional problems are considered in \citeasnoun{TiDu13} with the weight function $W(z)=|z|$ used for $z\in \Ro$. In the multi-dimensional case considered in \citeasnoun{dtty19ima}, the weight function $W(\zb)= |\zb|^2/|\zb|_1$ is chosen, where $|\cdot|_1$ denotes  the {$\ell_1$-norm} in the $\ddd$-dimensional vector space whereas $|\cdotsp|$ denotes the standard Euclidean norm. The finite difference scheme on the Cartesian grid is then given by
\begin{equation}\label{eq:FD_multi_linear}
	\mcLdh u(\xb_\ib)
	= 2\int_{\mathcal{B}_\del({\bm 0})}\mathcal{I}_h\biggl(\dfrac{(u(\zb+\xb_\ib)-u(\xb_\ib)}{W(\zb)}\biggr)W(\zb)\gamma_\del(|\zb | ) \D \zb,
\end{equation}
where $\mathcal{I}_h(\cdotsp)$ denotes the piecewise $\ddd$-multilinear interpolation operator in $\zb$ associated with the grid $\mcTh$.

\pagebreak %%20200331
There are two major features of the finite difference scheme \eqref{eq:FD_multi_linear}. The first is a uniform consistency result. In this regard, 
the following result about what we refer to as {\em quadratic exactness} is proved in \citeasnoun{dtty19ima}; it is a nonlocal analogue of its local counterpart. That is,  the fact that the centred difference approximation to the Laplacian is exact for quadratic {polynomials  holds in the nonlocal case as well}.

\paragraph{Quadratic exactness.} 
For any quadratic polynomial in $\Rd$ given as $u(\xb)=\xb\otimes\xb:M$, where $M=(m_{kj})$ denotes a constant matrix and $\otimes$ denotes the tensor product, we have  
\begin{equation}\label{err0}
          \mcLdh u(\xb_\ib)=\mcLd u(\xb_\ib)=\sum_k m_{kk} \quad\text{for all } \ib.
\end{equation}
Quadratic exactness plays a vital role in the analysis of the AC  property of quadrature-based finite difference schemes. It leads to the {\em uniform consistency result},
{which} says that \eqref{eq:FD_multi_linear} is an $O(h^2)$ approximation of the nonlocal diffusion operator independent of the parameter $\del$, which means that the truncation error is independent of $\del$ for small $\del$.

\paragraph{Uniform truncation error.} 
Assume that  $u\in C^4(\overline{\omg_\del})$. Then, for all $\xb_\ib \in \mcTh\cap \omg $, we {have}
\begin{equation}\label{err1}
		 |\mcLdh u(\xb_\ib)-\mcLd u(\xb_\ib)|\leq \CCC|D^4u|_\infty h^{2},
\end{equation}
where $\CCC$ denotes a constant independent of $\del$ and $h$.

With this observation, it then follows that $u_{\del,h}$ approximates $u_0$ at the rate $O(h^2+\del^2)$, once numerical stability is established. 

Another major feature of \eqref{eq:FD_multi_linear} is that it satisfies a discrete maximum principle, thus it is a stable numerical scheme. In fact, \eqref{eq:FD_multi_linear} can be rewritten as 
\begin{equation}\label{eq:FD_multi_linear_detailed}
	\mcLdh u(\xb_\ib) =
	\sum_{\xb_\jb\in {B}_\del(\xb_\ib)} a_{\ib,\jb} (u(\xb_\jb)-u(\xb_\ib)) ,
\end{equation}
where each $a_{\ib,\jb}$ is a nonnegative number given by
\begin{equation}\label{eq:entries}
a_{\ib,\jb}=\ffrac{2}{W(h(\jb-\ib))}\int_{ \mathcal{B}_\del( \bm{0} ) }
\varphi_\jb(\zb+h\ib) W(\zb) \gamma_\del(|\zb|) \D \zb ,
\end{equation}
with $\varphi_\jb$ being the piecewise multilinear basis function centred at $\xb_\jb=h \jb$; $\varphi_\jb$ satisfies $\varphi_\jb(\xb_\ib)=0$ when $\ib\neq \jb$ and
$\varphi_\jb(\xb_\jb)=1$. 
The discrete maximum principle of \eqref{eq:FD_multi_linear_detailed} can then be easily seen from the fact {that} $a_{\ib,\jb}\geq 0$ for each $\ib$ and $\jb$.

Equation \eqref{eq:FD_multi_linear_detailed} gives a second-order accurate AC finite difference scheme on uniform Cartesian grids. Higher-order finite difference schemes may also be constructed. However,
\pagebreak %%20200331
the discrete maximum principle is not satisfied for higher-order methods because $a_{\ib,\jb}$ could be negative for higher-order interpolants. In this case, one needs new techniques to study the stability of numerical schemes, which is still an open problem except for some specialized kernels; see \citeasnoun{LTTF19}.

We note that there have been other works concerning AC schemes for nonlocal models based on the strong form. For example,  {\citeasnoun{dzz18cicp}, \citeasnoun{you2019asymptotically} and \citeasnoun{ttd19mms} have studied AC schemes for coupled local and nonlocal models.
\citeasnoun{dhzz18sisc}, \citeasnoun{zyzd17cicp} {and} \citeasnoun{dzz18cicp}  presented implementation techniques and numerical experiments of AC schemes for problems defined in infinite domains    via the development of nonlocal artificial boundary conditions.} In these works, {time-dependent} nonlocal models are considered. 
Similar studies have also been made for models that are nonlocal in time and {space \cite{cdlz17chaos},} which are generalizations of models developed in \citeasnoun{dyz17dcdsb} \bo{and \citeasnoun{DuLZ:2020}}.
\roo{\citename{dyz17dcdsb}\ \citeyear{dyz17dcdsb,DuLZ:2020}} only considered nonlocal memory/history effects in time but the spatial interactions remained local.
{\citeasnoun{cdlz17chaos} also} replaced the local spatial differential operators  {with} nonlocal operators.

AC  difference approximations of nonlinear models have also been 
investigated, including approximations of
 nonlocal hyperbolic conservation laws.
From a modelling perspective, one can argue that conservation laws in 
integral forms
may be particularly more natural than their local counterpart, especially in the presence of singular solution behaviours. Important physics might be lost in the local formulation, {so} additional assumptions such as entropy conditions have to be re-introduced to maintain  validity. On the other hand, it is {possible} to introduce nonlocal conservation laws for which appropriate entropy conditions are automatically satisfied \cite{dhl17sinum}, thus retaining important physical features in the modelling process that are missing from  local models. The model {of} \citeasnoun{dhl17sinum} also
 improves the model studied {by} \citeasnoun{du12siap}. With AC discretization \cite{dh17jmra}, the numerical convergence has been demonstrated with or without singular solutions.

We also mention some other works concerning  difference approximations of nonlocal models. {For example, \citeasnoun{jha2019numerical} studied} the discretization of peridynamic models involving bond-softening.  {Regarding} the consistency between the  approximations of nonlocal models and their local limits, one can also find comparisons in \citeasnoun{BYASSX09} {and} \citeasnoun{du16cmame}.

{\citeasnoun{LTTF19} discretized the strong form of nonlocal diffusion operators by} using the reproducing kernel particle method (RKPM) \cite{WK:95}. The work shows the consistency and stability of the linear RK collocation scheme on quasi-uniform Cartesian grids. Define a quasi-uniform Cartesian grid on $\Rd$ as
\[
\mathcal{T}_\hb = \{ \xb_\jb :=  \hb \odot \jb  \colon  \jb \in \Mbzd \},
\]
where $\hb$ denotes a vector in $\Rd$ and $\odot$ denotes component-wise multiplication of vectors. We also assume that $\hb = h {\overline\hb}$, where $h\in \Ro$ and ${\overline\hb}$ is a fixed unit vector in $\Rd$, so that the convergence rate will again be given in terms of $h\in\Ro$.  RKPM provides a systematic means for generating basis functions from scattered particles such that polynomials can be exactly represented up to a certain order. Here we assume that the RK basis functions $\{ \Psi_{\jb} (\xb)\}_{\jb\in \Mbzd}$ are generated with respect to the quasi-uniform Cartesian grid $\mathcal{T}_\hb$. Then the RK interpolation of any continuous function $u$ is defined as 
\[
(\mathcal{R}_h u)(\bm{x}) := \sum\limits_{\jb \in \Mbzd} \Psi_\jb(\xb)u(\xb_\jb).
\]
{\citeasnoun{LTTF19} considered} the RK basis functions given by 
\begin{equation}\label{eq:RK_basis}
\Psi_\jb (\xb):=\prod^{d}_{k=1}\varphi\biggl(\ffrac{[\xb]_k- [\xb_\jb]_k}{[\bm a]_k}\biggr),
\end{equation}
where $[\xb]_k$ denotes the {$k$th} component of the vector $\xb\in\Rd$. The vector $\bm a \in \Rd$ is assumed to satisfy $\bm a = m \hb $ with $m>0$ being an even number. The function $\varphi$ is the cubic B-spline given by 
\[
\varphi(x) = 
\begin{cases}
\,\ffrac{2}{3}-4|x|^2+4|x|^3   & 0 \leq |x| \leq \ffrac{1}{2} , \\[9pt]
\,\ffrac{4}{3}(1-|x|)^3 & \ffrac{1}{2} \leq |x| \leq1 , \\[9pt]
\, 0 & \textnormal{otherwise}.
\end{cases}
\]
Those assumptions allow the RK approximation $\mathcal{R}_h u$ to represent multilinear polynomials exactly, and a special synchronized convergence property is satisfied \cite{dt-wkliu96}.

{Given the RK interpolation operator $\mathcal{R}_h$,} the RK collocation scheme is then defined as 
\begin{equation}\label{eq:RK_multi_linear}
\mcLdh u(\xb_\ib) = \mcL_\del (\mathcal{R}_h u) (\xb_\ib) ,
\end{equation}
for any $\xb_\ib \in \mathcal{T}_\hb \cap \omg$. One nice property about the interpolation operator $\mathcal{R}_h$ defined above is that although it only reproduces multilinear polynomials exactly, it shifts quadratic polynomials by a constant, so that $\mcL_{\del,h}$ defined in \eqref{eq:RK_multi_linear} actually satisfies the quadratic exactness result given earlier. This property and the synchronized convergence property become the crucial reasons that {cause} the scheme \eqref{eq:RK_multi_linear} to be AC. Indeed, {\citeasnoun{LTTF19} showed} that the collocation scheme \eqref{eq:RK_multi_linear} satisfies the same uniform consistency result presented {above}. The stability of the scheme, however, is a more tricky issue. In fact, all RK collocation schemes fail to satisfy the discrete maximum principle. But the linear RK basis function $\Psi_{\bf 0}$ given by \eqref{eq:RK_basis} has a special property that its Fourier transform, given by 
\begin{equation}\label{ftft}
\widehat{\Psi_{\bf 0}}(\bm{\xi})= \prod_{k=1}^{d} [\bm a]_k \cdot \widehat{\varphi}([\bm a]_k[\bm\xi]_k) = \prod_{k=1}^{d} \ffrac{[\bm a]_k}{2}\biggl(\ffrac{\textnormal{sin}([\bm a]_k[\bm\xi]_k/4)}{([\bm a]_k[\bm\xi]_k/4)}\biggr)^4 ,
\end{equation}
is always nonnegative. This result is a key observation used {by} \citeasnoun{LTTF19} to show that the strong form discretization  \eqref{eq:RK_multi_linear}
is comparable in terms of Fourier symbols with the Galerkin approximation \eqref{eq:wfscalarapprox}, with $V_{\del,h}$ being the span of RK basis functions \eqref{eq:RK_basis}. This immediately implies  the stability of the scheme \eqref{eq:RK_multi_linear}, because the standard Galerkin approximation is naturally stable. Stability issues again {prevent} us from discussing higher-order RK collocation methods because the higher-order RK basis functions fail to have purely nonnegative Fourier transforms such as one shown in \eqref{ftft}. More careful investigations are needed for the stability analysis of higher-order methods. As a last comment, we note that {\citeasnoun{LTTFNavier} have extended the analysis of \citeasnoun{LTTF19} to the peridynamics Navier system}.

\vspace{5pt}
\section[Numerical methods for the strong form of\\ fractional diffusion]{Numerical methods for the strong form of fractional diffusion}\label{sec:strf}

In contrast to Sections~\ref{sec:femintlap}, \ref{sec:femsl} and~\ref{sec:smf}, in which {we studied} the discretization of a weak formulation of a fractional diffusion model, in this section we consider discretizations of a strong formulation. We consider three such approaches: quadrature rule-based finite difference methods, Monte Carlo methods and radial basis function methods.

\subsection{{Quadrature rule-based} finite difference methods}\label{sec:strfq}

The most used class of methods for discretizing fractional diffusion models in bounded domains is {quadrature rule-based} finite difference methods, where their popularity stems from their simplicity. These {methods} directly discretize the integral fractional Laplacian that features a singular, non-integrable integrand. A common approach to {dealing} with the difficulties that singular integrals pose is to split the integral into the sum of two integrals, one isolating the singular part, for which some care should be exercised, {and the other having} a smooth integrand so that standard quadrature rules can be used to obtain accurate approximations. \bo{This approach has been widely used for solving the nonlocal diffusion models \cite{dtty19ima,TiDu13}, as reported in Section \ref{sec:FD-nonlocal}. \roo{In} the same spirit, some variants were also developed for fractional diffusion models.} 

We have the one-dimensional fractional equation and volume constraint
\begin{align}\label{eqn:fde5}
\begin{cases}
\displaystyle
(-\Delta)^{s} u(x) = \CCC_{1,s} \int_{-\infty}^\infty \ffrac{u(x) {-u(x+z)}}{|z|^{1+2s}}\D {z} =f(x)   &\text{for all } x \in \omg ,\\
 u(x)= 0 & \text{for all } x \in\omgii ,
\end{cases}
\end{align}
with $s\in(0,1)$, \bo{$\omg= (0,L)$} with $L>0$, and $\omgii= \Ro\setminus\omg$. \bo{To discretize the problem, we use the uniform Cartesian grid $\{x_j=jh \colon j\in\mathbb{Z} \}$
with $h=L/N$ for some integer $N$}. In \citeasnoun{HuangObreman:2014}, the fractional equation in \eqref{eqn:fde5} is evaluated at the grid points $x_j$ and then split into two integrals, \tie, {for $j=0,\ldots,N$, we have}
\begin{align}\label{eqn:split-1}
(-\Delta)^s u(x_j) 
&=\underbrace{\CCC_{1,s}\int_{|z|< h}   \frac{u(x_j)-u(x_j+z)}{|z|^{1+2s}} \D z}_{I_1}\notag 
\\[3pt]
&\qquad+
\underbrace{\CCC_{1,s} \int_{|z|>h} \frac{u(x_j)-u(x_j+z) }{|z|^{1+2s}} \D z}_{I_2} .
\end{align} 
For smooth $u$ 
the use of Taylor series and central difference scheme leads to
\begin{equation}\label{eqn:I1}
I_1  = - \frac{\CCC_{1,s} }{2-2s}\biggl(\frac{u(x_{j+1})-2u(x_j)+u(x_{j-1})}{h^{2s}}\biggr) + O(h^{4-2s}).
\end{equation}

The second term $I_2$ in \eqref{eqn:split-1} has a regular integrand so that term can be well {approximated using} standard quadrature rules. As an example, let $I_h^k$ denote the Lagrange interpolation operator for piecewise polynomials of order less than or equal to $k$. \bo{Then, for smooth $u$,}
\begin{equation}\label{eqn:I2}
\bo{I_2  = \CCC_{1,s} \int_{|z|>h} \frac{\mathcal{I}_h^k ( u(x_j)-u(x_j+z))}{|z|^{1+2s}} \D z +  O(h^{k+1-2s}).}
\end{equation}
Then \eqref{eqn:I1} and \eqref{eqn:I2} result in a discrete approximation of the fractional Laplacian:
\begin{equation}\label{eqn:I3}
(-\Delta)^s u(x_j) = D_h u(x_j) + O(h^{\min(4,k+1)-2s}),
\end{equation}
where the discrete operator $D_h$ is defined by
\begin{align}\label{eqn:I4}
D_h u(x_j) &:=  \ffrac{\CCC_{1,s} }{2-2s} \ffrac{-u(x_{j-1})+2u(x_j)-u(x_{j+1})}{h^{2s}} \notag \\*
&\quad\, + \CCC_{1,s} \int_{|z|>h} \frac{\mathcal{I}_h^k ( u(x_j)-u(x_j+z) )}{|z|^{1+2s}} \D z.
\end{align}
As a result, the corresponding quadrature-based finite difference scheme is given as follows: find $U_j$, $j\in \mathbb{A}= \{1,\ldots,N-1\}$, such that
\begin{equation}\label{eqn:I5}
D_h U_j  = f(x_j) \quad \text{with $ U_j = 0$ for all $  j\in \mathbb{Z}\backslash \mathbb{A}$.}
\end{equation}
\bo{Similar to discussions given in Section \ref{sec:FD-nonlocal}, \roo{a} maximum principle for the above scheme (and hence its stability) can be proved 
  \roo{(see \citeb{HuangObreman:2014} for $k=1$ and $2$)}},  by using a barrier function method, and therefore the error estimate follows immediately.
\bo{\citename{Endal:2018} \citeyear{Endal:2018,Endal:2019} extended the aforementioned approach to higher-dimensional problems, and improved the convergence rate by using adapted vanishing viscosity approximation.}

\bo{\citeasnoun{MindenYing:2018} discussed the discretization in two and three dimensions  by applying the `window' function $\www(z):=\www(|z|)$ such that $1-\www(z)=O(|z|^p)$ as $z\rightarrow0$ for some positive integer $p$.}  
Then the following splitting is applied:
\begin{align*}
&(-\Delta)^s u(x_j) 
  \\*
  &\quad = \cds \int_{\Rd} \ffrac{u(x)-u(y) + \www(|x-y|)\sum_{1\le |\beta|\le 3}D^\beta u(x)(y-x)^\beta/|\beta|!}{|x-y|^{\ddd+2s}}  \D y\\*
&\quad\quad\, -\cds\int_{\Rd} \ffrac{\www(|x-y|)\sum_{1\le |\beta|\le 3}D^\beta u(x)(y-x)^\beta/|\beta|!}{|x-y|^{\ddd+2s}}  \D y.
\end{align*}
Note that the second term is regular and {can thus} be easily approximated by standard quadrature rules. On the other hand, by using a window function, the first term can be evaluated by the trapezoidal rule. If the grid is uniform, then the resulting scheme provides a discrete operator that is {translation-invariant} and therefore an FFT can be used \cite[Section~3]{MindenYing:2018} for {an} efficient solution.

{\citeasnoun{Duo:JCP2018} and \citeasnoun{DuoZhang:2019} suggested another splitting approach,} namely
\[
 \ffrac{u(x) - u(y)}{|x-y|^{d+2s}} = (|x-y|^{-d+\beta-2s})\, ((u(x) - u(y))|x-y|^{-\beta})
\]
for an appropriate parameter $\beta\in(2s,2)$. The first term is singular but can be evaluated analytically, whereas the second term is regular so it can be approximated {using the trapezoidal rule, \feeg}.

The discrete systems obtained by the quadrature-based finite difference methods introduced above can be Toeplitz-like, which can be tackled by using an FFT, but only for shift-invariant discretizations that usually require the use of uniform grids.
As we mentioned in Section~\ref{fracreg}, the solution of the fractional Poisson equation \eqref{eqn:fde5} is weakly singular near the boundary, and hence graded meshes are preferred. In this case, {discretization} results in dense and unstructured matrices, and hence the storage complexity will be $O(N^2)$ whereas the computational complexity will be $O(N^3)$. One promising strategy for efficiently solving the nonlocal problem \eqref{eqn:fde5} with unstructured meshes is to apply hierarchical matrices so that only linear storage and computational complexity are required; see \citeasnoun{XuDarve:2018H}.

In one dimension, {problem} \eqref{eqn:fde5} is closely related to the fractional boundary problem involving the Riemann--Liouville fractional derivative or Riesz derivative  \cite{Podlubny:book,Ortigueira:2006}. For such problems, many
finite difference methods have been studied \cite{Meerschaert:2004,TMS:2006,CelikDuman:2012,Meerschaert:2006,ChenDeng:2014,Ding:2015JCP} and fast algorithms \cite{PanNgWang:2016,WangBasu:2012,WangDu:2013,WangDu:2014,ZhangSun:2015,LeiSun:2013}.  See also the comprehensive survey {by} \citeasnoun{LiChen:2018} and the monograph {by} \citeasnoun{George:Book2019} for details about these {kinds} of models and their numerical approximation.

\subsection{Monte Carlo method by Feynman--Kac formula}\label{sec:strfm}

Since its formulation in \citeasnoun{Kac:1949}, the Feynman--Kac formula has been a powerful tool for both theoretical reformulations and practical simulations of local PDEs. The link it establishes between PDEs and related stochastic processes can be exploited to develop Monte Carlo methods that are effective, especially for the numerical simulation of high-dimensional problems \cite{Curtiss:1950,Muller:1956,Sabelfeld:1994,YanCaiZeng:2013}.

{\citeasnoun{Kyprianou:2018} have developed a Monte Carlo algorithm} for approximately solving the fractional Poisson problem with inhomogeneous volume constraints given by
\begin{equation}\label{eqn:fde6}
\begin{cases}
(-\Delta)^{s} u(\xb) =f(\xb) &\text{for all } \xb\in\omg ,\\
 u(\xb)= g(\xb)  &\text{for all } \xb\in\omgii.
\end{cases}
\end{equation}
The solution representation of \eqref{eqn:fde6} has been extensively studied from both the analytical and probabilistic perspectives\footnote{\bo{The connection between stochastic processes and {\em fractional} equations has been extensively investigated in, \eg, \citeasnoun{Meerschaert2012book}. In \roo{contrast}, \citeasnoun{Burch2014exit-time}, \citeasnoun{Burch2015sym-stoch} and \citeasnoun{DElia2017nlcd} discuss the connection between finite-range jump processes and {\em nonlocal} (truncated) equations.}}; see \eg\  \citeasnoun{BogdanByczkowski:1999}, \citeasnoun{Bucur:2016}, \citeasnoun{Ros-Oton:2016} {and} \citeasnoun{Kyprianou:2018}. In particular, we assume that $\omg$ is bounded, that $f\in C^{2s+\epsilon}(\overline\omg)$ for some $\epsilon>0$, and that $g$ is continuous function belonging to $L_{2s}^1(\Rd\setminus\omg)$, \roo{\tie,}
\[
 \int_{\Rd\setminus\omg} \ffrac{|g(\xb)|}{1+|\xb|^{2s+\ddd}}\D x <\infty.
\]
Then the solution of the fractional Poisson problem \eqref{eqn:fde6} has the probabilistic representation \cite[Theorem~6.1]{Kyprianou:2018}
\begin{equation}\label{eqn:FK}
 u(\xb) = \mathbb{E}_\xb[g(X_{\sigma_\Omega})] + \mathbb{E}_\xb\biggl[ \int_0^{\sigma_\Omega}f(X_s) \D s \biggr],
\end{equation}
{with $\mathbb{E}_\xb(\cdotsp)$ denoting the expected value \bo{conditioned on $X_0=\xb$} and where $X=(X_t,t\ge0)$ is an isotropic $\alpha$-stable L\'evy process with index $\alpha=2s$ and $\sigma_\Omega = \inf\{ t>0:~X_t\notin\Omega \}$
 is the first exit time of $X_t$ from $\Omega$. In particular, let $B(0,1)\subset \mathbb{R}^d$ be a unit ball centred at the origin, and then 
for $|{\yb}|>1$ {we have} \cite[Theorem~A]{Blumenthal:1961}
\begin{equation}\label{eqn:P0}
\mathbb{P}_0(X_{\sigma_{B(0,1)}}\in  \D \yb) =  \ffrac{\Gamma(d/2)\,\sin(\pi\alpha/2)}{\pi^{d/2+1}}\, |1-|\yb|^2|^{-\alpha/2}|\yb|^{-d} \D \yb.
\end{equation}
This gives the distribution of the stable process that begins from the origin when it first exits the unit ball.}

The formula \eqref{eqn:FK} is a nonlocal analogue to the Feynman--Kac formula of the classical Laplacian, but the role of $\ddd$-dimensional Brownian motion
is replaced by an isotropic $\alpha$-stable L\'evy process in $\Rd$. One significant difference {to Brownian} motion
is that the $\alpha$-stable L\'evy process will exit the domain by a jump rather than hitting the boundary.

Then the solution to \eqref{eqn:fde6} can be computed numerically by applying a {Monte Carlo} sampling method, based on the Feynman--Kac formula \eqref{eqn:FK}.
That is to say, if $(X_t^i,t\le \sigma_\Omega^i)$ are
{independent and identically distributed} copies of $(X_t,t \le \sigma_\Omega)$ issued from ${\xb} \in \Omega$,
then by the law of large numbers, {we obtain}
\begin{equation}\label{eqn:MC}
 u({\xb}) = \lim_{n\rightarrow\infty} \ffrac1n\sum_{i=1}^n \biggl(g(X_{\sigma_\Omega^i}^i) +  \int_0^{\sigma_\Omega^i}f(X_s^i) \D s  \biggr),
\end{equation}
where sample paths of the stable process could be constructed by  \eqref{eqn:P0}.
However, the direct application of a Monte Carlo method based on the Feynman--Kac formula is inefficient,
because the evaluation of $u({\xb})$ requires the simulation of a very large number of paths of $X_t$, beginning from the single point ${\xb}\in \Omega$.

{One way to speed things up is to apply the {\em walk-on-spheres} strategy (WOS), that was first established in \citeasnoun{Muller:1956} for the Laplace equation
and then adopted in \citeasnoun{Kyprianou:2018} for the fractional Laplace equation. The WOS method
\pagebreak %%20200331
does not require a complete simulation of the entire path and takes advantage of
the distributional symmetry of the stochastic process $X_t$.
Next we briefly discuss the WOS method {for solving} the fractional diffusion problem \eqref{eqn:fde6} within a convex domain $\Omega\subset\mathbb{R}^d$, $d \geq 2$.} 

{With an arbitrary $\xb\in\Omega$, we let $\rhob_0 = \xb$ and define $r_1$  to be the radius of the largest sphere inscribed in $\Omega$ that is centred at $\rhob_0$.
  Then we set
  $  B_1 = \{\xb\in\mathbb{R}^d\colon |\xb - \rhob_{0}|< r_1\}  $
  and select $\rhob_1$ according to $X_{\sigma_{B_1}}$ under $\mathbb{P}_{\rhob_{0}}$, which is known from  \eqref{eqn:P0}.
Repeating the argument, we construct a sequence $\{\rhob_n\}_{\,n\geq 1}$ inductively.
The algorithm ends at the random index $N = \min\{n\geq 0\colon \rhob_n\notin \Omega\}$. 
As a result, the Feynman--Kac formula \eqref{eqn:MC} can be replaced by \cite[Corollary~6.4]{Kyprianou:2018}
\begin{equation*}
u(\xb) =\mathbb{E}_{\xb}[g(\rhob_{N})] 
+
\mathbb{E}_{\xb} \Biggl[
\sum_{n=0}^{N-1} r_n^{\alpha} \int_{|\yb|<1} f(\rhob_n + r_n \yb ) \,V_1(0, \dd \yb)\Biggr],
\end{equation*}
where $V_1(0, \dd \yb)$ denotes the expected occupation measure of the stable process prior to exiting a unit ball centred at the origin, which is given 
for $|\yb|<1$ by \cite[Corollary~4]{Blumenthal:1961}
\begin{equation*}
V_1(0, \dd \yb)= 2^{-\alpha}\,\pi^{-d/2}\, \ffrac{\Gamma(d/2)}{\Gamma(s)^{2}}\,            
			|\yb|^{2s -d}\, \biggl(\int_0^{|\yb|^{-2}-1}(u+1)^{-d/2}u^{s-1} \D u\biggr) \D \yb. 
\end{equation*}
We refer interested readers to \citeasnoun{Kyprianou:2018} for more details about convergence analysis, discussion of nonconvex domain and implementation of WOS.}

A multilevel Monte Carlo method based on the WOS algorithm was proposed {by} \citeasnoun{Shardlow:2019} for efficiently computing the solution for the entire domain $\Omega$. {\citeasnoun[Section~5]{Shardlow:2019} drew comparisons between} the  multilevel WOS Monte Carlo method and the adaptive finite element method of \citeasnoun{AG:2017}.

\subsection{Radial basis function methods}\label{sec:strfr}

The development, analysis and implementation of mesh-free methods using radial basis functions (RBFs) has been studied thoroughly for the integer order (local) PDEs; see {\citename{Buhmann:Acta} \citeyear{Buhmann:Acta,Buhmann:Book}} and the references cited therein. Compared to finite difference methods, the RBF method is seemingly easier to implement for approximately solving the nonlocal problem \eqref{eqn:fde5}, especially in high {dimensions}, because only minor modifications to existing algorithms for local PDEs are needed.

However, studies of RBF methods for {the} fractional diffusion model {are} fairly scarce.
{\citeasnoun{Rosenfeld:JCP2019} used RBF interpolants to} approximate
\pagebreak %%20200331
the fractional Laplacian of a given function through a mesh-free pseudo-spectral method.
Specifically, they used the  compactly supported Wendland functions $\Phi(\xb) =  \phi_{\ddd,k}(|\xb|)$
defined in
{\citename{Wendland:1995} \citeyear{Wendland:1995,Wendland:Book}}, where $k$ is a smoothness parameter. With the collection of points
\[ %%20200331
\roo{X=\{\xb_1,\xb_2\ldots,\xb_m \}\subset \omg,}
\]
the mesh size $h$ of $X$ is defined by
\[
h_X=\sup_{\xb\in\omg}\min_{\xb_j\in X} \| \xb-\xb_j\|_2.
\]
Then, for any function $f\in N_\Phi(\Rd)$, where $N_\phi(\Rd)$ denotes the native space of the RBFs \cite{Wendland:Book}, the interpolant of $f$ is defined by
\[
I_X f(\xb) = \sum_{j=1}^n w_j \Phi(  \xb - \xb_j  )\quad \text{such that $I_X f(\zb) = f(\zb)$ for all $ \zb\in X$.}
\]
Suppose that $\Phi \in C^{2k}(\Rd)$ is symmetric and strictly positive definite, and $f\in N_\Phi(\Rd)$ is compactly supported. Then, by estimating the inverse Fourier transform,  the following interpolation error result is proved  in \citeasnoun[Proposition~3.2]{Rosenfeld:JCP2019}:
\[
|(-\Delta)^s (I_X f)(x) -  (-\Delta)^s f(x) | \le (c h_X^{k-2\beta} \| f \|_{N_\phi(\Rd)}  + E),
\]
where the parameter $\beta\in\mathbb{N}$ satisfies $2\beta<k$ and $2\beta-2s>n$. Here $E$ denotes the residual error that can be further bounded by applying a smooth cut-off function \cite[Theorem~3.1]{Rosenfeld:JCP2019}. Numerical experiments showed no measurable difference in the resulting estimations when using the {cut-off function, compared} to when the {cut-off} function is not used \cite[Section~4.1]{Rosenfeld:JCP2019}. However, the optimal convergence rate of the method without {cut-off functions remains} an open question.

The above interpolation property inspires the future study of collocation methods (or Galerkin methods) using RBFs for approximately solving the fractional diffusion problem \eqref{eqn:fde5}. Many interesting questions, \feeg\ the selection of {suitable} RBFs, {the} stability of numerical methods, optimal convergence {rates and adaptive algorithms,}
are largely open and warrant further investigation. In addition, for high-dimensional problems, an additional challenge stems from the computation of the fractional Laplacian of the basis functions, which requires fast and accurate numerical approximation. See also \citeasnoun{Lehoucq:2016} {and} \citeasnoun{Lehoucq:2018} for a Galerkin method using RBFs for solving nonlocal diffusion with \bo{kernel functions that are both radial and integrable}, and \citeasnoun{PangChenFu:2015} {and} \citeasnoun{SunLiuZhangPang:2017} for collocation methods using RBFs for solving fractional diffusion with {Riemann--Liouville-type} fractional derivatives.

%%\vspace{5pt} %%20200331
\section{Conditioning and fast solvers}

The effectiveness of nonlocal modelling and simulations relies on the effective solution of the algebraic systems resulting from the discretization of nonlocal models. Thus a good theoretical understanding of the conditioning of the stiffness matrices for nonlocal problems is important. Results in this direction {have} been provided, using the  Fourier analysis of the point spectrum for nonlocal operators, as given in \citeasnoun{du10sinum}; see also \citeasnoun{AkMe10}, \citeasnoun{AkPa11}, \citeasnoun{AU14}, \citeasnoun{du12sirev} {and} \citeasnoun{SPGL09} for additional  discussions. 

It is {well known that} a typical local diffusion model yields a condition number of $\Bigoh (h^{-2})$ for a discretization having a meshing parameter $h$. The corresponding nonlocal models have condition numbers that depend, in general, on both $\delta$ and $h$. For example, {\citeasnoun{dz17mg} give sharp lower and upper bounds for} the condition number of the stiffness matrix corresponding to a finite element discretization of a nonlocal diffusion operator based on a quasi-uniform regular  triangulation. Typically, if a {fractional-type} {translation-invariant} kernel such as those satisfying \eqref{frackernel} is used to describe  nonlocal interactions, the condition number of the resulting discretized nonlocal diffusion model {is} of
$\Bigoh (h^{-2s} \delta^{2s-2})$,   where, as always, $\del$ denotes the size of the horizon parameter. In practice, both $h$ and $\del$ could be small parameters, {so the view that nonlocal models yield better conditioned  systems than  their local counterparts should be taken with due care}. Additional detailed estimates on the conditioning of nonlocal stiffness {matrices} can be found in \citeasnoun{AU14}.

For effective algebraic solvers of the resulting linear system, we refer to studies {of} the use of Toeplitz \cite{WaTi12,vollman} and multigrid  solvers \cite{dz17mg}.

\subsection{Fast algorithm for kernels with non-smooth truncation}

\begin{figure} %%fig12.1
  \centering
\subfigure[]{\includegraphics[width=110pt]{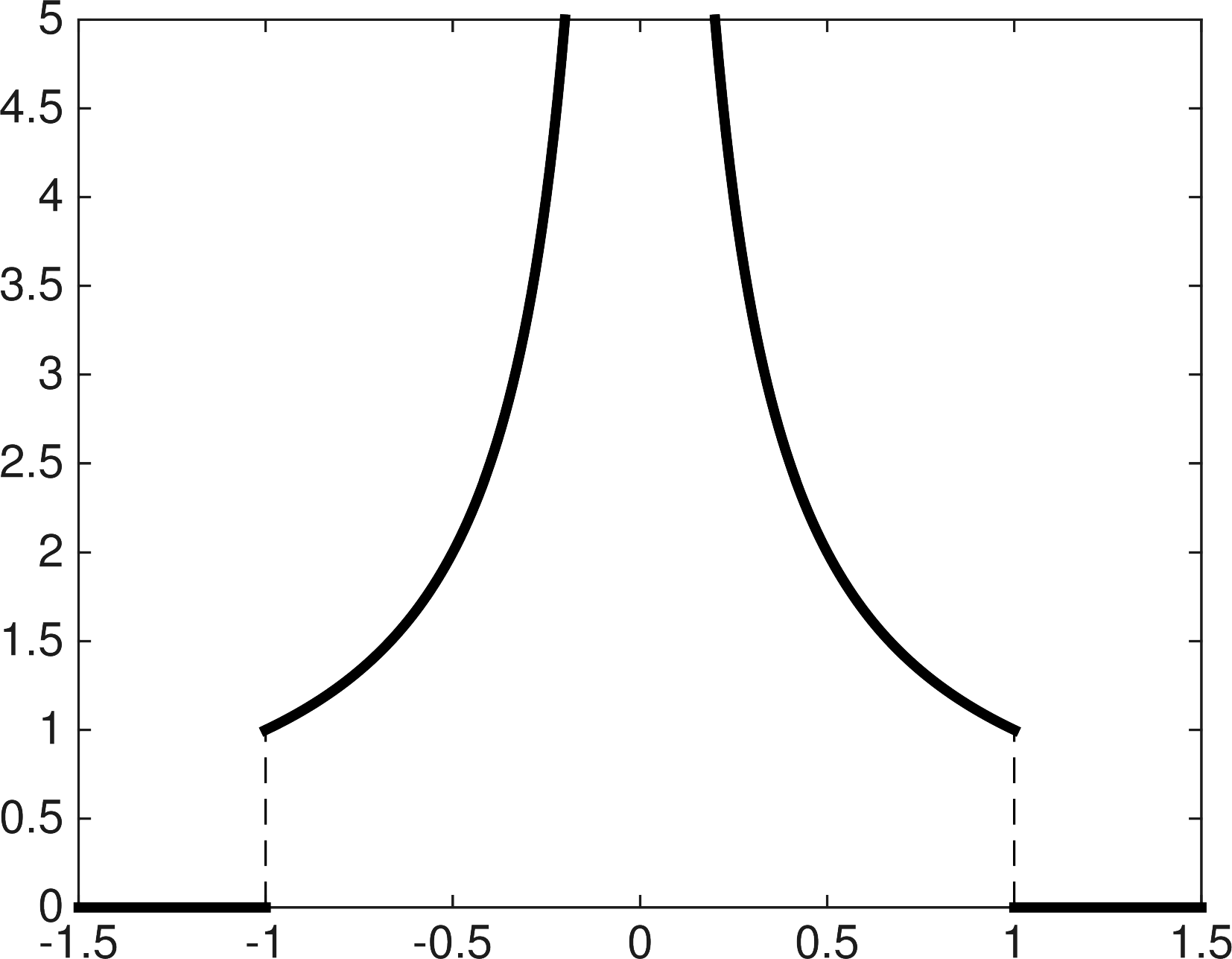}}\hfill %%kernel_ga
\subfigure[]{\includegraphics[width=110pt]{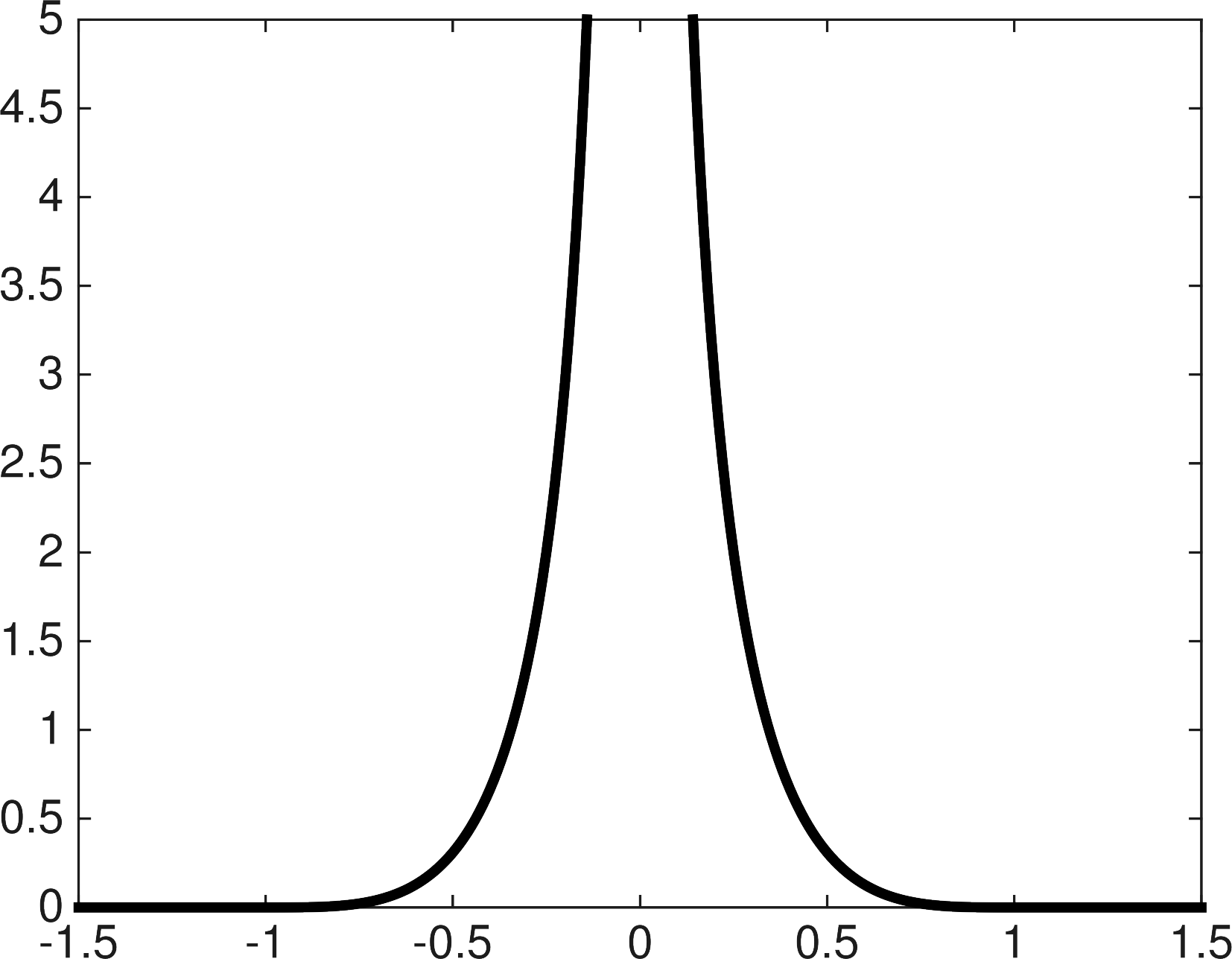}}\hfill %%kernel_ka
\subfigure[]{\includegraphics[width=110pt]{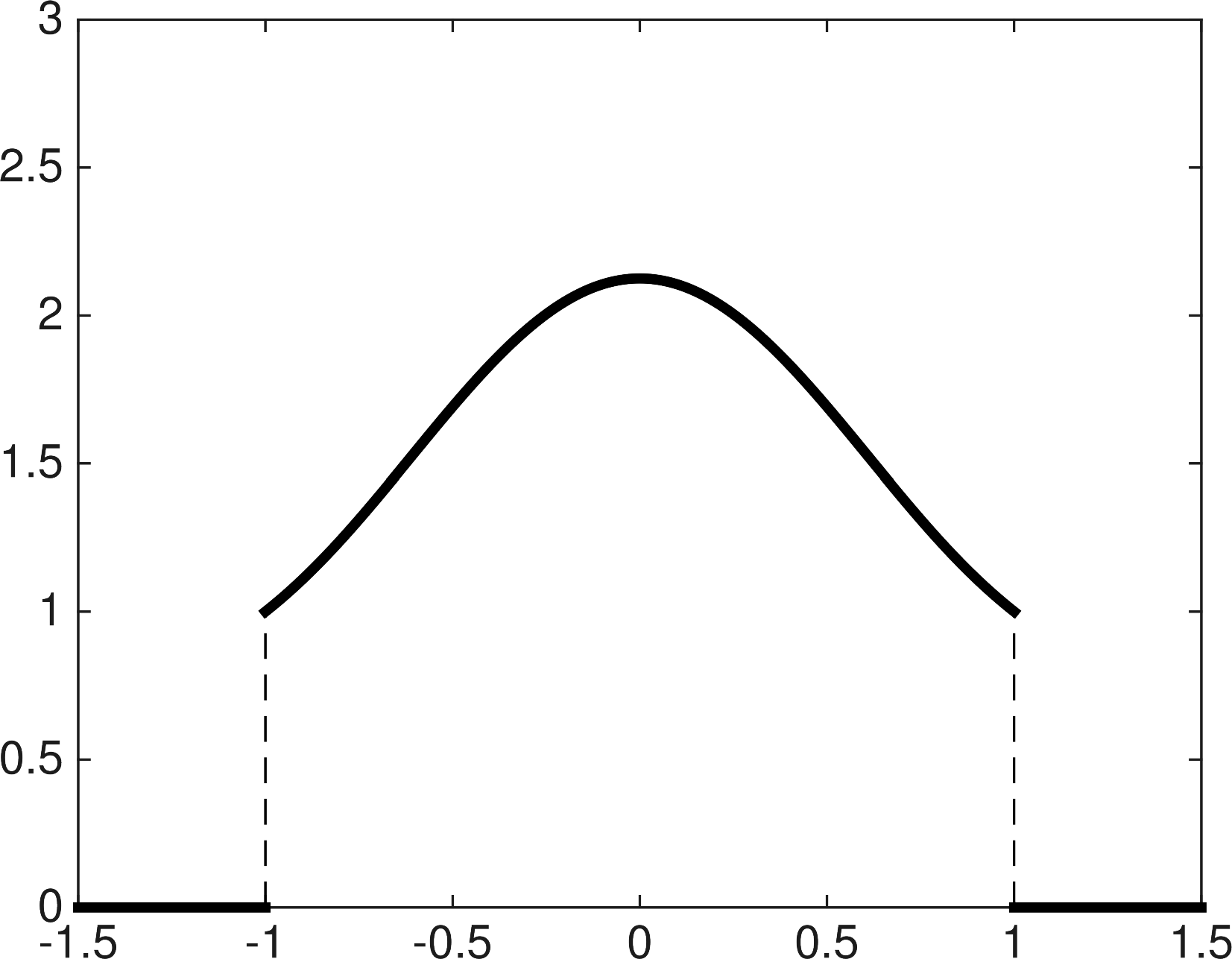}}\\ %%kernel_p
\parbox{290pt}{\caption{The kernel $\gamma(x)$~{(a)} splits into $\gamma_1(x)$~{(b)} and $\gamma_2(x)$~{(c)}.} \label{fig:tian2019fast_kernel}}
    \end{figure}

\begin{figure} %%fig12.2
\centering
\subfigure[]{\includegraphics[height=145pt]{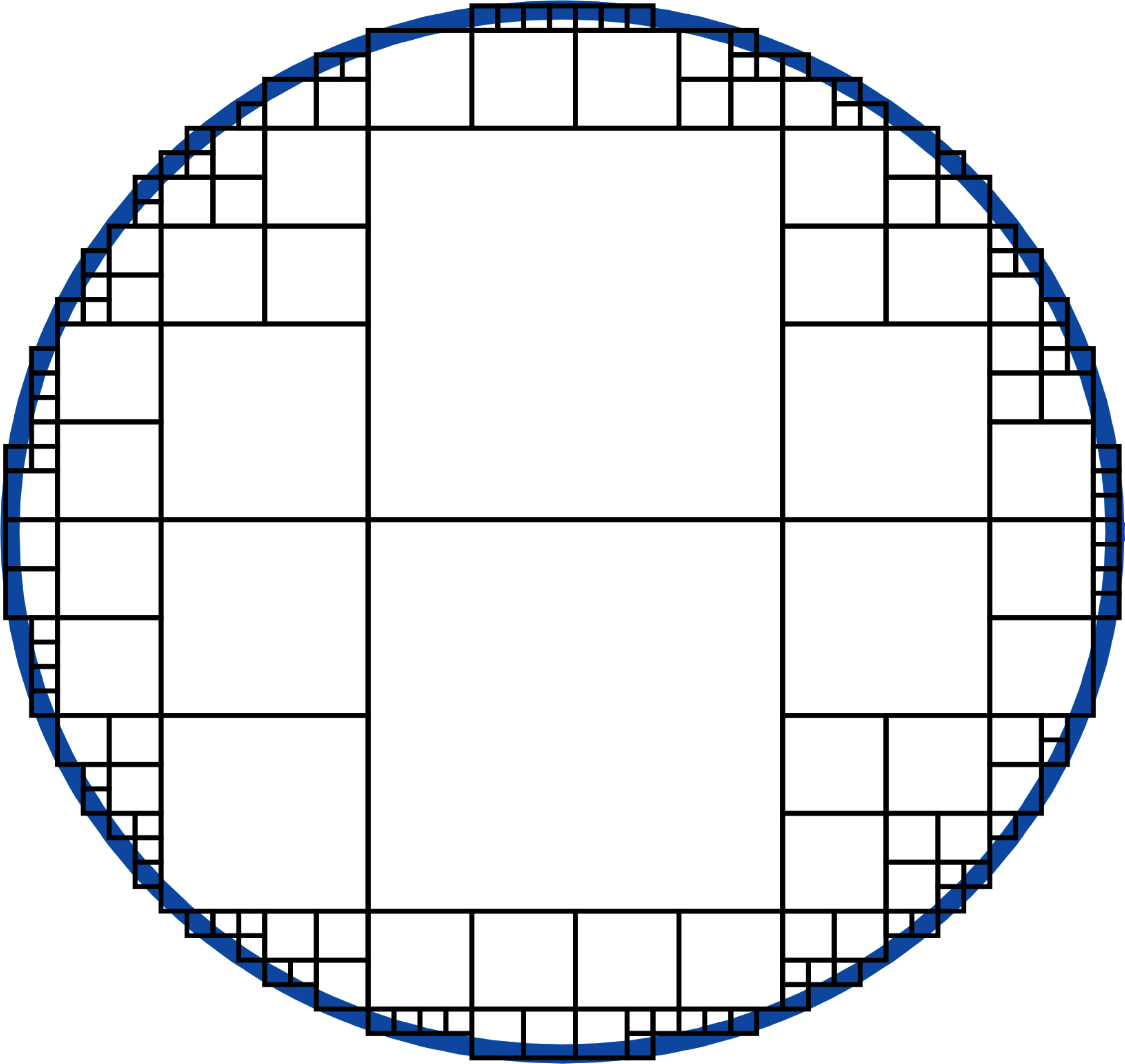}}\hspace{48pt} %%partition
\subfigure[]{\includegraphics[height=145pt]{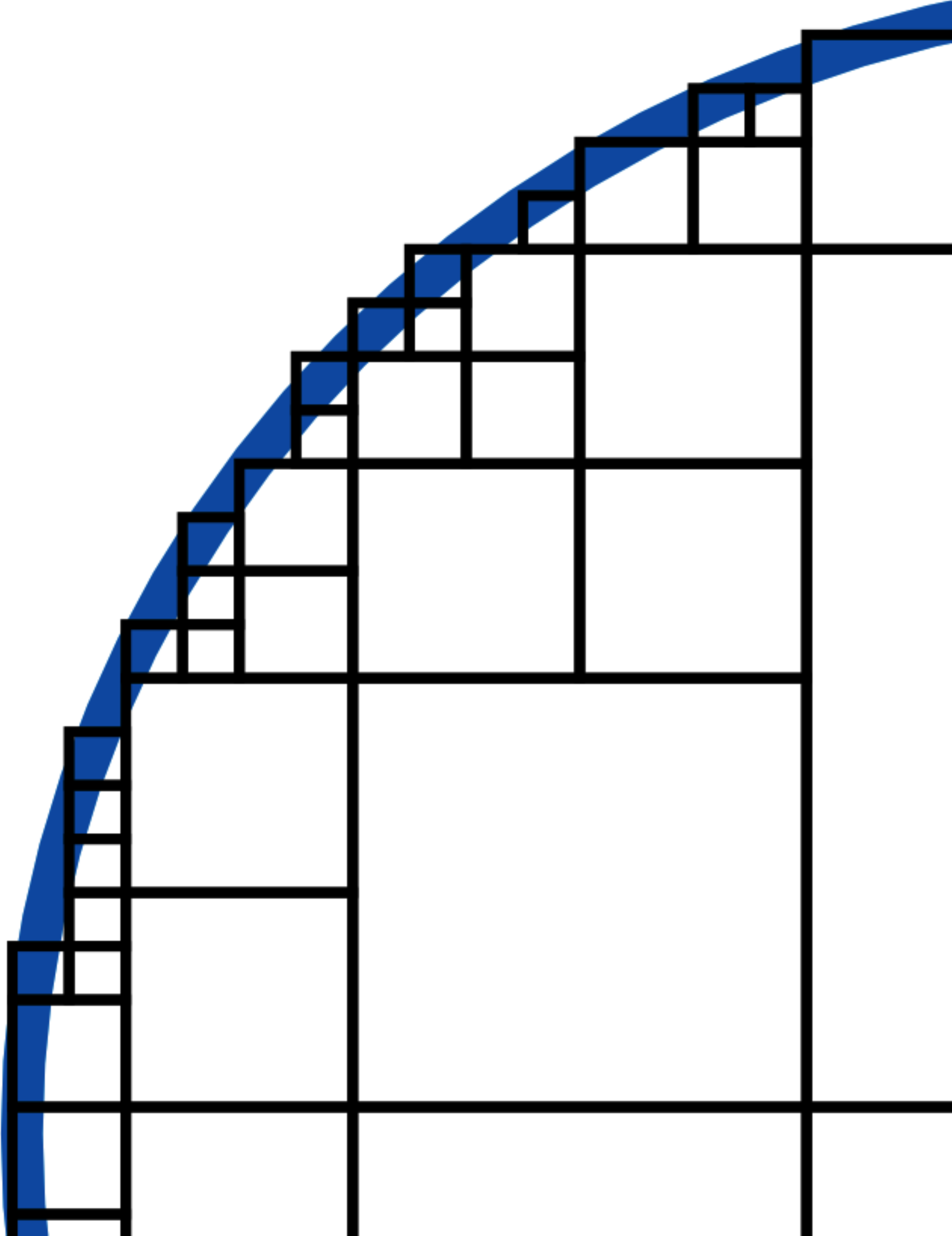}}\\ %%partition2
\caption{The circular region {represents} the interaction region corresponding to a point located at its centre. The decomposition of the circular region into hierarchical boxes is illustrated, with {plot~(b) being a close-up of plot~(a)}.}
 \label{fig:tian2019fast_decomposition}
\end{figure}

{\citeasnoun{tian2019fast} have developed methods} based on fast multipole methods (FMMs) and hierarchical matrix techniques.  {Their key observation 
  is} that the non-smooth transition of the kernel function typically used in peridynamics and nonlocal diffusion models {can reduce} the effectiveness of many standard fast solvers that are based on the compression of {far-field} interactions. A typical kernel  $\gamma$ used in practice, shown in \fo{Figure~\ref{fig:tian2019fast_kernel}}{(a)}, has a singularity at origin {and its} nonlocal interaction is truncated at a finite distance. The kernel is then decomposed into two parts, $\gamma_1$ and $\gamma_2$, as illustrated in  \fo{Figure~\ref{fig:tian2019fast_kernel}}{(b,\,c)}. The first part ($\gamma_1$) is smooth away from {the origin, so fast} solvers using FMMs or hierarchical matrix techniques for the compression of {far-field} interactions can be successfully applied.  The second part ($\gamma_2$) is smooth inside {the support.}  {\citeasnoun{tian2019fast} developed an FMM-type} algorithm for the fast evaluation of nonlocal operators with such a kernel function.  The key idea is to compress the nonlocal interaction away from {the} boundary of the support. This idea is depicted in \fo{Figure~\ref{fig:tian2019fast_decomposition}}, \roo{where the} geometric boundary of the interaction kernel centred at a certain point  is finely resolved by small boxes, whereas away from the {boundary large} boxes are used because there the kernel function is smooth. Because the number of small boxes needed to resolve co-dimension~$1$ surfaces increases with dimension, the complexity of the algorithm also increases with dimension. {\citeasnoun{tian2019fast} showed} that the optimal complexity $O(N\log N)$ can be achieved for $N$ unknowns in one dimension. In higher dimensions there is algebraic complexity $O(N^{2 -1/d})$, where $d$ is the spatial dimension of the problem.   
   
\subsection{Conditioning and solvers for finite element discretizations of the integral fractional Laplacian model}\label{sec:femslcon}

{\citeasnoun{AinsworthMcleanEtAl1999_ConditioningBoundaryElementEquations} gave the following results for} the stiffness matrices resulting from  finite element methods for the integral fractional Laplacian model.

  For $s<\ddd/2$ and a family of {shape-regular} triangulations $\mcTh$ with minimal and maximal element size $h_{\min}$ and $h$, the spectrum of the stiffness matrix $\mat{A}_{s}$ satisfies
 \[
    c N_h^{-2s/d}h_{\min}^{\ddd-2s}\mat{I} \leq \mat{A}_{s} \leq Ch^{\ddd-2s}\mat{I}
    \] and \[
    c N_h^{-2s/d}\mat{I} \leq  (\mat{D}^{s}) ^{-1/2}\mat{A}_{s}  (\mat{D}^{s}) ^{-1/2} \leq C\mat{I},
\]
  where $\mat{D}^{s}$ denotes the diagonal part of $\mat{A}_{s}$.
Moreover, the condition number of the stiffness matrix satisfies
  \begin{align*}
    \kappa (\mat{A}_{s})  = C \biggl(\ffrac{h}{h_{\min}}\biggr)^{\ddd-2s}N_h^{2s/d}
    \quad\text{and}\quad
    \kappa ( (\mat{D}^{s}) ^{-1}\mat{A}_{s})  = CN_h^{2s/d}.
  \end{align*}

  If an implicit time-stepping scheme is used {for the} fractional-order heat equation \roo{\eqref{eq:diffusion-z-RHS2}},
  systems having matrices of the form
\begin{align*}
  \mat{M}+\Delta t\mat{A}_{s}
\end{align*}
need to be solved, where $\mat{M}$ denotes the usual finite element mass matrix.
{\citeasnoun{AG:2017} have given} the following result about the condition number of $\mat{M}+\Delta t\mat{A}_{s}$.

  For a {shape-regular} and globally quasi-uniform family of triangulations $\mcTh$ and for a time step $\Delta t \leq 1$,
  \begin{align*}
    \kappa (\mat{M}+\Delta t\mat{A}_{s}) 
    &\leq C \biggl(1+\ffrac{\Delta t}{h^{2s}}\biggr).
  \end{align*}
  More generally, for a family of triangulations that is only locally quasi-uniform and $\Delta t\leq h_{\min}^{2s}N_h^{2s/d}$, {we have}
  \begin{align}
    \kappa (\mat{M}+\Delta t\mat{A}_{s})
    &\leq C \biggl(\ffrac{h}{h_{\min}}\biggr)^{\ddd}  \biggl(1+\ffrac{\Delta t}{h^{2s}}\biggr). 
  \end{align}
  If $\mat{D}^{0}$ is taken to be the diagonal part of the mass matrix and $\Delta t \leq h^{2s}N_h^{2s/d}$, then
  \begin{align}
    \kappa ( (\mat{D}^{0}) ^{-1} (\mat{M}+\Delta t \mat{A}_{s}) ) 
    &\leq C\biggl(1+\ffrac{\Delta t}{h_{\min}^{2s}}\biggr) .
  \end{align}
{These results show} that for small fractional order $s$, if the time step $\Delta t$ is chosen small enough with respect to the mesh size, the conjugate gradient method with diagonal preconditioning will converge in a fixed number of steps.
For larger fractional orders or for the steady-state problem, the number of iterations depends on the problem size.
It has been shown that by applying a multigrid solver, one can restore a uniform bound on the number of iterations; see \citeasnoun{SauterSchwab2010_BoundaryElementMethods},
{\citename{Hackbusch1985_MultiGridMethodsApplications} \citeyear{Hackbusch1985_MultiGridMethodsApplications,Hackbusch1994_IterativeSolutionLargeSparseSystemsEquations}},
\citeasnoun{AinsworthMclean2003_MultilevelDiagonalScalingPreconditioners} {and \citename{AG:2017} \citeyear{AG:2017,AinsworthGlusa2018_TowardsEfficientFiniteElement}}.

We also observe that the system matrix $\mat{A}_{s}$ is entirely dense, owing to the nonlocal interactions.
This means that for the efficient solution, efficient techniques for computing matrix--vector products with $\mat{A}_{s}$ need to be used.
In the literature, fast transforms and matrix compression \cite{AinsworthGlusa2018_TowardsEfficientFiniteElement,AG:2017} have been explored.
The drawback of the former is their limitation to uniform meshes, whereas the latter are more difficult to implement.
Both approaches lead to quasi-optimal complexity, \ie\ $\Bigoh(N_h\log N_h)$ operations to solve the system.

We note that a brief discussion about solvers for  extended fractional Laplacian problems is given in Remark~\ref{femslls}.

\vspace{5pt}
\part{Selected extensions}\label{part3}

So far in the article we have {mostly} focused on steady-state nonlocal diffusion models, including fractional  models. In this part we provide brief accounts {of} the extension of the models we have considered to other settings,
{along with} additional approaches for obtaining approximate solutions of the  {models and a} few applications of the models.

\vspace{5pt}
\section[Weakly coercive, indefinite and\\ non-self-adjoint problems]{Weakly coercive, \roo{indefinite and} non-self-adjoint problems} \label{sec:variations}

The discussion about discretization schemes given in Part~\ref{part2} dealt with problems that fall into the category of what are known as {\em Rayleigh--Ritz} or {\em strongly coercive} problems. As such, the Lax--Milgram theorem is a fundamental tool in proving {well-posedness} of both continuous and discrete problems. In this section we consider more general settings involving indefinite and non-self-adjoint problems, settings for which the Lax--Milgram theorem cannot be applied.

\subsection{Indefinite and non-self-adjoint problems}\label{weakcoercive}

Instead of the symmetric and coercive bilinear form \eqref{weak-bform} (see also \eqref{peri-bform}), we now consider bilinear forms that have neither of these properties. As such, the bilinear forms considered in this section can be used in many settings that cannot be treated using the bilinear form \eqref{weak-bform}. The discussion here largely follows the formulations given in \citeasnoun{TiDu19} for a more general class of parametrized problems.

The bilinear form $\mcB_\del(\cdot,\cdot)$ is defined on a trial space ${V_\del}$  and test space ${W_\del}$ and satisfies the following requirements.
\renewcommand{\leftmargini}{22pt}
\begin{enumerate}\renewcommand{\theenumi}{(\roman{enumi})}
  \setlength\itemsep{4pt}
\item $\mcB_\del$ is {\em bounded}: there exists a constant $C_b>0$ independent of $\del$ such that
\[
\mcB_\del (u, v) \leq C_b \|u\|_{{V_\del}} \|v \|_{{W_\del}}\quad\text{for all } u \in {V_\del},\;  v\in {W_\del}.
\]
\item 
{\em Inf-sup condition}:  
there exists a constant $C_c>0$ independent of $\delta$ such that
\[
\inf_{u\in {V_\del}} \sup_{v\in {W_\del}} \ffrac{ \mcB_\del(u, v) }{\| u\|_{{V_\del}} \|v \|_{{W_\del}}}\geq C_c >0 .
\]
\item  If $\mcB_\del(u, v)=0$ for all $u\in {V_\del}$, then $v=0$.
 \end{enumerate}
These conditions, as first shown {by \citename{necas1} \citeyear{necas2,necas1}}, guarantee that the problem 
\begin{equation}\label{petgal}
 \text{{Find} $u(\xb)\in V_\del$ such that $ \mcB_\del(u, v) = \lrang{f}{v}$ for all $ v\in W_\del$}
\end{equation}
 is {well-posed}, provided that the linear functional on the right-hand side is bounded. Note that there are no symmetry or self-adjoint conditions placed on the bilinear form in \eqref{petgal}. Note also that the discrete inf-sup condition  is automatically satisfied in the case of coercive self-adjoint problems with $V_{\delta,h}=W_{\delta,h}$. Otherwise, it has to be verified for the chosen
finite element spaces and the problem under consideration. Problems such as \eqref{petgal} that feature different test and trial spaces are often referred to as {\em Petrov--Galerkin} formulations. If $W_\del=V_\del$, then \eqref{petgal} is often referred to as a {\em {weakly coercive}} formulation.  
 
 Conforming discretizations are defined, as introduced in  \citeasnoun{TiDu19} for the parametrized setting under consideration here, by choosing approximation subspaces  $V_{\del,h}\subset V_\del$ and $W_{\del,h}\subset W_\del$ satisfying the following requirements.
\renewcommand{\leftmargini}{20pt}
 \begin{enumerate}\renewcommand{\theenumi}{(\arabic{enumi})}
   \setlength\itemsep{4pt}
 \item For a given $\delta \in (0,\delta_0)$, the family $\{W_{\delta, h},   h\in(0,h_0]\}$ of discrete subspaces of  $ {W_{\delta}}$ is  dense in $ {W_\del}$ as $h\to 0$.

\item {Discrete inf-sup condition}: there exists a constant $\widetilde C_c>0$, independent of $\delta$ and $h$, such that
\[
\inf_{u\in V_{\delta,h}} \sup_{v\in W_{\delta,h}} \ffrac{ \mcB_\del(u, v) }{\| u\|_{ V_\delta} \|v \|_{ W_\delta}}\geq \widetilde C_c >0 .
\]
\end{enumerate}
These conditions, as first shown in \citeasnoun{babuska1} {and} \citeasnoun{babuska2}, guarantee that for a given $\delta \in (0,\delta_0)$  the problem 
\begin{equation}\label{petgalh}
 \text{{Find} $u(\xb)\in V_{\del,h}$ such that $ \mcB_\del(u, v) = \lrang{f}{v}$ for all $ v\in W_{\del,h}$}
\end{equation}
is {well-posed}. {Note that the independence of the constants $C_b,C_c$ and $\widetilde C_c$ on the parameter $\delta$ is not required if one is only interested in solving the problem for a fixed $\delta$, as in the original theory of \citeasnoun{babuska1} {and} \citeasnoun{babuska2}. It is imposed here for the study of asymptotic compatibility discussed next in Section~\ref{ac-nsym} \cite{TiDu19}.}
 
{Problem} \eqref{petgal} includes several important settings such as nonlocal mixed \bo{and operator-splitting formulations}  (see Section~\ref{mixed}),  nonlocal convection--diffusion problems (see Section~\ref{sec:convection-diffusion}),
nonlocal diffusion models and {bond-based} peridynamic models involving both attractive and repulsive interactions \cite{du13dcdsb}, and nonlocal systems such as the nonlocal Stokes equation introduced in \citeasnoun{DuTi19fcm} {and} \citeasnoun{ld19ng}.

 \subsection{Asymptotically compatible schemes}\label{ac-nsym}
 
 The discrete problem \eqref{petgalh} involves the horizon parameter $\del$ and the {grid size} parameter $h$ so that the {asymptotic} compatibility of particular choices of finite element spaces should be questioned. {A} study of this question is given in \citeasnoun{TiDu19} based on an extension of the original AC framework presented in \citeasnoun{TiDu14}. The latter deals with only symmetric and coercive bilinear forms of the types similar to \eqref{weak-bform} and~\eqref{peri-bform}. 
  
 Concerning the discrete approximations, in addition to the requirements (1) and (2) listed in Section~\ref{weakcoercive}, we impose the following requirement \cite{TiDu19}.
 \begin{enumerate}\setlength\itemsep{4pt}
\item[(3)] The family of discrete subspaces of  $ \{W_{\del,h}\}$  is 
{\em asymptotically dense} in $ W_0$ as $\delta\to 0$ and $h\to 0$, in the sense {of definition} \eqref{eq:asdense}. Here $ W_0$ refers to the energy space for the {original continuum} local PDE problem corresponding to \eqref{petgal} and \eqref{petgalh}.
\end{enumerate}
Now any scheme that satisfies requirements (1), (2) and (3) is provably an AC scheme \cite{TiDu19}.

\subsection{\bo{Operator-splitting and mixed formulations}}\label{mixed}

PDE equations such as \bo{$-\nabla\cdot(\Dbf\nabla u)=f$} are often derived by first postulating a balance law $\nabla\cdot\wb = f$ and choosing a constitutive (Darcy, Fick, Fourier, Ohm, {\it etc.})\ law $\wb=-\Dbf\nabla u$. In some settings there are advantages to directly solving the two first-order equations instead of the single second-order equation, perhaps the most important and useful being that {well-posedness} can be proved for $u\in L^2(\omg)$ instead of $H^1(\omg)$, as is the case for the second-order equation. On the other hand, there are problems that are most often posed as a mixed formulation, the most common {being} the Stokes equations for incompressible flows.

\pagebreak %%20200331
    One can mimic the local setting by {recasting problem} \eqref{str-vcp} with $\mcLd u= -\mcDd\cdot(\mcDsd u)$ in $\omg$ and with, say, $u=0$ on $\omgid$,  into the equivalent \bo{operator-splitting}  formulation
\begin{equation}\label{mixedp}
\begin{cases}
\mcDd \nub\xyp = f(\xb)   & \text{for all } \xb \in\omg ,\\[2pt]
   \nub\xyp = \mcDsd u  & \text{for all } \xb, {\yb} \in\omg\cup\omgid ,
   \end{cases}
\end{equation}
along with $u(\xb)= 0$ for all $ \xb\in\omgid$.
\bo{A corresponding mixed weak formulation is given by}
\begin{equation}\label{mixedpw}
\begin{cases}
 {(\nub\xyp, \mub\xyp)_{L^2([\omg\cup\omgid]^2)}} - (\mcDd \mub\xyp, u(\xb) )=0,
\\[2pt]
 (\mcDd \nub\xyp, v(\xb) )  =  (f(\xb), v(\xb)), 
\end{cases}
\end{equation}
along with $u(\xb)= 0$ for all $ \xb\in\omgid$,
where \bo{$(\cdot,\cdot)_{L^2([\omg\cup\omgid]^2)}$ and $(\cdot,\cdot)$ denote $L^2$ inner products on the respective domains} and the pairs $\nub,\mub$ and $u,v$ belong to appropriate function spaces. 

{Problem} \eqref{mixedpw} can be treated as an \bo{extension} {of problem} \eqref{petgal} with the bilinear form
\begin{align*}
\mcB_\del ((\nub,u) ,\,(\mub,v)) &=
 (\nub\xyp, \mub\xyp)  - (\mcDd \mub\xyp, u(\xb) )  \\*[2pt]
&\quad\,  + (\mcDd \nub\xyp, v(\xb) ) .
\end{align*}
\bo{While one may attempt to adopt the general inf-sup theory in \citeasnoun{TiDu19} to study the mixed weak formulation,
  the mixed finite element \roo{approximations would} require the discretization of $\nub\xyp$ and $\mub\xyp$ in both $\xb$ and $\yb$. On the other hand,
finite} element discretizations of nonlocal models based on \boo{the operator-splitting} formulation can be found in {\citename*{ldjl19sc} \citeyear{ldjl19sc,djl19mc}} for \bo{interaction kernels that are both radial and \roo{integrable, and}}  that recover the local {discontinuous} Galerkin (LDG) discretization of the local PDE problem as $\del\to 0$. 

The nonlocal Stokes model is another example that can be formulated as a system in  mixed form; see the discussions {of \citeasnoun{DuTi19fcm} and \citeasnoun{ld19ng}, who also analysed spectral and finite difference approximations in} a  periodic boundary condition setting.

\vspace{5pt}
\section{Nonlocal convection--diffusion problems} \label{sec:convection-diffusion}

In this section we consider nonlocal analogues of the local convection--diffusion (also referred to as advection--diffusion) problem
\begin{equation}\label{lcodi}
\begin{cases}
    -\nabla\cdot(\Dbf \nabla u) + \nabla\cdotsp ({\bm U} u) = f &\text{for all } \xb\in\omg,
    \\[2pt]
    u=0 &\text{for all } \xb\in\partial\omg,
\end{cases}
\pagebreak %%20200331
\end{equation}
where ${\bm U}(\xb)$ denotes a given velocity field and $\Dbf(\xb)$ denotes a given symmetric, positive definite matrix. The most common approach towards defining a nonlocal convection--diffusion model is to replace the second-order diffusion term in \eqref{lcodi} {with} a nonlocal analogue $-\mcL u$ but to keep the convection term as it is in \eqref{lcodi}. Such classical convection--nonlocal diffusion problems have been investigated for fractional models, including their connection to L\'evy jump processes; see \eg\ \citeasnoun{Meerschaert2012book}.

However, we consider fully nonlocal analogues of \eqref{lcodi} in which the convection term in \eqref{lcodi} is also replaced by a nonlocal analogue. As a result, the nonlocal convection--diffusion models we consider feature non-symmetric kernels $\gamd\xyp$, {so they can also} be viewed as modelling non-symmetric diffusion. Among other {descriptions} of non-symmetric diffusion that are not necessarily related to stochastic processes, we mention \citeasnoun{bame:10}, \citeasnoun{mebe:99} {and} \citeasnoun{Meerschaert2012book} where the equations are set either in free space or in bounded domains,  {\citeasnoun{erro:07}  \bo{that treats} the same problem in a variational setting, \citeasnoun{Felsinger2015}  \bo{that \roo{analyses}} 
  a variational formulation of non-symmetric diffusion for \bo{integrable kernels with square-integrable symmetric parts and for} non-integrable kernels, and \citeasnoun{amrt:10}  \bo{that considers}  the strong form of non-symmetric diffusion equations for \bo{kernels that are positive, translation-invariant and integrable.} 

Based on \citeasnoun{DElia2017nlcd}, we consider the most general form of a nonlocal analogue of \eqref{lcodi}, treating the two nonlocal terms as separate phenomena. Note that \citeasnoun{DElia2017nlcd} {provide} a generalization of \citeasnoun{Du2013conv} to a more general class of kernels.

\subsection{Non-symmetric kernels and nonlocal convection--diffusion operators}

Let $\mcL_{cd,\del}$ denote the nonlocal convection--diffusion operator defined as 
\begin{equation}\label{non-symm-L}
-\mcL_{cd,\del} u(\xb) := \mcD_{d,\del_d}(\Theb \mcD^\ast_{d,\del_d}u)(\xb) +\mcD_{c,\del_c}(\mub u)(\xb),
\end{equation}
where $\mcD_{d,\del_d}$ and $\mcD_{c,\del_c}$ are nonlocal divergence operators associated with the anti-symmetric functions $\alpb_d\xyp$ and $\alpb_c\xyp$ and where $\mub\xyp=\mub(\yb,\xb)$ is a given function. Note that we allow for different horizons and different kernel functions for the diffusion and convection terms. We refer to the second-order tensor $\Theb\xyp$ as the nonlocal diffusion tensor (see Section~\ref{sec:mod-nlvc}) and to the vector $\mub\xyp$ as the nonlocal convection `velocity'. 
Specifically, from \eqref{str-gamd}, we have
\begin{align*}
&\mcD_{d,\del_d}(\Theb \mcD^\ast_{d,\del_d} u)(\xb)
\\*
& \quad= -2\int_{\Rd}  ( u(\yb) - u(\xb) )  \alpb_d\xyp \cdot (\Theb\xyp\alpb_d\xyp) {\mcX}_{B_{\del_d}(\xb)}(\yb) \dyb
\end{align*}
and similarly, from \eqref{str-div} with $\nub\xyp=\mub\xyp u(\xb)$, we have
\begin{align*}
\mcD_{c,\del_c}(\mub u)(\xb) &= \int_{\Rd} (\mub\xyp u(\xb)  +  \mub(\yb,\xb) u(\yb)  ) \cdot\alpb_c\xyp {\mcX}_{B_{\del_c}(\xb)}(\yb) \dyb \notag \\*
&= \int_{\Rd} (u(\xb)  +  u(\yb)  ) \mub\xyp \cdotsp\alpb_c\xyp {\mcX}_{B_{\del_c}(\xb)}(\yb) \dyb.
\end{align*}
Setting
\begin{align}\label{gamma_nonsym}
\gamma_{cd,\del}\xyp & :=
\underbrace{2\alpb_d\xyp\cdot (\Theb\xyp\alpb_d\xyp) \mcX_{B_{\del_d}(\xb)}(\yb)}_{\text{symmetric part} \,\, \gamma_{d,\del_d}\xyp} \notag \\*
&\quad\, \underbrace{- \mub\xyp\cdot\alpb_c\xyp\mcX_{B_{\del_c}(\xb)}(\yb)}_{\text{anti-symmetric part} \,\,\gamma_{c,\del_c}\xyp},
\end{align}
we can rewrite \eqref{non-symm-L} as
\begin{equation}\label{Lexplicit}
\mcL_{cd,\del} u(\xb) = \int_\Rd  (u(\yb)\gamma_{cd,\del}{(\yb,\xb)}-u(\xb)\gamma_{cd,\del}\xyp)  \dyb \quad\text{for all } \xb\in\Rd.
\end{equation}
In the case of 
 $\alpb_c(\xb,\yb)=\alpb_d(\xb,\yb)$ and $\del_d=\del_c$, this non-symmetric diffusion operator and the corresponding initial value problem are analysed \bo{for a special class of kernels} in \citeasnoun{Du2013conv}.

Note that with $\del=\max\{\del_d,\del_c\}$, we have
\begin{equation}\label{gamma-conds}
\gamma_{cd,\del}\xyp  = 0 \quad\text{for all } \yb \not\in \ballx.
\end{equation}
Also note that although $\alpb_d\xyp$ and $\alpb_c\xyp$ are often radial functions, in general $\Theb\xyp$ and $\mub\xyp$ are not, so that, also in general, $\gamma_{cd,\del}\xyp$ is not radial {or} translationally invariant. See Remark~2.3.

{\citename{du15cmame} \citeyear{du15cmame,tjd17cmame} and \citeasnoun{ld19sinum} considered AC discretizations of nonlocal convection--diffusion problems.  \citeasnoun{Du2013conv} considered} a convection--diffusion operator in one dimension, which turns out to be a special case of the general nonlocal {convection--diffusion} operator \eqref{non-symm-L} if {we choose} $\del_d=\del_c$, $\Theta(x,y)=\kappa=$ constant,
$ (\alpha_d(x,y)) ^2=\sigma_{d,\del}(|y-x|)$, $\mu(x,y)=U$ a constant, and $\alpha_c(x,y)=(y-x)\sigma_{c,\del}(|y-x|)$, with $\sigma_{d,\del}(\cdot)$ and $\sigma_{c,\del}(\cdot)$ being even functions having unit second moments. The resulting operator is given by
\begin{align*}
\label{nl-conv-diff}
\mcL_{cd,\del} u (x) & =  2 \kappa  \int_{\mathbb{R}}  (u(y)- u(x))  \sigma_{d,\del}(|y-x|) \D y \\*
& \quad\, - U\int_{R}   (u(y)+ u(x))  (y-x)\sigma_{c,\del}(|y-x|)  \D y .
\end{align*}
The local counterpart {of this} operator is $\mathcal{L}_0 u(x)=  \kappa u''(x) + U u'(x)$.
One may also connect this  nonlocal convection diffusion {model} with non-sym\-metric jump processes; see \citeasnoun{Du2013conv}. 

Other nonlocal convection--diffusion models are presented in \citeasnoun{tjd17cmame}, including some that are reminiscent of  state-based peridynamic models. For example, a conservative formulation is defined by the nonlocal convection--diffusion operator
\begin{align*}
  &\mcL_{cd,\del} u (\xb) \\*
  &= \intoi (A(\xb)+A(\yb)) (u(\yb)-u(\xb))  \sigma_{d,\del}(|\ymx|) (|\yb-\xb|)
  {\mcX}_c(\xb,\yb)  \dyb  \\*
  &   +  \intoi  ({\bm b}(\xb) u(\xb){\mcX}_c(\xb,\yb) + {\bm b}(\yb) u(\yb) {\mcX}_c(\yb,\xb)
) \cdot (\yb-\xb) \sigma_{c,\del}(|\ymx|)  \dyb,
\end{align*}
where we have the indicator function 
\[ 
  {\mcX}_{c}\xyp=
  \begin{cases} 1 & \text{if $ |\yb-\xb| < \delta$  and $ {\bm b}(\xb)\cdot (\yb-\xb)> 0$,} \\
    0 & \text{otherwise}, \end{cases}
\]
which is generically non-symmetric and is dependent on the velocity field ${\bm b}(\xb)$. This operator can {also be shown to be} a special case of operator given in \eqref{non-symm-L}. Its local counterpart  is
\[ %%20200331
\roo{\mathcal{L}_0 u(\xb)= \nabla\cdot(A(\xb) \nabla u(\xb))+\nabla\cdot({\bm b}(\xb) u(\xb)).}
\]

Although AC discretizations of the models using these operators have been discussed {by \citename{du15cmame} \citeyear{du15cmame,tjd17cmame} and \citeasnoun{ld19sinum}}, the attendant analyses are done using different techniques and for specialized kernels. One may apply the general framework given in \citeasnoun{TiDu19} to possibly offer a unified treatment of systems of non-self-adjoint problems.
Indeed,  additional studies {of} asymptotically compatible schemes may also shed new light on
improving the robustness of numerical methods based on various nonlocal smoothing approaches, that are in use for  local PDE models such {as smoothed particle hydrodynamics, as well as the} construction of well-posed nonlocal models such as the nonlocal Stokes equation considered in \citeasnoun{DuTi19fcm} {and} \citeasnoun{ld19sinum}.

\subsubsection{Steady-state nonlocal convection--diffusion problems}

The nonlocal analogue of \eqref{lcodi} is given by
\begin{equation}\label{non-symm-diff}
\begin{cases}
-\mcL_{cd,\del} u = f  & \xb\in\omg ,\\
u        = 0&  \xb\in{\omgid}.
\end{cases}
\end{equation}
Weak {formulations corresponding} to \eqref{non-symm-diff} can be defined in the usual way. Simplifying some notation, for $\ell=d$ or $\ell=c$, one can define the constrained energy space 
\[
V_\ell^0  
:= \{ v  \in L^2(\omgomgi)  \colon  |v|_{V_\ell^0} < \infty \ \text{and $  v=0$  on $\omgid$}\} 
\]
for which the semi-norm
\[
|v|_{V_\ell^0}   := 
\ffrac12\int_\omgomgi\int_\omgomgi|\mcD^\ast_\ell(v)\xyp|^2 
\gamma_{d,\del_\ell}\xyp\dydx 
\]
is a norm. We assume that the norm $|v|_{V_\ell^0}$ satisfies the nonlocal Poincar\'e inequality $\|v\|_{L^2(\omgomgi)}\leq C_p |v|_{V_\ell^0}$ for all \bo{$v\in V_\ell^0$} and that \bo{$V^0_d\subset V^0_c$} so that $|v|_{V_c^0} \le |v|_{V_d^0}$. The latter assumption implies that solution operators for nonlocal diffusion problems under consideration effect greater smoothing compared to those for nonlocal convection problems, as is the case for local partial differential operators. 

Let the bilinear form $\mcA_{cd,\del}(\cdot,\cdot)$ be defined, for all $u,v\in V_d^0$, by
\begin{align}\label{bilinear_a}
\mcA_{cd,\del}(u,v)& = \int\limits_\omgomgi \int\limits_\omgomgi \mcD^\ast_{d,\del}(u)(\xb,\yb)\cdot 
(\Theb \mcD^\ast_{d,\del} v)(\xb,\yb) \dydx
\notag \\*
& \quad\, -\int\limits_\omg \mcD_{c,\del}(\mub u)(\xb)\,v(\xb) \dxb ,
\end{align}
and let the linear functional $\lrang{f}{v}$ be defined, for all $v\in V_d^0$, by
\begin{equation}\label{linearG}
\lrang{f}{v} = \int_\omg f(\xb)v(\xb) \dxb.
\end{equation}
Then a weak formulation of \eqref{non-symm-diff} is given as follows: given $f$ belonging to the dual space of $V_d^0$, find $u \in V_d^0$ that satisfies 
\begin{equation}\label{nlcd-weak}
\mcA_{cd,\del}(u,v)=\lrang{f}{v} \quad\text{for all } v\in V_d^0.
\end{equation}
{\citeasnoun{DElia2017nlcd} proved the well-posedness of \eqref{nlcd-weak} using} three different approaches. Here we state results that mimic what is obtained for local convection--diffusion problems.  Specifically,
  \renewcommand{\leftmargini}{20pt}
\begin{itemize}\setlength\itemsep{4pt}
\item if $f(\xb)$ {belongs} to the dual space of $V_d^0$,

\item if $\Theb\xyp$ is such that there exist ${\underline\vartheta},{\overline\vartheta}>0$ satisfying
\begin{equation}\label{theta_cond}
0<{\underline\vartheta}\leq \inf_{\xb\in\Rd}(\min_i \theta_i)\quad\text{and}\quad \sup_{\xb\in\Rd}(\max_i \theta_i)\leq{\overline\vartheta}<\infty, 
\end{equation}
where $\theta_i$ denote the singular values of $\Theb$,

\item if $\mub$ is such that $C_p^2\|\mcD_c\mub\|_\infty\leq 2{\underline\vartheta}$ and $\| \,|\mub| \,\|_\infty\leq {\overline\mu}$, 
\end{itemize}
{then problem} \eqref{nlcd-weak} has a unique solution $u\in V_d^0$. Furthermore, that solution satisfies the {\it a~priori} estimate
\begin{equation}
| u|_{V_d^0}  \leq C \|f\|_{V'_d},
\end{equation}
where $C=\gfrac{1}{C_{\rm coer}}$ {and} 
$$C_{\rm coer}= {\underline\vartheta}-\frac12 C_p^2\|\mcD_c\mub\|_\infty$$
{denotes} the coercivity constant for the bilinear form $\mcA_{cd,\del}(\cdot,\cdot)$. Note that the above analysis is effective for {diffusion-dominated} problems. For the {convection-dominated case one} may use the formulation in \citeasnoun{TiDu19}, discussed in Section~\ref{weakcoercive}, to get results in more general cases.

\vspace{5pt}
\section{Time-dependent nonlocal problems}\label{sec:weak-time}

Although we do not consider nonlocal time-dependent problems other than in this subsection and briefly in some other sections, a brief discussion is warranted. We do {not 
delve} into nonlocality in time, for which there is a vast literature devoted to fractional time derivatives and other settings in which memory effects are present.

Weak formulations of time-dependent problems can be defined in the same manner as that for local PDE time-dependent problems, once one knows how to treat steady-state problems. Likewise, for the discretization of nonlocal time-dependent problems, spatial discretization can be effected using any discretization method for the corresponding steady-state problem, including those discussed in Part~\ref{part2}, and temporal discretization can be effected using any discretization method for the corresponding local PDE problem, \eg\ the backward-Euler or Crank--Nicolson method for \eqref{weak-timed} or, for \eqref{weak-timewa}, a leap-frog or other explicit method. Furthermore, the analysis of weak formulations and discretizations of nonlocal time-dependent problems, including the derivation of well-posedness results and error estimates, also follows the same paths as those for the corresponding local PDE problems.

Here we {only} consider a small sample of the differences between time-dependent local and nonlocal problems.

\subsection{Time-dependent nonlocal diffusion}

Using the notation used in \eqref{str-vcp}, we have the nonlocal time-dependent diffusion equation 
\begin{equation}\label{weak-timed}
\begin{cases}
  \rrr \ffrac{\partial u}{\partial t} = \mcLd u +  f(\xb,t) & \text{for all } \xb\in\omg\times(0,T] , \\[7pt]
    \mcV u = g(\xb,t) & \text{for all } \xb\in\omgid\times(0,T] , \\[2pt]
      u(\xb,0) = u_0(\xb) & \text{for all } \xb\in\omg,
\end{cases}
\end{equation}
where $\rrr(\xb)>0$ and $u_0(\xb)$ are given functions defined on $\omg$ and $f(\xb,t)$ and $g(\xb,t)$ are given functions defined on $\omg\times(0,T]$ and $\omgid\times[0,T]$, respectively.

\pagebreak %%20200331
Using the nonlocal Green's first identity, it is an easy matter to show, for $f=0$ and $g=0$, that the nonlocal diffusion equation \eqref{weak-timed} implies that
\[
\ffrac12 \ffrac{\dd}{\dd t}\int_\omg  u^2\dxb +
\int_{\omgomgi}\int_{\omgomgi} \mcDsd u\cdotsp(\Thebd\mcDsd u) \dydx = 0 ,
\]
so that
\begin{equation}\label{weak-decay}
\int_\omg  u^2(\xb,t)\dxb {\leq} \int_\omg  u_0^{2}(\xb)\dxb
\quad\text{for all } t>0,
\end{equation}
and for time-independent $\Thebd$
\begin{align}\label{weak-decay-ene}
&
\int_{\omgomgi}\int_{\omgomgi} \mcDsd u(\xb,t)\cdotsp(\Thebd\mcDsd u(\xb,t) \dydx \\*
&\quad\leq
\int_{\omgomgi}\int_{\omgomgi} \mcDsd u(\xb,0)\cdotsp(\Thebd\mcDsd u(\xb,0) \dydx \notag
\quad\text{for all } t>0.
\end{align}
These are decay characteristics of diffusive processes, \feeg, \eqref{weak-decay} and the local version of \eqref{weak-decay-ene} hold for parabolic PDEs. However, for \bo{kernel functions that are both radial and integrable}, although solutions of the nonlocal diffusion equation \eqref{weak-timed} satisfy the decay properties \eqref{weak-decay} and \eqref{weak-decay-ene}, unlike the case for parabolic PDEs, those solutions may not be much smoother than the data.
One can also consider various types of  nonlocal-in-time versions of nonlocal diffusion equations; see for example \citeasnoun{cdlz17chaos}.

We note that the {\em fractional heat equation} 
  \begin{equation}\label{eq:diffusion-z-RHS2}
\begin{cases}
u_t + (-\Delta)^s u  = f & \text{for all }\xb\in\omg,\;t\in(0,T) , \\
u = g & \text{for all }\xb\in\Rd\setminus\omg,\;t\in(0,T) ,\\
u(\cdot,0)=u_0(\xb) & \text{for all }\xb\in\omg.
\end{cases}
  \end{equation}
is perhaps of even greater interest compared to the steady-state case. Similarly, one can consider the fractional heat equation with the regional fractional Laplacian; \bo{see \citeasnoun{Gal2016}}.

\subsection{Time-dependent convection--diffusion problems}\label{sec:time-CD}

For $T>0$, we introduce the time-dependent function spaces
\begin{align*}
L^2(0,T;V_d^0) & =\{v(\cdot,t)\in V_d^0 \colon | v(\cdot,t)|_{V_d} \in L^2(0,T)\},\\*
L^2(0,T;V'_d) & =\{f(\cdot,t)\in V'_d \colon \|f(\cdot,t)\|_{V'_d}\in L^2(0,T)\}.
\end{align*}
The time-dependent  nonlocal {convection}{-}{diffusion} problem is then given by
\begin{equation}\label{unsteady-strong}
\begin{cases}
u_t-\mcL_{cd,\del} u   = f    &  \text{for all } \xb\in\omg,\; t\in(0,T] , \\
u(\xb,t)     = 0    &  \text{for all } \xb\in {\omgid},\; t\in(0,T] , \\
u(\xb,0)    = u_0(\xb)  & \text{for all } \xb\in\omg ,
\end{cases}
\end{equation}
\bo{where the operator $\mcL_{cd,\del}$ is defined in \eqref{non-symm-L}.}
The corresponding weak formulation {is as} follows: given $f\in L^2(0,T;V'_d)$ and $u_0\in V_d^0$, find $u\in L^2(0,T;V_d^0)$ that satisfies $u(\xb,0)=u_0(\xb)$ and, for all $v\in V_d^0$ and for almost every $t\in(0,T]$,
\begin{equation}\label{unsteady-weak-forms}
(u_t,v)_\omg + \mcA_{cd,\del}(u,v) = \lrang{f}{v},
\end{equation}
where $(\cdot,\cdot)_{\omg}$ denotes the $L^2$ inner product over $\omg$. The coercivity and continuity of $\mcA_{cd,\del}(\cdot,\cdot)$ and the continuity of $\lrang{\cdot}{\cdot}$ guarantee that the weak formulation \eqref{unsteady-weak-forms} is {well-posed}. However, as pointed out {by} \citeasnoun{DElia2017nlcd}, standard arguments of variational theory \cite{evans:98} imply that actually much weaker assumptions on $\mub$ are required for {well-posedness}, namely $\|\mcD_{c,\del}\mub\|_\infty<\infty$, \bo{where $\mub(\xb,\yb)$ denotes the nonlocal convection `velocity' introduced in \eqref{gamma_nonsym}.}

\subsection{Nonlocal wave equations}

One can also consider the nonlocal wave equation
\begin{equation}\label{weak-timewa}
\begin{cases}
  \rrr \ffrac{\partial^2 u}{\partial t^2} = \mcLd u +  f(\xb,t) & \text{for all } \xb\in\omg\times(0,T]  , 
\\[7pt]
  \mcV u = g(\xb,t) &\text{for all } \xb\in\omgid\times(0,T] , 
    \\[7pt]
    u(\xb,0) = u_0(\xb) & \text{for all } \xb\in\omg , 
    \\[7pt]
    \ffrac{\partial u}{\partial t}(\xb,0) = u_1(\xb) & \text{for all } \xb\in\omg,
\end{cases}
\end{equation}
where $\rrr(\xb)>0$, $u_0(\xb)$, $f(\xb,t)$ and $g(\xb,t)$ are defined as for \eqref{weak-timed} and $u_1(\xb)$ is a given function defined on $\omg$. 

For \eqref{weak-timewa} with $f=0$ and $g=0$, we have conservation of energy, \tie,
\[
\ffrac{\dd}{\dd t} \biggl(
\ffrac12 \int_\omg  \biggl(\ffrac{\dd u}{\dd t}\biggr)^2\dxb +
\int_{\omgomgi}\int_{\omgomgi} \mcDsd u\cdotsp(\Thebd\mcDsd u) \dydx
\biggr)
=0 ,
\]
{which} is a characteristic of wave processes such as the PDE wave equation.  One can find studies related to these nonlocal wave equations in \citeasnoun{guanwave}, \citeasnoun{dzz18cicp} {and} \citeasnoun{du16esiam}.

One of the stark differences between local and nonlocal models is in their dispersion relations for wave equations.
For the one-dimensional local PDE wave equation
\[
\ffrac{\partial^2 u}{\partial^2 t}
= c^2\ffrac{\partial^2 u}{\partial^2 x},
\pagebreak %%20200331
\]
where $c$ is a constant, by setting
\[
u(x,t) = \rme^{-\rmi\omega t + \rmi k x},
\]
{we obtain} the dispersion relation
\begin{equation}\label{ldisper}
      \omega^2 = c^2 k^2.
\end{equation}
In fact, this relation shows that there is no dispersion. The velocity of the wave is $\omega /k = \pm c$, which is independent of $\omega$ and $k$. 

{\citeasnoun{guanwave} showed that} the one-dimensional nonlocal wave equation
\[
\ffrac{\partial^2 u}{\partial^2 t}
 =
\ffrac{2-2s}{\del^{2-2s}}c^2\int_{x-\delta}^{x+\del}
\ffrac{u(y,t)-u(x,t)}{|y-x|^{1+2s}} \D y,\quad 0\leq s< 1/2,
\]
has the dispersion relation
\begin{equation}\label{nldisper}
\omega^2 =
\ffrac{2-2s}{\delta^{2-2s}}
c^2\int_{0}^{\delta}
\ffrac{2-2\cos(ky) }{y^{1+2s}} \D y.
\end{equation}
Observe that the wave velocity $\omega/k$ is a nonlinear function of $k$.
{\citeasnoun{guanwave} also showed that} as $\del\to0$, $\omega$ given by \eqref{nldisper} converges (quadratically with respect to $\del$) to the local $\omega$ of \eqref{ldisper}. Similar results are obtained there for the two-dimensional case. Similar dispersion relations have been discussed  for nonlocal operators; see \eg\ \citeasnoun{du10sinum}, \citeasnoun{du11esiam} {and} \citeasnoun{du19cbms}.

\vspace{5pt}
\section{Inverse problems}\label{sec:control}

Among the many challenges faced when dealing with nonlocal problems, we {find,
even more than} for local PDE problems, that mathematical models are not known with exactitude{;} \feeg, source terms, volume constraint data, coefficients, and even the functional form of the kernel itself may be unknown or subject to uncertainty. If there are experimental data or other {\it a~priori} information (that may be sparse and/or noisy) available about the state of the system or about an output of interest that depends on the state, one can then resort to control or optimization strategies to identify the unknown entities and thus define a data-driven mathematical model that is more faithful to the physics being considered. 

Here we consider inverse problems for the nonlocal diffusion problem \eqref{str-vcp}
in which the boundary operator $\mcV$ could correspond to Dirichlet, Neumann or Robin volume constraints. Let $V(\omg\cup\omgi)$ denote a function space for the state $u(\xb)$ and let $\WWW$ denote a set of controls that could consist of function spaces or parameter vectors or a combination of both. Then a general inverse problem for nonlocal diffusion is given as follows: 
\begin{equation}\label{eq:generalized-control}
\begin{minipage}[b]{18pc}
{Seek} $u(\xb)\in V(\omg\cup\omgi)$ and $\mumu\in \WWW$ such that
\[\bo{\min\limits_{\mumu\in W}  \;\mcJ(u;\mumu)}  \]
subject to \eqref{str-vcp} being satisfied, \bo{where}
\[
\bo{\mcJ(u;\mumu)=\mcQ(u;\mumu) + \mcR(\mumu)}.
\]
\end{minipage}
\end{equation}
In \eqref{eq:generalized-control}, \bo{the first term $\mcQ(\cdotsp;\cdotsp)$ in the objective functional $\mcJ(\cdotsp;\cdotsp)$} denotes a cost functional that depends on the state and control whereas $\mcR(\cdotsp)$ denotes a regularization functional that serves to guarantee the {well-posedness} of the problem. The control set $\WWW$ could contain data functions such as $f$, $g$ and $\Theb$, and also parameters appearing in the model definition such as the horizon $\del$ or the fractional exponent $s$ if \eqref{str-vcp} represents a fractional Laplacian problem. The functions may be parametrized, in which case $\WWW$ only contains a parameter vector. Additionally, the control set may contain constraints on the control; constraints on the state may be also be imposed. {See} Section~\ref{sec:vi} for an example of the latter. In some such cases, the regularization term in \eqref{eq:generalized-control} may not be needed because such constraints may be sufficient to guarantee {well-posedness}. 

The literature about the control and optimization of nonlocal problems is still limited; however, interest in such topics in the setting of nonlocal diffusion is quickly growing. Recent studies in this direction focus on the {well-posedness} and stability of the minimization problem \eqref{eq:generalized-control}, the asymptotic behaviour of its solution, and its numerical discretization. In particular, with respect to the latter, numerical convergence analyses, error estimates, and solver performance are of interest. In this section we provide brief reviews of selected contributions devoted to the control and optimization of nonlocal problems, including \bo{integral} fractional models, treating both control and identification problems.

\subsection{Inverse problems for nonlocal diffusion}\label{sec:inverse}

In this section we focus on operators of the form of \eqref{str-lap} with $\gamd\xyp$ given \bo{in terms of the kernel function $\phid\xyp$ and constitutive function $\thed\xyp$ by \eqref{gamd1}, \eqref{ealpbd} and \eqref{ethed},} \tie, we have
\begin{align}\label{eq:truncated-laplacian}
\mcLd u(\xb) &=  \mcDd (\Thebd\mcDs u) (\xb) \notag \\*
& = -2\int_{\omgomgi}  (u(\yb)-u(\xb))  \gamd\xyp \dyb \\*
&= -2\int_{\ballx} \thed\xyp (u(\yb)-u(\xb))  \phid\xyp \dyb \quad\text{for all } \xb\in\omg. \notag
\end{align}

\subsubsection{Distributed optimal control \bo{in nonlocal diffusion} for a matching functional} \label{sec:controlnonlocal}

We consider the minimization problem \eqref{eq:generalized-control} with the state space $V(\omg\cup\omgi)=\VVV$, the operator $\mcL$ now given by \eqref{eq:truncated-laplacian}, and with perhaps the most commonly used cost functional \bo{and regularization term} for both PDE and nonlocal optimal control problems, namely 
\begin{equation}\label{eq:matching-Tikhonov}
\mcJ(u;\mumu)=\dfrac{1}{2} \|u-\widehat u\|_\UUU^2 + \dfrac{\beta}{2}\|\mcM(\mumu)\|^2_\WWW,
\end{equation}
where the first term is usually referred to as a {\em matching functional}, \bo{the second term as {\em Tikhonov regularization}}, and the given function $\widehat u(\xb)$ for $\xb\in\omg$ as the {\em target function}. The operator $\mcM$ could be, \feeg, a local derivative or a nonlocal operator such as $\mcDsd$, that is chosen with the purpose of keeping $\mcM(\mumu)$ under control and to either guarantee {well-posedness} or improve the conditioning of the problem. Furthermore, $\widehat u$ need not belong to state space $\VVV$ and, in \eqref{eq:matching-Tikhonov}, we have norms $\|\cdotsp\|_\UUU$ and $\|\cdotsp\|_\WWW$ that are {well-defined} for $\widehat u\in U(\omg)$ and $\mcM(\mumu)$ with $\mumu\in\WWW$, respectively. Often target functions are not regular, {so a} reasonable choice is $\UUU=L^2(\omg)$.

Perhaps \citeasnoun{DElia2014DistControl} {were the first to analyse} this problem for \bo{square-integrable and also} non-integrable kernel functions $\phid\xyp$, albeit for \bo{$\theta_\delta=1$} and $\mcV u=u$, \ie\ for Dirichlet volume constraints. Specifically, {they considered the problem of finding} the optimal forcing term $f$ such that the nonlocal solution $u$ is as close as possible to a given target function $\widehat u$, \tie, we have that $\mumu=f$, $\mcM$ is the identity operator, and $\|\cdotsp\|_\WWW=\|\cdotsp\|_{L^2(\omg)}$. As a result, we have the functional
\begin{equation}\label{eq:J-delia2014}
\mcJ(u;f)=\dfrac{1}{2} \|u-\widehat u\|_{L^2(\Omega)}^2 + \ffrac{\beta}{2}\|f\|^2_{L^2(\Omega)}.
\end{equation}
There are no additional constraints on the solution, {so the} optimization is solely constrained by the nonlocal diffusion equation. The {well-posedness} of that equation is sufficient to guarantee the existence and uniqueness of an optimal pair $(u^*,f^*)$. Furthermore, {\citeasnoun{DElia2014DistControl} showed} that in the limit of vanishing nonlocality, \ie\ as $\del\to0$, the optimal nonlocal state and control converge to the optimal solution $(u_l^*,f_l^*)$ of the local counterpart of \eqref{eq:generalized-control} given {as follows:} 
\begin{equation}\label{eq:generalized-control-local}
\begin{minipage}[b]{19pc}
{Seek} $u_l(\xb)\in H^1(\omg)$ and $f_l\in L^2(\omg)$ such that
\[
\min\limits_{f_l\in L^2(\omg)}   \mcJ(u_l;f_l)  
\]
subject to 
$\begin{cases}
-\Delta u_l = f_l & \text{for all } \xb\in\omg ,\\
u_l=g & \text{for all } \xb\in\partial\omg.
\end{cases}$
\end{minipage}
\end{equation}

\pagebreak %%20200331
 {For finite element discretizations of the state and control variables, \citeasnoun{DElia2014DistControl} proved the convergence of both variables with respect to $\del$ and the mesh size $h$, along with} error estimates. Also, numerical results for discontinuous target functions show that nonlocal models, for which irregular solutions are admissible, allow one to match non-smooth functions in a much better way compared to local models.

An example is provided in \fo{Figure~\ref{fig:control-paper}} for the domains $\omg=(0,1)$ and $\omgid=(-\delta,0)\cup(1,1+\delta)$, and for a target function having a jump discontinuity at $x=0.5$. Discontinuous {piecewise linear} finite element {discretizations} are used for both the state and control. For comparison purposes, continuous {piecewise linear} finite element approximations of the local optimal control problem \eqref{eq:generalized-control-local} are also computed. In \fo{Figure~\ref{fig:control-paper}}{(a) we plot} the target function $\widehat u$, the optimal local state $u_l^*$ (the solution of \eqref{eq:generalized-control-local}) and nonlocal optimal state $u^*$ (for two values of $\delta$). Note that for a large horizon $\del$ the nonlocal solution perfectly matches the target, whereas for a small horizon the nonlocal optimal solution is visually identical to the local one. In \fo{Figure~\ref{fig:control-paper}}{(b) we plot the} corresponding optimal source terms $f^*$. Here, for a large horizon, the control has a smaller amplitude and a smaller {$L^2$-norm} (which can be viewed as indications of a smaller cost of control) even though that control does a better job {of} matching the target function. As explained in \citeasnoun{DElia2014DistControl}, this behaviour can be justified by the fact that the nonlocal model allows for discontinuous behaviour in the optimal state, and thus the optimal control has an `easier time' forcing a match between the optimal state and the non-smooth target function.
\begin{figure} %%fig16.1 %%wasfig15.1
\centering
\subfigure[]{\includegraphics[width=177pt,viewport=0 0 393 330,clip]{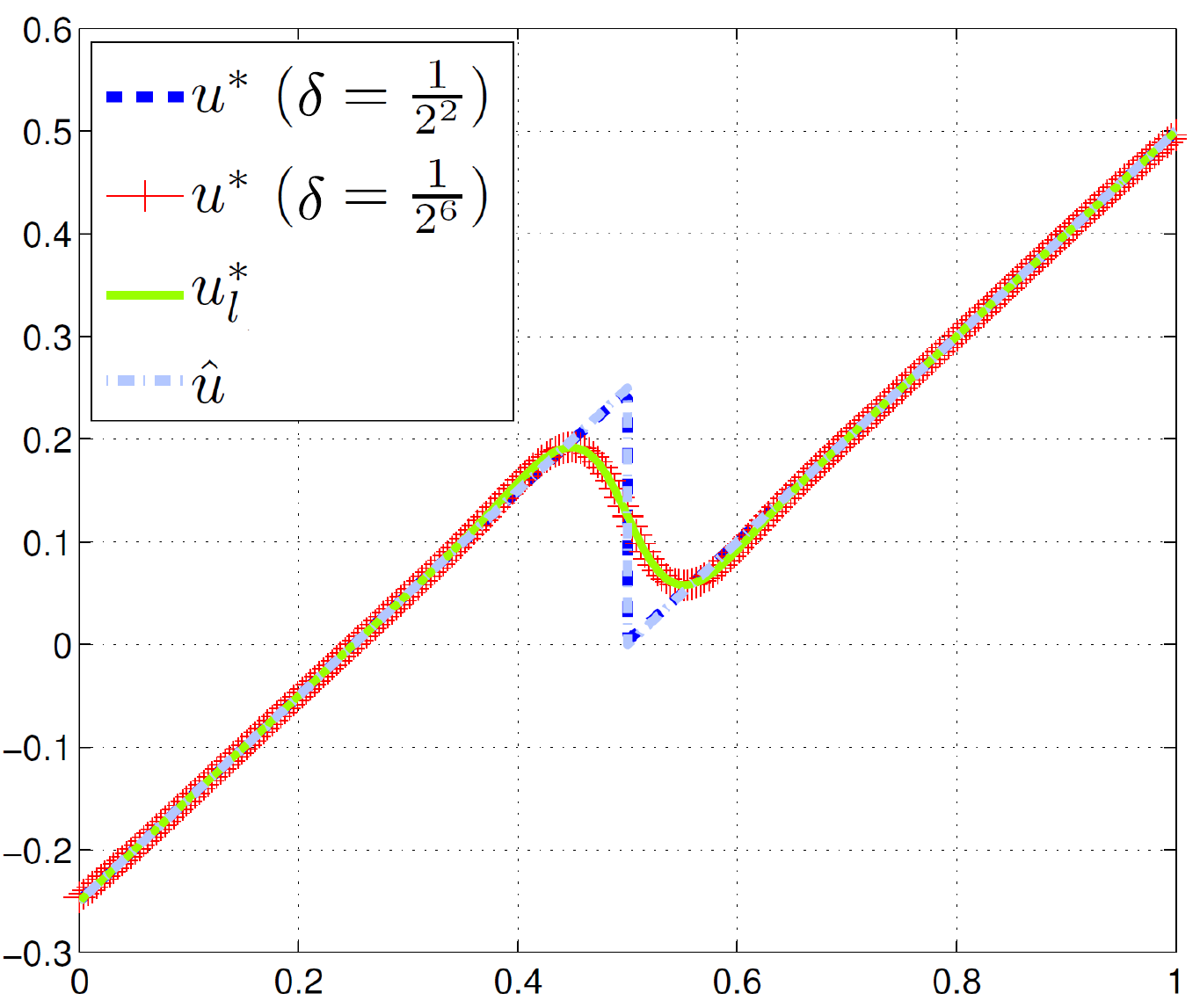}}\hfill %%control
\subfigure[]{\includegraphics[width=177pt,viewport=393 0 786 330,clip]{Figures/fig16_1.pdf}}\\ %%control
\caption{{(a)~The} optimal nonlocal state $u^*$ for $\delta=2^{-2}$ and $2^{-6}$, the target function $\widehat u$, and the optimal local state $u_l^*$. {(b)~The} corresponding optimal source terms.}
\label{fig:control-paper}
\end{figure}

\subsubsection{Coefficient identification}

The identification of kernel parameters or of the kernel function itself is one of the most important open problems in nonlocal modelling. In fact, in general, the choice of kernel {function} and its parameters is based on intuition or designed through heuristic techniques. Here we report on three approaches that tackle the kernel identification problem using different \bo{objective}  functionals. 

{\citeasnoun{Fuensanta2015OptimalDesign} considered the} problem of identifying a constitutive function in \eqref{eq:truncated-laplacian} having the form \bo{$\theta_\del \xyp=\thet(\xb)+\thet(\yb)$, that is, to identify the parameter function $\mumu=\thet(\xb)$.} In addition, $\mcV u=u$ and $\phi_\del\xyp=k\xyp|\ymx|^{-2}$ with $k\xyp\geq C|\ymx|^{2-\ddd-2s}$ {for $s\in(0,1)$} and ${\rm support} (k{\xyp}) =\ballx$. The set of admissible $\thet(\xb)$ is defined as
\[
\WWW=\{\thet(\xb)\in[\thet_{\min},\thet_{\max}], \; \thet = 0 \;{\rm on} \;\omgi, \; \int_\omg \thet(\xb) \D \xb=\overline\thet\}
\]
for positive constants $\thet_{\min}$, $\thet_{\max}$ and $\overline\thet$. The \bo{objective functional} in \citeasnoun{Fuensanta2015OptimalDesign} \bo{consists of a cost functional only and it is} referred to as a {\em compliance} functional; \bo{it is defined as}
\begin{equation}\label{eq:J-Fuensanta2015}
\mcJ(u)=
\int_{\Omega\cup\Omega_I}\int_{\Omega\cup\Omega_I}
F(\xb,\yb;u) \dydx
\end{equation}
and is then minimized over $\thet(\xb)\in\WWW$. First, $F$ is chosen such that $\mcJ(u)=\trip{u}^2$, where $\trip{\cdotsp}$
denotes the `energy' norm corresponding to the kernel function $\phi\xyp=k\xyp|\ymx|^{-2}$. For this choice of $F$ and with no other constraints applied on the state or control, the existence of a solution of the problem of minimizing the functional \eqref{eq:J-Fuensanta2015} is proved. Note that a regularizing term is not included in this functional because the box {constraints} on $\thet(\xb)$ included in the admissibility set $\WWW$ already guarantee the {well-posedness} of the minimization problem.  Extensions to more general functionals are then considered, including $F$ {only being} required to be measurable with respect to $u$ and lower semi-continuous with respect to $\xb$. Furthermore, for the compliance case, the convergence of the optimal nonlocal solution to its local counterpart \eqref{eq:generalized-control-local} for $\mumu=\thet$ is proved. In this work, neither discretizations {nor} numerical tests are provided.

{\citeasnoun{DElia2016ParamControl} also considered the problem of identifying the constitutive function \bo{$\theta_\delta\xyp$} in \eqref{eq:truncated-laplacian} for} \bo{square-integrable and also}
non-integrable kernel functions $\phid\xyp$ {-- specifically,} again $\mcV u=u$, and, in the functional \eqref{eq:matching-Tikhonov}, $\beta=0$ and $\|\cdotsp\|_\UUU=\|\cdotsp\|_{L^2(\omg)}$. The functional is minimized over the set of admissible controls given by
\[
\bo{
  \WWW\zz = \zz \{\theta_\delta \zz \in \zz  W^{1,\infty}(\omgomgi \zz \cup\omgomgi),\, \theta_\delta \zz \in \zz [\thet_{\min},\thet_{\max}],\, \|\theta\|_{1,\infty}\leq  \zz \theta_\delta^{\max} \zz < \zz \infty\}.}
\]
Again, the box constraints on \bo{$\theta_\delta$} included in the admissibility set suffice to prove {that problem} \eqref{eq:generalized-control} has at least one solution. {\citeasnoun{DElia2016ParamControl} developed a mixed finite element discretization of the state and control variable  and proved} the convergence of the discretization error as the mesh size is refined. Numerical tests are also provided that show that the approach taken there allows for the identification of both smooth and discontinuous diffusion coefficients for both \bo{square-integrable} and (truncated) fractional kernels. A sample result is shown in \fo{Figure~\ref{fig:identification-paper}}. For that figure, $\omg=(0,1)$, $\omgid=(\del,0)\cup(1,1+\del)$, $g(x)=0$, $f(x)=5$, and the spatial grid size used to discretize the state is $2^{-12}$. A surrogate for the target functional $\widehat u$ is a very {fine-grid} finite element approximation of the nonlocal diffusion problem with the data just listed and with
\begin{equation}\label{eq:parametercoeff}
  {\theta_\delta(x,y)} = \thet \biggl(\ffrac{x+y}2\biggr) 
  \quad\text{with }
  \thet(z) = 
  \begin{cases}
  1 & \text{for all } z\in (0,0.2),  \\[2pt]
  0.1 & \text{for all } z\in (0.2,0.6),  \\[2pt]
  1 & \text{for all } z\in (0.6,1).
  \end{cases}
\end{equation}
Thus the goal of the minimization problem is to identify this constitutive coefficient function. \fo{Figure~\ref{fig:identification-paper}} illustrates, for two values of the horizon $\del$, the convergence (with respect to the grid size $1/N_{\thet}$ used to approximate the coefficient function $\thet(z)$) of the approximation.
\begin{figure} %%fig16.2 %%wasfig15.2
\centering
\subfigure[]{\includegraphics[width=173pt,viewport=0 0 402 335,clip]{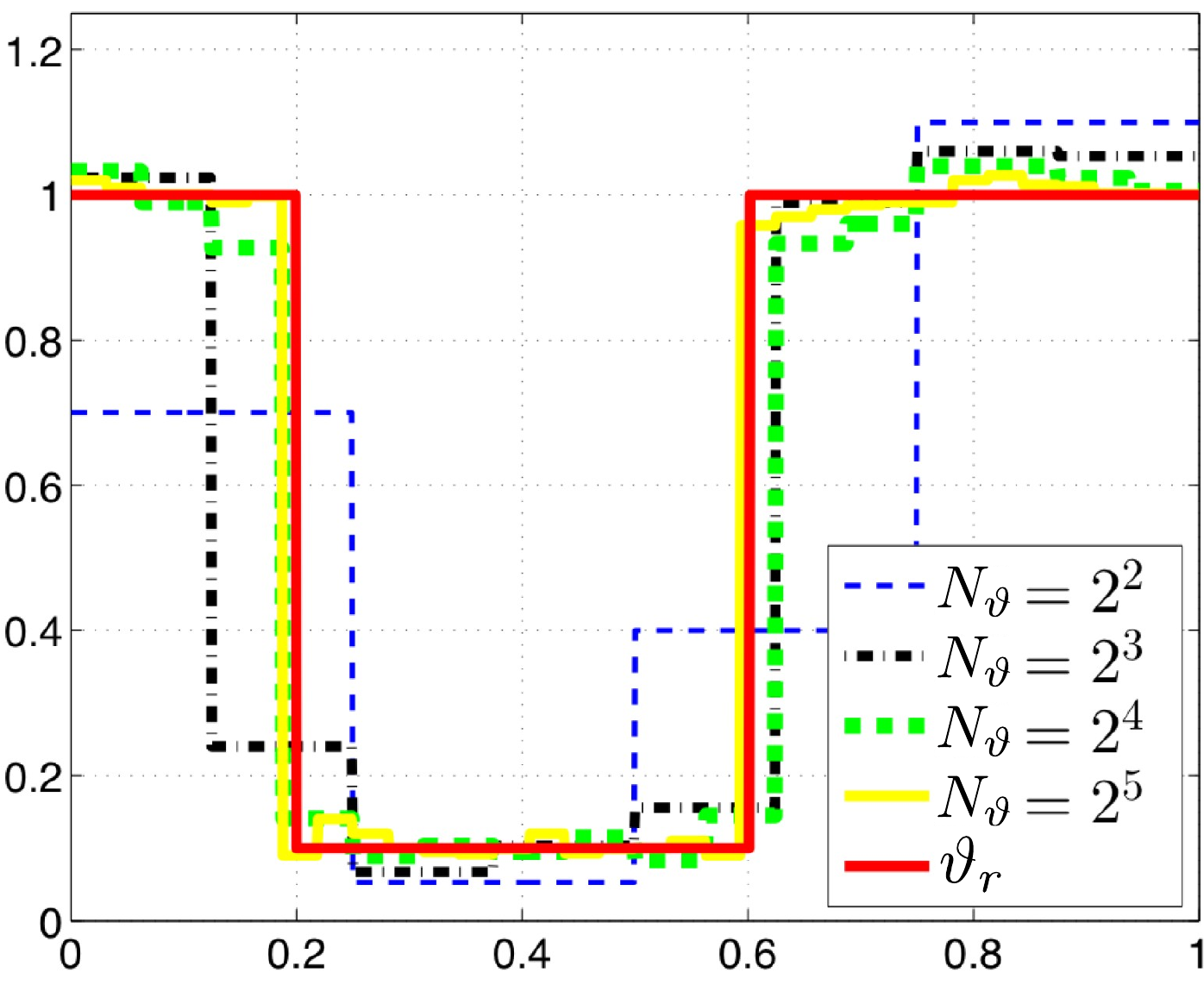}}\hfill %%identification
\subfigure[]{\includegraphics[width=173pt,viewport=432 0 834 335,clip]{Figures/fig16_2.pdf}}\\ %%identification
\caption{Optimal approximate coefficient functions for $\del=2^{-9}$~{(a)} and $2^{-4}$~{(b)} for different numbers of degrees of freedom $N_\vartheta$ in the discretization of the coefficient function $\thet(z)$ defined in \eqref{eq:parametercoeff}.}
\label{fig:identification-paper}
\end{figure}

We also mention the approach introduced in \citeasnoun{Pang2019nPINNs}, in which the task of parameter identification for truncated diffusion operators of fractional type is pursued by including the nonlocal diffusion equation in the \bo{objective} functional, \tie, we have 
\begin{equation}\label{eq:J-pang2019}
\mcJ(u;\mub)=\dfrac{1}{2} \|u-\widehat u\|_{L^2(\Omega)}^2 + \dfrac{\beta}{2}\|-\mcL u -f\|^2_{L^2(\Omega)}.
\end{equation}
Thus the state equations are only weakly prescribed through the minimization of the residual. {\citeasnoun{Pang2019nPINNs} approximated the} state and control variables using fully connected neural networks. Numerical tests on analytic solutions and turbulence models show the ability of this technique to recover kernel parameters such as the interaction radius $\del$ and the variable fractional order $s(\xb)$, \tie, we have $\mumu=\{\del,s(\xb)\}$.

\subsection{Inverse problems for fractional operators}\label{sec:truncated-control}

In this section we focus on control and optimization problems for the fractional Laplace operator. The formulations described in what follows can be easily extended to more general fractional operators (see \eg\ \citeb{Meerschaert2012book}), but at a cost of more complicated analyses.

\subsubsection{Distributed optimal control \bo{in fractional diffusion} for a matching functional}

The formulation presented in the previous section for distributed control can be used for fractional operators with (almost) no modification. We mention several works in the literature that analyse the theoretical and numerical aspects of both the elliptic and parabolic fractional problems. 

{\citeasnoun{DElia2019FracControl} considered the} same problem as that introduced in Section~\ref{sec:controlnonlocal}. Specifically, for ${\mcJ(u;f)}$ as in \eqref{eq:J-delia2014} and $\mcV u=u$, {they considered} the problem of minimizing ${\mcJ(u;f)}$ with respect to
$f\in\WWW=\{r\in L^2(\omg)\colon r \in[r_{\min},r_{\max}]\}$
subject to
\begin{equation}\label{eq:diffusion-z-RHS}
\begin{cases}
(-\Delta)^s u  = f & \text{for all }\xb\in\omg ,\\
u = g & \text{for all } \xb\in\Rd\setminus\omg .
\end{cases}
\end{equation}
 The {well-posedness} of the control problem is {proved}, optimality conditions are derived, and regularity estimates for the optimal variables are also {proved}.  Furthermore, based on an {\it a~priori} error analysis for the state equation, a semidiscrete scheme is constructed for which {\it a~priori} error estimates for the approximation of the control variable are derived. A fully discrete scheme is also considered for which state and control error estimates are derived. Several two-dimensional numerical illustrations of the theoretical results are also provided. 

 {\citeasnoun{Glusa2019ParabControl}, in a follow-up paper, considered} the fractional parabolic equation \bo{\eqref{eq:diffusion-z-RHS2}}.
The \bo{objective} functional, similar to that in \eqref{eq:J-delia2014}, now involves an integral over time \bo{for both the cost and regularization terms}, \tie,  
\begin{equation}\label{eq:J-glusa2019}
\mcJ(u;f)=\dfrac{1}{2} \int\limits_0^T 
 (\|u-\widehat u\|_{L^2(\Omega)}^2 + 
\beta\|f\|^2_{L^2(\Omega)})  \D t.
\end{equation}
    {Further,} the control variable $f$ is now also {time-dependent} and belongs to the admissibility set
    \[ %%20200331
\roo{    \WWW(t)=\{r(t)\in L^2(\Omega) \colon r(t)\in [r_{\min}(t),r_{\max}(t)]\ \text{for all}\ t\in(0,T]\}.}
      \]
      As in \citeasnoun{DElia2019FracControl}, the existence and uniqueness of optimal solutions are {proved,} and first-order necessary and sufficient optimality conditions are derived. Also derived are regularity estimates for the optimal state and control. Then discrete stability results and {\it a~priori} error estimates are derived for the discretized problem resulting from the standard backward Euler scheme for temporal discretization and a piecewise linear finite element spatial discretization. The theoretical findings are illustrated by one- and two-dimensional numerical experiments.

{\citeasnoun{Antil2019ExternalControlS} and \citeasnoun{Antil2019ExternalControl} chose the} control as the data $g$ in the volume constraint for the steady-state and time-dependent cases, respectively. In \eqref{eq:generalized-control}, nonlocal Dirichlet, Neumann and Robin volume constraints are considered for the operator $\mcV$. As an example, in the Robin case we have 
\[
\mcV u  =\kappa_N \mcN u + \kappa_D u = g,
\]
where $\mcN$  \bo{is} defined as in \eqref{str-flux} with the appropriate changes made to reflect that here the interaction domain is $\omgii=\Rd\setminus\omg$. Even though the theory is presented for the general functional 
\begin{equation}\label{eq:J-antil2019external}
\mcJ(u;g)= {\mcQ(u)} +
\dfrac{\beta}{2}\int\limits_0^T\|g\|^2_{L^2(\omg_{\mathcal I})} \D t
\end{equation}
with a convex \bo{cost} functional \bo{$\mcQ(u)$}, numerical experiments are performed for a matching functional in the usual form; see \eg\ \eqref{eq:J-glusa2019} for the time-depend\-ent version. In the more general formulation, for the time-dependent case, {\citeasnoun{Antil2019ExternalControl}  mostly focus} on nonlocal Dirichlet and Robin optimal control problems. {Well-posedness} and regularity are discussed and a discretization scheme is proposed. The theoretical results are illustrated by two-dimensional numerical experiments. The main contribution of \citeasnoun{Antil2019ExternalControl} is to show the ability of nonlocal models to take advantage of information outside the domain and not only on the boundary, which is one of the limitations of control problems for PDEs.

Even if not entirely focused on operators such as that in \eqref{intfl}, we mention that {\citeasnoun{Antil2019FracControl} consider control} problems for both a spectral fractional semilinear operator and for the integral fractional Laplacian. {They derive existence and regularity results for} the spectral case ({which} is not treated in this section).
{They also state that} the results obtained for the spectral case can be extended, after small modifications, to the integral definition of the operator. 

\setcounter{rem}{1}
\medskip\noi{\bf Remark 15.1 (kernel identification).}
For kernel identification in the setting of fractional operators, {we should mention}  that {\citeasnoun{Pang2019fPINNs} studied an} algorithm for parameter identification based on physics-informed neural networks. As such, this method is a special instance of the algorithm presented in \citeasnoun{Pang2019nPINNs}. Specifically, this work is focused on the estimation of the fractional power of the integral fractional Laplacian. Other works on the identification of kernel parameters, more specifically, of the fractional power $s$, only deal with the spectral definition of $(-\Delta)^s u$.
\medskip

\setcounter{rem}{2}
\noi\bo{{\bf Remark 15.2 (control for the spectral fractional Laplacian).}
We should also mention that control of equations involving the spectral definition of the fractional Laplacian has been \roo{analysed; see \eg}\  \citename*{Antil2016SpaceTimeControl} \citeyear{Antil2016SpaceTimeControl,Antil2018sSpectralControl} and \citeasnoun{Otarola2015SparseControl}.}

\vspace{5pt}
\section{Variational inequalities and obstacle problems}\label{sec:vi}

We consider the nonlocal obstacle problem
\begin{equation}\label{obsprob}
\begin{cases}
 -\mcL u \geq f   &    \text{for all } \xb\in \omg , \\ 
 u\geq\psi          &  \text{for all } \xb\in \omg ,\\
(-\mcL u-f)(u-\psi)=0 & \text{for all } \xb\in \omg ,\\
u=0 & \text{for all } \xb\in \omgii,
\end{cases}
\end{equation}
{where $\psi(\xb)$} denotes the obstacle function. Nonlocal obstacle problems such as this one are used in studying the deformation of elastic membranes, in contact mechanics and in finance, \feeg\ the pricing of American put options in L\'evy jump-diffusion models. Clearly, \eqref{obsprob} is a nonlocal analogue of the local PDE obstacle problem
\[
\begin{cases}
 -\Delta u \geq f    &  \text{for all } \xb\in \omg , \\ 
 u \geq\psi       &    \text{for all } \xb\in \omg ,\\
(-\Delta  u-f)(u-\psi)=0 & \text{for all } \xb\in \omg ,\\
u=0 & \text{for all } \xb\in \partial\omg.
\end{cases}
\]
The well-posedness analysis of nonlocal obstacle problems {needs} less smooth obstacles {than the corresponding} local PDE obstacle problems. Moreover, as is illustrated at the end of this section, the behaviours and properties of  solutions in the local and nonlocal setting can be quite different. 

There is an ever-growing literature on the analysis and approximation of nonlocal obstacle problems, especially in the fractional setting. Here we give a brief account {of} an approach used {by \citeasnoun{burko1}, who considered regularity estimates, well-posedness analyses, and finite element methods and their numerical analysis}. Using a different approach, {\citeasnoun{guan2017} established the well-posedness of nonlocal obstacle {problems} and the convergence of the finite element approximation for} fractional Laplacian \bo{kernels and for kernels that are both radial and integrable}. With respect to other works about obstacle problems for the fractional Laplacian, we mention {\citeasnoun{servadei}, who obtained Lewy--Stampacchia-type estimates similar to those obtained by \citeasnoun{burko1}}, but with restrictions on the fractional exponent and requiring greater smoothness of the obstacle.
The regularity of the obstacle problem measured in H\"older and Lipschitz spaces {was studied by \citeasnoun{silvestre2007} and \citeasnoun{caffarelli2017}, \feeg. \citeasnoun{borthagaray_preprint_18} studied a finite element approximation of the obstacle problem for the fractional Laplacian, and proved} error estimates. Regularity results for the solution are derived in weighted Sobolev spaces  under additional regularity assumptions on the right-hand side (H\"older continuity) and the obstacle. {\citeasnoun{bonito_preprint_18} studied the regularity of the obstacle problem involving integro-differential operators, with the fractional Laplacian as the nonlocal term, and proposed and analysed a finite element-based} discretization. In the purely nonlocal case, the same regularity for the solution is {proved} as in \citeasnoun{burko1}, but for more restricted cases. Other than \citeasnoun{burko1}, none of these works on the obstacle problem for the fractional Laplacian treat truncated kernels.

{\citeasnoun{burko1} obtained the well-posedness and regularity results for the nonlocal obstacle problem \eqref{obsprob},
   and used} the mixed formulation 
\begin{equation}\label{mixform}
\begin{cases}
\mcA(u,v)-\mcB(\laa,v)=\lrang{f}{v} &\text{for all } v\in V ,
\\
\mcB(\eta-\lambda,u-\psi) \geq 0 &\text{for all } \eta\in M\subset V_d
\end{cases}
\end{equation}
{to} define, analyse and apply finite element methods. In \eqref{mixform}, $\mcA(\cdotsp,\cdotsp)\colon V\times V\to\Ro$ is the usual bilinear form corresponding to the nonlocal operator $\mcL$, $V$ is the energy space  associated with that bilinear form and the homogeneous volume constraint, $V_d$ is the dual space for $V$, $\mcB(\cdotsp,\cdotsp)\colon V_d\times V\to\Ro$ is defined as $\mcB(\eta,v)=\lrang{\eta}{v}$, and $M$ denotes the closed convex dual cone
$$M:=\{\eta\in V_d\colon \lrang{\eta}{v}\geq 0\ \text{for all}\  v\in V, v\geq 0 \}.$$
Of course, $\mcA(\cdotsp,\cdotsp)$ is continuous and coercive in $V\times V$, and {\citeasnoun{burko1} showed that} $\mcB(\cdotsp,\cdotsp)$ is  continuous and inf-sup stable on $V_d\times V$, \tie, we {have}
\[
\inf_{\eta\in V_d}\sup_{v\in V}\ffrac{\mcB(\eta,v)}{\|\eta\|_{V_d}\|v\|_{V}}\geq\beta_0>0.
\]
So, the task at hand is to find, for  $f\in V_d$ and $\psi\in V $, $u\in V$ and $\laa\in M$ satisfying \eqref{mixform}.

Under the assumptions $f\in L^2(\omg)$ and $(-\mcLd\psi-f)_+|_\omg\in L^2(\omg)$, the improved regularity results
\[
\mcLd u\in L^2(\omg), \quad\laa\in L^2(\omg),\quad \laa\leq(-\mcLd \psi-f)_+
\]
are derived. For the fractional Laplacian obstacle problem, the improved regularity results obtained are given by
\[
\laa\in L^2(\omg)\quad\text{and}\quad u\in H^{s+\beta}_\omg(\omg)
\]
{with $ \beta=\min\{s,1/2-\varepsilon\}$, $\varepsilon>0$, $ s\in(0,1)$.}
Note that these regularity results hold for all $s\in(0,1)$.

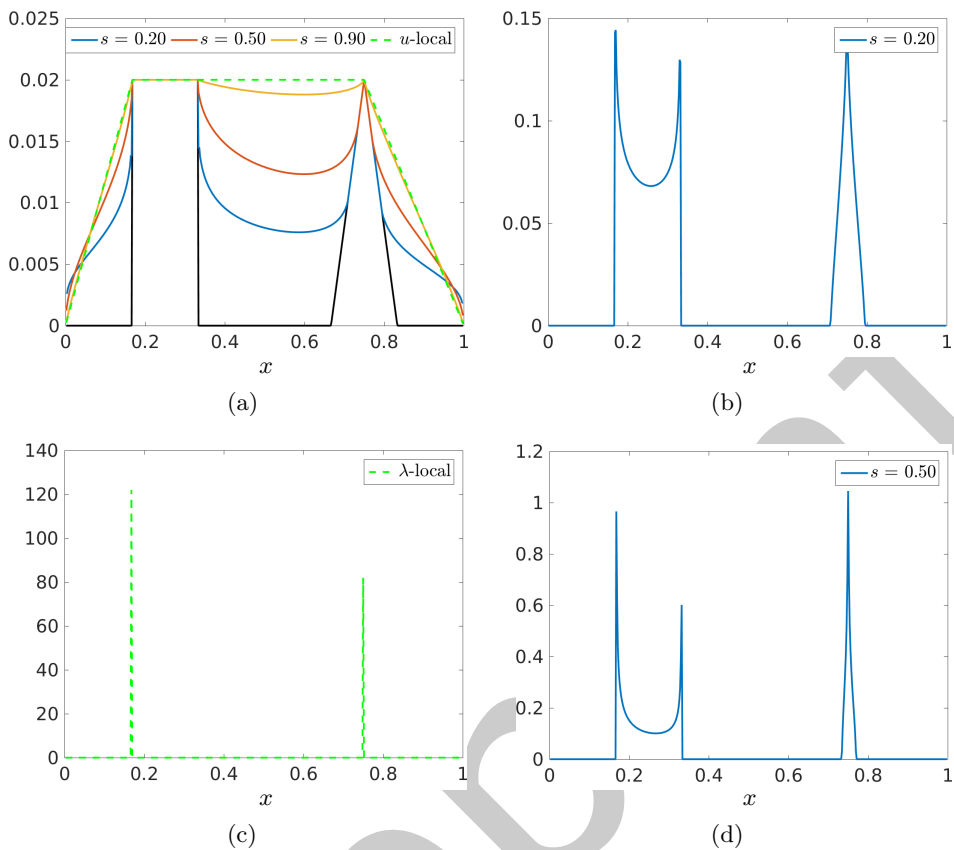
\begin{figure} %%fig17.1 %%wasfig16.1
\centering
\setlength{\unitlength}{1pt}
\begin{picture}(360,320)
  \put(0,186){\includegraphics[width=173pt,viewport=0 30 1075 835,clip]{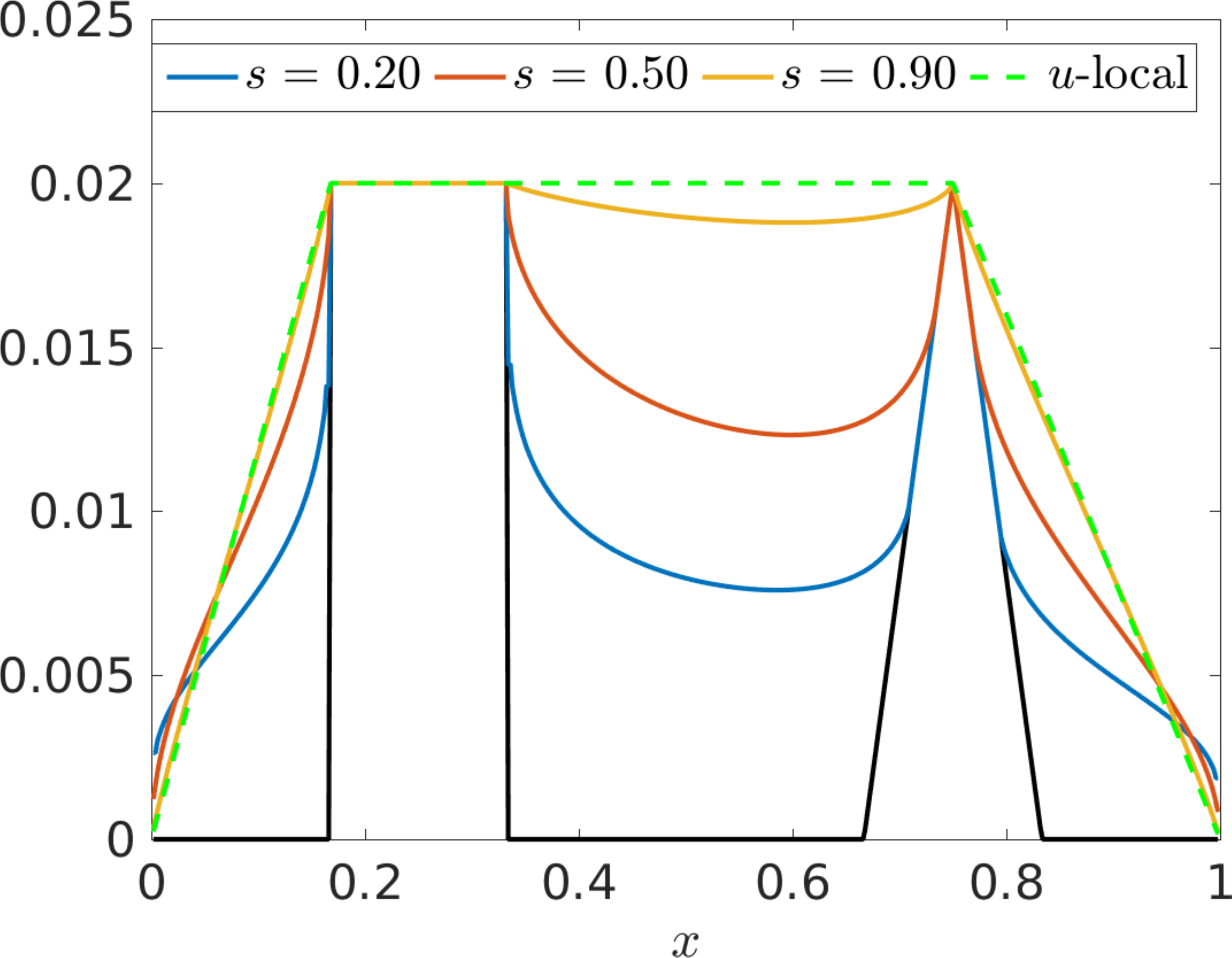}}
  \put(186,186){\includegraphics[width=173pt,viewport=0 30 1075 835,clip]{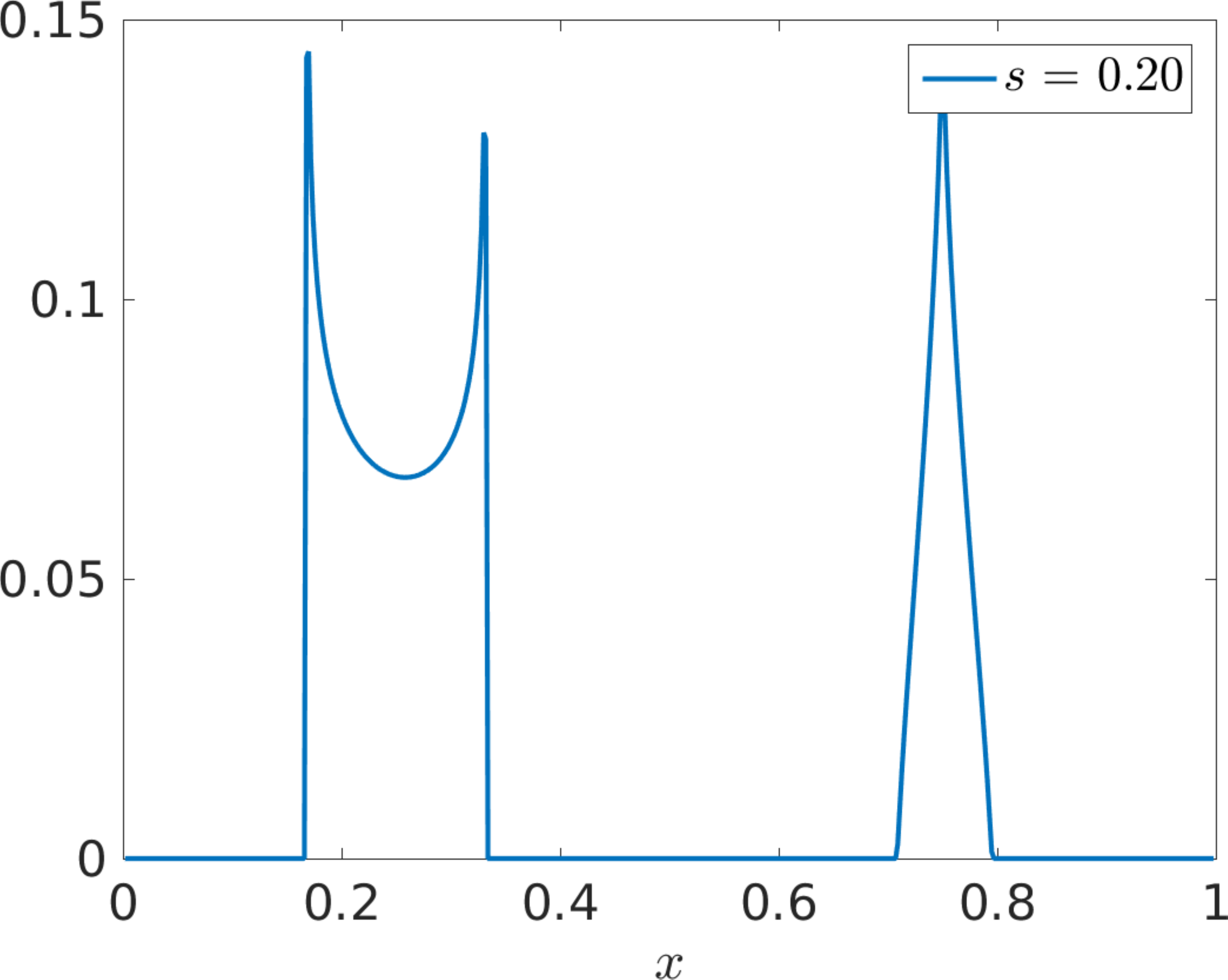}}
  \put(0,24){\includegraphics[width=173pt,viewport=-40 30 1035 835,clip]{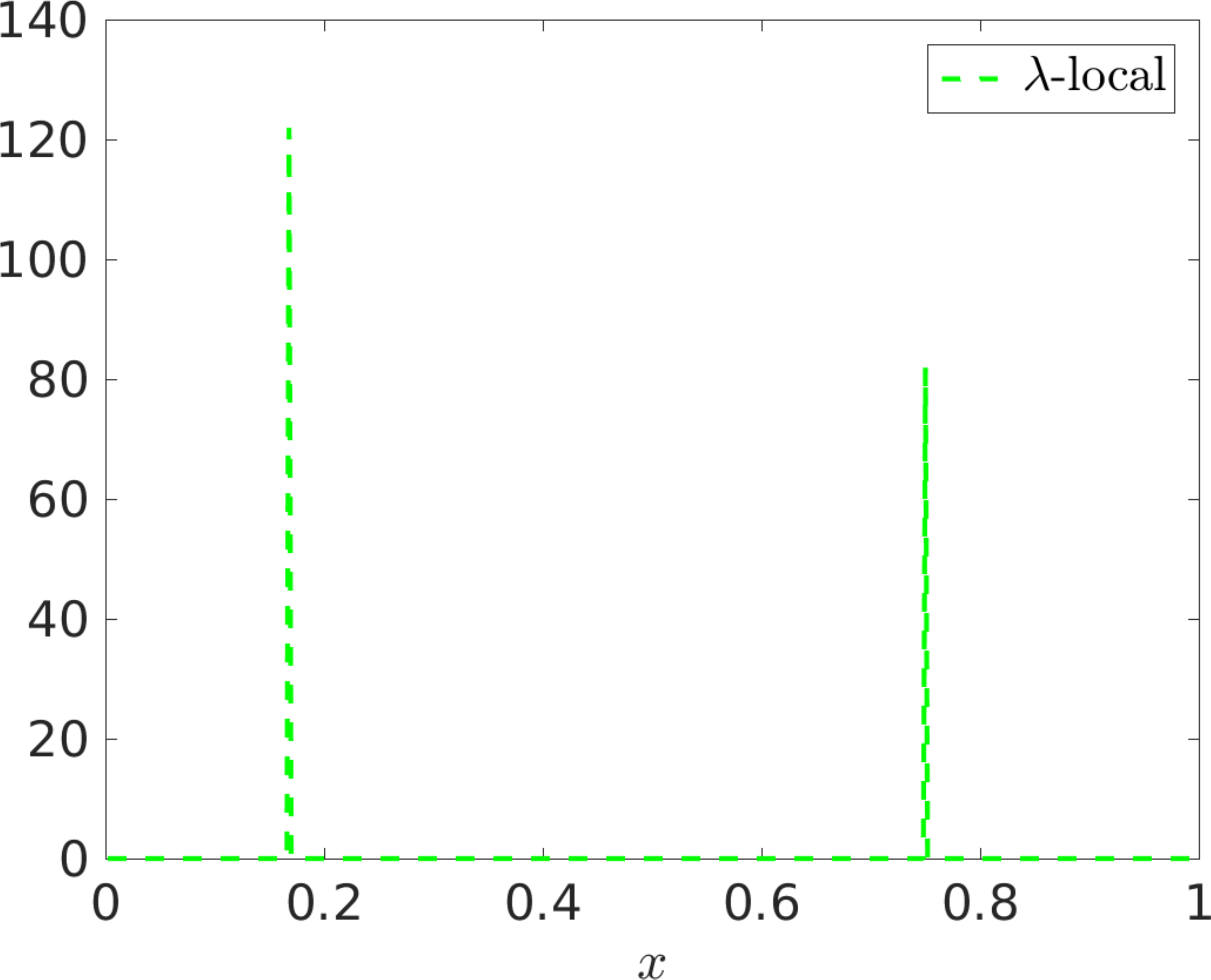}}
  \put(186,24){\includegraphics[width=173pt,viewport=-30 12 1045 820,clip]{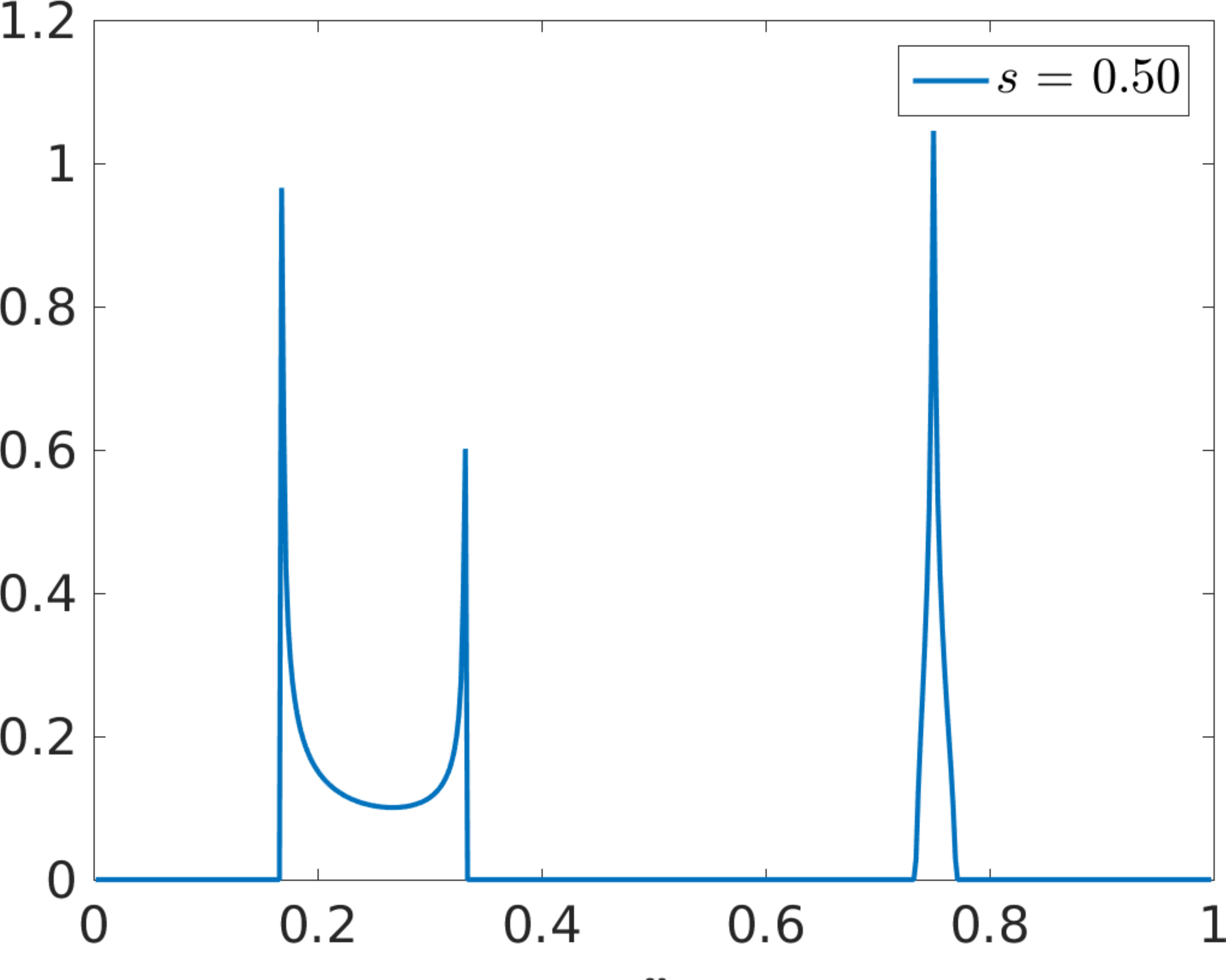}}
  \put(94,180){\fns{$x$}}
  \put(276,180){\fns{$x$}}
  \put(94,18){\fns{$x$}}
  \put(276,18){\fns{$x$}}
  \put(82,166){\fns{(a)}}
  \put(264,166){\fns{(b)}}
  \put(82,4){\fns{(c)}}
  \put(264,4){\fns{(d)}}
\end{picture}\\
\caption{{(a)~Primal} nonlocal solutions for three values of the fractional exponent and the primal local solution. {(b)~The} local Lagrange multiplier.
  {(c,\,d)~The} nonlocal Lagrange multiplier for $s=0.2$~{(c)} and $s=0.5$~{(b)}. Note the different ordinate scales in the three Lagrange multiplier plots.}\label{fig:obstex}
\end{figure}

Finite element approximations are defined for the space $V^h\subset V$ of {piecewise linear} continuous polynomials and the space $V_d^h\subset V_d$ of discontinuous linear polynomials
for which locally bi-orthogonal basis functions can be constructed, \tie,  for any basis functions $\xi_j(\xb)\in V_d^h$ and $\phi_{j'}(\xb)\in V^h$ and for any finite element $K$, {we have}
\[
\int_{K}\xi_j\phi_{j'}=\delta_{jj'}\int_{K}\phi_j\geq 0.
\]
The inf-sup stability with respect to $V_d^h\times V^h$ is {proved, so the} discrete problem is {well-posed}. 

%%fig 17.1 original position

An example numerical result is provided in \fo{Figure~\ref{fig:obstex}}. Note the differences in the primal solution $u(\xb)$ and in the Lagrange multiplier $\laa$ obtained using nonlocal and local models. In particular, note that for the nonlocal case $\laa\in L^2(\omg)$, whereas for the local case $\laa$ consists of Dirac delta functions.

%% fig16.1 original position

Note that all the results in \citeasnoun{burko1} pertinent to the fractional Laplacian were also proved for the truncated fractional Laplacian kernel introduced in Section~\ref{truncted-flap}. Moreover, the convergence of the solution of the obstacle problem with the truncated kernel is shown to converge to the solution of (the un-truncated) fractional Laplacian obstacle problem.

\vspace{5pt}
\section{Reduced-order modelling}\label{sec:rom}

Reduced-order modelling (ROM) is the task of constructing a very low-dim\-ensional discretization for {parametrized} problems that, in order to achieve a desired fidelity, are usually approximated by a high-dimensional discretization (HDD). The construction of a ROM usually requires an off-line cost incurred by having to do runs of the HDD {for relatively} few parameter choices. Once the ROM is constructed, it can be used {on-line} instead of the HDD by, \feeg, doing additional simulations for many more parameter choices at much lower cost than if the HDD {were} used instead. For example, in uncertainty quantification settings, one may need to obtain many simulations in order to obtain good statistical information, {so using} the ROM instead of the HHD can result in huge computational savings. The huge body of literature {on} ROMs for PDEs {attests} to their usefulness. 

PDE discretizations ubiquitously involve large, sparse linear or nonlinear discrete systems. Analogous discretizations of nonlocal models usually involve discrete systems of similar size but, due to nonlocality, having much less sparsity. In some settings, such as those involving fractional Laplacians, the discrete systems may even be full. The double curse of high-dimensionality and lack of sparsity of the discrete systems means that{:}
\begin{equation*}
  \begin{minipage}{23pc}
{Reduced-order} modelling is needed much more for nonlocal
models than for corresponding local PDE models.
  \end{minipage}
\end{equation*}

The setting we use to discuss ROMs is parametrized nonlocal diffusion. Thus we assume that the constitutive function $\thed(\xb,\yb;\pbt)$ in \eqref{gamd1} depends on the components of a parameter vector $\pbt\in\Gamt\subset\RNt$, and the kernel function $\phid(\xb,\yb;\pbp)$ in \eqref{gamd1} depends on the components of a parameter vector $\pbp\in\Gamp\subset\RNp$. We refer to $\Gamt$ and $\Gamp$ as {\em parameter domains.} Then, in the weak formulation \eqref{weak-weakf} of the nonlocal diffusion problem, we have the bilinear form
\begin{align*}
&\mcA_{\pb}(u,v)
  \\*
  &
=\int_{\omgomgi}\int_{\ballx}  (v(\yb)-v(\xb))  (u(\yb)-u(\xb))  \phid(\xb,\yb;\pbp)\thed(\xb,\yb;\pbt)
 \dydx,
\end{align*}
where we have used \eqref{gamd1}, \ie\ $\gamd\xyp= \phid(\xb,\yb;\pbp) \thed(\xb,\yb;\pbt) \mcX_{\ballx}(\yb)$. We also assume that we have in hand a finite element or other variational discretization\footnote{ROMs for discretizations of the strong form of nonlocal models are also of interest.} of \eqref{weak-weakf} with the bilinear form $\mcA_{\pb}(\cdotsp,\cdotsp)$, \tie, we have the nonlocal discrete problem for $\udhp(\xb)\in V^h_c(\omgomgi)$ given by
\begin{equation}\label{weak-weakfph}
   \mcA_{\pb}(\udhp,v_h) = \langle f,v^h \rangle \quad\text{for all } v_h\in V^h_c(\omgomgi),
\end{equation}
where, for economy of exposition, we impose the homogeneous Dirichlet volume constraint $\udhp(\xb)=0$ on $\omgid$ and where $V^h_c(\omgomgi)$ denotes
{a subspace of $\VVc$, \eg\ one generated via finite elements.} Thus, given any parameter vectors $\pbp\in\Gamp$ and $\pbt\in\Gamt$, \eqref{weak-weakfph} can be solved to determine the approximation ${\udhp}(\xb)$ of the exact solution $u_{\del,\pb}(\xb)$ of \eqref{weak-weakf} with $\mcA(\cdot,\cdot)$ replaced by $\mcA_{\pb}(\cdot,\cdot)$. We let $N_h$ denote the dimension of  $V^h_c(\omgomgi)$, \ie\ of the HDD.

Assume we have a ROM basis $V_{N_{\text{rom}}}:=\{v_{n,\text{rom}}(\xb)\}_{n=1}^{N_{\text{rom}}}$ consisting of $N_{\text{rom}}$ functions $v_{n,\text{rom}}(\xb)\in V^h_c(\omgomgi)$. Then, for any $\pbp\in\Gamp$ and $\pbt\in\Gamt$, the {\em ROM-Galerkin approximation} $u_{\pb,\text{rom}}(\xb)\in V_{N_{\text{rom}}}$ of the solution $\udhp(\xb)$ of \eqref{weak-weakfph} is determined by solving the problem
\begin{equation}\label{weak-rom}
   \mcA_{\pb }(u_{\pb,\text{rom}},v) = \fv \quad\text{for all } v\in V_{N_{\text{rom}}}.
\end{equation}
Thus the task at hand is to construct a ROM basis such that $N_{\text{rom}}\ll \NNh$, {which} results in ROM solutions $u_{\pb,\text{rom}}(\xb)$ that are acceptably accurate approximations of $\udhp(\xb)$ for any $\pbp\in\Gamp$ and $\pbt\in\Gamt$, or at least for subsets of the parameter domains that are of interest.

Two examples of the use of reduced-order models for nonlocal diffusion are given {by \citeasnoun{guan}, who use a greedy reduced basis (GRB) approach, and \citeasnoun{witman}, who use a} proper orthogonal decomposition (POD) approach.\footnote{Although both ROMs we consider, {like} most others, involve a `reduced basis', \ie\ a basis of lower dimension than that of the parent HDD model, {\em reduced-basis methods} usually refer to ROMs in which the basis consists of solutions of \eqref{weak-rom}. In contrast, a POD basis consists of linear combinations of such solutions.} Both approaches involve only the parameter vector $\pbt$. The basis construction and application processes for nonlocal and local models are the same, excepting of course that in the local case approximate solutions of PDEs are involved whereas the nonlocal case involves approximate solutions of nonlocal models such {as} \eqref{weak-rom}. Thus, here, we simply refer to \citeasnoun{guan} and \citeasnoun{witman} for detailed descriptions of the GRB and POD approaches, respectively, in the context of nonlocal diffusion models. {\citeasnoun{witman} considered a time-dependent problem. \citeasnoun{guan} considered random parameter vectors  and compared the GRB surrogates to} sparse-grid surrogates. In addition to the analyses of the errors incurred by the ROMs, \citeasnoun{guan} and \citeasnoun{witman} provide numerical examples that illustrate the usefulness of the ROMs considered therein. 

{\citeasnoun{burko2} considered the} construction and analysis of ROMs for nonlocal models parametrized by the horizon $\del$ and, {for} fractional models, by the fractional exponent $s$, so that now the parameter vector $\pbp$ is present in the kernel function. The parameter domains for both $\del$ and $s$ are intervals bounded away from the origin and infinity, and greedy algorithms are used to construct the reduced bases. Also provided in \citeasnoun{burko1} are illustrative numerical examples. 

\begin{figure} %%fig18.1 %wasfig17.1
\centering
\subfigure[]{\includegraphics[width=175pt]{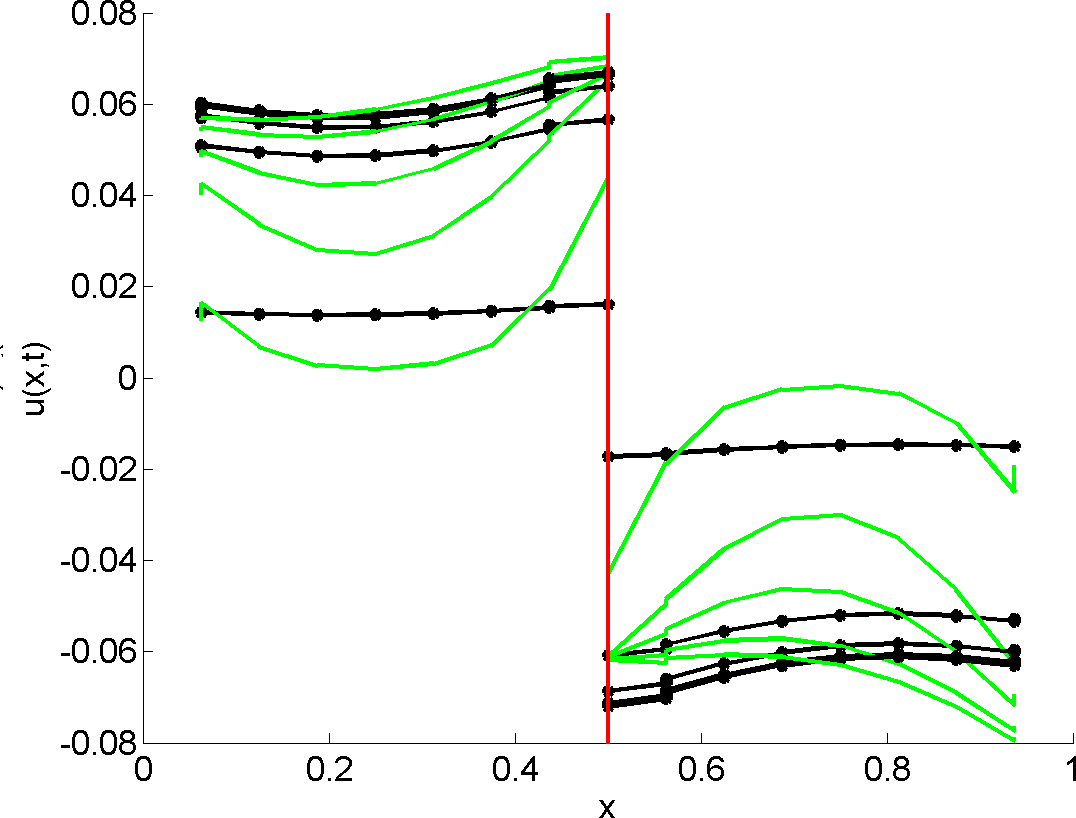}}\hfill %%witman1
\subfigure[]{\includegraphics[width=175pt]{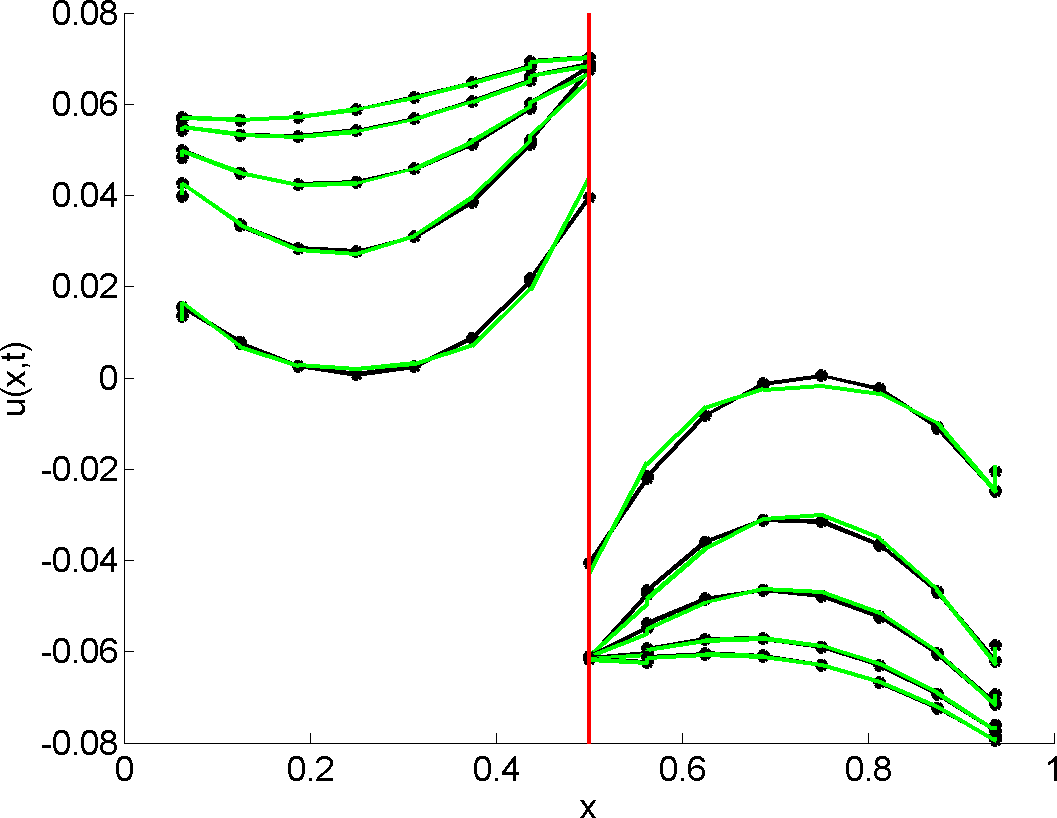}}\\ %%witman2
\caption{Comparison at different times between a fine-grid discontinuous Galerkin finite element approximation ({green line}) and POD approximations ({black line} with black dots) with one POD basis function~{(a)} and six POD basis functions~{(b)}.}
\label{fig:rom}
\end{figure}

A very simple one-dimensional, time-dependent illustration of the effectiveness of ROM in the nonlocal setting is given in \fo{Figure~\ref{fig:rom}}. The kernel function is proportional to $1/|y-x|$, {so it} is singular and non-integrable. A manufactured solution that has a jump discontinuity at all times is used to set the data of the problem. A discontinuous finite element Galerkin method is used {both for} generating the snapshots for the POD basis construction and for comparison purposes. Plots for five time instants are provided. From the figure, it is evident that a one-dimensional POD basis does a poor job of approximating the HDD finite element solution whereas a six-dimensional POD basis does {an} excellent job.

%%fig18.1 original position

\begin{rem}\hspz% 
As is the case for the local PDE setting, obtaining {\em solutions} of nonlocal ROM models such as \eqref{weak-rom} {incurs} costs that depend only on the dimension $N_{\text{rom}}$  of the ROM model, whereas the naive {\em assembly} of the ROM models involves steps whose cost depends on the dimension $N_h$ of the HDD models used to generate the ROM bases. In the PDE setting, several strategies have been developed to overcome this obstacle, all of which can be applied to the nonlocal setting. \hfill\actaqed
\end{rem}

\vspace{2pt} %%\vspace{5pt} %%20200331
\section{A turbulent flow application}\label{sec:richarson}

A fundamental problem of fluid flows is to ascertain information about the paths taken by two initially closely positioned particles, in particular what happens, over a {long time} period,
to the separation between a pair of particles. Such knowledge is crucial, \feeg, to {predicting} how pollutants spread. Intuitively, over a long time period, for {\em laminar flows} one expects that the separation between initially close particle pairs remains `small' and, in fact, the mean-square separation is proportion to $t$. Also intuitively, over a long time period, for {\em turbulent flows} one expects that at least some initially close particle pairs may become widely separated. What happens to the separation between pairs of particles in turbulent flows is the issue addressed in the classic paper {by} \citeasnoun{richardson}. This subject, often referred to {as} {\em Richardson pair dispersion}, {has remained of interest to the present day.}
Reviews on the subject of Richardson pair dispersion are provided in \citeasnoun{salazar} {and} \citeasnoun{swaford}.

Using dimensional analysis arguments, Richardson predicted that the mean-square separation in turbulent flows is proportional to $t^3$. Although it is generally accepted that the mean-square separation is {\em not} proportional to $t$ and that it grows faster than that, there is some controversy about the exponent $3$, referred to as the {\em Richardson constant}. For example, a $t^2$ dependence is advocated {by} \citeasnoun{Bourgoin}. Dispersion faster than $t$ is often referred to as {\em superdiffusion}. As discussed in Section~\ref{sec:mod-frac}, the fractional Laplace operator, when substituted for the classical  Laplacian, is known to result in superdiffusive spread in `parabolic' equations. {\citeasnoun{GJX} studied a} modification of the Navier--Stokes equations introduced in \citeasnoun{Chen}, where the modification is that a fractional Laplacian viscous term is {\em added} to the standard Laplacian viscosity term. 

The model considered is given by
\begin{equation}\label{rich:mns}
\begin{cases}
\ub_t + (\ub\cdot\nabla) \ub -\nu \Delta \ub+ \sigma (-\Delta )^s \ub+\nabla p = \fb & \text{in }  (0,T] \times \omg , \\
 \nabla \cdot \ub =0  & \text{in }  (0,T] \times \omg,
\end{cases}
\end{equation}
where $s\in(0,1)$, $\omg\subset \Rd$ denotes a bounded, open domain, $\sigma$ denotes a constant and $[0,T]$ denotes a temporal interval of interest.

{\citeasnoun{GJX} have investigated the energy spectrum of the modified Navier Stokes equations \eqref{rich:mns}. They have shown that} for the special value of the fractional exponent $s=1/3$, the corresponding power law of the energy spectrum in the inertial range has a deviation from the well-known Kolmogorov $-5/3$ scaling, \tie, instead of a $k^{-\gfrac53}$ decay in the spectrum, one has $k^{-\gfrac53 +\beta}$ for a constant $\beta$, {where $k$} denotes the frequency. For other values of $s\in(0,1)$, the power law of the energy spectrum is consistent {with Kolmogorov's} theory.

An IMEX scheme can be used to discretize the system \eqref{rich:mns}; a backward-Euler scheme of that type is considered in \citeasnoun{GJX} and is proved to be unconditionally stable and first-order accurate. Both the usual viscous term $\nu \Delta \ub$ and the added fractional viscous term $\sigma (-\Delta )^s \ub$ are treated implicitly. Because of the implicit treatment of the fractional Laplacian term, the scheme requires the solution of a {\em dense} linear system at each time step. Having to also handle the standard Navier--Stokes terms makes for an even greater computational challenge. Thus it is tempting to lag the fractional term to the previous time step. Unfortunately, this leads to serious stability issues, {so that} term has to be treated implicitly. To mitigate these challenges, the following two-stage operator splitting algorithm is introduced in \citeasnoun{GJX}.

\paragraph{Stage 1.}  Given $ u^n$, find $w^{n+1}$ satisfying}
\[
\begin{cases}
\ffrac{{\wb}^{n+1}-\ub^n}{\Delta t}+(\ub^n\cdot\nabla)\wb^{n+1}
 -\nu \Delta  \wb^{n+1} +\nabla p^{n+1} =\fb^{n+1}-\sigma {(-\Delta)^\alpha {\ub^n}} , \\[9pt]
\nabla \cdot \wb^{n+1} =0 .
\end{cases}
\]
 
\paragraph{Stage 2.}  Given $\ub^n$ and $\wb^{n+1}$, find $\ub^{n+1}$ satisfying
\[
2\Delta t\sigma{(-\Delta )^\alpha} ({\ub^{n+1}}-u^n)
+ \ub^{n+1}-\wb^{n+1}  =0 .
\]
The Stage 1 problem can be solved using a legacy Navier--Stokes code with the only modification necessary being in the construction of the right-hand side. In the second  {stage we solve} a nonlocal Poisson problem for the fractional Laplace operator which involves a symmetric, positive definite, albeit dense linear system. This two-stage algorithm, although involving two linear system solves per time step, when compared to the algorithm in which the two viscous terms are treated monolithically, requires much less coding effort and introduces efficiencies not possible for the monolithic scheme.

Of course, we have glossed over the fact \roo{that} the fractional Laplacian involves an integral over all of $\Rd$. Moreover, we also have to choose a spatial discretization scheme. For the first of these, we can,
\pagebreak %%20200331
\feeg, use the truncated approximation discussed in Section~\ref{truncted-flap}, whereas for the latter we can use a finite element method based {on the Taylor--Hood finite element pair \cite{Girault}, \feeg}. Such fully discrete schemes involving the two-stage time-stepping algorithm {have proved} to be unconditionally stable and satisfy the error estimate
\[
\| u_\infty(t^n)-u^n_{h,{\del}}\|_{L^2(\omg)} \leq C \biggl(\ffrac{1}{\del^{2s}}+h^3+\Delta t\biggr),
\]
where $u(t^n)$ denotes the exact solution evaluated at time $t=t_n$ and with no truncation of the fractional Laplace operator involved, $u^n_{h,\del}$ denotes the approximate solution, and $h$ and $\Delta t$ denote spatial and temporal grid size parameters, respectively. Note that for $s=\gfrac13$ we have that 
$
\gfrac{1}{\del^{2s}} = \gfrac{1}{\del^{2/3}}.
$

\begin{rem}[another fractional turbulence model]\hspz%
We point out that a similar fractional model for turbulence has been introduced in \citeasnoun{Song2018}; here, the authors identify a universal form of fractional order that holds for any Reynolds number. The same model is further investigated in \citeasnoun{Pang2019nPINNs} in the context of fractional parameter identification\roo{;} see Section~\ref{sec:inverse}. \hfill\actabox \end{rem}

\vspace{5pt}
\section{Peridynamics models for solid mechanics}\label{sec:mod-peri}

In {Sections~\ref{sec:mod-strong}--\ref{sec:mod-frac} we considered} nonlocal models for scalar-valued functions that are appropriate for {modelling  \bo{nonlocal} diffusion problems, \feeg}. To provide an example of a nonlocal model for vector-valued functions, in this section we consider a {\em peridynamics model} for the vector-valued displacement function in solid mechanics. Peridynamics is a nonlocal, continuum model that has been shown to provide an effective means for nucleating and propagating defects such as fractures \cite{HHB12,lipton1,lipton2,Bobaru15,dtt18je}. 

The particular model we consider is the linear {\em state-based peridynamics model} for solid mechanics, introduced and analysed in \citeasnoun{du14je}, that is a nonlocal analogue of the classical Navier equations of linear elasticity. The model of \citeasnoun{du14je} is a generalization of the peridynamics model introduced in \citeasnoun{SEWXA07} {and} \citeasnoun{Silling10} when specialized to linear constitutive laws. The model of \citeasnoun{du14je} is defined {via} an energy minimization principle for which the corresponding Euler--Lagrange equation provides a weak formulation of the problem. It should be noted that although the model of \citeasnoun{du14je} generalizes the models in \citeasnoun{SEWXA07} {and} \citeasnoun{Silling10}, those earlier models {are ubiquitous} in peridynamic computations. Additionally, it should be noted that the models of \citeasnoun{SEWXA07} {and} \citeasnoun{Silling10} cannot be cast in terms of fractional kernels.

To define the model considered, we need to introduce the operators\footnote{In \eqref{str-divs}, the operator $\mcDsd$ is defined {via} its action on a scalar-valued function $u(\xb)$, resulting in $\mcDsd u$ being a vector-valued function of $\xb$ and $\yb$. In \eqref{weak-divs}, $\mcDsd$ is defined {via} its action on a {\em vector-valued function} $\ub(\xb)$, resulting in $\mcDsd \ub$ being a tensor-valued function of $\xb$ and $\yb$. This is entirely analogous to the local case in which $\nabla u$ is a vector and $\nabla\ub$ is a tensor.} $\mcDsd$ and $\mcDdsw $, {which} are respectively defined {via} their action on a vector-valued function $\ub(\xb)$ by
\begin{equation}\label{weak-divs}
{
(\mcDs\ub) \xyp  := - (\ub(\yb) - \ub(\xb))\otimes \ffrac{{\yb}- {\xb}}{|{\yb} -  {\xb}|}}
\end{equation}
and
\[
{
(\mcDdsw \ub) ({\xb}):= \int_{\omg} (\mcDs \ub)\xyp \omega_\del\xyp\dyb,}
\]
where  
\[
  \omega_\del\xyp:=   
 \ffrac{d}{{\mmm_\del}(\xb)}  \gamma_\delta(|\yb-\xb|)|\yb-\xb|
\]
{and}
\[
{ \mmm_\del(\xb) := \int_{\omg} \gamma_\del({|\yb-\xb|})|\ymx|^2 \dyb .}
\]

We define the energy space
\[
\VVpd= \{ \ub\in [ L^2(\omgomgi)]^\ddd \colon 
\|\ub\|_{\VVpd} < \infty, \; \ub\mid_{\omgid}=0\}
\]
equipped with the norm
\[
\|\ub\|_\VVpd = \int_\omgomgi\int_\omgomgi  \gamma_\del(|\ymx|)  (\text{Tr}(\mcDs \ub)\xyp) ^2
\dydx,
\]
where $\text{Tr}(\cdotsp)$ denotes the trace operator acting on a tensor. 
{Also, for functions $\ub(\xb)\in{\VVpd}$ and $\vb(\xb)\in{\VVpd}$, we define} the symmetric bilinear form
\begin{align}\label{peri-bform}
\mcPd(\ub, \vb) & :=\int_{\omg} \bigg(
\biggl(  k - \ffrac{\alpha \,m_\del(\xb)}{\ddd^{2}}\biggr)\text{Tr}(\mcDdsw \ub)(\xb)\text{Tr}(\mcDdsw \vb)(\xb) \notag
\\*
 & \quad\, + {\alpha
\int_{\omg}
  \gamma_\del({|\yb-\xb|})
  \text{Tr}(\mcDs \ub)\xyp \text{Tr}(\mcDs {\vb})\xyp
  \dyb\bigg) \dxb,}
\end{align}
{where $k$ and $\alpha$ denote scalar constants} related to the bulk and shear modulus of the material.

\pagebreak %%20200331
With the notation established, the weak formulation of the peridynamics problem considered is given as follows{:} 
\begin{equation}\label{peri-weak}
  \begin{minipage}[b]{20pc}
{Seek} $\ub(\xb)\in{\VVpd}$ such that $\ub(\xb)=\zzb$ on $\omgid$ and
\[
\mcPd(\ub, \vb)= (\fb,\vb) \ \text{for all } \vb\in {\VVpd}
\]
with $\vb(\xb)=\zzb$ on $\omgid$,
\end{minipage}
\end{equation}
where $(\cdotsp,\cdotsp)$ denotes the $[L^2(\omg)]^\ddd$ inner product. For simplicity, we have imposed a homogeneous volume constraint on the displacement $\ub$. Nonlocal traction constraints and inhomogeneous volume constraints can also be treated.

{\citeasnoun{du14je} established the} coercivity and continuity of the bilinear form $\mcPd(\cdotsp,\cdotsp)$ with respect to the energy space ${\VVpd}$, {so the well-posedness of problem} \eqref{peri-weak} can be rigorously established via the Lax--Milgram theorem. Additional discussions, including studies with other types of nonlocal constraints, can be found in
{\citename{md15non} \citeyear{md15non,MeDu16} and} \citeasnoun{du19cbms}.

\begin{rem}\label{rem:jump}\hspz%
  Some choices for \boo{the  kernel} function $\gamma_\del(|\zb|)$ 
are those {with} bounded second moments, including those  
that are comparable to $|\zb|^{-\ddd-2s}$ in the sense that {they} have the same singular behaviour at the origin as that of the fractional kernel \eqref{fraker}. Of particular interest in the peridynamics setting are \bo{kernel functions that are both radial and integrable and also} fractional kernels with $s<\gfrac12$ because the corresponding function spaces in which the peridynamics model is {well-posed} admit functions with jump discontinuities. As such, discontinuities in the displacement $\ub(\xb)$ can be viewed as fracture. Functions comparable to $|\zb|^{-1}$ are integrable for $\ddd\ge2$ and are the ubiquitous choice made in peridynamics modelling, starting {with the works of \citeasnoun{Silling00}, \citeasnoun{SEWXA07} and \citeasnoun{Silling10},} in which peridynamics was first introduced.  \hfill\actaqed
\end{rem}

The strong formulation corresponding to the weak formulation \eqref{weak-weakf} {has also been derived by} \citeasnoun{du14je}.
In the local limit, \ie\ as $\del\to0$, the \bo{weak} solutions of {the} nonlocal state-based peridynamic model \bo{with an appropriately chosen kernel} converge to that of the classical Navier--Cauchy equation of linear elasticity.

\begin{rem}[bond breaking]\label{rem:bondb}\hspace{-6.3pt}%  
  There is more to peridynamics than just the weak and strong formulations introduced above. Certainly, {as discussed in Remark~\ref{rem:jump}, those models admit} solutions having jump discontinuities \bo{over \roo{co-dimension~$1$} interfaces}. For example, for \bo{kernels $\gamd(\zb)$ that are both radial and integrable,}   
  solutions are no smoother than the data so that, in general, discontinuities in the boundary data or in the source term will result in discontinuities in the solution. However, in the {time-dependent} setting, fracture and other defects can arise, \feeg\ due to large tensile loads, even if the initial condition data and other data are smooth. Thus, to model fracture nucleation, the peridynamics model equations, whether in weak or strong form, have {been supplemented by a {\em bond-breaking rule} in some cases}. Basically, such rules state that two points $\xb$ and $\yb$ that are initially bonded, \ie\ {within} a distance $\del$ from each other, become un-bonded if at a later time those points become separated by a distance greater than $\del$. Bond-breaking rules and how {they} effect
  the nucleation of cracks are more complicated than the simplistic description just given, but further discussion is beyond the scope of this article. Detailed explanations may be found {in \citeasnoun{Silling00},
    {\citename{lipton1} \citeyear{lipton1,lipton2}}, \citeasnoun{lipton3} and \citeasnoun{dtt18je}, \feeg}.  \hfill\actaqed
\end{rem}
 
For additional detailed discussions about the material presented in this section, see \eg\ \citeasnoun{du12sirev}, \citeasnoun{MeDu16} {and} \citeasnoun{du19cbms}\ro, as well as \citeasnoun{du14je}.

Numerical methods for both the strong and weak formulations of the peridynamics model have been developed.
For example, consider the finite element methods discussed in Section~\ref{femnl}.
Let $\{ V_{p,\del, h} \} \subset  \VVpd$
denote a family of finite element subspaces, where
$h$ characterizes the mesh size, and
 for any $\vb\in \VVpd$ we have a family of elements
$\{\vb_h\in V_{p,\del, h}\}$ such that $
\| \vb_h-\vb\|_{ \VVpd}\to0$
as $h\to0$ for any fixed $\del>0$. 
For the finite element discretization of the state-based linear peridynamic model \eqref{peri-weak},  {we merely replace $\VVpd$ with $ V_{p,\del, h}$ in that equation to arrive at the following discrete problem:} 
\begin{align}\label{eq:wfstateapprox}
 \text{{Find} $ \ub_{\del,h}\in V_{p,\del, h} $ such that }
      {\mcPd}(\ub_{\del,h}, \vb) =({\bm f}, \vb)_{L^2} \ \text{for all 
      $\vb\in V_{p,\del, h}$.}
\end{align}
The analysis of solutions of this problem can be formulated within the general framework of the AC schemes originally
presented in \citeasnoun{TiDu14}.
It is rigorously shown in \citeasnoun{TiDu14} that, for linear multidimensional state-based peridynamic systems associated with \eqref{peri-weak},  all conforming Galerkin approximations of the nonlocal models containing continuous piecewise linear functions are automatically AC. This means that they can recover the correct local limit as long as both $\del$ and $h$ are diminishing, even if the nonlocality measure $\del$ is reduced at a much faster pace than the mesh spacing $h$. 

As was the case for the nonlocal diffusion case, {if} $h=o(\del)$  as $\del\to0$, then we are able to obtain the correct local limit even for discontinuous piecewise constant finite element approximations when they are of the conforming type. Practically speaking, this implies that a mild growth of bandwidth in the finite element stiffness matrix is needed as the mesh is refined, in order to recover the correct local limit for {Riemann sum-like} quadratures or piecewise constant finite element
schemes. On the other hand, if a constant bandwidth is kept as the mesh is refined, as is often  \roo{practised} in the peridynamics community, approximations may converge to an incorrect local limit.

\vspace{5pt}
\section{Image denoising}\label{sec:imaging}

We describe the use of nonlocal diffusion operators in a variational setting for image denoising. In this context, by considering intensity patterns in neighbourhoods of points surrounding a pixel, nonlocal filters allow for simultaneous conservation of structures (patterns) and textures. Classical methods, which use differential operators, do not necessarily guarantee feature preservation because, by definition, they only consider infinitesimal neighbourhoods around points. 

Even though nonlocal-type filters have been used for decades, we refer to \citeasnoun{Buades2010} as the first foundational work that can be related to a nonlocal diffusion equation as presented in this article. In that paper, the {\it nonlocal-means} (NL-means) filter for image denoising is proposed. Given a blurred (or noisy) image $f$ defined by its intensity in the image domain $\omg$, the NL-means filter is defined as 
\begin{equation}\label{eq:NLmeans}
N\!L[f](\xb) = \ffrac{1}{c(\xb)}\int\limits_\omg \rme^{-d_a(f(\xb),f(\yb))/D^2} f(\yb)\dyb \quad\text{for all } \xb\in\omg  
\end{equation}
where
\[
d_a (f(\xb),f(\yb))  = \int\limits_\omg G_a(\zb)|f(\xb+\zb)
-f(\yb+\zb)| \D \zb,
\]
with $G_a(\zb)$ a Gaussian with standard deviation $a$ and with $c(\xb)$ a normalizing factor. Note that the nonlocal filter is directly applied to the blurred image, \tie, the nonlocal operator is not associated with a diffusion equation. 

The variational interpretation of such a filter {was introduced by \citeasnoun{Kindermann2005}, who interpreted the NL-means and other neighbourhood filters} as regularizing functionals. Specifically, they {formulated} the denoising problem as an unconstrained optimization problem where the objective is to minimize the nonlocal functional
\begin{equation}\label{eq:J-regularization}
\mcJ_{R} u = \int\limits_\Omega \int\limits_\Omega 
g\biggl(\dfrac{|u(\xb)-u(\yb)|}{D^2}\biggr) {\nu}(|\xb-\yb|) \D \yb + F(u,f),
\end{equation}
where {$\nu$} is a symmetric window that determines the extent of the nonlocal interactions and $g$ determines the type of filtering. Here, $F(u,f)$ is a fidelity term, usually a measure of the distance between the reconstructed and noisy images. The outcome of the optimization is the reconstructed image $u$. {For} NL-means filtering, \eqref{eq:J-regularization} becomes
\begin{equation}\label{eq:J-NLmeans}
\mcJ_{N\!L} u = \int\limits_\Omega \int\limits_\Omega 
 (1-\rme^{-d_a(u(\xb),u(\yb))/D^2})  {\nu}(|\xb-\yb|) \D \yb + F(u,f).
\end{equation}
Such a minimization problem cannot be directly related to the nonlocal diffusion theory reported {in} this article because of the potential non-convexity of the functionals. For the same reasons, {\citeasnoun{Gilboa2007} introduced a modified convex functional (still based on nonlocal filters) that} resembles common functionals used in PDE-based imaging approaches such as total variation {functionals}
\cite{Rudin1992}. In Section~\ref{eq:J-regularization2}, we describe methods based on this concept and highlight their connection to the solution of nonlocal diffusion problems, \bo{including integral fractional problems}. In all the examples below, unless otherwise stated, the relation between clean and noisy image is the {\em additive noise} model, \ie\ $f=u+\eta$, where $u$ denotes the clean image, $f$ the noisy {image} and $\eta$ the noise.

\subsection{Image deblurring via minimization of a nonlocal functional}\label{eq:J-regularization2}

{\citeasnoun{Gilboa2007} tackled the image deblurring problem by} solving an unconstrained optimization problem with \bo{objective} functional
\begin{equation}\label{eq:J-Gilboa}
\mcJ(u;f) = \dfrac{1}{2}{|u|_{V_c(\Omega)}^2} + \dfrac{w}{2}\|u-f\|^2_{L^2(\Omega)},
\end{equation}
where {$|\cdotsp|_{V_c}$ denotes the energy seminorm defined in \eqref{weak-enormc}} associated with an appropriate kernel $\gamd$, whose definition depends on the type of filtering (either NL-means or other neighbourhood filters) and the fidelity term controls the difference between the reconstructed image and the noisy {image}. Note that the weight $w\in\Ro^+$ could be replaced by a {spatially dependent} function that is included within the {$L^2$-norm}. 

The Euler--Lagrange equation that determines the necessary conditions for optimality is given by the nonlocal reaction--diffusion equation
\begin{equation}\label{eq:euler-lagrange}
-\mcL u + w(u-f) = 0,
\end{equation}
whose {well-posedness} can be determined as described in Section~\ref{weak-wposed} based on the properties of the kernel $\gamd$. 

\begin{figure} %%fig21.1 %%wasfig22.1
\centering
\subfigure[clean]{\includegraphics[width=166pt,viewport=15 371 232 681,clip]{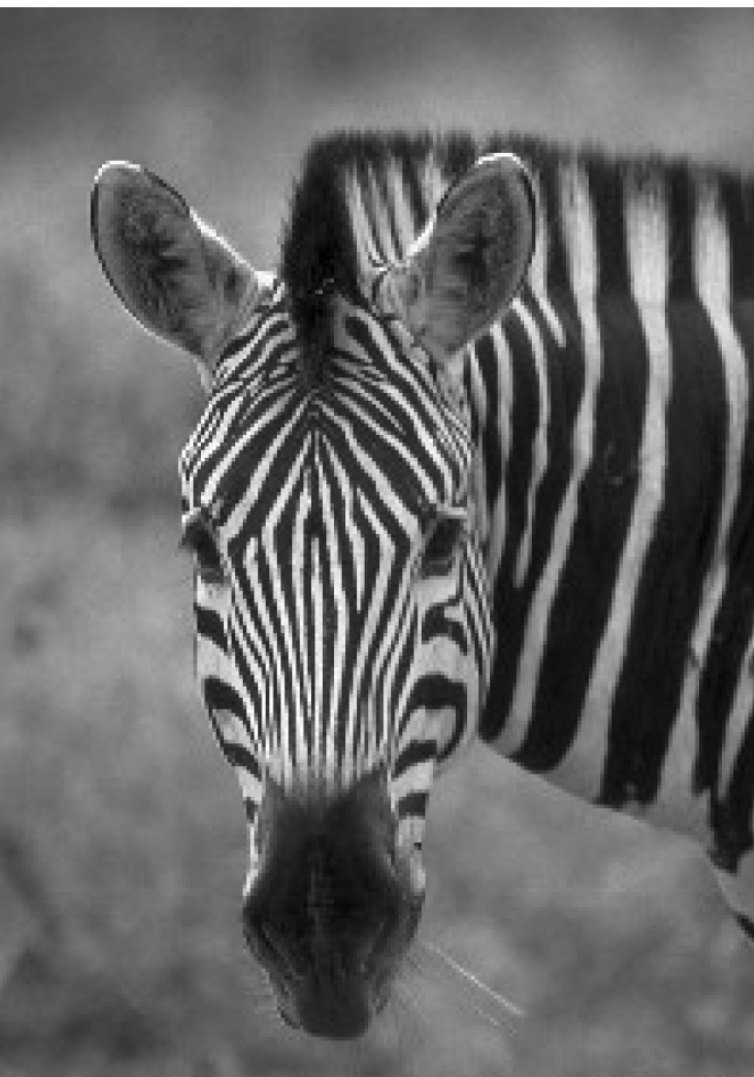}}\hfill %%imaging
\subfigure[noisy]{\includegraphics[width=166pt,viewport=272 371 489 681,clip]{Figures/fig21_1.pdf}}\\ %%imaging
\subfigure[nonlocal method]{\includegraphics[width=166pt,viewport=16 7 233 312,clip]{Figures/fig21_1.pdf}}\hfill %%imaging
\subfigure[local method (TV)]{\includegraphics[width=166pt,viewport=272 7 489 312,clip]{Figures/fig21_1.pdf}}\\ %%imaging
\caption{Clean, noisy, and nonlocal and local reconstructions of a noisy image featuring high-contrast features.}
\label{fig:imaging}
\end{figure}

\fo{Figure~\ref{fig:imaging}}, {taken from \citeasnoun{Gilboa2007},} provides an example of a nonlocal deblurring reconstruction and also a comparison with a classical method local method, namely the {total variation} (TV) method of \citeasnoun{Rudin1992}. The clean and noisy images are given in {Figure~\ref{fig:imaging}(a,\,b)
 and the nonlocal and local reconstructions are given in Figure~\ref{fig:imaging}(c,\,d)}. These results show the {superior ability of nonlocal filters to capture} high-contrast features compared to classical methods.

{\citeasnoun{Gilboa2008}, in a follow-up paper, provided} additional insights about the theoretical aspects of the method of \citeasnoun{Gilboa2007} and about other types of filters, and also {provided} comparisons with other methods.

\begin{rem}\hspz% 
A very similar formulation can be used for image deconvolution \cite{dt-LZOB10} where the relation between {a} clean and noisy image, $f=\mcK u + \eta$, involves a convolution operator $\mcK$. This renders the denoising problem more complicated. In this setting, the \bo{objective}  functional of the unconstrained optimization problem is given by
\begin{displaymath}
\mcJ_{C}(u) = \dfrac{1}{2}{|u|_{V_c(\Omega)}^2} +
\dfrac{w}{2}\int\limits_\Omega (k(\xb)\ast u(\xb)-f(\xb))^2 \D \xb
\end{displaymath}
and the Euler--Lagrange equation is the nonlocal diffusion equation
\begin{displaymath}
-\mcL u = w\, {\widetilde k \ast (f-k\ast u)},
\end{displaymath}
where \bo{$\widetilde k$ denotes the adjoint of $k$}.  \hfill\actaqed
\end{rem}

\begin{rem}\hspz% 
Similar operators to those discussed in this section have also been used for image segmentation \cite{Gilboa2007}, filtering \cite{Darbon2008} and image and video recovery \cite{Buades2008}. \hfill\actaqed
\end{rem}

\begin{rem}\hspz% 
In the literature, {denoising} nonlocal operators are not necessarily associated with nonlocal filters such as NL-means, but can be described by fractional operators, specifically by the spectral fractional Laplacian{;} see \eg\ \citeasnoun{Antil2017imagingSpectral}.  \hfill\actaqed
\end{rem}

\subsubsection{Optimization of the denoising parameters}

The quality of the reconstruction strongly depends on model parameters (\eg\ parameters in the kernel) and on the weight parameter $w$. However, model parameters are often unknown (see Section~\ref{sec:control}) and the selection of $w$ is not a trivial task.
{\citeasnoun{DElia2019imaging} considered a bilevel optimization approach for} the identification of kernel parameters ({for} integrable kernels, including NL-means) and of the weight function $w\colon \omg\to\Ro^+$. Because the estimation of kernel parameters is discussed in Section~\ref{sec:control}, here we only consider the formulation for the identification of $w$, in its simplest setting. 

Given a clean image and the corresponding blurred image, $(\hat u,\hat f)$, the bilevel optimization problem is formulated as 
\begin{align}\label{eq:lambda-opt}
\min_{w\in\WWW}
	& \;\; \biggl(\mcJ_w(u,w)= \|u-\hat u\|^2_{L^2(\Omega)} + \ffrac{\beta}{2} \|w\|^2_{H^1(\Omega)} \Biggr)\notag \\*
	\text{such that}
	& \;\; u = \arg\min_{u\in V_c} \biggl( \mcJ_u(u,w)=
	\dfrac{1}{2}{|u|_{V_c(\Omega)}^2} +\int_\Omega w (u-\hat f)^2 \D \xb \Biggr),
\end{align}
where the feasible set is $\WWW = \{w\in H^1(\Omega)\colon 0\leq {w(\xb) \leq \mathrm{w}}\}$, \bo{for some \roo{$\mathrm{w}<\infty$}}. Note that the constraint in \eqref{eq:lambda-opt} is equivalent to the diffusion--reaction equation in \eqref{eq:euler-lagrange}. The {well-posedness} of the bilevel optimization problem {has been proved by \citeasnoun{DElia2019imaging}. They also introduce a} second-order optimization algorithm for its solution, and {give insights into} implementation aspects and numerical performance. Also, {they illustrate the theory and advantages of using nonlocal models via} several numerical tests on standard benchmark images.

\section*{Acknowledgements}

{The research of MD, CG and (partially) MG is supported by Sandia\linebreak %%20200331
  National Laboratories (SNL). SNL is a multimission laboratory managed and operated by National Technology and Engineering Solutions of Sandia, LLC, a wholly owned subsidiary of Honeywell International, Inc., for the {US} Department of Energy's National Nuclear Security Administration under contract DE-NA-0003525. Specifically, this work was supported by the SNL Laboratory-directed Research and Development (LDRD) {programme}.
  Note that this paper describes objective technical results and analysis. Any subjective views or opinions that might be expressed in the paper do not necessarily represent the views of the {US} Department of Energy or the United States Government.}

{The research of QD is supported in part by the US National Foundation grants NSF DMS-1719699 and NSF CCF-1704833, the US Air Force Office of Scientific Research MURI Center for Material Failure Prediction Through Peridynamics, and by the US Army Research Office MURI grant W911NF-15-1-0562.} 

{The research of XT is supported in part by the {US} NSF grant DMS-1819233.}

{The research of ZZ is supported by a start-up grant from the  Hong Kong Polytechnic University  and Hong Kong RGC grant 25300818.}

Several figures in the article are reproduced by permission. Specifically, 
Figure~\ref{fig:stretch} is reproduced from \citeasnoun{xu2016multiscale},
Figures~\ref{fig:tian2019fast_kernel} and
\ref{fig:tian2019fast_decomposition} from
\citeasnoun{tian2019fast},
Figure~\ref{fig:control-paper} from  \citeasnoun{DElia2014DistControl},
Figure~\ref{fig:identification-paper} from
\citeasnoun{DElia2016ParamControl},
Figure~\ref{fig:obstex} from
\citeasnoun{burko1},
Figure~\ref{fig:rom} from
\citeasnoun{witman} 
and
Figure~\ref{fig:imaging} from
\citeasnoun{Gilboa2007}.

\newpage  %%20200331
%GS%\bibliography{nonlocal}

\end{document}